\pdfoutput=1
\documentclass[10pt]{amsart} 
\usepackage[final]{showlabels}
\usepackage{amsmath, accents}
\usepackage{amssymb}
\usepackage[toc,page]{appendix}
\usepackage{array}
\usepackage{bm, bbm}
\usepackage{calligra}
\usepackage{dsfont}
\usepackage{enumerate}
\usepackage{enumitem}
\usepackage{fullpage}
\usepackage{graphicx}
\usepackage{color}
\usepackage{xcolor}
\usepackage{mathMacros}
\usepackage{mathtools}
\usepackage{multirow}
\usepackage{pifont}
\usepackage{caption}
\usepackage{subcaption}
\usepackage{stmaryrd}
\usepackage{tablefootnote}

\usepackage[backend=biber,citestyle=numeric, maxbibnames=99, hyperref=true,backref=true]{biblatex}

\usepackage{hyperref}
\usepackage{nameref}
\usepackage{cleveref}
\usepackage{autonum}
\cslet{blx@noerroretextools}\empty 

\usepackage{hyperref}
\hypersetup{
    colorlinks=true,
    linkcolor={red!90!black},
    citecolor={red!90!black},
    filecolor={red!90!black},      
    urlcolor={red!90!black},
    }
\colorlet{linkequation}{red!90!black}
\numberwithin{equation}{section}
\theoremstyle{plain}
\newtheorem{theorem}{Theorem}[section]

\newtheorem{claim}[theorem]{Claim}
\newtheorem{corollary}[theorem]{Corollary}
\newtheorem{lemma}[theorem]{Lemma}
\newtheorem{proposition}[theorem]{Proposition}
\newtheorem{definition}[theorem]{Definition}
\newtheorem{assumption}[theorem]{Assumption}
\newtheorem{remark}[theorem]{Remark}

\theoremstyle{remark}

\renewcommand{\epsilon}{\varepsilon}
\def\scA{{\rmH}}
\def\hatH{{\widehat{\scA}^{(\vtau)}}}
\def\hatHx{{\widehat{\scA}^{(\vtau,\, x)}}}
\def\ratio{{\gamma}} 
\def\qr{{\mathfrak{g}}} 
\def\benrate{{\mathrm{h}}} 
\def\error{{\mathfrak{\xi}}} 
\def\radius{{\mathrm{r}_{\star}}} 
\def\resolvent{{\rmG}}
\def\eres{{\ermG}}
\def\Vonehigh{{\gV_{1}^{(\geq \tau_{1})}}}
\def\Vtwolow{{\gV_{2}^{(\leq \tau^{-}_{2})}}}
\def\Vtwohigh{{\gV_{2}^{(\geq \tau^{+}_{2})}}}
\def\pruneG{{\gG^{(\vtau)}}}
\def\pruneV{{\gV^{(\vtau)}}}
\def\hatvtau{{\widehat{\rvv}^{(\vtau)}}}
\newcommand{\pruneR}[1]{{\mathrm{r}_{#1, \vtau}}}

\newcommand*\circled[1]{\tikz[baseline=(char.base)]{
            \node[shape=circle,draw,inner sep=1pt] (char) {#1};}}

\DeclareMathAlphabet{\mathcalligra}{T1}{calligra}{m}{n}
\DeclareMathAlphabet{\mathpzc}{OT1}{pzc}{m}{it}

\DeclareRobustCommand*\circnum[1]{%
  \tikz[baseline=(num.base)] 
    \node[draw,circle,inner sep=1pt] (num) {\small #1};%
}

\usepackage{todonotes}
\makeatletter
\define@key{todonotes}{ID}[]{%
    \setkeys{todonotes}{author=ID, inline, color=magenta}}%
\makeatother
\makeatletter
\define@key{todonotes}{HW}[]{%
    \setkeys{todonotes}{author=HW, inline, color=orange}}%
\makeatother       
\makeatletter
\define@key{todonotes}{ZW}[]{%
    \setkeys{todonotes}{author=ZW, inline, color=green}}%
\makeatother 

\makeatletter
\define@key{todonotes}{YZ}[]{%
    \setkeys{todonotes}{author=YZ, inline, color=pink}}%
\makeatother

\newcommand{\toromanlower}[1]{\ifcase#1 \or i\or ii\or iii\or iv\or v\or vi\or vii\or viii\or ix\or x\fi}
\newcommand{\circromlower}[1]{%
  \tikz[baseline=(char.base)] \node[shape=circle,draw,inner sep=1pt] (char) {\normalfont\small\toromanlower{#1}};%
}

\usepackage{orcidlink}

\title{Singular values of sparse random rectangular matrices: Emergence of outliers at criticality}

\begin{document}

\author{Ioana Dumitriu}
\address{Department of Mathematics, University of California, San Diego, La Jolla, CA 92093}
\email{idumitriu@ucsd.edu}

\author{Hai-Xiao Wang \orcidlink{0000-0003-2730-1439}} 
\address{Department of Mathematics, University of California, San Diego, La Jolla, CA 92093}
\email{h9wang@ucsd.edu}

\author{Zhichao Wang}
\address{International Computer Science Institute and
       Department of Statistics,
       University of California, Berkeley,
       Berkeley, CA 94720}
\email{zhichao.wang@berkeley.edu}

\author{Yizhe Zhu}
\address{Department of Mathematics, University of Southern California,  Los Angeles, CA 90007}
\email{yizhezhu@usc.edu}


\date{\today}

\begin{abstract}
Consider the random bipartite Erd\H{o}s-R\'{e}nyi graph $\gG(n, m, p)$, where each edge with one vertex in $\gV_{1}=[n]$ and the other vertex in $\gV_{2} =[m]$ is connected with probability $p$, and $n=\lfloor \gamma  m\rfloor$ for a constant aspect ratio $\gamma \geq 1$. It is well known that the empirical spectral measure of its centered and normalized adjacency matrix converges to the Mar\v{c}enko-Pastur (MP) distribution. However, the largest and smallest singular values may not converge to the right and left edges, respectively, especially when $p = o(1)$. Notably, it was proved by \cite{dumitriu2024extreme} that there are almost surely no singular values outside the compact support of the MP law when $np = \omega(\log(n))$.
In this paper, we consider the critical sparsity regime where  $p = b\log(n)/\sqrt{mn}$ for some constant $b>0$. We quantitatively characterize the emergence of outlier singular values as follows. For explicit $b_{\star}$ and $b^{\star}$ as functions of $\ratio$, we prove that when $b > b_{\star}$, there is no outlier outside the bulk; when $b^{\star}< b < b_{\star}$, outliers are present only outside the right edge of the MP law; and when $b < b^{\star}$, outliers are present on both sides, all with high probability. Moreover, the locations of those outliers are precisely characterized by a function depending on the largest and smallest degree vertices of the random graph. We estimate the number of outliers as well. Our results follow the path forged by \cite{alt2021extremal}, and can be extended to sparse random rectangular matrices with bounded entries.
\end{abstract}


\maketitle

\tableofcontents

\section{Introduction}\label{sec:intro}
In real-world applications, large matrices encountered in science, engineering, and mathematics are often \emph{sparse}, that is, containing predominantly zero entries. The study of sparse matrices resides at the intersection of numerical linear algebra, graph theory, and probability theory. A simple yet effective model, that bridges these three areas and has attracted great interest since the 1950s, is the \emph{Erd\H{o}s-R\'enyi} graph. Initially, it was introduced to demonstrate the existence of graphs with specific properties. Since then, its simple but mathematically rich structure has provided a bevy of theoretical uses in computer science, physics, and engineering: as the stepping stone or building block model in the study of networks, as an expander, as a basic paradigm for studying particle interaction systems. Precisely, let $\gG(N,p)$ denote the Erd\H{o}s-R\'enyi graph with $N$ vertices, where each edge $\{j, k\}\subset [N]$ with $j\neq k$ is included independently with probability $p$; its adjacency matrix $\rmA\in \{0, 1\}^{N\times N}$ is sparse when $p< 1$, but it is particularly interesting when $p \ll 1$.

This paper focuses on a \emph{bipartite} variant of the Erd\H{o}s-R\'enyi graph, denoted as $\gG(n, m, p)$, where the vertex set $\gV = [N]$ is a disjoint union of $\gV_1= \{1, \ldots, n\} = [n]$ and $\gV_2 = \{n+1, \ldots, n+m\} = [N]\setminus [n]$ with $N = n + m$. Each edge of $\gG(n, m, p)$, containing one vertex in $\gV_1$ and the other one in $\gV_2$, is sampled with probability $p$, independently of all other edges. Let $\widetilde{\rmA}\in \{0, 1\}^{n\times m}$ denote the \emph{bi-adjacency} matrix of the $\gG(n, m, p)$, where $\widetilde{\ermA}_{jl} \sim \text{Bernoulli}(p)$ independently for each $j\in\gV_{1}$ and $l\in\gV_{2}$. Then the adjacency matrix admits the following block structure:
\begin{align}\label{eqn:adj}
     \rmA \coloneqq \begin{bmatrix}
    \bzero & \widetilde{\rmA}\\
    \widetilde{\rmA}^{\sT}  & \bzero
    \end{bmatrix}. 
\end{align}

\begin{figure}[h]
    \centering
    \begin{minipage}[h]{0.5\linewidth}
\centering
{\includegraphics[height=0.5\textwidth]{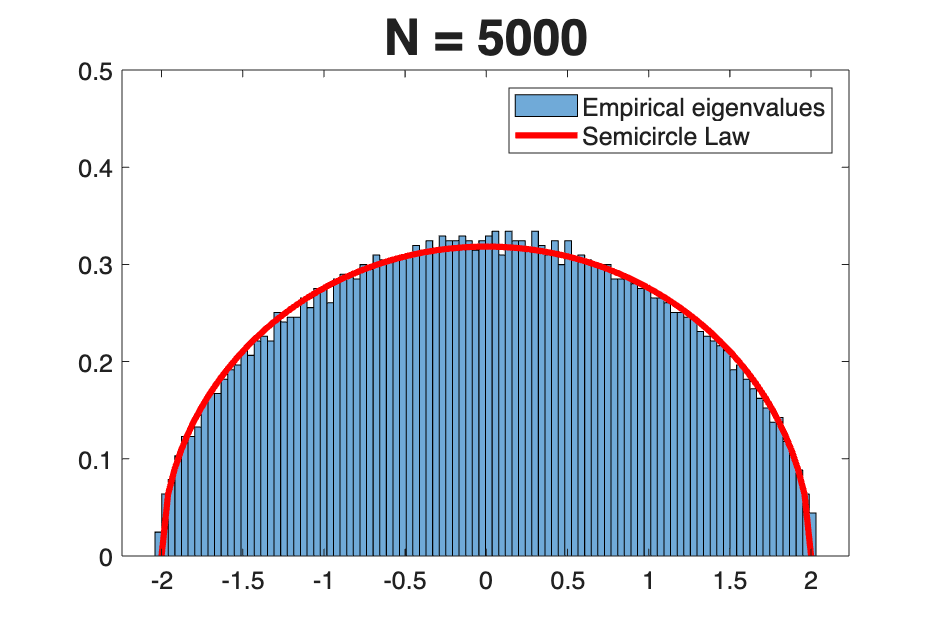}}  \\
\small (a) ESD of $\scA$ under $\gG(N, p)$.
\end{minipage}%
    \begin{minipage}[h]{0.5\linewidth}
\centering
{\includegraphics[height=0.5\textwidth]{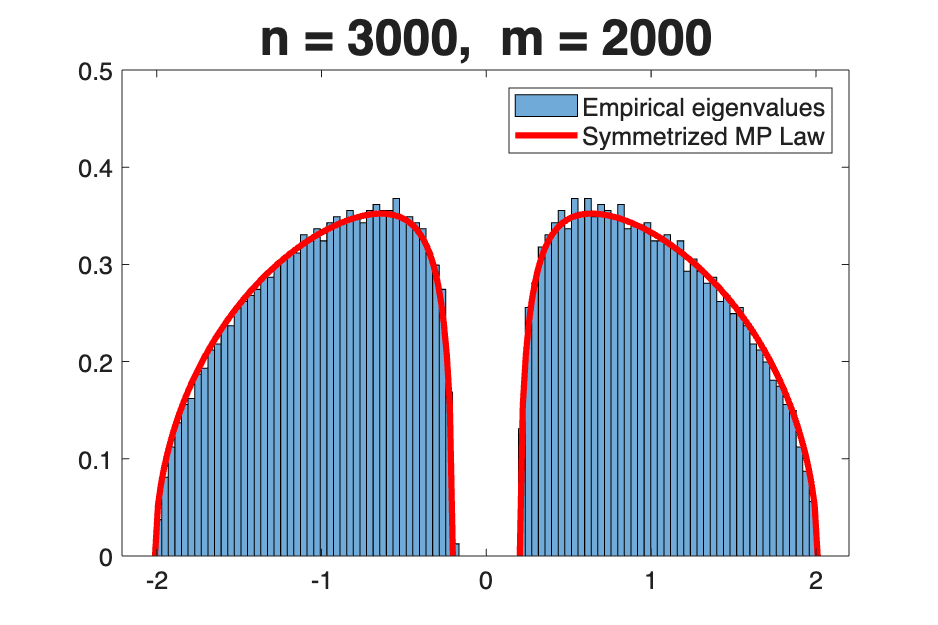}}  \\
\small (b) ESD of $\scA$ under $\gG(n, m, p)$.
\end{minipage}%
    \caption{The left plot (a) shows the semicircle law \eqref{eqn:SC_law}, while the right plot (b) shows a symmetrized version of the MP law. Plot (b) should also contain an atom at $0$ with weight $\frac{n-m}{n+m}=1/5$, which for practical purpose we didn't display.}
    \label{fig:SCandMP}
\end{figure}

Among the first properties of the (bipartite) Erd\H{o}s-R\'enyi graph that were studied were its connectivity and its spectral properties; both its adjacency matrix and its Laplacian operator (scaled or combinatorial) were subject to intense scrutiny, and their spectral properties were linked to the connectivity of the graph itself. Precisely for $\gG(N, p)$, let $\scA \coloneqq (\rmA - \E \rmA)/\sqrt{pN}$ be the centered and scaled
adjacency matrix, with $\lambda_{1}(\scA) \geq \ldots \geq \lambda_{N}(\scA)$ denoting its eigenvalues. It has long been known (e.g., \cite[Corollary 1.2]{ding2010spectral}, illustrated in Figure~\ref{fig:SCandMP}(a)) that if $Np(1-p) \to \infty$, the \emph{Empirical Spectral Distribution} (ESD) of $\scA$, defined as $\mu_{N}(\scA)\coloneqq\frac{1}{N} \sum_{j=1}^{N} \delta_{\lambda_{j}(\scA)}$, converges to the semicircle law $\mu_{\rm{SC}}$, where
\begin{align}
    d\mu_{\rm{SC}}(x)\coloneqq \frac{1}{2\pi} \sqrt{4 - x^{2}} \,\indi{x\in [-2, 2]}\,dx. \label{eqn:SC_law}
\end{align}

In the case of a bipartite Erd\H{o}s-R\'enyi graph, where we denote $\rmX = (\widetilde{\rmA} - \E \widetilde{\rmA})/\sqrt{d}$ with $d = pn$, its \emph{centered and scaled adjacency} matrix is defined by
\begin{align}\label{eqn:scAdj}
    \scA \coloneqq \frac{1}{\sqrt{d}} (\rmA -\E \rmA) = 
    \begin{bmatrix}
    \bzero &  \rmX \\
    \rmX^{\sT}  & \bzero
    \end{bmatrix}.
\end{align}
One can verify that $\rmH$ has a zero eigenvalue of multiplicity at least $n-m$, and that the nonzero eigenvalues of $\rmH$ come in pairs $\pm\lambda$, with $\lambda$ being a positive singular value of $\rmX$. 

Let $\mu_{m}(\rmX^{\sT}\rmX) = \frac{1}{m}\sum_{j=1}^{m}\delta_{\lambda_{j}(\rmX^{\sT}\rmX)}$ denote the ESD of  $\rmX^{\sT}\rmX \in \R^{m\times m}$, where $\lambda_{1}(\rmX^{\sT}\rmX) \geq \ldots \geq \lambda_{m}(\rmX^{\sT}\rmX)$ are the eigenvalues. In the dense case, when $n/m = \ratio$ for some $\ratio \geq 1$ as $n, m\to \infty$, the ESD will converge weakly to the Mar\v{c}enko-Pastur (MP) law $\mu_{\rm{MP}}$  \cite{marchenko1967distribution,tran2020local,nadutkina2024marchenko}, where
\begin{align}
    d\mu_{\rm{MP}}(x)= {\frac{\ratio\sqrt{(\lambda_{+}-x)(x-\lambda _{-})}}{2\pi x}}\,\indi{x\in [\lambda _{-},\lambda _{+}]}\,dx~.\label{eqn:MP_law_original}
\end{align}
Here $\lambda_{+} = (1 + \ratio^{-1/2})^{2}$ and $\lambda_{-} = (1 - \ratio^{-1/2})^{2}$ denote the right and left edges of the support, respectively. It is a quick exercise to check that the ESD of $\scA$ converges to a symmetrized version of the MP law, as shown in Figure~\ref{fig:SCandMP}(b).

\begin{remark}\label{remark:singular-eigenvalue}
   Given the correspondence established above, for the remainder of the paper we will often work with $\scA$ in place of $\rmX$, and we will use the terms ``singular values of $\rmX$'' and ``eigenvalues of $\scA$'' interchangeably to ease the presentation.
\end{remark}
Notably, the convergence of ESDs does not imply the convergence of extremal eigenvalues to the edges of the limiting distributions (both semicircle and MP law), especially in the sparse regime. 
The difficulties in the regime we study here stem from two aspects. The first is the sparsity of the model: for matrix entries of $\rmA/\sqrt{p}$, all their moments higher than $3$ diverge, leading to the failure of direct applications of classical Wishart/Wigner results.  The second is the complexity of the  least singular value: when the aspect ratio $\ratio \neq 1$, the left edge $\lambda_{-} = (1 - \ratio^{-1/2})^2$ of the MP law \eqref{eqn:MP_law_original} is bounded away from the origin, as shown in \Cref{fig:SCandMP}. The gap between $0$ and $\lambda_{-}$, present in rectangular Wishart matrices, is absent in the case of Wigner matrices, which makes it notoriously hard to compute the asymptotics of the smallest singular value in that case (see the literature review in Section~\ref{sec:literature}).

\subsection{Formulation of the problems}
Consider the Erd\H{o}s-R\'enyi graph $\gG(N, p)$, where $p = d/N$ with $d = b\log(N)$ for some constant $b>0$. Alt, Ducatez and Knowles \cite{alt2021extremal} proved that there exists a sharp threshold $b_{\star}= \frac{1}{\ln(4) - 1}$, such that the following holds with high probability for the eigenvalues of $\scA$:
\begin{itemize}
    \item $b > b_{\star}$: no eigenvalue appears outside the compact support of the semicircle law;
    \item $b < b_{\star}$: outlier eigenvalues emerge.
\end{itemize}

For the bipartite Erd\H{o}s-R\'enyi graph $\gG(n, m, p)$, Dumitriu and Zhu \cite{dumitriu2024extreme} demonstrated that, almost surely, no outliers exist outside the compact support of the MP law when $np \gg \log(N)$. This paper focuses on the regime where $np \asymp \log(N)$, known as the \emph{critical sparsity regime}. The following assumption on the model parameters is maintained throughout this work.
\begin{assumption}[Proportional regime at criticality]\label{ass:proportional_regime}
    Assume 
    \begin{align}
      n=\lfloor \gamma m\rfloor  \quad \textnormal{and} \quad \qr \coloneqq \ratio^{1/4}, \label{eqn:ratio}
    \end{align}
    where $\ratio \geq 1$ is some absolute constant. Additionally, the connection probability can be written as
    \begin{align}
        p\coloneqq \frac{d}{\sqrt{mn}} \quad \textnormal{and} \quad d \coloneqq b\log(N), \label{eqn:critical_regime}
    \end{align}
    where $b>0$ is some constant independent of $n,m$.
\end{assumption}

We aim to investigate the emergence of outlier singular values in this critical sparsity regime. In particular, we are interested in understanding if there exist $\ratio$-dependent \emph{thresholds} for the emergence of left, respectively right, outlier singular values in the critical sparse regime. Furthermore, we want to understand the mechanisms by which these outliers appear. 

Our main results shed light on these issues, as follows:
\begin{enumerate}
    \item Such thresholds for the parameter $b$ exist, and they depend only on $\ratio$. The thresholds for the emergence of left, respectively right, outlier singular values are generally different.
    \item The mechanism by which outlier singular values on the right are created is from the loss of degree concentration.  It relies on the existence of ``high-degree" vertices; as $b$ decreases, these vertices start appearing in $\gV_2$, and as $b$ gets smaller, they also start appearing in $\gV_1$.
    \item The mechanisms by which outlier singular values are created on the left are more complex. In part, they are generated by the emergence of ``low-degree" vertices in $\gV_2$ (perhaps surprisingly, ``low-degree" vertices in $\gV_1$ do not in fact generate outliers). Due to our current method, at lower values of $\ratio$, we can't eliminate the possibility of other mechanisms that create outliers.
    
    \item The locations of those outliers are precisely characterized by a function depending purely on the degrees of the vertices.
    \item We extend our results to more general sparse matrices resulting from a Hadamard product of the bi-adjacency matrices with i.i.d. rectangular matrices with mean 0, variance 1 entries and bounded support. This model corresponds to bipartite Erd\H{o}s-R\'enyi graphs with independent (bounded) edge weights. 
\end{enumerate}

\subsection{Emergence of outlier singular values}
Following \Cref{ass:proportional_regime} and \eqref{eqn:MP_law_original}, the singular value distribution $\mu_m(\rmX)=\frac{1}{m}\sum_{j=1}^m\delta_{\lambda_{j}(\rmX)}$ converges weakly to a distribution with density function:
\begin{equation}\label{eqn:MP_law_singular_value}
    f(s)
\coloneqq \frac{\qr^2\sqrt{(s^2-\lambda_-)(\lambda_+ - s^2)}}{\pi\, s}\cdot\indi{[\sqrt{\lambda_-},\sqrt{\lambda_+}]}
\end{equation}
where $\lambda_\pm\coloneqq  (\qr\pm\qr^{-1})^2.$ Here, the support of this limiting singular value distribution (which we will often refer to as the \emph{bulk}) is $[\qr-\qr^{-1},\,\qr+\qr^{-1}]$. We call $(\qr+\qr^{-1})$ the \textit{right edge} and $(\qr-\qr^{-1})$ the \textit{left edge} of the MP law, respectively. 

Throughout the paper, we work on the following three threshold functions depending only on the aspect ratio $\ratio=\qr^{4}\ge 1$:
        \begin{align}
            r^{\star}_2 \coloneqq &\, [(\qr^{2} + 1) \log(1 + \qr^{-2}) - 1]^{-1}, \label{eqn:right_threshold_2}\\
            \quad r^{\star}_1 \coloneqq&\, [(\qr^{-2} + 1) \log(1 + \qr^{2}) - 1]^{-1}\, ,\label{eqn:right_threshold_1}\\   l^{\star}_2 \coloneqq&\, [(\qr^{2} - 1) \log(1 - \qr^{-2}) + 1]^{-1}.\label{eqn:left_threshold}
        \end{align}
We first describe two phase transitions regarding the largest and smallest singular values of $\scA$, which can be seen as the analog of the Bai-Yin law \cite{yin1988limit, bai1993limit} for sparse random matrices at the critical sparsity level. Here, outlier singular values correspond to vertices with extreme degrees. This so-called \emph{vertex-to-outlier} mechanism that induces outliers at criticality will be illustrated in full detail in \Cref{sec:vertex_to_outlier}.
\begin{theorem}[Emergence of right outliers]\label{thm:right_outlier_emergence}
Suppose that \Cref{ass:proportional_regime} is true. Then there exists  $\nu >0$ such that the following holds with probability at least $1 - N^{-\nu}$ for sufficiently large $N$:
\begin{enumerate}
    \item (No outlier) When $b > r^{\star}_2$, \[\lambda_1(\scA)=\qr+\qr^{-1}+o(1).\]
    \item (Emergence) When $b < r^{\star}_2$, there exists a $\nu$-dependent constant $\epsilon >0$ such that, 
    \[
    \lambda_1(\scA)>\qr+\qr^{-1}+\epsilon.
    \]
\end{enumerate}
\end{theorem}

\begin{figure}[h]
    \centering
    \begin{minipage}[h]{0.5\linewidth}
\centering
{\includegraphics[height=0.75\textwidth]{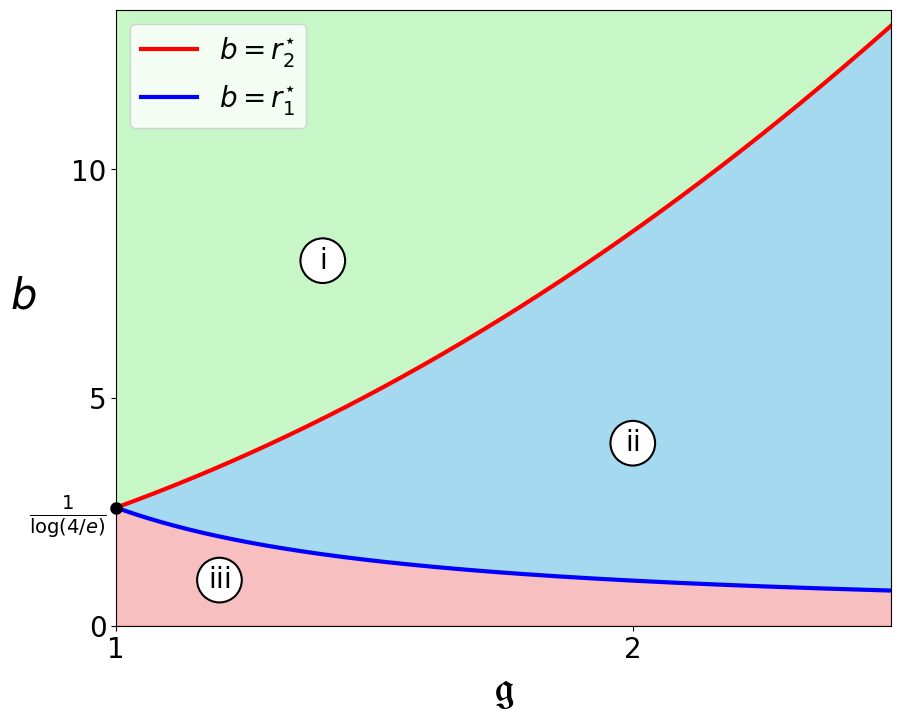}}  \\
\small (a) $r^{\star}_{2}$ and $r^{\star}_{1}$ are defined in \eqref{eqn:right_threshold_2} and \eqref{eqn:right_threshold_1}. 
\end{minipage}%
    \begin{minipage}[h]{0.5\linewidth}
\centering
{\includegraphics[height=0.75\textwidth]{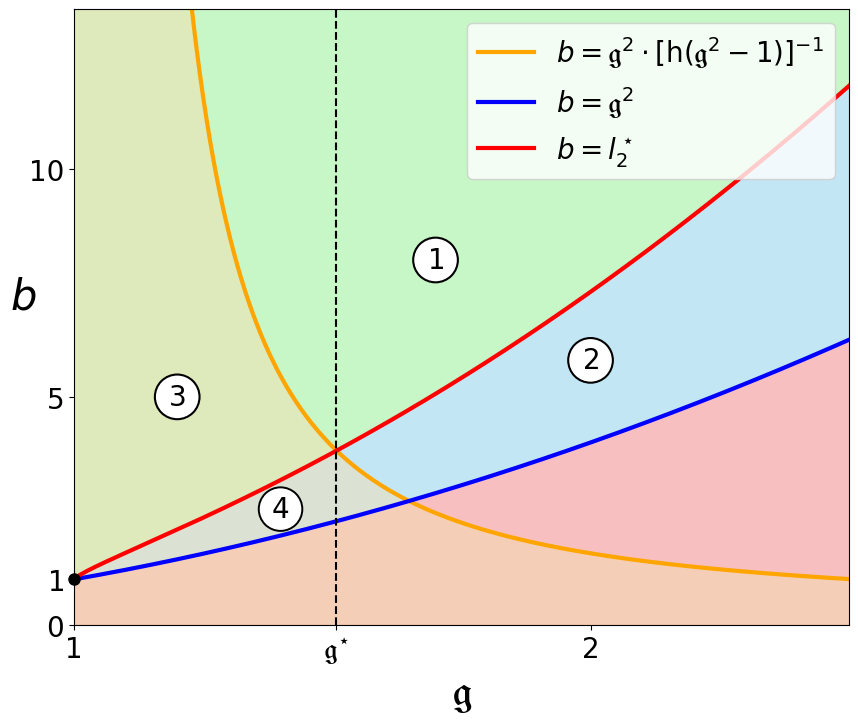}}  \\
\small (b) $l^{\star}_{2}$ is defined in \eqref{eqn:left_threshold}. 
\end{minipage}%
    \caption{Thresholds for the emergence of left (a) and right (b) outlier singular values, with respect to the sparsity parameter $b$ and different ratio parameters $\qr$. Here, $\qr^{\star}$ is the solution to the equation $l_2^\star=\qr^2\cdot [\benrate(\qr^{2} - 1)]^{-1}$, approximately $\qr^\star\approx 1.5084747$.}
    \label{fig:threshold}
\end{figure}

\begin{remark} Figure \ref{fig:threshold}(a) divides the parameter space $(b, \qr)$ into three slices. In region \circromlower{1},
$b>r_2^{\star}$ and there are no outliers on the right. In region \circromlower{2}, $r_1^{\star}<b<r_2^{\star}$; there exist outliers, and they are in correspondence with high-degree vertices, all of which are in $\gV_2$. Finally, in region \circromlower{3}, there are outliers in correspondence with high-degree vertices both in $\gV_1$ and in $\gV_2$.  We give a precise description of the correspondence in \Cref{thm:right_edge_behavior}, and an estimation of the number of outliers in Theorem \ref{thm:number_outliers}.
\end{remark}

To describe the left outliers, we need to introduce two additional assumptions. One (Assumption \ref{ass:critical_connectivity}) is natural; in the absence of graph connectivity, there will almost certainly be outliers to the left of the bulk (and zero singular values), generated by small components (respectively, isolated vertices). The second one (Assumption \ref{ass:Ihara_Bass}) is a technical one induced by our method. We explain more below. 

\begin{assumption}[Above critical connectivity]\label{ass:critical_connectivity}
We assume that $b > \qr^2$.
\end{assumption}
According to \Cref{lem:connectivity_bipartite}, Assumption~\ref{ass:critical_connectivity} ensures that the sampled bipartite graph is connected.

\begin{assumption}[Positive definiteness of $\rmI - d^{-1}\rmD^{(1)}$]\label{ass:Ihara_Bass}
Define the function $\benrate(u) \coloneqq (1+u)\log(1 + u) - u$. We assume
\begin{align}
    b > \qr^2\cdot [\benrate(\qr^{2} - 1)]^{-1}.
\end{align}
\end{assumption}
Assumption~\ref{ass:Ihara_Bass} guarantees $\rmI - d^{-1}\rmD^{(1)}$ is positive definite, where $\rmD^{(1)}$ is defined by~\eqref{eqn:degree_matrix}; see \Cref{lem:invertiblility_D1} for more details. 
We now present our second main result under  Assumptions \ref{ass:critical_connectivity} and \ref{ass:Ihara_Bass}.

\begin{theorem}[Emergence of left outliers]\label{thm:left_outlier_emergence}
Recall $l_2^\star$ in \eqref{eqn:left_threshold}. Under Assumptions \ref{ass:critical_connectivity} and \ref{ass:Ihara_Bass}, there exists $\nu >0$ such that the following holds with probability at least $1 - N^{-\nu}$ for sufficiently large $N$:
\begin{enumerate}
    \item  (No outlier) When $b> \max\{l_{2}^{\star}, \qr^2[\benrate(\qr^{2} - 1)]^{-1}\}$, \[\lambda_m(\scA) =\qr-\qr^{-1} - o(1).\]
    \item  (Emergence) When $\max\{\qr^{2}, \qr^2[\benrate(\qr^{2} - 1)]^{-1}\} < b < l_{2}^{\star}$, there exists a $\nu$-dependent constant $\epsilon>0$ such that
    \[\lambda_m(\scA) \leq \qr-\qr^{-1} -\epsilon.
    \]
    \end{enumerate}
\end{theorem}

\begin{remark} Similarly, Figure \ref{fig:threshold}(b) divides the parameter space $(b, \qr)$ into six slices. In the unnumbered regions below the curve $b=\qr^2$, we can certify that there are outliers coming from small components and low-degree vertices in the giant component. In region \circnum{1} where $b>l_2^{\star}$, there are no outliers on the left. In region \circnum{2}, there exist outliers, and they are in one-to-one correspondence with low-degree vertices, all of which are in $\gV_2$; an estimation of their number is given in Theorem \ref{thm:number_outliers}. For regions \circnum{3} and \circnum{4}, due to limitations of our method, we cannot eliminate the possibility of other mechanisms that induce outliers. 
\end{remark}

To illustrate the emergence of right and left outliers, we present a numerical experiment in Figure~\ref{fig:bipartite_outliers} below. 
\begin{figure}[h]
     \centering
     \begin{subfigure}{0.32\textwidth}
         \centering
    \includegraphics[height = 3cm, width=\textwidth]{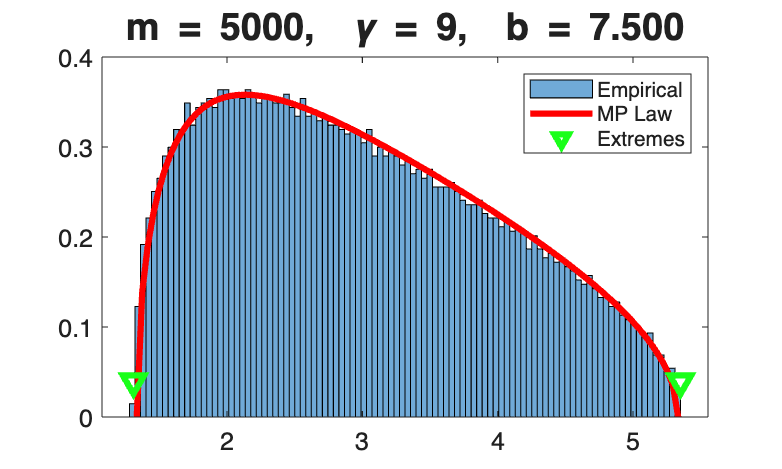}
    \subcaption{$b = 7.500$: no outlier both side.}
     \end{subfigure}
     \begin{subfigure}{0.32\textwidth}
         \centering
    \includegraphics[height = 3cm, width=\textwidth]{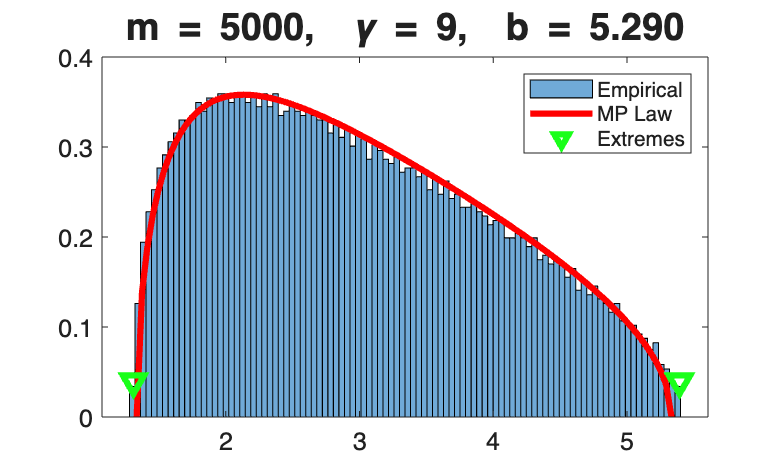}
    \subcaption{$b = 5.290$: right outliers appear.}
     \end{subfigure}
     \begin{subfigure}{0.32\textwidth}
         \centering
    \includegraphics[height = 3cm, width=\textwidth]{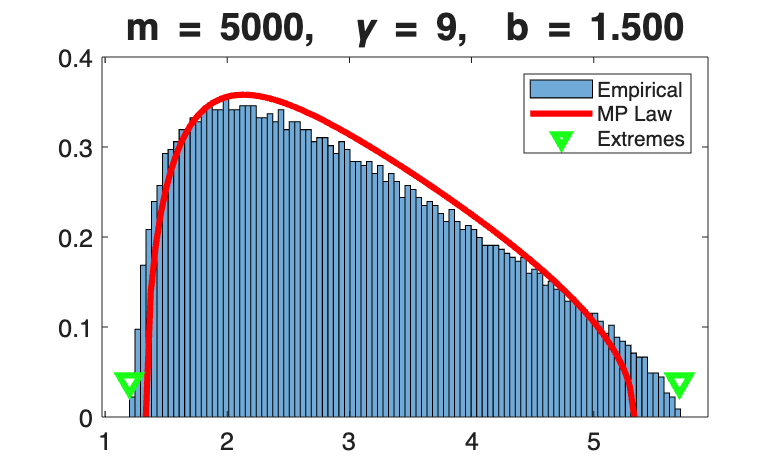}
    \subcaption{$b = 1.500$: left outliers appear.}
     \end{subfigure}
     \caption{\small{Emergence of outliers when $\ratio = 9$, where $\qr =\sqrt{3}$, $r_{2}^{\star} \approx 6.634$, $l_{2}^{\star} \approx 5.289$, $r_{1}^{\star} \approx 1.179$.}}
    \label{fig:bipartite_outliers}
\end{figure}

Note that $r^{\star}_2 \geq l^{\star}_2$. As $b$ decreases, the right outliers induced by the large-degree vertices of $\gV_{2}$ appear first; the left outliers induced by low-degree vertices of $\gV_{2}$ appear subsequently. This indicates that the left edge is more robust in this sparse random matrix model, which parallels a similar phenomenon observed in \cite{tikhomirov2015limit,tikhomirov2016smallest,koltchinskii2015bounding} for heavy-tailed random matrices.

\subsection{Quantitative characterization of outlier locations}
Before presenting the results, we introduce the necessary notations first. For vertex $x\in \gV$, let $\rD_{x} \coloneqq \sum_{y\in \gV}\ermA_{xy}$ denote its degree. The degree matrix $\rmD$ is a diagonal matrix with the degrees, which admits the following block structure:
\begin{align}\label{eqn:degree_matrix}
    \rmD \coloneqq \begin{bmatrix}
    \rmD^{(1)} & \bzero\\
    \bzero & \rmD^{(2)}
    \end{bmatrix}, \textnormal{ where } \rmD^{(1)}\coloneqq \text{diag} \{\rD_x\}_{x\in\gV_1} \,\,\textnormal{and}\,\, \rmD^{(2)}\coloneqq \text{diag} \{\rD_x\}_{x\in\gV_2}.
\end{align}
Furthermore, with $d = p\sqrt{mn}$ in \Cref{ass:proportional_regime}, we define \emph{normalized degree} as
\begin{align}
    \alpha_{x} \coloneqq \rD_{x}/ d, \label{eqn:normalized_degree}
\end{align}
Let $d_1$ and $d_2$ denote the expected degrees of vertices in $\gV_1$, respectively, in $\gV_2$. Obviously, $d_1 = m \cdot d/\sqrt{mn} = d/\qr^{2}$ and $d_2 = d\qr^{2}$. We also have the following relation:
\begin{align}
    d=\sqrt{d_1d_2}.
\end{align}
Define the deterministic map
\begin{align}
    \Lambda_{\qr}(t)\coloneqq\sqrt{t+\qr^{-2}+\frac{1}{t-\qr^{2}}}.\label{eqn:Lambda_qr}
\end{align}
Note that $\Lambda_{\qr}(t)$ is non-decreasing for $t \in [0, \qr^{2} - 1) \cup [\qr^{2}+1, \infty)$. Meanwhile, $0 \leq \Lambda_{\qr}(t) \leq \qr-\qr^{-1}$ when $t\in [0, \qr^{2} - 1)$ and $ \Lambda_{\qr}(t)\geq \qr + \qr^{-1}$ when $t \in [\qr^{2}+1, \infty)$.
Similarly, the map
\begin{align}
    \Lambda_{\qr^{-1}}(t)\coloneqq\sqrt{t+\qr^{2}+\frac{1}{t-\qr^{-2}}}.\label{eqn:Lambda_qr_inverse}
\end{align}
is non-decreasing and $\Lambda_{\qr^{-1}}(t)\geq \qr + \qr^{-1}$ for $t \in [\qr^{-2}+1, \infty)$.

We introduce the following \emph{error control} parameter
\begin{align}
    & \error \coloneqq \sqrt{\frac{\log(d)}{d}}\label{eqn:error_parameters}
\end{align}
and define the following two sets of vertices via $\alpha_{x}$ in \eqref{eqn:normalized_degree}
\begin{subequations}
    \begin{align}
            \sR_2 \coloneqq &\, \{x\in \gV_2| \Lambda_{\qr}(\alpha_x) \geq \qr + \qr^{-1} + \error^{1/4} \},\label{eqn:R_2}\\
        \sR_1 \coloneqq &\, \{x\in \gV_1| \Lambda_{\qr^{-1}}(\alpha_x) \geq \qr + \qr^{-1} + \error^{1/4} \}.\label{eqn:R_1}
    \end{align}
\end{subequations}
Denote $\Lambda^{\sR_1}:=\{\Lambda_{\qr}(\alpha_{x}), x \in \sR_1\}$  and $\Lambda^{\sR_2}:=\{\Lambda_{\qr^{-1}}(\alpha_{x}), x \in \sR_2\}$. We arrange all the elements in $\Lambda^{\sR_1} \cup \Lambda^{\sR_2}$ in a non-increasing order, denoted by
\begin{align}
    \Lambda_1 \geq \Lambda_2\geq \cdots \geq \Lambda_{|\sR_1| + |\sR_2|}.
\end{align}

\Cref{thm:right_edge_behavior} shows that \emph{all} the right outlier eigenvalues are $\Omega(\error)$ away from their predicted locations $\Lambda_1, \ldots, \Lambda_{|\sR_1|+|\sR_2|}$.

\begin{theorem}[Right edge behavior]\label{thm:right_edge_behavior}
    Under \Cref{ass:proportional_regime}, for any $\nu > 0$, there is a  constant $\const >0$ such that the following holds with probability at least $1 - N^{-\nu}$ for sufficiently large $N$:
\begin{enumerate}
    \item For all $1\leq j\leq |\sR_1| + |\sR_2|$, with $\error$ defined in \eqref{eqn:error_parameters}, we have
      \begin{align}
          |\lambda_j(\scA) - \Lambda_{j}| \leq \const \error.
      \end{align}
      \item For all $|\sR_1| + |\sR_2| + 1\leq j\leq m$,
      \begin{align}
          \lambda_j(\scA) \leq \qr+\qr^{-1} + \error^{1/4}.
      \end{align}
\end{enumerate}
\end{theorem}

We now turn to the left edge.  Consider the following set of \emph{low-degree} vertices.
\begin{align}\label{eqn:L2_set}
    \sL_2 \coloneqq &\, \{x\in \gV_2| \Lambda_{\qr}(\alpha_x) \leq \qr - \qr^{-1} - \error^{1/4} \}.
\end{align}

Let $\pi: [m]\to [m]$ be the permutation such that $\alpha_{\pi(1)}\geq \alpha_{\pi(2)}\geq \cdots \geq \alpha_{\pi(m)}$. 

\begin{theorem}[Left edge behavior]\label{thm:left_edge_behavior}
Suppose that Assumptions \ref{ass:critical_connectivity} and \ref{ass:Ihara_Bass} are true. There exist constants $\nu > 0$ and $\const >0$ such that the following holds with probability at least $1 - N^{-\nu}$ for sufficiently large $N$:
\begin{enumerate}
    \item For all $m - |\sL_2| + 1 \leq j \leq m$, with $\error$ defined in \eqref{eqn:error_parameters}, 
      \begin{align}
          |\lambda_j(\scA) - \Lambda_{\qr}(\alpha_{\pi(j)})| \leq \const \error.
      \end{align}
    \item For all $1 \leq j \leq m - |\sL_2|$,
      \begin{align}
          \lambda_j(\scA) \geq \qr - \qr^{-1} - \error^{1/4}.
      \end{align}
\end{enumerate}
\end{theorem}
\begin{remark}
The existence of left outliers can be established only when $\qr > 1$, since the left edge of the \emph{MP} law is zero when $\qr = 1$.
\end{remark}
Similarly, \Cref{thm:left_edge_behavior} shows that the locations of \emph{all} of the left outlier eigenvalues are most $\Omega(\error)$ from their predicted locations $\{\Lambda_{\qr}(\alpha_x)| x\in \gV_2\}$.

\subsection{Estimating the number of outliers}\label{sec:phase_diagram}
Based Theorems \ref{thm:right_edge_behavior}, \ref{thm:left_edge_behavior}, left and right outliers, which are at least $\Omega(\error^{1/4})$ away from the bulk, are induced by vertices in $\sR_1\cup \sR_2$ and $\sL_2$, respectively. Therefore, we can characterize the phase transition behavior for the source of outliers depending on when $\sR_1,\sR_2,$ and $\sL_2$ become nonempty as the sparsity parameter $b$ changes. This is made precise in the following theorem using the threshold functions $r^{\star}_1, r^{\star}_2$, and $ l^{\star}_2$ defined in \eqref{eqn:right_threshold_1}, \eqref{eqn:right_threshold_2}, and \eqref{eqn:left_threshold}.

\begin{theorem}[Number of outliers] \label{thm:number_outliers}
Under \Cref{ass:proportional_regime}, the following holds with high probability:
\begin{enumerate}
    \item   When $b > r^{\star}_2$, $\sR_2=\emptyset$, and when $b< r^{\star}_2$, $|\sR_2|= m\cdot N^{-b/r_2^\star+o(1)}$.
    \item When $b > r^{\star}_1$, $\sR_1=\emptyset$, and when $b< r^{\star}_1$, $|\sR_1| = n\cdot N^{- b/r_1^\star+o(1)} $.
    \item Under additional Assumptions~\ref{ass:critical_connectivity} and \ref{ass:Ihara_Bass},  when $b>\max\{l_{2}^{\star},\qr^2[\benrate(\qr^{2} - 1)]^{-1}\}$,
    $\sL_2 =\emptyset$ and when $\max\{\qr^{2}, \qr^2[\benrate(\qr^{2} - 1)]^{-1}\} < b < l_{2}^{\star}$, $|\sL_2|= m\cdot N^{-b/l_2^\star+o(1)}$.
\end{enumerate}
\end{theorem}

\subsection{Extension to general sparse random rectangular matrices}
Let $\widetilde{\rmA}\in\R^{n\times m}$ be the \emph{bi-adjacency} matrix of $\gG$ sampled from $\gG(n, m,d/\sqrt{mn})$. Let $\widetilde{\rmW} \in \C^{n \times m}$ be a rectangular random matrix whose entries $\widetilde{\ermW}_{xy}$ are independent real or complex-valued random variables that satisfy the following requirements.
\begin{assumption}\label{ass:homogeneous_variance}
    $\E \widetilde{\ermW}_{xy} = 0$, $\E\widetilde{\ermW}_{xy}^{2} = 1$ and $|\widetilde{\ermW}_{xy}| \leq \gK$ almost surely for some constant $\gK > 0$.
\end{assumption}
Let $\widetilde{\rmM} \coloneqq \widetilde{\rmW} \circ \widetilde{\rmA}$ denote the Hadamard product between $\widetilde{\rmW}$ and $\widetilde{\rmA}$. Consider
\begin{align}\label{eqn:sparse_wigner}
    \scA \coloneqq \frac{1}{\sqrt{d}} \rmM\,, \quad \textnormal{where} \quad \rmM = \begin{bmatrix}
        \bzero & \widetilde{\rmM}\\
        \widetilde{\rmM}^{*} & \bzero
    \end{bmatrix}.
\end{align}
Similarly to \eqref{eqn:normalized_degree}, we define the normalized degree $\alpha_{x}$ as
\begin{align}
    \alpha_{x} = \frac{1}{d}\sum_{y\in[N]} |\ermM_{xy}|^{2} = \frac{1}{d}\sum_{y\in[N]} |\ermW_{xy}|^{2} \cdot \ermA_{xy}.\label{eqn:new_degree}
\end{align}
We define $\sR_1,~ \sR_2,~\sL_2$ analogously using the normalized degree above. \Cref{thm:sparse_rectangular} generalizes the results in Theorems \ref{thm:right_edge_behavior} and \ref{thm:left_edge_behavior}.

\begin{theorem}\label{thm:sparse_rectangular}
Under Assumptions \ref{ass:proportional_regime} and \ref{ass:homogeneous_variance}, consider $\scA$ in \eqref{eqn:sparse_wigner} and $\alpha_x$ in \eqref{eqn:new_degree}. 
\begin{enumerate}
\item 
 For any $\nu > 0$, there is a  constant $\const >0$ such that the following holds with probability at least $1 - N^{-\nu}$ for sufficiently large $N$:
\begin{enumerate}
    \item For all $1\leq l\leq |\sR_1| + |\sR_2|$, 
    we have
$ |\lambda_l(\scA) - \Lambda_{l}| \leq \const \error.$
      \item For all $|\sR_1| + |\sR_2| + 1\leq l\leq m$, $          \lambda_l(\scA) \leq \qr+\qr^{-1} + \error^{1/4}.$
\end{enumerate}
\item Under additional Assumptions \ref{ass:critical_connectivity} and \ref{ass:Ihara_Bass}, there exist constants $\nu > 0$ and  $\const >0$ such that  with probability $1-N^{-\nu}$ for sufficiently large $N$:
\begin{enumerate}
    \item For all $m - |\sL_2| + 1 \leq l \leq m$,  $ |\lambda_l(\scA) - \Lambda_{\qr}(\alpha_{\pi(l)})| \leq \const \error.$
    \item For all $1 \leq l\leq m - |\sL_2|$, $ \lambda_l(\scA) \geq \qr - \qr^{-1} - \error^{1/4}.$
\end{enumerate}
\end{enumerate}
\end{theorem}

The matrix $\rmM$ in \eqref{eqn:sparse_wigner} corresponds to the adjacency matrix of a sparse weighted random graph, a model that has been extensively studied in the large deviation principle; see \cite{ganguly2022large,ganguly2024spectral,augeri2025large,augeri2023large}.  It is possible that the boundedness condition in \Cref{ass:homogeneous_variance} can be removed; see, for example, the local law for unbounded sparse random matrices in \cite{he2024sparse}. Moreover, when $\qr = 1$, a lower bound for the least singular value of $\rmM$, where $\rmM$ has i.i.d. mean $0$ and variance $1$ sub-gaussian entries, was recently established in \cite{sah2025spielman}.

\subsection{Related work}\label{sec:literature}

For the classical, dense Wishart/Wigner matrix models, investigating the existence of spectra outliers has generally been approached through different methods, depending on whether they were right-edge or, in the case of Wishart ensembles, left-edge. For a comprehensive list of references, we refer the reader to \cite{dumitriu2024extreme}.

The classical Bai-Yin theorem \cite{yin1988limit,bai1993limit} guarantees that for dense Wishart matrices whose entries' distributions have 4th moments, almost surely there are no outliers in the spectrum, neither to the left nor to the right. Since then, this theorem has been extended for dense matrices to include heavy-tailed random matrices \cite{adamczak2010quantitative,tikhomirov2015limit,tikhomirov2016smallest,koltchinskii2015bounding,tikhomirov2018sample,bao2024phase, bao2024phasebs, han2024smallest} and anisotropic sample covariance matrices \cite{koltchinskii2017concentration,zhivotovskiy2024dimension,bandeira2024matrix}.

For sparse ensembles, the level of sparsity turns out to determine both the approach and the strength of the results. Starting with the seminal paper \cite{furedi1981eigenvalues} and subsequently moving in a graph-related direction, through \cite{vu2007spectral}, which let go of the i.i.d. entry assumption and reached down to $np = \mathrm{polylog}(n)$ sparsity, the regimes of investigation split into supercritical ($np = \omega(\log n)$), critical ($np \asymp \log n$), and subcritical ($np= o(\log n)$). 
In the case of supercritical ensembles, a sharp analysis of the second largest eigenvalue of sparse random graphs and Wigner matrices has been provided in \cite{latala2018dimension,bandeira2016sharp, bandeira2021matrix,brailovskaya2024extremal, bandeira2024matrix} via a broad universality principle, and in \cite{benaych2020spectral} via graph methodology, including the use of the non-backtracking operator.
Furthermore, many concentration inequalities for the spectral norm of sparse random matrices have been proved down to $np=\Omega(\log n)$, up to a constant factor \cite{feige2005spectral,le2017concentration,benaych2020spectral}.  
Additionally, the extreme singular values of a random matrix with a deterministic sparsity pattern were investigated in \cite{bandeira2016sharp,altschuler2024spectral}.

In the critical and subcritical cases,  \cite{khorunzhy2001sparse} identified the phase transition of the order of the spectral norm as happening when $np\asymp \log n$. For the first time, the largest eigenvalue of sparse matrices 
was investigated in \cite{krivelevich2003largest} based on the highest degree vertex. The subcritical behavior of the largest eigenvalues was described up to a constant in  
\cite{benaych2019largest} and more recently in \cite{alt2023localized}. At criticality, for homogeneous ensembles, the second largest eigenvalue (and therefore the existence of outliers in the spectrum) was characterized by \cite{tikhomirov2021outliers,alt2021extremal,alt2021delocalization}. Finally, in subcritical and very sparse cases, the locations of the largest eigenvalues were proved in \cite{hiesmayr2023spectral}.

Additionally, in the setting of $d$‑regular sparse random graphs, the behavior of the second eigenvalue has been a central topic in the theory of sparse random matrices for decades \cite{friedman1989second,broder1998optimal,friedman1995second}, and a variety of techniques have now led to a comprehensive understanding of the spectral gap across all sparsity regimes \cite{bordenave2020new,cook2018size,tikhomirov2019spectral,brailovskaya2022universality,sarid2023spectral,he2024spectral,chen2024new}. Most recently, for random $d$ regular graphs with constant $d$, \cite{huang2024ramanujan} proved that the fluctuation of the second eigenvalue follows the Tracy–Widom law.  For random bipartite biregular graphs, the spectral gap has been analyzed in \cite{brito2022spectral,zhu2022second}. 

For the bipartite sparse random graph, its bi-adjacency matrix is a sparse random rectangular matrix, whose extreme singular values have also been investigated recently, with applications to spectral clustering and community detection \cite{florescu2016spectral, zhou2019analysis}.  A two-sided spectral norm bound for sparse random rectangular matrices was obtained in \cite{gotze2022largest} when $np=\mathrm{polylog}(n)$ up to a constant. For bipartite Erd\H{o}s-R\'enyi graphs $\gG (m,n,p)$, \cite{dumitriu2024extreme} showed that no outliers exist outside the MP law when $np=\omega(\log n)$, and their results cover inhomogeneous random bipartite graphs. The general results in \cite{brailovskaya2022universality}  imply no outliers for a more general inhomogeneous random bipartite graph model when $np=\omega(\log^4 n)$. Very recently, \cite{tropp2025comparison} studied a broad class of random positive semidefinite matrices and, as a corollary, obtained a nontrivial lower bound (up to a constant) on the smallest singular value of sparse rectangular matrices in the regime $np=\Omega(\log n)$. 

\subsection{Discussion of the technical condition}
The current technical assumption \ref{ass:Ihara_Bass} is due to the condition in Bennett's inequality (\Cref{lem:Bennett}). We believe that this is a limitation of the proof technique, not an intrinsic requirement for the sparse random bipartite graphs. We leave this investigation for future work.

\subsection{Notation}\label{sec:notation}
We denote the set of positive numbers by $\N = \{1, 2, \ldots\}$ and define $\N_0 = \N \cup \{0\}$. For any $n\in \N$, $[n] \coloneqq \{1, \ldots, n\}$, and $\llbracket n \rrbracket \coloneqq \{0, 1, \ldots, n\}$. 
For any two sequences of numbers $\{a_N\}$, $\{b_N\}$, denote $a_N = O(b_N)$ or $a_N \lesssim b_N$ (resp. $a_N = \Omega(b_N)$ or $a_N \gtrsim b_N$) if there exist some constants $\const$ and $N_0$ such that $a_N \leq \const\, b_N$ (resp. $a_N \geq \const\, b_N$) for all $N \geq n_0$. Denote $a_N = \Theta(b_N)$ or $a_N \asymp b_N$ if both $a_N \lesssim b_N$ and $a_N \gtrsim b_N$. Denote $a_N = o(b_N)$ or $a_N \ll b_N$ (resp. $a_N = \omega(b_N)$ or $a_N \gg b_N$) if for any $\epsilon >0$, there exists $N_0\in \N$ s.t. $a_N < \epsilon b_N$ (resp. $a_N > \epsilon b_N$) for all $N\geq N_0$.

\subsubsection*{Matrix notation}
The lowercase letters (e.g., $a, b$), lowercase boldface letters (e.g., $\rvx, \rvy$), and uppercase boldface letters (e.g., $\rvx, \rmA$) are used to denote scalars, vectors and matrices, respectively. For vector $\rvx\in \R^N$, let $\|\rvx\|_{p} = (\sum_{i=1}^{N}|\ervx_i|^p)^{1/p}$ denote the $\ell_p$ norm of $\rvx$ for $p \geq 1$. For matrix $\rmX \in \R^{m \times n}$, let $\ermX_{ij}$, $\rmX_{i:}$ and $\rmX_{:j}$ denote its $(i, j)$th entry, $i$-th row and $j$-th column respectively. Let $\|\rmX\| $ and $\|\rmX\|_{\frob}$ denote the operator and Frobenius norm respectively. 
For any subset $\setS \subset [N]$, define the matrix $\rmA_{\setS} \in \R^{|\setS| \times |\setS|}$ and the family $\rmA_{(\setS)}$ through
\begin{align}
    \rmA_{\setS} \coloneqq \ermA_{ij} \cdot \indi{i\in\setS} \cdot \indi{j\in \setS}\,, \quad \rmA_{(\setS)} \coloneqq \ermA_{ij} \cdot \big( \indi{i\in\setS} + \indi{j\in \setS} \big). \label{eqn:Asubsets}
\end{align}
We write $\rmA_{(x)}$ instead of $\rmA_{(\{x\})}$ when $\layer = \{x\}\subset [N]$. Let $\rmA\in \R^{N \times N}$ be a Hermitian matrix. Its eigenvalues are sorted in decreasing order and denoted by $\lambda_1(\rmA)\geq \ldots \geq \lambda_N(\rmA)$. Moreover, for matrices $\rmA, \rmB \in \R^{N \times N}$, we write $\rmA \succ \rmB$ ($\rmA \succeq \rmB$) if the matrix $\rmA - \rmB$ is positive (semi-) definite, i.e., $\lambda(\rmA - \rmB) >0$ ($\geq 0$).

Let $\ones_{j} \in \{0, 1\}^{N}$ denote the column vector, where the $j$-th entry is the only non-zero one. For any subset $\layer \subset [N]$, define $\ones_{\layer} \coloneqq \sum_{x\in \layer} \ones_x$. For convenience, we denote by $\ones_{N} \coloneqq \sum_{v\in \gV} \ones_v$ and its normalized version $\rve_{N} \coloneqq N^{-1/2}\ones_N$.

\subsubsection*{Graph notation} 
For graph $\gG = (\gV, \gE)$, let $\gV$, $\gE$ denote the set of vertices and edges, respectively. For two graphs $\gH$ and $\gG$, we write $\gH \subset \gG$ if $\gV(\gH) \subset \gV(\gG)$ and $\gE(\gH) \subset \gE(\gG)$. Let $\gG$ be the graph on the set of vertices $[N]$ and $\rmA\in \R^{N \times N}$ be the adjacency matrix of $\gG$. For any $\gV \subset [N]$, let $\gG|_{\gV}$ denote the induced subgraph of $\gG$ on the set of vertices $\gV$, and $\rmA_{\gV} \in \R^{|\gV|\times |\gV|}$ denote the adjacency matrix of $\gG|_{\gV}$. The degree of vertex $x$ on $\gG$ is denoted by $\rD_{x}^{\gG} \coloneqq \sum_{y\in [N]} \ermA_{xy}$. The distance between two vertices $x$ and $y$ on $\gG$ (the shortest path from $x$ to $y$) is denoted by
\begin{align}
    d^{\gG}(u, v) = \min\{k\in \N_0|(\rmA^{k})_{uv} \neq 0\}. \label{eqn:graph_distance}
\end{align}
For vertex $v\in \gG$, we define the $j$-layer $\layer_{j}^{\gG}(v)$ and the $r$-ball $\ball_{r}^{\gG}(v)$ centered at vertex $v$ by
\begin{align}\label{eqn:LiBi(x)}
    \layer_{j}^{\gG}(v) \coloneqq \{u\in \gV: d^{\gG}(u, v)  = j\}, \quad \ball_{r}^{\gG}(v) \coloneqq \cup_{j=0}^{0} \,\, \layer_{k}^{\gG}(v).
\end{align}
By convention, we denote $\layer_0^{\gG}(x) = \{x\}$ and $\layer_{-1}^{\gG}(x) = \emptyset$.

In the remainder of this paper, let $\gG$ denote the random bipartite Erd\H{o}s-R\'{e}nyi graph sampled from $\gG(n, m,d/\sqrt{mn})$ with the adjacency matrix $\rmA \in \R^{N \times N}$, where $N = n+ m$. For convenience, we denote $\rD_v$, $d(u, v)$, $\layer_{j}(v)$, $\ball_{j}(v)$ instead of $\rD_{v}^{\gG}$, $d^{\gG}(u, v)$, $\layer_{j}^{\gG}(v)$, $\ball_{j}^{\gG}(v)$.

\subsubsection*{Probability notation} 
Throughout the paper, the term \emph{very high probability} refers to the following.
\begin{definition}[Very high probability]
Let $\Xi \equiv \Xi_{N, \nu}$ be a family of events parametrized by $N\in \N$ and $\nu >0$. We say that $\Xi$ holds with very high probability if for every $\nu >0$, there exists some $\nu$-dependent constant $\const_{\nu} >0$ such that for all $N\in \N$,
\begin{align}
    \P(\Xi_{N, \nu}) \geq 1 - \const_{\nu} N^{-\nu}.
\end{align}
Moreover, the bound $|\rX|\lesssim \rY$ with very high probability, or equivalently $\rX = O(\rY)$, refers to the fact that for every $\nu>0$, there exists some $\nu$-dependent constants $\const_{\nu}, c_{\nu}>0$ such that for all $N\in \N$
    \begin{align}
        \P(|\rX|\leq c_{\nu} \rY)\geq 1- \const_{\nu} N^{-\nu}.
    \end{align}
Here, $\rX$ and $\rY$ are allowed to depend on $N$.
\end{definition}

\subsection{Organization}
The rest of the paper is organized as follows.
Proofs of Theorems~\ref{thm:right_outlier_emergence}, \ref{thm:left_outlier_emergence} and \ref{thm:number_outliers} are presented in \Cref{sec:proof_emergence}. We summarize the proof ideas of Theorems \ref{thm:right_edge_behavior} and \ref{thm:left_edge_behavior} in \Cref{sec:outline}. Theorems \ref{thm:right_edge_behavior} (1) and \ref{thm:left_edge_behavior} (1) are proved in \Cref{sec:outlier_locations}, and Theorems \ref{thm:right_edge_behavior} (2) and \ref{thm:left_edge_behavior} (2) are proved in~\Cref{sec:bulk_boundness}. The proof of \Cref{thm:sparse_rectangular} is deferred to \Cref{sec:sparse_rectangular}.

\section{Proofs of Theorems~\ref{thm:right_outlier_emergence}, \ref{thm:left_outlier_emergence} and \ref{thm:number_outliers} } \label{sec:proof_emergence}

We first prove Theorems \ref{thm:right_outlier_emergence}, \ref{thm:left_outlier_emergence}, and \ref{thm:number_outliers}, based on Theorems~\ref{thm:right_edge_behavior} and \ref{thm:left_edge_behavior}, and the following two lemmas. The proof of Lemma~\ref{lem:approxBinomProb} and Lemma~\ref{lem:degreeProbApprox} can be found in Appendix \ref{sec:proof_emergence_lemma}. In the following, let $\rD_{x}^{(j)}$ denote the degree of the vertex $x\in \gV_{j}$ for $j = 1, 2$.

\begin{lemma}\label{lem:approxBinomProb}
Under \Cref{ass:proportional_regime}, the following holds:
\begin{enumerate}
    \item When $\alpha \geq \qr^{-2} + 1$, $\P(\rD_x^{(1)} \geq \alpha d) \asymp \P(\rD_x^{(1)} = \alpha d)$ for any $x\in \gV_1$.
    \item When $\alpha \geq \qr^{2} + 1$, $\P(\rD_x^{(2)} \geq \alpha d ) \asymp \P(\rD_x^{(2)} = \alpha d)$ for any $x\in \gV_2$.
    \item When $\alpha \leq \qr^{2} - 1$, $\P(\rD_x^{(2)} \leq \alpha d) \asymp \P(\rD_x^{(2)} = \alpha d)$ for any $x\in \gV_2$.
\end{enumerate}
\end{lemma}

\begin{lemma}\label{lem:degreeProbApprox}
    Define the following function 
    \begin{align}
        f_{\qr, d} (\alpha) = d\Big( \alpha \log(\qr^{-2}\alpha) - \alpha + \qr^{2}\Big) + \frac{1}{2}\log(2\pi \alpha d).\label{eqn:stirling_rate}
    \end{align}
    Then $\P(\rD_x^{(1)} = \alpha d) = \exp(- f_{\qr^{-1}, d}(\alpha))$ and $\P(\rD_x^{(2)} = \alpha d) = \exp(- f_{\qr, d}(\alpha))$.
\end{lemma}

\begin{proof}[Proof of \Cref{thm:right_outlier_emergence}]

Let $\alpha^{(j)}_{x}$ denote the normalized degree of $x\in \gV_{j}$ for $j = 1, 2$. Without loss of generality, arrange $\alpha^{(2)}_{1} \geq \ldots \geq \alpha^{(2)}_{m}$ in non-increasing order. We further define the counting function $\rN^{(j)}_{t} \coloneqq \sum_{x\in \gV_j} \indi{\rD_{x}^{(j)} \geq t}$.

For part (1) of \Cref{thm:right_outlier_emergence}, by taking $t=\alpha d$, we have
\begin{align}
   \E \rN^{(j)}_{\alpha d} =&\, (n\cdot \indi{j=1} + m\cdot \indi{j=2}) \cdot \P(\rD^{(j)}_x \geq \alpha d).
\end{align}
Meanwhile, according to Lemmas \ref{lem:approxBinomProb} and \ref{lem:degreeProbApprox}, the following holds when $\alpha \geq \qr^{2} + 1$, 
\begin{align}
    \E \rN^{(2)}_{\alpha d} \overset{\textnormal{(\Cref{lem:approxBinomProb})}}{\asymp} m \cdot \P(\rD_x = \alpha d)
    \overset{\textnormal{(\Cref{lem:degreeProbApprox})}}{=} &\, m  \cdot \exp(- f_{\qr, d}(\alpha)).
\end{align}
Then by Markov, for any $l\in [m]$ and $\alpha \geq \qr^{2} + 1$, the probability of the event $\{\alpha^{(2)}_{l} \geq \alpha\}$ is bounded by
\begin{align}
    \P(\alpha^{(2)}_{l} \geq \alpha) = &\, \P(\rN^{(2)}_{\alpha d} \geq l) \leq \frac{\E \rN^{(2)}_{\alpha d}}{l} \\
    \lesssim &\, \frac{m}{l} \exp\Big(- f_{\qr, d}(\alpha) \Big) = \exp\Big( \log(m/l) - f_{\qr, d}(\alpha) \Big).\label{eqn:lth_prob_upper}
\end{align}
We choose $l = 1$ and $\alpha = \qr^{2} + 1$. Since $b > r_{2}^{\star}$ defined in \eqref{eqn:right_threshold_2}, we have $\P(\alpha^{(2)}_{1} \geq \qr^{2}+ 1) \leq  N^{-\nu}$ for some constant $\nu>0$, consequently $\sR_{2} = \emptyset$.
Then $\lambda_1(\scA) \leq \qr+\qr^{-1} +o(1)$ follows from part (2) of \Cref{thm:right_edge_behavior}. 

Part (2) of \Cref{thm:right_outlier_emergence} is proved by a standard second-moment argument. For any $l$ with $1\leq l \leq \E \rN^{(2)}_{\alpha d}/2$, according to Chebyshev's inequality, the following holds for some constant $\const > 0$,
\begin{align}
    \P(\alpha^{(2)}_{l} \geq \alpha) &\, = \P( \rN^{(2)}_{\alpha d} - \E \rN^{(2)}_{\alpha d} \geq l - \E \rN^{(2)}_{\alpha d} ) \geq 1 - \P( |\rN^{(2)}_{\alpha d} - \E \rN^{(2)}_{\alpha d}| \geq |l - \E \rN^{(2)}_{\alpha d}| ) \\
    &\, \geq 1 - \frac{\Var(\rN^{(2)}_{\alpha d})}{(l - \E \rN^{(2)}_{\alpha d})^2} \geq 1 - \frac{4\Var(\rN^{(2)}_{\alpha d})}{(\E \rN^{(2)}_{\alpha d})^2} \geq 1 -\frac{\const}{m\cdot \exp(- f_{\qr, d}(\alpha))}\,,\label{eqn:lth_prob_lower}
\end{align}
where the last inequality is due to the fact $ \Var(\rN^{(2)}_{\alpha d})\leq \const (\E \rN^{(2)}_{\alpha d})^2$ (see, for example, \cite[Lemma 3.11]{bollobas1998random}). When $b < r_{2}^{*}$ and we choose $l = 1$, $\alpha = \qr^{2} + 1$, there exists some $\epsilon>0$ such that $\P(\alpha^{(2)}_{1} \geq \qr^{2} + 1 + \epsilon)\geq 1-N^{-\nu}$, hence $\sR_{2} \neq \emptyset$. By part (1) of \Cref{thm:right_edge_behavior}, the location of $\lambda_1(\scA)$ is characterized by $\Lambda_{\qr}(\alpha^{(2)}_{1})$ in \eqref{eqn:Lambda_qr} up to some tiny correction $\error = o(1)$. Since $\Lambda_{\qr}(t)$ in \eqref{eqn:Lambda_qr} increases monotonically when $t\geq \qr^{2} + 1$, we find $\Lambda_{\qr}(\alpha^{(2)}_{1})\geq \qr+\qr^{-1} +\widetilde{\epsilon}$ for a constant $\widetilde{\epsilon}>0$. Consequently,  
 $\lambda_1(\scA)=\Lambda_{\qr}(\alpha^{(2)}_{1}) + \error > \qr + \qr^{-1} +\widetilde{\epsilon}/2$ for sufficiently large $N$, which completes the proof.
\end{proof}
The proof of \Cref{thm:left_outlier_emergence} follows the same strategy as above. However, we use part (3) of \Cref{lem:approxBinomProb} instead. We sketch the main differences below.
\begin{proof}[Proof of \Cref{thm:left_outlier_emergence}]
For (1), an argument similar to the analysis above shows that $\alpha^{(2)}_{m} \geq \qr^{2} - 1$ with probability at least $1 - N^{-\nu}$, hence $\sL_{2} = \emptyset$. We conclude the proof by using part (2) of \Cref{thm:left_edge_behavior}.

For (2), one can show that $\P(\alpha^{(2)}_{m} \leq \qr^{2} - 1 - \widetilde{\epsilon})\geq 1 - N^{-\nu}$ for some constant $\widetilde{\epsilon}>0$, hence $\sL_{2} \neq \emptyset$ . By monotonicity, $\Lambda_{\qr}(\alpha^{(2)}_{m}) \leq \qr -\qr^{-1} -\epsilon$ for some constant $\epsilon>0$.  According to part (1) of \Cref{thm:left_edge_behavior}, $\lambda_m(\scA) \leq \qr -\qr^{-1}-\epsilon/2$ for sufficiently large $N$, which completes the proof.
\end{proof}

\begin{proof}[Proof of \Cref{thm:number_outliers}]
    For (1), the proof for $b> r_{2}^{\star}$ follows directly according to \Cref{thm:right_outlier_emergence} (1). When $b< r_{2}^{\star}$, we perform the calculations in \eqref{eqn:lth_prob_upper} and \eqref{eqn:lth_prob_lower} again. In particular, for fixed $b$, we have
    \begin{align}
        \log(m/l) = \big( b / r_{2}^{\star} + o(1)\big) \cdot \log(N).
    \end{align}
Here, with high probability, there are at least $l-1$ vertices in $\gV_{2}$ with normalized degree no smaller than $\qr^{2} + 1$. Then we write $l = m \cdot N^{-b/r_{2}^{\star} + o(1)}$, which completes the proof.

The proofs of (2) and (3) follow by similar calculations.
\end{proof}

\section{Outline for proofs of Theorems~\ref{thm:right_edge_behavior} and \ref{thm:left_edge_behavior}}\label{sec:outline}

The proofs of Theorems~\ref{thm:right_edge_behavior} and \ref{thm:left_edge_behavior} follow the three-step strategy below.
\begin{enumerate}
    \item[(1)] ``Vertex-to-outlier'' mechanism: one large (resp. small) degree vertex induces two right (resp. left) outliers, one positive and one negative.
    \item[(2)] Sufficiency: different large (resp. small) degree vertices induce distinct right (resp. left) outliers. 
    \item[(3)] Necessity: no other eigenvalues lie outside the left and right edges.
\end{enumerate}
Here, (1) and (2) together prove \Cref{thm:right_edge_behavior} (1) and \Cref{thm:left_edge_behavior} (1), while (3) establishes the proofs of \Cref{thm:right_edge_behavior} (2) and \Cref{thm:left_edge_behavior} (2). We illustrate the main idea of each step below.

\subsection{Vertex-to-outlier}\label{sec:vertex_to_outlier}
For the triple $\vtau = (\tau_{1}, \tau_{2}^{+}, \tau_{2}^{-})$ with $\tau_1 = \qr^{-2} + 1 + \error^{1/2}$, $\tau^{+}_{2} = \qr^{2} + 1 + \error^{1/2}$, and $\tau^{-}_{2} = \qr^{2} - 1 - \error^{1/2}$, we consider the following three sets of vertices.
\begin{align}
    \Vonehigh \coloneqq \{x\in \gV_{1}| \alpha_{x} \geq \tau_1 \}, \,\, \Vtwohigh \coloneqq \{x\in \gV_{2}| \alpha_{x} \geq \tau^{+}_{2} \}, \,\, \Vtwolow \coloneqq \{x\in \gV_{2}| \alpha_{x} \leq \tau^{-}_{2} \}. \label{eqn:atypical_vertices}
\end{align}
For convenience, we denote their union by
\begin{align}
    \pruneV \coloneqq \Vonehigh \cup \Vtwohigh \cup \Vtwolow.\label{eqn:pruneV}
\end{align}
To understand the mechanism, it is insightful to analyze the local structure of $\gG$ in the neighborhood of some vertex $x \in \pruneV$. Consider the tree produced by the \emph{Breadth-First-Search} algorithm and let $x$ be the root vertex. Recall that $\layer_{j}^{\gG}(x)$, defined in \eqref{eqn:LiBi(x)}, denotes the $j$-th layer on $\gG$ rooted at $x$, and $\ball_{r}^{\gG}(x)$ denotes the ball on $\gG$ within radius $r$. With the proofs deferred to \Cref{sec:approx_eigenpairs} and \Cref{sec:proofs_of_approximate_eigen_pairs}, there exists some $r_{x} \in \N$ such that
\begin{enumerate}
    \item[(a)] The radius $r_{x}$ tends to infinity as $N$ grows, as defined in \eqref{eqn:defrx}.
    \item[(b)] The ratio $|\layer_{j}^{\gG}(x)|/|\layer_{j}^{\gG}(x)|$ concentrates around $d_{1}$ or $d_{2}$ for $1 \leq j \leq r_{x}$ (see \Cref{lem:concentrationSi}).
    \item[(c)] The subgraph $\gG |_{\ball_{r_x}^{\gG}(x)}$ is a tree up to a bounded number of edges (see \Cref{lem:few_even_cycle_prob}).
\end{enumerate}
We note that this ``locally tree-like'' property was first utilized for spectral analysis of random $d$-regular graphs in \cite{dumitriu2012sparse}. Due to (b) and (c), it is natural to study the spectral properties of the following idealized deterministic tree.
\begin{definition}\label{def:biregular_tree}
Let $x$ be the root vertex of the tree $\tree$ and $\layer_{0}^{\tree}(x) = \rD_{x}$. For $1\leq j \leq \lfloor r_x/2 \rfloor$, 
\begin{enumerate}
    \item $x\in \gV_1$: vertices in $\layer_{2j - 1}^{\tree}(x)$ have $d_2$ children, while vertices in $\layer_{2j}^{\tree}(x)$ have $d_1$ children.
    \item $x\in \gV_2$: vertices in $\layer_{2j - 1}^{\tree}(x)$ have $d_1$ children, while vertices in $\layer_{2j}^{\tree}(x)$ have $d_2$ children.
\end{enumerate}
\end{definition}

Let $\rmA^{\tree}$ denote the adjacency matrix of $\tree$. In the following, we first construct an approximate eigenpair of $\rmA^{\tree}/\sqrt{d}$ up to some orthogonal transformation, which then motivates the construction of an approximate eigenvector of $\scA$.

Let $\rvs_0^{\tree}, \rvs_1^{\tree}, \rvs_2^{\tree}, \ldots, \rvs_{r_{x}}^{\tree}$ be the vectors obtained from $\ones_{x}$, $(\rmA^{\tree})\ones_x$, $(\rmA^{\tree})^2\ones_x$, $\ldots$, $(\rmA^{\tree})^{r_{x}}\ones_x$ through Gram–Schmidt orthonormalization and $\rvs_{r_{x} + 1}^{\tree}, \ldots, \rvs_{N-1}^{\tree}$ be any completion to an orthonormal basis of $\R^{N}$. Define the matrix $\rmS^{\tree} \coloneqq [\rvs_0^{\tree}, \rvs_1^{\tree}, \ldots, \rvs_{N-1}^{\tree}]$ and let
\begin{align}
    \rmZ^{\tree} \coloneqq d^{-1/2} (\rmS^{\tree})^{\sT} \, \rmA^{\tree}\, \rmS^{\tree}
\end{align}
denote the representation of $\rmA^{\tree}/\sqrt{d}$ under this basis. Note that $\rmZ^{\tree}$ and $\rmA^{\tree}/\sqrt{d}$ share the same spectrum. The upper-left $(r_{x} + 1)\times (r_{x} + 1)$ principal submatrix $\rmZ_{\llbracket r_{x} \rrbracket}^{\tree}$ of $\rmZ^{\tree}$ can be written as
\begin{align}
    \rmZ_{\llbracket r_{x} \rrbracket}^{\tree}(\alpha_{x}) = \indi{x\in \gV_1} \rmZ^{(1)}(\alpha_x) + \indi{x\in \gV_2} \rmZ^{(2)}(\alpha_x),
\end{align}
where $\rmZ^{(1)}(\alpha_{x})$ and $\rmZ^{(2)}(\alpha_{x})$ are tridiagonal matrices due to the tree structure, defined by
\begin{align}\label{eqn:tridiag_Z_outline}
\rmZ^{(1)}(\alpha_{x}) \coloneqq \begin{bmatrix}
0  & \sqrt{\alpha_{x}} & \\
\sqrt{\alpha_{x}} & 0 & \qr\\
   &  \qr &  0 &  \qr^{-1}\\
    & &\qr^{-1} & 0 & \ddots\\
    & &  &  \ddots &  \ddots & \qr \\
    & &  & &   \qr & 0
\end{bmatrix}, \quad
\rmZ^{(2)}(\alpha_{x})=\begin{bmatrix}
0  & \sqrt{\alpha_{x}} & \\
\sqrt{\alpha} & 0 & \qr^{-1}\\
   &  \qr^{-1} &  0 &  \qr\\
    & &\qr & 0 & \ddots\\
    & &  &  \ddots &  \ddots & \qr^{-1} \\
    & &  & &   \qr^{-1} & 0
\end{bmatrix}.
\end{align}
With the detailed calculations deferred to \Cref{sec:tri_regular}, one can show that as $r_{x} \to\infty$,
\begin{enumerate}
    \item $\alpha_{x} > \qr^{-2} + 1$: the largest eigenvalue and smallest eigenvalue of $\rmZ^{(1)}(\alpha_{x})$ converge to $\Lambda_{\qr^{-1}}(\alpha_{x})$ and $ -\Lambda_{\qr^{-1}}(\alpha_{x})$, respectively, where $\Lambda_{\qr^{-1}}(t)$ is defined in \eqref{eqn:Lambda_qr_inverse}.
    \item $\alpha_{x} > \qr^{2} + 1$: the largest eigenvalue and smallest eigenvalue of $\rmZ^{(2)}(\alpha_{x})$ converge to $\Lambda_{\qr}(\alpha_{x})$ and $ -\Lambda_{\qr}(\alpha_{x})$, respectively, where $\Lambda_{\qr}(t)$ is defined in \eqref{eqn:Lambda_qr}.
    \item $0< \alpha_{x} <\qr^2-1$: the smallest positive eigenvalue and largest negative eigenvalue of $\rmZ^{(2)}(\alpha_{x})$ converge to $\Lambda_{\qr}(\alpha_{x})$ and $-\Lambda_{\qr}(\alpha_{x})$, respectively.
\end{enumerate}
    
Concurrently, one is able to construct the approximate eigenvectors of $\rmZ^{(1)}(\alpha_{x})$ and $\rmZ^{(2)}(\alpha_{x})$ as follows. 

For $x\in \Vonehigh$, let $\rvu_{+} = \{\ervu_{j}\}_{j=0}^{r_{x}}$ and $\rvu_{-} = \{(-1)^{j}\ervu_{j}\}_{j=0}^{r_{x}}$ have the components
\begin{subequations}
\begin{align}
    \ervu_0 \in \R \setminus \{0\}, \quad \ervu_1 &\,= \alpha_{x}^{-1/2} \Lambda_{\qr^{-1}} (\alpha_{x})\ervu_0, \quad \ervu_2 = \alpha_{x}^{1/2}\qr(\alpha_{x} - \qr^{-2})^{-1} \ervu_0\,,\label{eqn:u012V1}\\
    \ervu_{2j+1} &\, = (\alpha_{x} - \qr^{-2})^{-j} \ervu_1, \quad \ervu_{2j+2}= (\alpha_{x} -\qr^{-2})^{-j} \ervu_2,\label{eqn:ujsV1}
\end{align}
\end{subequations}
for $1 \leq j \leq \lfloor (r_{x}-3)/2 \rfloor$. Here, $\ervu_0$ is chosen to ensure $\|\rvu_{+}\|_2 =  1$. One can easily verify that $\rvu_{+} \perp \rvu_{-}$, and  as $r_{x} \to \infty$, $\rvu_{+}$ and $\rvu_{-}$ are the approximate eigenvectors of $\rmZ^{(1)}(\alpha_{x})$ corresponding to $\Lambda_{\qr^{-1}}(\alpha_{x})$ and $ -\Lambda_{\qr^{-1}}(\alpha_{x})$, respectively. 

For $x\in \Vtwohigh \cup \Vtwolow$, let $\rvu_{+} = \{\ervu_{j}\}_{j=0}^{r_{x}}$ and $\rvu_{-} = \{(-1)^{j}\ervu_{j}\}_{j=0}^{r_{x}}$ have the components
\begin{subequations}
    \begin{align}
        \ervu_0 \in \R \setminus \{0\}, \quad \ervu_1 &\,= \alpha^{-1/2}_{x} \Lambda_{\qr} (\alpha_{x})\ervu_0, \quad \ervu_2 = \alpha_{x}^{1/2}\qr^{-1}(\alpha_{x} - \qr^{2})^{-1} \ervu_0\,, \label{eqn:u012V2}\\
        \ervu_{2j+1} &\,= (\alpha_{x} - \qr^{2})^{-j} \ervu_1, \quad \ervu_{2j+2} = (\alpha_{x} -\qr^{2})^{-j} \ervu_2, \label{eqn:ujsV2}
    \end{align}
\end{subequations}
for $1 \leq j \leq \lfloor (r_x-3)/2 \rfloor$. Similarly, one can show that $\rvu_{+}$ and $\rvu_{-}$ are the approximate eigenvectors of $\rmZ^{(2)}(\alpha_{x})$ corresponding to $\Lambda_{\qr}(\alpha_{x})$ and $-\Lambda_{\qr}(\alpha_{x})$ as $r_{x} \to \infty$, respectively.

With all of the ingredients above, we can build approximate eigenpairs of $\scA$ in \eqref{eqn:scAdj}. Let $\ones_{y}\in \{0, 1\}^{N}$ denote the vector where only its $y$-entry is $1$. Let $\ones_{\layer_{j}^{\gG}(x)} = \sum_{y\in \layer_{j}^{\gG}(x)} \ones_{y}$ and $\rvs_{j}^{\gG} = |\layer_{j}^{\gG}(x)|^{-1/2}\ones_{\layer_{j}^{\gG} (x)}$ for $0\leq j\leq r_{x}$. Consider the following two vectors in $\R^{N}$
\begin{align}
    \rvv_{+}(x) \coloneqq \sum_{j=0}^{r_{x}} \ervu_{j} \rvs_{j}^{\gG}, \quad \rvv_{-}(x) \coloneqq \sum_{j=0}^{r_{x}} (-1)^{j} \ervu_{j} \rvs_{j}^{\gG}, \label{eqn:vx_plus_minus}
\end{align}
where coefficients $\ervu_{j}$ are determined by \eqref{eqn:u012V1} and \eqref{eqn:ujsV1} if $x\in \Vonehigh$, otherwise by \eqref{eqn:u012V2} and \eqref{eqn:ujsV2} if $x\in \Vtwohigh \cup \Vtwolow$. It is not difficult to verify that $\rvv_{+}(x) \perp \rvv_{-}(x)$ and $\|\rvv_{+}(x) \| = \|\rvv_{-}(x)\| = 1$. As will be shown in \Cref{prop:approximate_eigenvalues}, the construction of $\rvv_{+}(x)$ and $\rvv_{-}(x)$ locates one positive and one negative eigenvalue of $\rmH$ outside the bulk, namely, the outlier eigenvalues corresponding to $x$.

\subsection{Sufficiency}
Let $\lambda > 0$ be a right outlier of $\scA$. \Cref{prop:approximate_eigenvalues} proves that there exists some vertex $x\in \pruneV$ such that $\lambda$ can be approximated by $\Lambda_{\qr}(\alpha_{x})$ if $x\in \Vtwohigh \cup \Vtwolow$, or $\Lambda_{\qr^{-1}}(\alpha_{x})$  if $x\in \Vonehigh$. However, the largest outlier of $\scA$ may not necessarily correspond to the vertex with the highest normalized degree $\alpha_x$ in $\pruneV$. Comparisons between outlier eigenvalues are required to establish the one-to-one correspondence between outliers of $\scA$ and vertices in $\pruneV$, as shown in \Cref{thm:right_edge_behavior} (1) and \Cref{thm:left_edge_behavior} (1). This is accomplished by establishing a lower (resp. upper) bound on the $l$-th largest (resp. smallest) eigenvalue of $\scA$, in terms of the $l$-th largest (resp. smallest) degree of the vertices in $\gG$.

To that end, we construct the pruned graph $\pruneG$, which unites the local trees and is used to distinguish the eigenvalues. As will be shown in \Cref{lem:existence_pruned_graph}, $\pruneG$ is a subgraph of $\gG$ such that with high probability:
\begin{enumerate}
    \item $\rmA$ is well approximated by $\rmA^{(\vtau)}$, where $\rmA^{(\vtau)}$ denotes the adjacency matrix of $\pruneG$.
    \item for distinct vertices $x, y \in \pruneV$, the balls $\ball_{\radius}^{\pruneG}(x)$ and $\ball_{\radius}^{\pruneG}(y)$ are disjoint. 
\end{enumerate}
Consequently, $\rvv_{+}(x)$ and $\rvv_{+}(y)$, defined in \eqref{eqn:vx_plus_minus}, are two orthogonal approximate eigenvectors of $\scA$, which subsequently induce two distinct eigenvalues, or the same eigenvalue with multiplicity at least two. Detailed proofs are deferred to \Cref{sec:outlier_locations}.

\subsection{Necessity}
The only thing left is to show that ``vertex-to-outlier'' is the unique mechanism that produces outliers in the critical regime \eqref{eqn:critical_regime}. In other words, we need to prove that the remaining eigenvalues lie in the bulk $[\qr - \qr^{-1}, \qr + \qr^{-1}]$ up to some vanishing error. We sketch the proof for the lower bound, and the proof for the upper bound follows similarly. Details can be found in \Cref{sec:bulk_boundness}. 

Let $(\lambda,\rvw)$ be an eigenpair such that $\lambda>0$ and $\rvw$ is perpendicular to all the eigenvectors corresponding to outliers characterized above. We denote
\begin{align}
    \rvw = \begin{bmatrix}
    \rvw^{(1)}\\
    \rvw^{(2)}
    \end{bmatrix},
\end{align}
where $\rvw^{(1)}=[\ervw_{x}]_{x\in\gV_1}$ and $\rvw^{(2)}=[\ervw_{x}]_{x\in\gV_2}$, since $\scA$ in \eqref{eqn:scAdj} and $\scA^{(\vtau)}$ in \eqref{eqn:def_H_tau} both admit a block structure. Furthermore, since $\|\rvw\|_{2}^2 = 1$, it is not difficult to show that 
    \begin{align}
         \|\rvw^{(1)}\|^2= \|\rvw^{(2)}\|^2=\frac{1}{2}.
    \end{align}
At the same time, with $\|\rmE\| = o(1)$, we obtain the following Loewner order in \Cref{prop:lower_bound_H}:
\begin{align}
    {\rmX}^{*}(\rmI - d^{-1}\rmD^{(1)})^{-1} {\rmX}\succeq d^{-1}\rmD^{(2)}- \rmI - \rmE.
\end{align}
We use $\rvw^{(2)}$ as a test vector on both sides of the inequality. For the left hand, since $\lambda \rvw^{(1)} = \rmX \rvw^{(2)}$ and $\|\rvw^{(1)}\|^2 = 1/2$, we can obtain that
\begin{align}
   \< \rmX\rvw^{(2)}, (\rmI - d^{-1}\rmD^{(1)})^{-1} \rmX\rvw^{(2)}\> = \lambda^{2} \< \rvw^{(1)}, (\rmI - d^{-1}\rmD^{(1)})^{-1} \rvw^{(1)}\> \leq \lambda^{2} \big(1 - \qr^{-2} - o(1) \big)^{-1}/2.
\end{align}
For the right hand side, we apply the delocalization estimate \Cref{prop:delocalization_vec}, which then gives us
\begin{align}
    \< \rvw^{(2)}, (d^{-1}\rmD^{(2)}- \rmI - \rmE) \rvw^{(2)}\> \geq (\qr^{2} - 1 - o(1))/2.
\end{align}
Combining all the estimates above, one can conclude that with high probability, $\lambda \geq \qr - \qr^{-1} - o(1)$.

\section{Existence of Approximate eigenpairs for outliers}\label{sec:approx_eigenpairs}

Throughout the section, unless otherwise stated, we denote $\ball_{j}$, $\layer_{j}$ instead of $\ball_{j}^{\gG}(x)$, $\layer_{j}^{\gG}(x)$ for convenience.

\subsection{Approximate eigenpairs for outliers}
Recall that $\pruneV$ in \eqref{eqn:pruneV} denotes the vertex set with extreme degrees. As illustrated in \Cref{sec:vertex_to_outlier}, vertex $x \in \pruneV$ induces two outliers, one positive and one negative, with the corresponding approximate eigenvectors defined in \eqref{eqn:vx_plus_minus}. In the following, we focus on positive outliers, since dealing with negative outliers requires straightforward modifications.

For $x\in \pruneV$, define the associated radius by
\begin{align}
    r_x \coloneqq \left\lfloor \frac{d}{6\log \rD_{x}} \right\rfloor. \label{eqn:defrx}
\end{align}
Denote $\rvs_j = |\layer_{j}|^{-1/2}\ones_{\layer_{j}}$ for $0\leq j \leq r_{x}$. As in \eqref{eqn:vx_plus_minus}, we consider the following approximate eigenvector   
\begin{align}
    \rvv = \sum_{j=0}^{r_{x}} \ervu_{j} \rvs_{j},\label{eqn:approx_eigenvector}
\end{align}
where are coefficients $\ervu_{j}$'s are determined by \eqref{eqn:u012V1} and \eqref{eqn:ujsV1} if $x\in \Vonehigh$, or by \eqref{eqn:u012V2} and \eqref{eqn:ujsV2} if $x\in \Vtwohigh \cup \Vtwolow$. 

With the proof deferred to \Cref{sec:proof_prop_approx_eigenvalues}, \Cref{prop:approximate_eigenvalues}  shows that $\rvv$ is an approximate eigenvector of $\scA$ corresponding to some positive outlier, which can be approximated by either $\Lambda_{\qr^{-1}}(\alpha_x)$ or $\Lambda_{\qr}(\alpha_x)$.

\begin{proposition}[Approximate eigenpairs] \label{prop:approximate_eigenvalues}
The following holds with very high probability:
\begin{enumerate}
    \item For $x\in \Vtwohigh$ with $\tau^{+}_{2} = \qr^{2} + 1 + \error^{1/2}$, under \Cref{ass:proportional_regime},  for any $r\in [\log d, r_x]$
\begin{align}
    \| (\scA-\Lambda_{\qr}(\alpha_x) ) \rvv\|\leq \const \error.
\end{align}
    \item For $x\in \Vtwolow$ with $\tau^{-}_{2} = \qr^{2} - 1 - \error^{1/2}$,  under Assumptions \ref{ass:critical_connectivity} and \ref{ass:Ihara_Bass}, for any $r\in [\log d, r_x]$, 
    \begin{align}
         \| (\scA-\Lambda_{\qr}(\alpha_x) ) \rvv\|\leq \const \error.
    \end{align}
\item For $x\in \Vonehigh$ with $\tau_1 = \qr^{-2} + 1 + \error^{1/2}$, under \Cref{ass:proportional_regime}, for any for any $r\in [\log d, r_x]$,
\begin{align}
    \| (\scA-\Lambda_{\qr^{-1}}(\alpha_x)) \rvv\| \leq \const \error.
\end{align}
\end{enumerate}
\end{proposition}

\subsection{Bounds on approximation errors} We present several Lemmas in this section, which are crucial to conclude \Cref{prop:approximate_eigenvalues}. 
With the proof deferred to \Cref{sec:five_term_decomp}, \Cref{lem:five_term_decomp} splits the entire error into five different terms. Each term measures a specific way in which the real graph $\gG$ locally departs from the ideal biregular tree $\tree$ in \Cref{def:biregular_tree}.

\begin{lemma}[Error decomposition]\label{lem:five_term_decomp}
Recall $\scA$ in \eqref{eqn:scAdj} and $\rvv$ in \eqref{eqn:approx_eigenvector}. Let $\NS{j}{y}$ denote the set of edges which starts at $y$ and ends at some vertex in $\layer_{j} \coloneqq \layer_{j}^{\gG}(x)$. Its cardinality can be then calculated by
\begin{align}
    |\NS{j}{y}|\coloneqq &\, \<\ones_y, \rmA \ones_{\layer_{j}}\> = |\layer_{j}(x) \cap \layer_{1}(y)| \\
    =&\, \indi{ y\in \layer_{j-1} } (\rD_y - \indi{j \geq 2 }) + \indi{y\in \layer_{j+1}}. \label{eqn:Nj(y)}
\end{align}

\begin{enumerate}
    \item[\rm{(i)}] For $x\in \gV_2$, with $\Lambda_{\qr}(t)$ defined in \eqref{eqn:Lambda_qr}, the following decomposition holds
    \begin{align}
         \big(\scA - \Lambda_{\qr}(\alpha_x) \big) \rvv = \rvw_0 + \rvw_1 + \rvw_2 + \rvw_3 + \rvw_4, \label{eqn:errors_V2}
    \end{align}
where for $r\leq r_{x}$, the error terms $\rvw_0, \rvw_1, \rvw_2, \rvw_3, \rvw_4$ can be written as
\begin{subequations}
\begin{align}
    \rvw_0 &\, \coloneqq - \frac{1}{\sqrt{d}} (\E \rmA) \, \rvv\,,\\
    \rvw_1 &\, \coloneqq \frac{1}{\sqrt{d}}\sum_{j=0}^{r} \frac{\ervu_{j}}{\sqrt{|\layer_{j}|}}\Bigg( \sum_{y\in \layer_{j}} |\NS{j}{y}| \ones_y + \sum_{y\in \layer_{j+1}}(|\NS{j}{y}|-1)\ones_y \Bigg)\,,\\
     \rvw_2 &\, \coloneqq \frac{1}{\sqrt{d}}\sum_{j=1}^{r} \frac{\ervu_{j}}{\sqrt{|\layer_{j}|}} \sum_{y\in \layer_{j-1}}\bigg(|\NS{j}{y}|-\frac{|\layer_{j}|}{|\layer_{j-1}|} \bigg)\ones_y\,,\\
    \rvw_3 &\, \coloneqq \ervu_2 \Big( \frac{\sqrt{|\layer_2|}}{\sqrt{d}\cdot \sqrt{|\layer_1|}} - \qr^{-1} \Big)\rvs_1 \\
    &\, \quad \quad + \sum_{j=1}^{\lfloor r/2\rfloor - 1 } \Bigg( \ervu_{2j-1}\Big(\frac{\sqrt{|\layer_{2j}|}}{ \sqrt{d} \cdot \sqrt{|\layer_{2j-1}|}} - \qr^{-1} \Big) + \ervu_{2j+1}\bigg(\frac{\sqrt{|\layer_{2j+1}|}}{\sqrt{d}\cdot \sqrt{|\layer_{2j}|}} - \qr \bigg) \Bigg) \rvs_{2j}\\
    &\, \quad \quad + \sum_{j=1}^{\lfloor r/2\rfloor - 1} \Bigg( \ervu_{2j}\Big(\frac{\sqrt{|\layer_{2j+1}|}}{\sqrt{d} \cdot \sqrt{|\layer_{2j}|}} - \qr \Big) + \ervu_{2j+2}\bigg(\frac{\sqrt{|\layer_{2j+2}|}}{\sqrt{d}\cdot \sqrt{|\layer_{2j+1}|}} - \qr^{-1} \bigg) \Bigg) \rvs_{2j+1}\,, \notag\\
    \rvw_4 &\, \coloneqq \bigg( \frac{\ervu_{r-1}}{\sqrt{d}} \frac{\sqrt{|\layer_{r}|}}{\sqrt{|\layer_{r - 1}|}} - \indi{r \textnormal{ even}} \big( \qr^{-1} \ervu_{r - 1} + \qr \ervu_{r + 1} \big) - \indi{ r \textnormal{ odd}} \big( \qr \ervu_{r - 1} + \qr^{-1} \ervu_{r + 1}\big)\bigg) \, \rvs_{r} \\
    &\, \quad \quad + \frac{\ervu_{r}}{\sqrt{d}}\frac{\sqrt{|\layer_{r + 1}|}}{\sqrt{|\layer_{r}|}} \, \rvs_{r + 1}. \notag
\end{align}
\end{subequations}

\item[\rm{(ii)}] For $x\in \gV_{1}$, with $\Lambda_{\qr^{-1}}(t)$ defined in \eqref{eqn:Lambda_qr_inverse}, the following decomposition holds
\begin{align}
    \big(\scA - \Lambda_{\qr^{-1}}(\alpha_x) \big) \rvv = \rvw_0 + \rvw_1 + \rvw_2 + \rvw_3 + \rvw_4, \label{eqn:errors_V1}
\end{align}
where in the definitions of $\rvw_{0}, \rvw_{1}, \rvw_{2}, \rvw_{3}, \rvw_{4}$, $\qr$ is substituted by $\qr^{-1}$.
\end{enumerate}
\end{lemma}
Following the framework in \cite{alt2021extremal, alt2021delocalization}, we present an informal description of each term below.
\begin{itemize}
\item $\rvw_0$ -- effect of the expectation. The largest eigenvalue of $\E \rmA$ corresponds to the completely delocalized vector $\mathbf{e} =N^{-1/2}\sum_{x\in [N]}\ones_{x}$. However, $\rvv$ is localized near $x$ by construction. Consequently, the overlap between $\rvv$ and $\mathbf{e}$ is small, leading to the smallness of $\|\rvw_0\|$.

\item  $\rvw_1$ -- local deviation from a tree. According to \eqref{eqn:Nj(y)}, $\rvw_1$ vanishes \emph{exactly} when $\gG |_{\ball_{r_{x}}(x)}$ contains no cycle. The error $\|\rvw_1\|$ turns out to be small with high probability, since it is controlled by the number of cycles $\ball_{r_{x}}(x)$, which is bounded according to \Cref{lem:few_even_cycle_prob}.

\item $\rvw_2$ -- deviation from the \emph{same-number-of-children} requirement for vertices in the same layer. The vector $\rvw_{2}$ quantifies the extent to which $\gG |_{\ball_{r_{x}}(x)}$ deviates from a tree with the property that for $1\leq j\leq r_{x}$, all the vertices in $\layer_{j}(x)$ have the same number of children in $\layer_{j+1}(x)$. With high probability, $\|\rvw_{2}\|$ is small since $\NS{j}{y}$ concentrates around $|\layer_j|/| \layer_{j-1}|$, and the number of cycles in $\gG |_{\ball_{r_{x}}(x)}$ is bounded.

\item  $\rvw_3$ -- deviation from layerwise growth ratios. The vector $\rvw_{2}$ quantifies the extent to which $\gG |_{\ball_{r_{x}}(x)}$ deviates from the idealized tree $\tree$ in \Cref{def:biregular_tree}. With high probability, $|\layer_{j+1}|/|\layer_j|$ concentrates around either~$d_{1}$ or $d_{2}$ depending on the community that $x$ belongs to, which in turn ensures the smallness of $\|\rvw_{3}\|$.

\item $\rvw_4$ -- leakage through the tree boundary $\layer_r$. The vector $\rvw_4$ quantifies the error arising from edges that connects $\ball_{r}(x)$, and to the rest of the graph $[N]\setminus \ball_r(x)$. The smallness of $\|\rvw_4\|$ is ensured by the exponential decay of the coefficients ~$\ervu_{r-1}$, $\ervu_{r}$ and $\ervu_{r+1}$, as defined in \eqref{eqn:ujsV1} and \eqref{eqn:ujsV2}.
\end{itemize}

Furthermore, errors are bounded quantitatively in \Cref{lem:w0tow5}, whose proof is deferred to \Cref{sec:w0tow5}.
\begin{lemma}[Approximation error bounds]\label{lem:w0tow5}
 Let $\epsilon > 0$ be some absolute constant. For $x\in \pruneV$, define
\begin{align}
    \Xi_{x}(\epsilon) \coloneqq  \{\epsilon d \leq \rD_{x} \leq \sqrt{N} (2d)^{-r_x}\} \label{eqn:Dx_Event},
\end{align}
where $\sqrt{N} (2d)^{-r_{x}} \lesssim N^{1/4}$ according to the choice of $r_x$ in \eqref{eqn:defrx}. Conditioned on event \eqref{eqn:Dx_Event}, the following hold with very high probability for any $r\le r_x$:

\begin{enumerate}
    \item For $x \in \Vtwohigh$, the following hold under \Cref{ass:proportional_regime},
\begin{subequations}
    \begin{align}
        \|\rvw_0\| \lesssim &\, \sqrt{d} N^{-1/4},\label{eqn:norm_w0}\\
        \|\rvw_1\| \lesssim &\, (d\rD_{x})^{-1/2},\label{eqn:norm_w1}\\
        \|\rvw_2\| \lesssim &\, (\log d)^{1/2}d^{-1/2}\big(1 + \log(N)/\rD_x \big)^{1/2}, \label{eqn:norm_w2}\\
        \|\rvw_3\| \lesssim &\,  (\log N)^{-1/2}(d\rD_x)^{-1/2},\label{eqn:norm_w3}\\
        \|\rvw_{4}\| \lesssim &\, |\alpha_{x} - \qr^{2}|^{-\lceil(r-1)/2\rceil}. \label{eqn:norm_w4}
    \end{align}
\end{subequations}

    \item For $x \in \Vtwolow$, under Assumptions \ref{ass:critical_connectivity} and \ref{ass:Ihara_Bass}, \eqref{eqn:norm_w0}, \eqref{eqn:norm_w1}, \eqref{eqn:norm_w2}, \eqref{eqn:norm_w3} and \eqref{eqn:norm_w4} hold.

    \item For $x \in \Vonehigh$, under \Cref{ass:proportional_regime}, \eqref{eqn:norm_w0}, \eqref{eqn:norm_w1}, \eqref{eqn:norm_w2}, \eqref{eqn:norm_w3} and \eqref{eqn:norm_w4_V1} hold.
    \begin{align}
        \|\rvw_{4}\| \lesssim &\, |\alpha_{x} - \qr^{-2}|^{-\lceil(r-1)/2\rceil}. \label{eqn:norm_w4_V1}\,\,
    \end{align}
\end{enumerate}
\end{lemma}


\subsection{Proof of \Cref{prop:approximate_eigenvalues}}\label{sec:proof_prop_approx_eigenvalues}
We prove (1) first. For $x\in \Vtwohigh$, the event \eqref{eqn:Dx_Event} occurs with very high probability under \Cref{ass:proportional_regime}. Recall $\error$ in \eqref{eqn:error_parameters}, by \Cref{lem:w0tow5}, 
\begin{align}
   \|\rvw_{0}\| \lesssim &\, \sqrt{d} N^{-1/4} \ll \error, \quad \,\, \|\rvw_{1}\| \lesssim (\log N)^{-1} \ll \error, \\
   \|\rvw_{2}\| \lesssim &\, \error, \quad \quad \quad \quad \quad \quad \quad \|\rvw_{3}\| \lesssim (\log N)^{-1/2} \ll \error. 
\end{align}
Note that $\tau^{+}_{2} = \qr^{2} + 1 + \error^{1/2}$, as a consequence of \eqref{eqn:norm_w4},  
\begin{align}
    \|\rvw_{4}\| \lesssim (1 + \error^{1/2})^{-\lceil(r_{x} -1)/2\rceil} \ll \error.
\end{align}
The proof is then concluded by \Cref{lem:five_term_decomp}, triangle inequality and the estimates above, since
\begin{align}
    \| (\scA-\Lambda_{\qr}(\alpha_x) ) \rvv\| \lesssim \|\rvw_{2}\| \leq \error.
\end{align}

For (3), the proof for $x\in \Vonehigh$ with $\tau_1 = \qr^{-2} + 1 + \error^{1/2}$ follows in a similar way under \Cref{ass:proportional_regime}.

For (2), where $x\in \Vtwolow$ with $\tau^{-}_{2} = \qr^{2} - 1 - \error^{1/2}$, the event \eqref{eqn:Dx_Event} occurs with high probability under Assumptions \ref{ass:critical_connectivity} and \ref{ass:Ihara_Bass}. The rest of the proof follows verbatim.

\section{Locations of outlier singular values}\label{sec:outlier_locations}
The proofs of \Cref{thm:right_edge_behavior} (1) and \Cref{thm:left_edge_behavior} (1) follow the same strategy as in {\cite[Section 7]{alt2021extremal}}.

\subsection{Existence of a pruned graph}
Let \eqref{eqn:bennett_rate} denote the rate function in Bennett inequality (\Cref{lem:Bennett})
\begin{align}
    \benrate(\alpha) \coloneqq (1 + \alpha) \log(1 + \alpha) - \alpha, \label{eqn:bennett_rate}.
\end{align}
For $\qr \geq 1$ and $\tau >0$, define the following functions through \eqref{eqn:bennett_rate}:
\begin{subequations}
    \begin{align}
        \benrate_{\qr^{-1}}(\tau) \coloneqq \qr^{-2} \cdot \benrate \Big(\frac{|\tau - \qr^{-2}|}{2\qr^{-2}}\Big), \quad &\,\benrate_{\qr}(\tau) \coloneqq \qr^{2} \cdot \benrate \Big(\frac{|\tau - \qr^{2}|}{2\qr^{2}}\Big), \label{eqn:htau} \\
        \mathrm{r}_{\qr^{-1}}(\tau) \coloneqq  \frac{d}{2\log d} \cdot \benrate_{\qr^{-1}}(\tau) - 1, \quad &\, \mathrm{r}_{\qr}(\tau) \coloneqq \frac{d}{2\log d} \cdot \benrate_{\qr}(\tau) - 1. \label{eqn:rtau}
    \end{align}
\end{subequations}
Recall $\Vonehigh$, $\Vtwohigh$, $\Vtwolow$ in \eqref{eqn:atypical_vertices}, $\pruneV$ in \eqref{eqn:pruneV} and $\error$ in \eqref{eqn:error_parameters}. In this section, we choose
 \begin{align}
   \tau_{1} = \qr^{-2} + 1 + \error^{1/2}, \quad \tau_{2}^{+} = \qr^{2} + 1 + \error^{1/2}, \quad \tau_{2}^{-} = \qr^{2} - 1 - \error^{1/2}.
\end{align}
For each $x\in \pruneV$, define the following radius associated to $x$ and $\vtau =(\tau_{1}, \tau_{2}^{+}, \tau_{2}^{-})$
\begin{align}
    \pruneR{x} \coloneqq \Big(\frac{1}{4} r_x \Big) \wedge \bigg[\frac{1}{2}\mathrm{r}_{\qr^{-1}}(\tau_{1}) \cdot \indi{x\in \Vonehigh} + \frac{1}{2}\mathrm{r}_{\qr}(\tau^{+}_{2}) \cdot \indi{x\in \Vtwohigh} + \frac{1}{2}\mathrm{r}_{\qr}(\tau^{-}_{2}) \cdot \indi{x\in \Vtwolow} \bigg].
\end{align}
We define the following radius
\begin{equation}\label{eqn:radius} 
     \radius = \left\lfloor \frac{1}{4\qr^2}\sqrt{\frac{d}{\log(d)}}\right\rfloor.
\end{equation}
Since $\benrate(\alpha)\geq \frac{\alpha^2}{2(1+\alpha/3)}$ for $\alpha\geq 0$, the following holds for sufficiently large $d$
\begin{align}
    \radius \leq  1 + 2\cdot \inf_{x\in \pruneV} \,\, \Big\lfloor \frac{1}{2} \pruneR{x} \Big\rfloor.
    \label{eqn:radius0}
\end{align}

The existence of the pruned graph $\pruneG$ is guaranteed by \Cref{lem:existence_pruned_graph}. For convenience, we denote
\begin{align}\label{eqn:S_j_tau}
    \layer_{j}^{(\vtau)}(x) \coloneqq \layer_{j}^{\pruneG}(x), \quad \ball_{j}^{(\vtau)}(x) \coloneqq \ball_{j}^{\pruneG}(x).  
\end{align}

\begin{lemma}[Existence of the pruned graph]\label{lem:existence_pruned_graph}
For the triple $\vtau = (\tau_{1}, \tau_{2}^{-}, \tau_{2}^{+})$, with very high probability, there exists a subgraph $\pruneG$ of $\gG$ that satisfies the following properties:
\begin{enumerate}
    \item If two vertices $x, y\in \pruneV$ are connected by a path $\ell$ in $\pruneG$, then the length of $\ell$ is at least $2\radius + 1$. In particular, the balls $\ball_{\radius}^{(\vtau)}$ for $x\in \pruneV$ are disjoint.  \label{item:one}
    \item The induced subgraph $\pruneG|_{\ball_{2\radius}^{(\vtau)}}$ is a tree for each $x\in \pruneV$. \label{item:two}
    \item For each edge in $\gG \setminus \pruneG$, there is at least one vertex in $\pruneV$ incident to it. \label{item:three}
    \item For each $x \in \pruneV$ and each $j\in [1, \radius]$, we have $\layer^{(\vtau)}_{j}(x) \subset \layer^{\gG}_{j}(x)$, and for all $y \in \ball^{(\vtau)}_{\radius}\setminus \{x\}$ \label{item:four}
        \begin{align}
            \layer^{(\vtau)}_{1}(y) \cap \layer^{(\vtau)}_{j}(x) = \layer^{\gG}_{1}(y) \cap \layer^{\gG}_{j}(x).
        \end{align}
        
    \item The degrees of $\gV$ induced on $\gG \setminus \pruneG$ satisfy \label{item:five}
        \begin{align}
            \max_{x\in \gV} \rD^{\gG \setminus \pruneG}_{x} \lesssim 1 + \frac{\log(N)}{d} \cdot \Big( \min\Big\{\benrate_{\qr^{-1}}(\tau_{1}),\,\, \benrate_{\qr}(\tau^{-}_{2}),\,\,\benrate_{\qr}(\tau^{+}_{2}) \Big\}\Big)^{-1}.
        \end{align}
    \item For each $x\in \pruneV$, we have the following bound for each $j\in [2, \radius]$: \label{item:six} 
        \begin{align}
            |\layer^{\gG}_{j}(x) \setminus \layer^{(\vtau)}_{j}(x)| \lesssim \rD_{x}^{\gG \setminus \pruneG} \cdot d^{j-1}.
        \end{align}
\end{enumerate}
\end{lemma}
The proof of Lemma~\ref{lem:existence_pruned_graph} is deferred to \Cref{sec:proof_existence_prunded_graph}. 

\subsection{Proofs of \Cref{thm:right_edge_behavior} (1) and \Cref{thm:left_edge_behavior} (1)}
\begin{proof}[Proof of \Cref{thm:right_edge_behavior} (1)]

We first present several estimates which are crucial to conclude the proof.

Let $\rmA^{(\vtau)}$ denote the adjacency matrix of the pruned graph $\pruneG$ defined in Lemma~\ref{lem:existence_pruned_graph}. Let $\chi^{(\vtau)}$ denote the orthogonal projection onto the space spanned by $\{\ones_{y}: y\not\in \cup_{x\in \pruneV} \ball^{(\vtau)}_{2\radius}(x)\}$. Define
\begin{align}\label{eqn:def_H_tau}
    \scA^{(\vtau)} \coloneqq \frac{1}{\sqrt{d}} \Big(\rmA^{(\vtau)} - \chi^{(\vtau)}(\E\rmA^{(\vtau)})\chi^{(\vtau)} \Big) = \begin{bmatrix}
    \bzero &  \rmX^{(\vtau)} \\
    (\rmX^{(\vtau)})^{\sT}  & \bzero
    \end{bmatrix}.
\end{align}
The matrix $\scA^{(\vtau)}$ is defined in a way such that
\begin{enumerate}
    \item $\scA^{(\vtau)}$ is close to $\scA$ given that $\rmA^{(\vtau)}$ is close to $\rmA$, since the kernel of $\chi^{(\vtau)}$ has a relatively low dimension.
    \item When restricted to vertices with distance at most $2\radius$ from $\pruneV$, $\scA^{(\vtau)}$ coincides with $\rmA^{(\vtau)}/\sqrt{d}$, meaning that $\scA^{(\vtau)}$ inherits the local structure of the matrix $\rmA/\sqrt{d}$.
\end{enumerate}
It is not difficult to verify claim (1), since by applying \eqref{item:five} of \Cref{lem:existence_pruned_graph}, we have
\begin{align} \label{eqn:approx_H_Htau}
  \|\scA - \scA^{(\vtau)}\| \leq d^{-1/2} \max_{x\in \gV} \rD_{x}^{\gG \setminus \pruneG} \lesssim d^{-1/2}  
\end{align}

Below, we present a quantitative version of claim (2) in \eqref{eqn:normalized_approx}. According to \Cref{prop:approximate_eigenvalues} (1), for any $l\in \Vtwohigh$, there exists a unit vector $\rvv_{l}$, defined in \eqref{eqn:approx_eigenvector}, such that $(\Lambda_{\qr}(\alpha_{l}), \rvv_{l})$ is an approximate eigenpair of $\scA$. Let $\rvv^{(\vtau)}_{l}$ below denote the vector with entries restricted to the ball $\ball^{(\vtau)}_{\radius}(l)$,
\begin{align}
    \rvv^{(\vtau)}_{l} \coloneqq \rvv_{l} |_{\ball^{(\vtau)}_{\radius}(l)} \label{def:vltau}
\end{align}
where the entries not in the ball $\ball^{(\vtau)}_{\radius}(l)$ are zeroed out, hence $\|\rvv^{(\vtau)}_{l}\|_{2} \leq 1$. According to \Cref{lem:concentrationSi} and \eqref{item:five} \eqref{item:six} in \Cref{lem:existence_pruned_graph}, and the  fact that $\sum_{j = 0}^{\radius} \ervu_{j}^{2} \leq 1$, the following holds with high probability
\begin{align}
    \|\rvv^{(\vtau)}_{l} - \rvv_{l} \|_{2}^{2} = &\, \sum_{j = 0}^{\radius} \ervu_{j}^2 \frac{|\layer^{\gG}_j \setminus \layer^{(\vtau)}_{j}|}{|\layer^{\gG}_j|} \lesssim \frac{\rD_{x}^{\gG \setminus \pruneG} \cdot d^{j-1}}{\rD^{\gG}_{x} \cdot d^{j-1}} \\
    \lesssim &\, \frac{1}{\rD_{x}^{\gG}} + \frac{\log(N)}{d\rD_{x}} \cdot \Big( \min \big\{\benrate_{\qr^{-1}}(\tau_{1}),\,\, \benrate_{\qr}(\tau^{-}_{2}), \,\,\benrate_{\qr}(\tau^{+}_{2})\big\} \Big)^{-1} \lesssim d^{-1}, \label{eqn:vl_tau}
\end{align}
where the last inequality follows since $l\in \Vtwohigh$ and $\tau_{2}^{+} \geq \qr^{2} + 1$ for vertices inducing outliers. Let 
\begin{align}
    \widehat{\rvv}^{(\vtau)}_{l} \coloneqq \frac{1}{\|\rvv^{(\vtau)}_{l}\|_{2}} \rvv^{(\vtau)}_{l}
\end{align}
denote the normalized version of $\rvv^{(\vtau)}_{l}$. According to \Cref{prop:approximate_eigenvalues} and \eqref{eqn:vl_tau}, we have
\begin{align}
    \| \big( \scA - \Lambda_{\qr}(\alpha_{l})\, \big)\widehat{\rvv}^{(\vtau)}_{l}\| \leq &\, \frac{1}{\|\rvv^{(\vtau)}_{l}\|_{2}}\left(\| \big( \scA - \Lambda_{\qr}(\alpha_{l})\, \big) \rvv_{l} \| +  \| \big( \scA - \Lambda_{\qr}(\alpha_{l})\, \big)\cdot (\rvv^{(\vtau)}_{l} - \rvv_{l})\| \right)\\
    \lesssim &\, \error + d^{-1/2} \leq 2\error. \label{eqn:normalized_approx}
\end{align}

Furthermore, by applying the triangle inequality, together with \eqref{eqn:approx_H_Htau} and \eqref{eqn:normalized_approx}, the following holds for some constant $\const >0$
\begin{align}
    \<\widehat{\rvv}^{(\vtau)}_{l}, \scA^{(\vtau)}\widehat{\rvv}^{(\vtau)}_{l}\>  = &\, \Lambda_{\qr}(\alpha_{l}) + \<\widehat{\rvv}^{(\vtau)}_{l}, \big( \scA^{(\vtau)} - \Lambda_{\qr}(\alpha_{l}) \big)\widehat{\rvv}^{(\vtau)}_{l}\>  \\
    \geq &\, \Lambda_{\qr}(\alpha_{l}) - \|\big( \scA - \Lambda_{\qr}(\alpha_{l}) \big) \widehat{\rvv}^{(\vtau)}_{l}\| - \| \scA - \scA^{(\vtau)}\|\\
    \geq & \Lambda_{\qr}(\alpha_{l})- \const \error. \label{eqn:vltau_lowerbound}
\end{align}

We now present the proof of the main result. Define the space $\widehat{\sW}^{(\vtau)} = \mathrm{Span}\{\widehat{\rvv}^{(\vtau)}_{j}: j\in [l]\}$, where $\{\widehat{\rvv}^{(\vtau)}_{j}\}_{j=1}^{l}$ are pairwise orthogonal since they have disjoint supports due to \eqref{item:one} \eqref{item:two} of \Cref{lem:existence_pruned_graph}. For any linear subspace $\sW \subset \R^{N}$, let $\S(\sW)$ denote the unit sphere with respect to the Euclidean norm. According to the max-min principle, as well as the estimates \eqref{eqn:approx_H_Htau} and \eqref{eqn:vltau_lowerbound}, we have
    \begin{align}
        \lambda_{l}(\scA) = &\, \max_{\dim \sW = l} \,\, \min_{\rvw \in \S(\sW)} \<\rvw, \scA \rvw\> \\
        \geq &\, \min_{\rvw \in \S(\widehat{\sW}^{(\vtau)})} \<\rvw, \scA \rvw\> \\
        \geq  &\, \min_{\rvw \in \S(\widehat{\sW}^{(\vtau)})} \<\rvw, \scA^{(\vtau)} \rvw\> -\|\scA-\scA^{(\vtau)}\|\\
        \geq &  \Lambda_{\qr}(\alpha_{l}) - \const \error.\label{eqn:lower_bound_lth}
    \end{align}
Note that if $\rmM$ is a Hermitian matrix and $\|\rvv\|_{2} = 1$ such that $\|\rmM \rvv\| \leq \epsilon$, then $\rmM$ has an eigenvalue in $[-\epsilon, \epsilon]$. Consequently, there exists an eigenvalue of $\scA$ that lies in the interval $[\Lambda_{\qr}(\alpha_{l}) - \const\error,\,\, \Lambda_{\qr}(\alpha_{l}) + \const\error]$, due to \Cref{prop:approximate_eigenvalues}, this eigenvalue must be the $l$-th eigenvalue of $\scA$ according to \eqref{eqn:lower_bound_lth}. Therefore, $ |\lambda_{l}(\scA)- \Lambda_{\qr}(\alpha_{l})| \leq \const\error$, which concludes the proof.

The proof for the case $l\in \Vonehigh$ follows similarly by applying \Cref{prop:approximate_eigenvalues} (3).
\end{proof}

\begin{proof}[Proof of \Cref{thm:left_edge_behavior} (1)]
    The proof is similar to the proof of \Cref{thm:right_edge_behavior} (1) under \Cref{ass:Ihara_Bass}, where we apply \Cref{prop:approximate_eigenvalues} (2) for $x\in \Vtwolow$.

    Without loss of generality, we assume $\alpha_1\geq \alpha_2\geq \cdots \geq\alpha_{m}$.
For each $l$ such that $m - |\sL_2| + 1 \leq l \leq m$, we can define $\widehat{\rvv}_{l}^{(\vtau)}$ using \eqref{def:vltau}.
Define $\widehat{\sW}^{(\vtau)} = \mathrm{Span}\{\widehat{\rvv}^{(\vtau)}_{j}: l\leq j\leq m\}$, which is a subspace of dimension $m-l+1$.  Following  the same argument as the proof of \eqref{eqn:vltau_lowerbound}, we have 
\begin{align}
     \<\widehat{\rvv}^{(\vtau)}_{l}, \scA^{(\vtau)}\widehat{\rvv}^{(\vtau)}_{l}\>  \leq  \Lambda_{\qr}(\alpha_{l}) + \|\big( \scA - \Lambda_{\qr}(\alpha_{l}) \big) \widehat{\rvv}^{(\vtau)}_{l}\| + \| \scA - \scA^{(\vtau)}\| \leq \Lambda_{\qr}(\alpha_{l})+C\error.
\end{align}
According to the max-min principle and \eqref{eqn:approx_H_Htau}, the following holds for some constant $\const >0$
    \begin{align}
        \lambda_{l}(\scA) = &\, \min_{\dim \sW = m - l + 1} \,\, \max_{\rvw \in \S(\sW)} \<\rvw, \scA \rvw\> \\
        \leq &\, \max_{\rvw \in \S({\widehat{\sW}^{(\vtau)}})} \<\rvw, \scA \rvw\> \\
        \leq &\max_{\rvw \in \S({\widehat{\sW}^{(\vtau)}})} \<\rvw, \scA^{(\vtau)} \rvw\> + \|\scA - \scA^{(\vtau)} \| \\
         \leq & \Lambda_{\qr}(\alpha_{l}) + \const \error.
    \end{align}
    Similarly, due to  \Cref{prop:approximate_eigenvalues}~(2), $|\lambda_{l}(\scA)- \Lambda_{\qr}(\alpha_{l})| \leq \const\error$ for $m - |\sL_2| + 1 \leq l \leq m$.
\end{proof}

\section{Bounding the remaining singular values}\label{sec:bulk_boundness}

We present some preliminaries in Sections \ref{sec:Loewner-order} and \ref{sec:delocalization}, where proofs of lemmas are deferred to \Cref{sec:proof_bulk_boundness}.

\subsection{Loewner order via the non-backtracking operator}\label{sec:Loewner-order}
Below, we present two inequalities in Propositions \ref{prop:upper_bound_H} and \ref{prop:lower_bound_H}, which are crucial for proofs of Theorems \ref{thm:right_edge_behavior} (2) and \ref{thm:left_edge_behavior} (2), respectively.

\begin{proposition}\label{prop:upper_bound_H}
Consider the centered and normalized adjacency matrix $\scA$ in \eqref{eqn:scAdj}. For $4 \leq d \leq (mn)^{1/13}$ and some positive constants $c, \epsilon > 0$, the following holds with probability at least $1 - N^{3 - c\sqrt{d} \log(1 + \epsilon)}$,
    \begin{align}
         \scA \preceq \id_{N} + d^{-1}\rmD + \rmE,\label{eqn:upper_bound_H}
    \end{align}
    where $\rmD = \mathrm{diag}\{\rD_{x}\}_{ x\in \gV}$. Under this event, the error matrix $\rmE$ satisfies
    \begin{align}
        \|\rmE\| \lesssim d^{-3/2}(\Delta+ d) \label{eqn:error_upper_bound_norm}
    \end{align}
    where $\Delta = \max_{x\in \gV} \rD_{x}$ denotes the maximum degree of the graph.
\end{proposition}

\begin{proposition}\label{prop:lower_bound_H}
Consider the centered and normalized adjacency matrix $\scA$ in \eqref{eqn:scAdj} with $\rmX = (\widetilde{\rmA} - \E\widetilde{\rmA})/\sqrt{d}$. Under \Cref{ass:Ihara_Bass}, for $\log(N) \leq d \leq N^{2/13}$ and some positive constants $c, \epsilon > 0$, the following holds with probability at least $1 - N^{3 - c\sqrt{d} \log(1 + \epsilon)}$,
    \begin{align}
        {\rmX}^{*}(\rmI - d^{-1}\rmD^{(1)})^{-1} {\rmX}\succeq d^{-1}\rmD^{(2)}- \rmI - \rmE \label{eqn:lower_bound_H}
    \end{align}
where $\rmD^{(1)} = \mathrm{diag}\{\rD_{x}\}_{x \in \gV_{1}}$ and $\rmD^{(2)} = \mathrm{diag}\{\rD_{x}\}_{x \in \gV_{2}}$, and the error matrix $\rmE$ satisfies $\|\rmE\| \lesssim d^{-1/2}$.
\end{proposition}
For the sake of simplicity, we only present the proof of \Cref{prop:lower_bound_H} below. The proof of \Cref{prop:upper_bound_H} follows similarly, and is deferred to \Cref{sec:Loewner-order-app}.

We first introduce the following non-backtracking operator.
\begin{definition}[Non-backtracking operator of a matrix]
For $\rmH \in \mathbb{M}_{N} (\C)$, let $\rmB$ denote its non-backtracking operator, which is an $N^2\times N^2$ matrix with its entry defined as
\begin{align}\label{eqn:defB}
    \ermB_{ef} \coloneqq \ermH_{kl}\indi{j=k}\indi{i\not=l},
\end{align}
where $e =\{ i, j\} \subset [N]^{2}$ and $f = \{k, l\}\subset [N]^{2}$.
\end{definition}

\Cref{lem:Ihara-Bass} and \Cref{lem:block_Ihara} present the Ihara-Bass formula and its generalization to bipartite block matrices, which are crucial for the proof of \Cref{prop:lower_bound_H}. 
\begin{lemma}[{\cite[Lemma 4.1]{benaych2020spectral}}]\label{lem:Ihara-Bass}
Let $\rmH \in \mathbb{M}_{N} (\C)$ associated with the non-backtracking matrix $\rmB$ and let $\lambda \in \C$ satisfy $\lambda^2 \neq \ermH_{jl} \ermH_{lj}$ for all $j, l \in [N]$. Define the matrices $\rmH(\lambda)$ and $\rmM(\lambda) = \mathrm{diag}\{ \ermM_{jj}(\lambda)\}_{j\in [N]}$ through
\begin{align}
    \ermH_{jl}(\lambda) \coloneqq \frac{\lambda \ermH_{jl}}{\lambda^2 - \ermH_{jl} \ermH_{lj}}, \quad \ermM_{jj}(\lambda) \coloneqq 1 + \sum_{l\in [N]} \frac{\ermH_{jl} \ermH_{lj}}{\lambda^2 - \ermH_{jl} \ermH_{lj}}
\end{align}
Then $\lambda \in \mathrm{Spec}(\rmB)$ if and only if $\det(\rmM(\lambda) - \rmH(\lambda)) = 0$.
\end{lemma}

\begin{lemma}[{\cite[Lemma 3.2]{dumitriu2024extreme}}]\label{lem:block_Ihara}
Let $\rmX$ be an $n\times m$ complex matrix and define the $N \times N$ matrix $\rmH = \begin{bmatrix}
\bzero & \rmX\\
\rmX^{*} &\bzero 
\end{bmatrix}$
where $N = n + m$. Let $\rmB$ be the non-backtracking operator associated with $\rmH$. Define an $n\times m$ matrix $\rmX(\lambda)$, and two diagonal matrices $\rmM^{(1)}(\lambda)=\mathrm{diag}\{\ermM^{(1)}_{jj}(\lambda) \}_{j\in [n]}$,  $\rmM^{(2)}(\lambda)=\mathrm{diag}\{\ermM^{(2)}_{ll}(\lambda) \}_{l\in [m]}$ as follows:
\begin{align}\label{eqn:XM_entries}
    \ermX_{jl}(\lambda) = \frac{\lambda \ermX_{jl}}{\lambda^{2}- |\ermX_{jl}|^{2}},\quad 
    \ermM^{(1)}_{jj}(\lambda) = 1 + \sum_{l\in [m]}\frac{|\ermX_{jl}|^{2}}{\lambda^{2}- |\ermX_{jl}|^{2}}, \quad 
     \ermM^{(2)}_{ll}(\lambda) = 1+\sum_{j\in [n]}\frac{|\ermX_{jl}|^{2}}{\lambda^{2} - |\ermX_{jl}|^{2}}.
\end{align}
Let $\lambda \in \C$ satisfy $\lambda \neq |\ermX_{jl}|$ for all $j\in [n]$, $l\in [m]$. 
Assume that $\rmM^{(1)}(\lambda)$ is non-singular, then $\lambda \in \mathrm{Spec}(\rmB)$ if and only if 
\begin{align}
    \mathrm{det}\big(\rmM^{(2)}(\lambda) - \rmX^{*}(\lambda) [\rmM^{(1)}(\lambda)]^{-1}\rmX(\lambda) \big) = 0.
\end{align}
Here, $[\rmX^{*}(\lambda)]_{jl} = \frac{\lambda \ermX_{jl}^{*}}{\lambda^{2}- |\ermX_{jl}|^{2}}$.
\end{lemma}

As shown below, for a Hermitian random matrix with a block structure, with high probability, the spectral radius of its non-backtracking operator can be bounded from above.
\begin{lemma}[Modified version of Theorem 4.1 in {\cite{dumitriu2024extreme}}] \label{lem:upper_bound_rhoB}
Let $\rmX$ be an $n\times m$ complex matrix, where the entries $\{\ermX_{jl}\}_{j\in [n], l \in [m]}$ are independent mean-zero random variables. Consider the $(n+m) \times (n+m)$ random Hermitian matrix
$\rmH = \begin{bmatrix}
    \bzero & \rmX \\
    \rmX^{\star} & \bzero
\end{bmatrix}$.
Let $\rmB$ denote the non-backtracking matrix associated with $\rmH$. Suppose $\ratio = n/m \geq 1$. Define $\rmH(\lambda)$ and $\rmM(\lambda)$ through $\rmH$ as in \Cref{lem:Ihara-Bass}. Let $\rho(\rmB)$ denote the spectral radius of $\rmB$. Moreover, assume there exists some $q>0$ and $t \geq 1$ such that
\begin{align}
    \max_{l\in [m]} \sum_{j\in [n]} \E |\ermX_{jl}|^2 \leq \qr^{2}, \quad \max_{j\in [n]} \sum_{l\in [m]} \E |\ermX_{jl}|^2 \leq 1/\qr^{2}, \quad \max_{j, l} \E  |\ermX_{jl}|^2 \leq \frac{t}{n}, \quad \max_{j, l} |\ermX_{jl}| \leq \frac{1}{q} \quad \textnormal{a.s.}
\end{align}
Then for $\epsilon \geq 0$ and $1 \vee q \leq (n + m)^{1/10} t^{-1/9}$, there exists some universal constants $c>0$ such that
\begin{align}
    \P(\rho(\rmB) \geq  1 + \epsilon ) \lesssim n^{3 - cq \log(1 + \epsilon)},
\end{align}
where $\lesssim$ only hides some universal constant.
\end{lemma}

With the proof deferred to \Cref{sec:Loewner-order-app}, \Cref{lem:invertiblility_D1} shows that with high probability, $\id_{n} - d^{-1} \rmD^{(1)}$ is positive definite.
\begin{lemma}\label{lem:invertiblility_D1}
 Under \Cref{ass:Ihara_Bass}, there exist a constant $\epsilon>0$ and $\nu>0$ such that the following holds with probability at least $1-N^{-\nu}$
\begin{align}
    \|(\id_{n} - d^{-1} \rmD^{(1)})^{-1}\|\leq \epsilon^{-1}.
\end{align}
\end{lemma} 

With all the ingredients above, we now present the proof of \Cref{prop:lower_bound_H}. To that end, we make use of the imaginary eigenvalues of the non-backtracking operator, inspired by \cite{brito2022spectral, dumitriu2024extreme}, and conduct a more detailed analysis to establish the Loewner order.

\begin{proof}[Proof of \Cref{prop:lower_bound_H}]
Recall $\scA$ in \eqref{eqn:scAdj} where $\rmX = (\widetilde{\rmA} - \E\widetilde{\rmA})/\sqrt{d}$. The matrices $\rmX(\lambda)$, $\rmM^{(1)}(\lambda)$ and $\rmM^{(2)}(\lambda)$ are defined in \eqref{eqn:XM_entries} through $\rmX \in \R^{n \times m}$. Define 
\begin{subequations}
    \begin{align}
        \rmH^{(2)} \coloneqq &\, {\rmX}^{*}(\rmI - d^{-1}\rmD^{(1)})^{-1} {\rmX}\label{eqn:H2}\\
        \rmH^{(2)}(\lambda)\coloneqq &\, \rmX^{*}(\lambda) [\rmM^{(1)}(\lambda)]^{-1}\rmX(\lambda). \label{eqn:H2lambda}
    \end{align}
\end{subequations}
Let $\lambda = \ii \theta$ with $\theta \in \R$, then $\rmM^{(1)}(\lambda)$ and $\rmM^{(2)}(\lambda)$ are real diagonal matrices, and $\rmH^{(2)}(\lambda)$ is Hermitian. Note that as $\theta \to \infty$, $\rmM^{(2)}(\lambda) - \rmH^{(2)}(\lambda) = \id_{m} + O(\theta^{-2})$. Hence, $\rmM^{(2)}(\lambda) - \rmH^{(2)}(\lambda)$ is strictly positive definite for sufficiently large $\theta$. Define $\theta_{\star}$ by
 \begin{align}
    \theta_{\star} \coloneqq \inf\{t >0: \lambda = \ii \theta, \rmM^{(2)}(\lambda) - \rmH^{(2)}(\lambda) \succ \bzero \textnormal{ for all } \theta > t\}.
\end{align}
Let $\lambda_{\star} = \ii \theta_{\star}$, then by continuity, the smallest eigenvalue of $\rmM^{(2)}(\lambda_{\star}) - \rmH^{(2)}(\lambda_{\star})$ is zero. Let $\rmB$ denote the non-backtracking operator associated with $\scA$ in \eqref{eqn:scAdj}. According to \Cref{lem:block_Ihara}, $\lambda \in \mathrm{Spec}(\rmB)$ if and only if $\mathrm{det}\big(\rmM^{(2)}(\lambda) - \rmH^{(2)}(\lambda) \big) = 0$, hence $\lambda_{\star} \in \mathrm{Spec}(\rmB)$ and $|\lambda_{\star}| \leq \rho(\rmB)$ since $\rho(\rmB)=\max_{\lambda \in \mathrm{Spec}(\rmB)}|\lambda|$. Consequently, for any $|\lambda| \geq \rho(\rmB)$, it implies $|\lambda| \geq |\lambda_{\star}|$, hence $\rmM^{(2)}(\lambda) \succeq \rmH^{(2)}(\lambda)$. We confirm that $\scA$ in \eqref{eqn:scAdj} satisfies the assumptions of \Cref{lem:upper_bound_rhoB} with $q = \sqrt{d}$. Choose $\lambda = \ii \theta$ with $\theta = 1 + \epsilon$, one concludes
\begin{align}
    \P\big(\rmM^{(2)}(1 + \epsilon) - \rmH^{(2)}(1 + \epsilon) \succeq \bzero \big) \geq \P(\rho(\rmB) \leq 1 + \epsilon ) \geq 1 - N^{3 - c\sqrt{d} \log(1 + \epsilon)}.
\end{align}    
Then with probability at least $1 - N^{3 - c\sqrt{d} \log(1 + \epsilon)}$, the following holds
\begin{align}
    d^{-1}\rmD^{(2)}- \rmI \preceq &\, -( \rmI -d^{-1}\rmD^{(2)} ) + (-\lambda^{2})\cdot [\rmM^{(2)}(\lambda) - \rmH^{(2)}(\lambda)] - \rmH^{(2)} + \rmH^{(2)}.\\
    = &\, \rmH^{(2)} + \big[(-\lambda^{2}) \cdot \rmM^{(2)} - \rmI + d^{-1}\rmD^{(2)} \big] + \big[ \lambda^{2}\rmH^{(2)}(\lambda) - \rmH^{(2)}\big].
\end{align}
where we choose $\lambda = \ii \theta$ with $\theta = 1 + \epsilon$. The proof of \eqref{eqn:upper_bound_H} then follows by the two inequalities below
\begin{align}
    \|(-\lambda^{2}) \cdot \rmM^{(2)}(\lambda) - \rmI + d^{-1}\rmD^{(2)}\| \lesssim d^{-1}, \quad \|\rmH^{(2)} - \lambda^{2}\rmH^{(2)}(\lambda)\| \lesssim d^{-1/2}. \label{eqn:upper_bound_M2H2}
\end{align}

For the first inequality in \eqref{eqn:upper_bound_M2H2}, note that $\rmM(\lambda)$ and $\rmD$ are diagonal matrices, then by definition
    \begin{align}
        &\, \|(-\lambda^{2}) \cdot \rmM^{(2)}(\lambda) - \rmI + d^{-1}\rmD^{(2)}\| = \max_{l\in \gV_2}\bigg| (-\lambda^{2}) \cdot \ermM^{(2)}_{ll}(\lambda) - 1 - \frac{1}{\lambda d}\sum_{j\in \gV_1}\ermA_{jl}  \bigg|.
    \end{align}
For $l\in \gV_{2}$, using definition in \eqref{eqn:XM_entries} and triangle inequality, we have
\begin{align}
 (-\lambda^{2})\cdot \ermM^{(2)}_{ll}(\lambda) - 1 - \frac{1}{\lambda d}\sum_{j\in \gV_1} \ermA_{jl} \leq &\, \sum_{j\in \gV_1} \frac{\ermX^{2}_{jl}}{\lambda^{2} - \ermX^{2}_{jl}} - \frac{1}{\lambda^2 d}\sum_{j\in \gV_1}\ermA_{jl}^{2} \\
 =&\, \sum_{j\in \gV_1} \bigg( \frac{\ermX_{jl}^4}{\lambda^2 d(\lambda^2 - \ermX^{2}_{jl})} - \frac{2\cdot\ermX_{uv}}{\lambda^2 \sqrt{mn}} - \frac{d}{\lambda^2 mn} \bigg).
\end{align}
Note that $\sqrt{mn} \gtrsim d^{13/2} \asymp \lambda^{13}$, and the last two terms inside the summation are relatively small compared with the first one. The desired result follows since $|\ermX_{jl}| \leq 1/\sqrt{d}$ and $|\ermX_{jl}|^2 \leq \lambda^2/2$.

For the second inequality in \eqref{eqn:upper_bound_M2H2}, by triangle inequality, we have
\begin{align}
         &\,  \|\rmH^{(2)} - \lambda^{2}\rmH^{(2)}(\lambda)\| \\
    \leq &\, \|\rmH^{(2)} - \rmX^{*}[\rmM^{(1)}(\lambda)]^{-1}\rmX \| + \|\rmX^{*}[\rmM^{(1)}(\lambda)]^{-1}\rmX - \lambda^{2}\rmH^{(2)}(\lambda)\| \\
    \leq &\, \|\rmH^{(2)} - \rmX^{*}[\rmM^{(1)}(\lambda)]^{-1}\rmX\| + \|\rmX^{*}[\rmM^{(1)}(\lambda)]^{-1}\rmX - \lambda \rmX^{*}(\lambda)[\rmM^{(1)}(\lambda)]^{-1} \rmX\| \label{eqn:H2_lambdaH2_bound} \\
    &\, \quad \quad \quad \quad \quad \quad \quad \quad \quad \quad \quad \,\,\,\, + \|\lambda \rmX^{*}(\lambda)[\rmM^{(1)}(\lambda)]^{-1} \rmX - \lambda^{2}\rmH^{(2)}(\lambda)\| \\
    \lesssim &\, (\log N)^{ 1/2} \cdot d^{-1} + d^{-1/2} + d^{-1/2} \lesssim d^{-1/2}.
\end{align}
where the last inequality is due to the following bounds  with very high probability, with the proof are deferred later.  We claim that
\begin{align}
    \|\rmX\| \lesssim 1, \quad \|[\rmM^{(1)}(\lambda)]^{-1}\|\lesssim 1, \quad \| (\rmI - d^{-1}\rmD^{(1)})^{-1}\| \lesssim [\log(N)]^{ 1/2}, \quad \|\rmX - \lambda \rmX(\lambda)\| \lesssim d^{-1/2}. \label{eqn:M1_inverse_bounds} 
\end{align}

For the first term in \eqref{eqn:H2_lambdaH2_bound}, we plug in the definition in \eqref{eqn:H2}, then by triangle inequality
\begin{align}
    &\, \|\rmH^{(2)} - \rmX^{*}[\rmM^{(1)}(\lambda)]^{-1}\rmX\| = \|{\rmX}^{*}(\rmI - d^{-1}\rmD^{(1)})^{-1} {\rmX} - \rmX^{*}[\rmM^{(1)}(\lambda)]^{-1}\rmX\|\\
    \leq &\, \|\rmX\|\cdot \Big( \| (\rmI - d^{-1}\rmD^{(1)})^{-1} - [\rmM^{(1)}(\lambda)]^{-1}\|\Big) \cdot \|\rmX\| \quad (\rmA^{-1} - \rmB^{-1} = \rmA^{-1}(\rmB - \rmA)\rmB^{-1} )\\
    \leq &\,  \|\rmX\|\cdot \|(\rmI - d^{-1}\rmD^{(1)})^{-1}\| \cdot \| \rmM^{(1)}(\lambda) - \rmI + d^{-1}\rmD^{(1)}\|\cdot \|[\rmM^{(1)}(\lambda)]^{-1}\| \cdot \|\rmX\|\\
    \lesssim &\, (\log N)^{ 1/2} \cdot d^{-1}.
\end{align}
where $\| \rmM^{(1)}(\lambda) - \rmI + d^{-1}\rmD^{(1)}\| \lesssim d^{-1}$ follows similarly to the proof of the first inequality in \eqref{eqn:upper_bound_M2H2}.

For the second term in \eqref{eqn:H2_lambdaH2_bound}, by triangle inequality
\begin{align}
    &\, \|\rmX^{*}[\rmM^{(1)}(\lambda)]^{-1}\rmX - \lambda \rmX^{*}(\lambda)[\rmM^{(1)}(\lambda)]^{-1} \rmX\| \leq \|\rmX^{*} - \lambda \rmX^{*}(\lambda)\| \cdot \|[\rmM^{(1)}(\lambda)]^{-1}\| \cdot \|\rmX\| \lesssim d^{-1/2}
\end{align}

For the third term in \eqref{eqn:H2_lambdaH2_bound}, we plug in the definition in \eqref{eqn:H2lambda}, then by triangle inequality
\begin{align}
    &\, \|\lambda \rmX^{*}(\lambda)[\rmM^{(1)}(\lambda)]^{-1} \rmX - \lambda^{2}\rmH^{(2)}(\lambda)\| = \|\lambda \rmX^{*}(\lambda)[\rmM^{(1)}(\lambda)]^{-1} \rmX - \lambda^{2}\rmX^{*}(\lambda) [\rmM^{(1)}(\lambda)]^{-1}\rmX(\lambda)\|\\
    \leq & \|\lambda \rmX^{*}(\lambda)\| \cdot \| [\rmM^{(1)}(\lambda)]^{-1}\| \cdot \|\rmX - \lambda \rmX(\lambda)\| \\
    \leq &\, \| [\rmM^{(1)}(\lambda)]^{-1}\| \cdot \|\rmX - \lambda \rmX(\lambda)\| \cdot \big( \|\rmX - \lambda \rmX(\lambda)\| + \|\rmX\| \big)
    \lesssim \, d^{-1/2}.
\end{align}

Now it remains to prove \eqref{eqn:M1_inverse_bounds}.
    First, with very high probability, $\|\rmX\| \leq \const$ for some absolute constant $\const$ follows directly from \Cref{thm:right_edge_behavior}.

    Second, we take $\lambda = \ii \theta$ with $\theta_{\star} \leq \theta = 1 + \epsilon$. Note that $p/\sqrt{d} \leq |\ermX_{jl}| \leq 1/\sqrt{d} \ll \theta$, then with very high probability $\sum_{l\in [m]}|\ermX_{jl}|^{2} = m \cdot \Var(\ermX_{12}) =(1 + o(1)) m p/d = \qr^{-2}$. Then by \eqref{eqn:XM_entries}, with very high probability,  
    \begin{align}
        \ermM^{(1)}_{jj}(\lambda) = 1 - \sum_{l\in [m]}\frac{|\ermX_{jl}|^{2}}{\theta^{2} + |\ermX_{jl}|^{2}} \geq 1 - \frac{1}{\theta^{2}} \cdot \qr^{-2} \geq c >0
    \end{align}
    under \Cref{ass:Ihara_Bass} where $\qr > 1$. Therefore, $\|[\rmM^{(1)}(\lambda)]^{-1}\| \leq \const$ for some absolute constant $\const$. 

    Third, $\| (\rmI - d^{-1}\rmD^{(1)})^{-1}\| \lesssim [\log(N)]^{1/2}$ follows by \Cref{lem:invertiblility_D1}, with very high probability.

    For the last claim in \eqref{eqn:M1_inverse_bounds}, we calculate the entrywise difference of two matrices, using $\lambda = \ii \theta$ with $\theta = 1 + \epsilon$,
    \begin{align}
        | \lambda [\ermX(\lambda)]_{jl} - \ermX_{jl} | = \Big| \frac{\lambda^{2} \ermX_{jl}}{\lambda^{2} - |\ermX_{jl}|^{2}} - \ermX_{jl} \Big| \leq \frac{|\ermX_{jl}|^{3}}{\theta^{2} + |\ermX_{jl}|^{2}} \leq \frac{1}{\theta^{2}\sqrt{d}} |\ermX_{jl}|^{2},
    \end{align}
    where the last inequality holds since $p/\sqrt{d} \leq |\ermX_{jl}| \leq 1/\sqrt{d} \ll \theta$. 
    
    Note that $\sum_{j\in \gV_{1}}\E |\ermX_{jl}|^{2} \lesssim 1$, $\sum_{l\in \gV_{2}}\E |\ermX_{jl}|^{2} \lesssim 1$. The desired upper bound then follows by the Schur test in \Cref{lem:schur_test} and Chebyshev's inequality.
\end{proof}

\subsection{Eigenvector delocalization estimate}\label{sec:delocalization}
Consider the triple $\vtau = (\tau_{1}, \tau_{2}^{+}, \tau_{2}^{-})$ and let
\begin{align} \label{eqn:tau_constrain}
    \tau_{1} = \qr^{-2} + \frac{1}{4\qr}\error^{1/4}, \quad \tau_{2}^{+} = \qr^{2} + \frac{\qr}{4}\error^{1/4}, \quad \tau_{2}^{-} = \qr^{2} - \frac{\qr}{4}\error^{1/4},
\end{align}
where $\error$ denotes the error control parameter in \eqref{eqn:error_parameters}. We consider a subset of $\pruneV$ in \eqref{eqn:pruneV}, defined as
    \begin{align}\label{def:W_set}
  \gV^{(\vtau^*)}\subset\pruneV,\text{ with } \vtau^{\star} \coloneq (1+\qr^{-2}+\gC\error^{1/8},\, 1 +\qr^{2}+\gC\error^{1/8},\, \qr^{2}-1 - \gC \error^{1/8}),
\end{align}
where $\gC>0$ is some constant. 

Recall the approximate eigenvectors $\rvv_{+}(x)$ and $\rvv_{-}(x)$ in \eqref{eqn:vx_plus_minus}, and the pruned graph $\gG^{(\vtau)}$ in \Cref{lem:existence_pruned_graph}. For $\sigma = \pm 1$, we construct the following approximate eigenvector on $\gG^{(\vtau)}$
\begin{align}
    \hatvtau_{\sigma}(x)\coloneqq  \sum_{j=0}^{\radius}\sigma^{j} \ervu_{j}(x)\frac{1}{|\layer_{j}^{\gG^{(\vtau)}} (x)|^{1/2}} \,\, \ones_{\layer_{j}^{\gG^{(\vtau)}}(x)}.
\end{align}
We introduce the block diagonal approximation of the pruned graph below.
\begin{definition}
    Define the following two orthogonal projection matrices
    \begin{align}
        \Pi^{(\vtau)} \coloneqq&\,  \sum_{x\in \gV^{(\vtau^*)}}\sum_{\sigma = \pm } \hatvtau_{\sigma}(x) [\hatvtau_{\sigma}(x)]^{*}, \quad \quad \quad \quad \overline{\Pi^{(\vtau)}} \coloneqq \id_{N} - \Pi^{(\vtau)}
    \end{align}
For $\rmH^{(\vtau)}$ in \eqref{eqn:def_H_tau}, define its block diagonal approximation as
\begin{align}
    \hatH \coloneqq &\, \sum_{x\in \gV^{(\vtau^*)}}\sum_{\sigma = \pm } \sigma \cdot (\Lambda_\qr(\alpha_x)\indi{x\in\gV_2}+\Lambda_{\qr^{-1}}(\alpha_x)\indi{x\in\gV_1})  \hatvtau_{\sigma}(x) [\hatvtau_{\sigma}(x)]^{*}\\
    &\, + \overline{\Pi^{(\vtau)}}\,\, \scA^{(\vtau)}\,\,\overline{\Pi^{(\vtau)}}. \label{eqn:hatH}
\end{align}
\end{definition}

Recall the notation in \eqref{eqn:Asubsets}, where $\rmM|_{\setS}$ denotes a submatrix of $\rmM$ with both ends in $\setS \subset [N]$.

For $x\in \pruneV\setminus \gV^{(\vtau^*)}$, denote $\hatHx \coloneqq \hatH |_{\ball^{(\vtau)}_{2\radius} (x)}$. On the other hand, for $x\in \gV^{(\vtau^*)}$, we consider 
\begin{align}
    \hatHx \coloneqq \hatH |_{\ball^{(\vtau)}_{2\radius} (x)} - \sum_{\sigma = \pm } \sigma \cdot \big( \Lambda_\qr(\alpha_x)\indi{x\in\gV_2}+\Lambda_{\qr^{-1}}(\alpha_x)\indi{x\in\gV_1} \big)\, \hatvtau_{\sigma}(x) [\hatvtau_{\sigma}(x)]^{*} \label{eqn:hatHx}
\end{align}
Here, $\hatHx$ is the restriction of $\hatH$ on the ball $\ball^{(\vtau)}_{2\radius} (x)$, eliminating the components corresponding to the outlier eigenvalues.

We now establish the delocalization estimate for eigenvectors.

\begin{proposition}\label{prop:delocalization_vec}
    Let $\vtau$ and $\vtau^{\star}$  satisfy \eqref{eqn:tau_constrain} and \eqref{def:W_set}. Let $(\lambda, \rvw)$ be an eigenpair of $\hatH$ defined in \eqref{eqn:hatH}, where $\|\rvw\| = 1$ and $\rvw \perp \hatvtau_{\pm}(x)$ for all $x\in \gV^{(\vtau^{\star})}$. The following holds with very high probability.
     \begin{enumerate}
        \item  If $\lambda > \qr +\qr^{-1}  +  \frac{1}{2} \error^{1/2}$, then 
        \begin{align}
            \sum_{y\in \pruneV}\ervw_{y}^{2} \leq \Big(\frac{\lambda }{\lambda - \qr -\qr^{-1} -  \frac{1}{2}\error^{1/2}}\Big)^{4} \Big( \frac{\qr +\qr^{-1}  + \frac{1}{2}\error^{1/2}}{\lambda} \Big)^{2\radius}. \label{eqn:right_delocalization}
        \end{align}
        \item If $0 < \lambda < \qr - \qr^{-1} - \frac{1}{2} \error^{1/2}$, then
        \begin{align}
            \sum_{y\in \pruneV}\ervw_{y}^{2} \leq \Big(\frac{ \qr - \qr^{-1} -  \frac{1}{2}\error^{1/2} }{ \qr - \qr^{-1} - \frac{1}{2}\error^{1/2} - \lambda} \Big)^4\Big( \frac{\lambda}{\qr - \qr^{-1} -  \frac{1}{2}\error^{1/2}} \Big)^{2\radius}. \label{eqn:left_delocalization}
        \end{align}
     \end{enumerate}
     Here $\error$ is defined in \eqref{eqn:error_parameters}, and $\radius$ is defined in \eqref{eqn:radius}. 
\end{proposition}
The proof of (1) follows the same strategy as {\cite[Proposition 3.14]{alt2021delocalization}}. To ensure that the paper is self-contained, we present details below, which will also enlighten the proof of Proposition~\ref{prop:delocalization_vec} (2).

\begin{proof}[Proof of Proposition~\ref{prop:delocalization_vec} (1)]
To prove \eqref{eqn:right_delocalization}, it suffices to show that if $x\in \gV^{(\vtau^{\star})}$ and $\rvw \perp \hatvtau_{\pm}(x)$ or if $x\in \pruneV \setminus \gV^{(\vtau^{\star})}$, the following holds
\begin{align}
    \frac{|\ervw_{x}|}{\|\rvw|_{\ball_{2 \radius}^{(\vtau)}(x)}\|} \leq \frac{\lambda^{2}}{(\lambda - \qr -\qr^{-1}  - \frac{1}{2} \error^{1/2})^{2}} \Big( \frac{\qr + \qr^{-1} + \frac{1}{2}\error^{1/2}}{\lambda} \Big)^{\radius}. \label{eqn:wx_ratio_right}
\end{align}
 This is because that with \eqref{eqn:wx_ratio_right}, we have  
\begin{align}
    \sum_{y\in \pruneV} \ervw_{y}^{2} \leq &\, \sum_{y\in \pruneV} \|\rvw_{\ball^{(\vtau)}_{2\radius}(y)}\|^{2} \frac{\lambda^{4}}{(\lambda - \qr -\qr^{-1}  -  \frac{1}{2}\error^{1/2})^{4}} \Big( \frac{\qr+\qr^{-1}  +  \frac{1}{2}\error^{1/2}}{\lambda} \Big)^{2\radius}\\
    \leq &\, \frac{\lambda^{4}}{(\lambda - \qr - \qr^{-1} - \frac{1}{2} \error^{1/2})^{4}} \Big( \frac{\qr + \qr^{-1} + \frac{1}{2}\error^{1/2}}{\lambda} \Big)^{2\radius}.
\end{align}
where we apply $\sum_{y\in \pruneV} \|\rvw_{\ball^{(\vtau)}_{2\radius}(y)}\|^{2} \leq \|\rvw\|^{2} = 1$ in the last step, since the balls $\{\ball^{(\vtau)}_{2\radius}(x): x\in \pruneV\}$ are disjoint according to \eqref{item:one} of \Cref{lem:existence_pruned_graph}.

The proof of \eqref{eqn:wx_ratio_right} consists of three steps:
    \begin{enumerate}
        \item There exists some constant $C>0$ such that with very high probability,
        \begin{equation}
            \|\hatHx\|\leq  \qr + \qr^{-1} + \frac{1}{2}\error^{1/2}.\label{eqn:upperbound_Hhatx}
        \end{equation}
        \item Let $\{\rvg_{j}\}_{j=0}^{\radius}$ be the Gram-Schmidt orthonormalization of the vectors $\{[\hatHx]^{j} \ones_{x}\}_{j=0}^{\radius}$. The following holds for $j = 0, 1,\ldots, \radius$.
        \begin{align}
            \mathrm{supp}\,(\rvg_{j}) \subset \ball^{(\vtau)}_{\radius + j}(x). \label{eqn:supp_in_ball}
        \end{align}
        \item Conduct the entrywise analysis on matrix resolvent to conclude \eqref{eqn:wx_ratio_right}.
    \end{enumerate}

For simplicity, the proofs of \eqref{eqn:upperbound_Hhatx} and \eqref{eqn:supp_in_ball} are deferred to \Cref{sec:delocalization_app}. Let $\rmZ = [\ermZ_{jl}]_{j,l=0}^{\radius}$ be the tridiagonal representation of $\hatHx$ up to radius $\radius$, where $\ermZ_{jl}\coloneqq \<\rvg_{j}, \hatHx \rvg_{l}\>$. We set $\ervu_{j} = \<\rvg_{j}, \rvw\>$ for any $0\leq j \leq \radius$. From \eqref{eqn:upperbound_Hhatx},
\begin{align}
   \|\rmZ\| \leq \|\hatHx\| \leq \qr + \qr^{-1} + \error^{1/2} < \lambda. \label{eqn:Z_right_edge_bound}
\end{align}
Note that $(\lambda, \rvw)$ is an eigenpair of $\hatH$, and $\rvw \perp \hatvtau_{\pm}(x)$. Then for any $0 \leq j \leq \radius$, \eqref{eqn:supp_in_ball} implies
\begin{align}
    \lambda \ervu_{j} = &\, \<\rvg_{j}, \hatH \rvw\> = \<\rvg_{j}, \big( \hatH - \sum_{\sigma = \pm} \sigma \Lambda_{\qr}(\alpha_{x}) \hatvtau_{\sigma}(x) [\hatvtau_{\sigma}(x)]^{*} \big)\rvw\>\\
    =&\, \<\hatHx \rvg_{j}, \rvw\> = \<\ermZ_{j,j} \rvg_{j} + \ermZ_{j,j+1} \rvg_{j+1} + \ermZ_{j,j-1} \rvg_{j-1}, \rvw\>\\
    =&\,\ermZ_{j,j} \ervu_{j} + \ermZ_{j,j+1} \ervu_{j+1} + \ermZ_{j,j-1} \ervu_{j-1},
\end{align}
with the conventions $\ervu_{-1} =0$ and $\ermZ_{0,-1} = 0$. According to \eqref{eqn:Z_right_edge_bound}, $\lambda - \rmZ$ is invertible. Let $\resolvent(\lambda) \coloneqq (\lambda - \rmZ)^{-1}$ be the resolvent of $\rmZ$ at $\lambda$. Using the fact $[(\lambda - \rmZ)\resolvent(\lambda) ]_{j, \radius} = 0$ for $j < \radius$, we find
\begin{align}
    \lambda \eres_{j,\radius}(\lambda) = \ermZ_{j,j} \eres_{j,\radius}(\lambda) + \ermZ_{j,j+1} \eres_{j+1,\radius}(\lambda) + \ermZ_{j,j-1} \eres_{j-1,\radius}(\lambda).
\end{align}
Thus $\{\eres_{j,\radius}(\lambda)\}_{j=0}^{\radius}$ and $\{\ervu_{j}\}_{j=0}^{\radius}$ satisfy the same linear recursive equation when the elements of $\rmZ$ are viewed as coefficients. By solving them recursively, it yields for all $j\leq \radius$,
\begin{align}
    \frac{\eres_{j,\radius}(\lambda)}{\eres_{\radius,\radius}(\lambda)} = \frac{\ervu_{j}}{\ervu_{\radius}}. \label{eqn:equivalence_G_u}
\end{align}

We consider the Neumann series $\resolvent(\lambda) = \lambda^{-1}\sum_{k=0}^{\infty}(\rmZ/\lambda)^{k}$, which converges since $\lambda > \|\rmZ\|$ by \eqref{eqn:Z_right_edge_bound}. For off-diagonal entries of the resolvent $\resolvent(\lambda)$, it admits the form
\begin{align}
    \eres_{0, \radius}(\lambda) = \frac{1}{\lambda} \sum_{k=0}^{\infty} \Big[(\rmZ/\lambda)^{k} \Big]_{0, \radius}.
\end{align}
Note that $\ermZ_{0, \radius} = \<\ones_{0}, \rmZ \ones_{\radius}\> \leq \max_{\|\rvx\| = 1}\<\rvx, \rmZ \rvx \> =\|\rmZ\|$, and $[(\rmZ/\lambda)^{k}]_{0, \radius} = 0$ for $k < \radius$ since $\rmZ$ is tridiagonal, then the following holds
\begin{align}
    |\eres_{0, \radius}(\lambda)| = \frac{1}{\lambda} \sum_{k=\radius}^{\infty} \Big[(\rmZ/\lambda)^{k} \Big]_{0, \radius} \leq \frac{1}{\lambda} \cdot \sum_{k=\radius}^{\infty} \bigg(\frac{\|\rmZ\|}{\lambda}\bigg)^{k} = \frac{1}{\lambda - \|\rmZ\|} \bigg(\frac{\|\rmZ\|}{\lambda}\bigg)^{\radius}.
\end{align}

For diagonal entries, by splitting the summation over $k$ into even and odd values,
\begin{align}
    |\eres_{\radius, \radius}(\lambda)| =&\, \frac{1}{\lambda} \sum_{k=0}^{\infty} \Big[(\rmZ/\lambda)^{k} \Big]_{\radius, \radius} = \frac{1}{\lambda} \sum_{k=0}^{\infty} \Big[(\rmZ/\lambda)^{k}(\id + \rmZ/\lambda)(\rmZ/\lambda)^{k} \Big]_{\radius, \radius}\\
    \geq &\, \frac{1}{\lambda} (\id + \rmZ/\lambda)_{\radius, \radius} \geq \frac{1}{\lambda} \bigg(1 - \frac{\|\rmZ\|}{\lambda}\bigg),
\end{align}
where in the second line, we only keep the term $k=0$ to obtain a lower bound since $\id + \rmZ/\lambda \succeq \bzero$ due to $\lambda > \|\rmZ\|$. By definition of $\ervu_{j} = \<\rvg_{j}, \rvw\>$, then \eqref{eqn:equivalence_G_u} implies
\begin{align}
    \frac{|\ervw_{x}|}{\|\rvw_{\ball^{(\vtau)}_{2\radius}(x)}\|} \leq \frac{|\ervu_{0}|}{(\sum_{j=0}^{\radius}\ervu_{j}^{2})^{1/2}} \leq \frac{|\ervu_{0}|}{|\ervu_{\radius}|} =  \frac{|\eres_{j,\radius}(\lambda)|}{|\eres_{\radius,\radius}(\lambda)|} \leq \bigg(\frac{\lambda}{\lambda - \|\rmZ\|} \bigg)^{2}\bigg(\frac{\|\rmZ\|}{\lambda}\bigg)^{\radius}.
\end{align}
The proof of \eqref{eqn:wx_ratio_right} for $x\in \gV^{(\vtau^{\star})}$ is concluded by filling in the upper bound of $\|\rmZ\|$ in \eqref{eqn:Z_right_edge_bound}. The case for $x\in \pruneV \setminus \gV^{(\vtau^{\star})}$ follows similarly by verifying \eqref{eqn:Z_right_edge_bound}. 
\end{proof}

\begin{proof}[Proof of \Cref{prop:delocalization_vec} (2)]
Similarly as in the proof of \Cref{prop:delocalization_vec} (1), to prove \eqref{eqn:left_delocalization}, it suffices to show that if $x\in \gV^{(\vtau^{\star})}$ and $\rvw \perp \hatvtau_{\pm}(x)$ or if $x\in \pruneV \setminus \gV^{(\vtau^{\star})}$, the following holds
    \begin{align}
        \frac{|\ervw_{x}|}{\|\rvw|_{\ball_{\radius}^{(\vtau)}(x)}\|} \leq \Big(\frac{ \qr - \qr^{-1} - \frac{1}{2}\error^{1/2} }{ \qr - \qr^{-1} - \frac{1}{2}\error^{1/2} - \lambda} \Big)^2\Big( \frac{\lambda}{\qr - \qr^{-1} - \frac{1}{2}\error^{1/2}} \Big)^{ \radius}. \label{eqn:wx_ratio_left}
    \end{align}

Let $\{\widehat{\rvg}_{j}\}_{j=\radius}^{0}$ be the Gram-Schmidt orthonormalization of the vectors $\{[\hatHx]^{-j} \ones_{\layer_{\radius}^{(\vtau)}}\}_{j=0}^{\radius}$, where we note that $\hatHx$ in \eqref{eqn:hatHx} and $\rmA^{(\vtau)}|_{\ball^{(\vtau)}_{2\radius} (x)}$ are invertible; see \Cref{sec:tri_regular}. Due to \eqref{eqn:AoneSl}, for $x\in \gV_2$ and $l\in [\radius -1]\setminus \{1\}$, we have
\begin{align}
    \ones_{\layer_{l}^{(\vtau)}(x)} = [\rmA^{(\vtau)}|_{\ball^{(\vtau)}_{2\radius}}]^{-1} \Big(\ones_{\layer_{l+1}^{(\vtau)}(x)} + \indi{ l \textnormal{ even} } d_1\ones_{\layer_{l-1}^{(\vtau)}(x)} + \indi{ l \textnormal{ odd} }d_2\ones_{\layer_{l-1}^{(\vtau)}(x)}\Big),
\end{align}
in the ball $\ball^{(\vtau)}_{2\radius} (x)$. Thus $\{[\rmA^{(\vtau)}|_{\ball^{(\vtau)}_{2\radius} (x)}]^{-j}\}_{j=0}^{\radius}$ is still a linear combination of $\{\ones_{\layer_{j}^{(\vtau)}}\}_{j=0}^{\radius}$. One can prove by induction that \eqref{eqn:supp_in_ball} holds for $\{\widehat{\rvg}_{j}\}_{j=\radius}^{0}$.

Let $\rmY = [\ermY_{jl}]_{j,l=0}^{\radius}$, where $\ermY_{jl}\coloneqq \<\rvg_{j}, [\hatHx]^{-1} \rvg_{l}\>$, be the tridiagonal representation of $[\hatHx]^{-1}$ up to radius $\radius$. With the proof deferred to \Cref{sec:delocalization_app}, we claim that with very high probability
    \begin{align}
        \|\rmY\| \leq \|[\hatHx]^{-1}\| \leq \frac{1}{\qr - \qr^{-1} - \frac{1}{2}\error^{1/2}} < \frac{1}{\lambda}. \label{eqn:Y_left_edge_bound}
    \end{align}
    
Then, we can apply the third step in the proof of \Cref{prop:delocalization_vec} (1) by considering $\resolvent(\lambda) = (\frac{1}{\lambda}-\rmY)^{-1}$. Note that in the proof of \eqref{eqn:wx_ratio_left}, \eqref{eqn:Y_left_edge_bound} plays the same role as \eqref{eqn:Z_right_edge_bound} in the proof of \eqref{eqn:wx_ratio_right}. By leveraging the fact that $(\lambda^{-1}, \rvw)$ is an eigenpair of $\rmM^{-1}$ iff $(\lambda, \rvw)$ is an eigenpair of $\rmM$, and following the resolvent analysis in \eqref{eqn:equivalence_G_u}, \eqref{eqn:wx_ratio_left} can be established similarly, which subsequently completes the proof of \eqref{eqn:left_delocalization}.
\end{proof}

\begin{proof}[Proof of \eqref{eqn:Y_left_edge_bound}]
Due to the symmetry of the spectrum, we focus on the positive eigenvalues. It suffices to prove that the smallest positive eigenvalue of $\hatHx$ is no smaller than $\qr-\qr^{-1}- \frac{1}{2}\error^{1/2}$ with very high probability. We only need to consider the case $x\in \Vtwolow$. 
 
Let $\widehat{\rmZ}$ be the tridiagonal representation of $\hatH\big |_{\ball^{(\vtau)}_{2\radius} (x)}$ with $\hatH$ defined in \eqref{eqn:hatH}, and $\rmZ(\alpha_x)$ in \eqref{eqn:tridiag_Z_outline} denote the tridiagonal representation of the idealized tree $\tree$ in \Cref{def:biregular_tree}. For simplicity, both matrices are of size $(\radius + 1)\times (\radius + 1)$. According to \Cref{lem:eigenZ1Z2}, when $0 < \alpha_{x} < \qr^{2} - 1$, the smallest positive eigenvalue of $\rmZ(\alpha_x)$ converges to $\Lambda_{\qr }(\alpha_x)$, while the remaining eigenvalues are in the interval $[\qr - \qr^{-1}, \qr + \qr^{-1}]$.

Note that component corresponding to $\Lambda_{\qr }(\alpha_x)$ is eliminated during the construction of $\hatHx$ in \eqref{eqn:hatHx}. By applying the Weyl inequality (\Cref{lem:weyl}), we have
\begin{align}
    \lambda_{\min}(\hatHx) =&\, \lambda_{\min} \Big(\hatHx + \Lambda_{\qr }(\alpha_x) \hatvtau_{\sigma}(x) [\hatvtau_{\sigma}(x)]^{*} - \Lambda_{\qr }(\alpha_x) \hatvtau_{\sigma}(x) [\hatvtau_{\sigma}(x)]^{*}\Big)\\
    =&\,  \lambda_{\min} \Big( \hatH |_{\ball^{(\vtau)}_{2\radius} (x)} - \rmZ(\alpha_x) + \rmZ(\alpha_x) - \Lambda_{\qr }(\alpha_x) \hatvtau_{\sigma}(x) [\hatvtau_{\sigma}(x)]^{*} \Big)\\
    \geq &\, \lambda_{\min} \Big(\rmZ(\alpha_x) - \Lambda_{\qr }(\alpha_x) \hatvtau_{\sigma}(x) [\hatvtau_{\sigma}(x)]^{*} \Big) - \|\widehat{\rmZ} - \rmZ(\alpha_x)\| \\
    \geq &\, \qr-\qr^{-1} - \frac{1}{2} \error^{1/2}
\end{align}
where $\|\widehat{\rmZ} - \rmZ(\alpha_x)\| \leq \frac{1}{2} \error^{1/2}$ is proved in \Cref{lem:Z-approx}.
\end{proof}

 \begin{lemma}\label{lem:Z-approx}
Consider the tridiagonal representation of $\hatH |_{\ball^{(\vtau)}_{2\radius} (x)}$ with $\hatH$ defined in \eqref{eqn:hatH}, where $x\in \pruneV$. Let $\widehat{\rmZ}$ denote its upper left $(\radius + 1)\times (\radius + 1)$ principal submatrix. With very high probability,
 \begin{equation}
    \|\widehat{\rmZ}-\rmZ(\alpha_x)\|\leq \frac{1}{2} \error^{1/2},
 \end{equation}
 where $\rmZ(\alpha_x)$ is defined in in \eqref{eqn:tridiag_Z_outline}.
\end{lemma}

\begin{proof}[Proof of \Cref{lem:Z-approx}]
From \Cref{lem:existence_pruned_graph}, $\gG^{(\vtau)} |_{\ball^{(\vtau)}_{2\radius}(x)}$ is a tree whose root is $x$ which has $\alpha_x d$ children. Let $\layer_{j}^{(\vtau)}(x)$ denote the $j$-th layer in $\gG^{(\vtau)}$. Below, we first show that 
  \begin{subequations}
    \begin{align}
    & \bigg|\frac{1}{d }\cdot\frac{|\layer^{(\vtau)}_{2j}|}{|\layer^{(\vtau)}_{2j - 1}|}-\qr^{-2}  \bigg| \lesssim \frac{1}{\sqrt{d}},    \label{eqn:layer_concentrate1}\\
        &\bigg| \frac{1}{d }\cdot \frac{|\layer^{(\vtau)}_{2j+1}|}{|\layer^{(\vtau)}_{2j}|} - \qr^2\bigg|  \lesssim \frac{1}{\sqrt{d}}  \label{eqn:layer_concentrate2}
     \end{align}
  \end{subequations}
  for every \(1\le j\leq  \radius\). Due to \Cref{lem:existence_pruned_graph} (5) (6) and the choices of $\vtau$ in \eqref{eqn:tau_constrain}, with very high probability, the following holds
\begin{align}
\frac{|\layer^{(\vtau)}_{j}|}{|\layer_j^{\gG}|}
= 1 - \frac{|\layer_j^{\gG}\setminus \layer^{(\vtau)}_{j}|}{|\layer_j^{\gG}|}
\geq &\, 1 - O\bigg(\frac{1}{d}\Big(1 + \frac{\log(N)}{d} \cdot \Big( \min\Big\{\benrate_{\qr^{-1}}(\tau_{1}),\,\, \benrate_{\qr}(\tau^{-}_{2}),\,\,\benrate_{\qr}(\tau^{+}_{2}) \Big\}\Big)^{-1}\Big)\bigg)\\
= &\, 1 - O\bigg(\frac{1}{d}\Big(1 + \frac{\log(N)}{d} \error^{-1/2} \Big) \bigg) = 1 - O\bigg(\frac{1}{(d^{3}\, \log(d))^{1/4}}\bigg),
\end{align}
where in the last step, we use the facts that $\error = \sqrt{\frac{\log d}{d}}$ in \eqref{eqn:error_parameters} and $\benrate (\alpha) \asymp \alpha^{2}$ for $\alpha = o(1)$ with $\benrate (\alpha)$ defined in \eqref{eqn:bennett_rate}. By applying \eqref{eqn:S2j} \eqref{eqn:S2j+1} in \Cref{lem:concentrationSi}, we then conclude \eqref{eqn:layer_concentrate2} since for $1\leq j \leq  \radius$,
\begin{equation} 
\frac{|\layer^{(\vtau)}_{2j+1} |}{d\,|\layer^{(\vtau)}_{2j}|}
= \frac{|\layer^{\gG}_{2j+1}|}{d\,|\layer^{\gG}_{2j}|}\frac{|\layer^{\gG}_{2j}|}{|\layer^{(\vtau)}_{2j}|}\frac{|\layer^{(\vtau)}_{2j+1} |}{|\layer^{\gG}_{2j+1}|}
=  \qr^2 + O\bigg(\frac{1}{\sqrt{d}}\bigg).
\end{equation}
The proof of \eqref{eqn:layer_concentrate1} follows similarly.

Let $\Zhat$ be the tridiagonal representation of $\hatHx$. As proved in \Cref{sec:tri_regular}, the orthogonal basis obtained after Gram-Schmidt can be written as $\rvg_j=|\layer^{(\vtau)}_{j}|^{-1/2}\ones_{\layer^{(\vtau)}_{j}}$ for \(0\le j\leq \radius\). Since $\hatHx$ corresponds to a tree, we have $\widehat{\ermZ}_{jj}=0$. Since $x\in \pruneV$,  according to \Cref{lem:existence_pruned_graph} (5) and choices of $\vtau$ in \eqref{eqn:tau_constrain}, we have
\begin{align}
    \widehat{\ermZ}_{0,1} =&\, \sqrt{\alpha_x} + O\bigg(\frac{1}{\sqrt{d}}\Big(1 + \frac{\log(N)}{d} \cdot \Big( \min\Big\{\benrate_{\qr^{-1}}(\tau_{1}),\,\, \benrate_{\qr}(\tau^{-}_{2}),\,\,\benrate_{\qr}(\tau^{+}_{2}) \Big\}\Big)^{-1}\Big)\bigg)\\
    = &\, \sqrt{\alpha_x} + O\bigg(\frac{1}{\sqrt{d}}\Big(1 + \frac{\log(N)}{d} \error^{-1/2} \Big) \bigg) = \sqrt{\alpha_x} + O\bigg(\frac{1}{[d\, \log(d)]^{1/4}}\bigg),
\end{align}
where we again use the fact $\benrate (\alpha) \asymp \alpha^{2}$ for $\alpha = o(1)$. Furthermore, due to \eqref{eqn:layer_concentrate1}, \eqref{eqn:layer_concentrate2},
\begin{align}
   \widehat{\ermZ}_{j,j+1}
      \;=\;\langle \rvg_j,\hatHx \rvg_{j+1}\rangle
      \;=\;\sqrt{\frac{|\layer^{(\vtau)}_{j+1}|}{d|\layer^{(\vtau)}_j|}} = \qr \indi{j \textnormal{ even}} + \qr^{-1} \indi{j \textnormal{ odd}} + O(d^{-1/4}), \,\, 1\leq j \leq \radius.
\end{align}
Recall $\rmZ(\alpha_x)$ be the matrix in \eqref{eqn:tridiag_Z_outline}. We conclude the proof using the entrywise estimates above, since
\begin{align}
    \|\widehat{\rmZ} - \rmZ(\alpha_x)\| \leq &\,
   2\,\max_{0\leq j<\radius}\bigl|[\widehat{\rmZ} - \rmZ(\alpha_x)]_{j,j+1}\bigr|\lesssim [d\, \log(d)]^{-1/4} + d^{-1/4} \ll \frac{1}{2} \error^{1/2}. 
\end{align}
\end{proof}

\subsection{Proof of \Cref{thm:right_edge_behavior} (2)}\label{sec:proof_bulk_bound_maintext}
We first present the following estimates.
 \begin{lemma}\label{lem:hatHtau_norm}
    Assume \eqref{eqn:tau_constrain}. For some constant $\const >0$, with very high probability
    \begin{align}
       \|\scA^{(\vtau)}-\rmH\|\le \frac{1}{16} \error^{1/2}, \quad \|\hatH -\scA^{(\vtau)}\|\leq \const \error.
    \end{align}
\end{lemma}
The proof of \Cref{lem:hatHtau_norm} is deferred to \Cref{sec:proof_bulk_bound_app}. \Cref{lem:project_norm} presents an operator norm upper bound of the projected matrix defined in \eqref{eqn:hatH}.

\begin{lemma}\label{lem:project_norm}
    Assume that $\vtau$ satisfies \eqref{eqn:tau_constrain}. Then with very high probability,
    \[ 
    \Big\|\overline{\Pi^{(\vtau)}} \scA^{(\vtau)}\overline{\Pi^{(\vtau)}}\Big\|
   \le  \qr  +\qr^{-1} + \frac{1}{2}\error^{1/4}.
    \]
\end{lemma}
\begin{proof}[Proof of \Cref{lem:project_norm}]
Assume that there exists an eigenpair $(\mu,\rvw)$ of
$\overline{\Pi^{(\vtau)}} \scA^{(\vtau)}\overline{\Pi^{(\vtau)}}$ with
\begin{align}
    \mu>\qr +\qr^{-1} +\frac{1}{2}\error^{1/4}. \label{eq:lem612_contradiction}
\end{align}
We will derive a contradiction below using \Cref{prop:upper_bound_H} and \Cref{prop:delocalization_vec} (1).

First, from Proposition~\ref{prop:upper_bound_H} and Lemma~\ref{lem:hatHtau_norm}, we have
\begin{align}
\scA^{(\vtau)}\preceq~&  \id_{N} + d^{-1}\rmD + \rmE+ \|\rmH-\scA^{(\vtau)} \|\\
\preceq~& \id_N+d^{-1}\rmD+\frac{1}{8}\cdot\qr\error^{1/2},\label{eqn:hatH_psd_order}
\end{align}with very high probability, where we also apply Lemma~\ref{lem:deg_bound} for $\rmE$ defined by \eqref{eqn:error_upper_bound_norm} such that $\|\rmE\|\ll \error$.

Note that $\scA$ in \eqref{eqn:scAdj} and $\scA^{(\vtau)}$ in \eqref{eqn:def_H_tau} both admit a block structure. Let $\rvw^{(1)}=[\ervw_{x}]_{x\in\gV_1}$ and $\rvw^{(2)}=[\ervw_{x}]_{x\in\gV_2}$, and denote 
\begin{align}
    \rvw = \begin{bmatrix}
    \rvw^{(1)}\\
    \rvw^{(2)}
\end{bmatrix}.\label{eqn:w_block_structure}
\end{align}
Furthermore, since $\|\rvw\|_{2}^2 = 1$, by the linear equation $\overline{\Pi^{(\vtau)}} \scA^{(\vtau)}\overline{\Pi^{(\vtau)}} \rvw=\mu \rvw$, we have  
    \begin{align}
         \|\rvw^{(1)}\|^2= \|\rvw^{(2)}\|^2=\frac{1}{2}.\label{eqn:norm_half} 
    \end{align}
  By the definition of $\Pi^{(\vtau)}$, the eigenvector $\rvw\in\R^N$ of $\overline{\Pi^{(\vtau)}} \scA^{(\vtau)}\overline{\Pi^{(\vtau)}}$ corresponding to $\mu$ satisfies $\rvw\perp \hatvtau_{\sigma}(x)$ for all $x\in\gV^{(\vtau^*)}$. Since $\rvw$ is orthogonal to the range of $\Pi^{(\vtau)}$, we know that 
\begin{align}
    \mu (\|\rvw^{(1)}\|^2+\|\rvw^{(2)}\|^2)=  \langle \rvw,\overline{\Pi^{(\vtau)}} \scA^{(\vtau)}\overline{\Pi^{(\vtau)}}\rvw\rangle =\langle \rvw,  \scA^{(\vtau)} \rvw\rangle =2\langle \rvw^{(1)},  \rmX^{(\vtau)} \rvw^{(2)}\rangle.\label{eqn:eigen_equation}
\end{align}
Let us then consider a vector $$\rvy=\begin{bmatrix}
    \qr \rvw^{(1)}\\
    \rvw^{(2)}
    \end{bmatrix}.$$
Thus, due to \eqref{eqn:hatH_psd_order}, \eqref{eqn:eigen_equation}, and \eqref{eqn:norm_half}, by choosing $\vtau$ satisfying \eqref{eqn:tau_constrain},
we conclude
\begin{align}
\qr \mu & = 2\qr \langle \rvw^{(1)},  \rmX^{(\vtau)}\rvw^{(2)}\rangle =\langle \rvy,\scA^{(\vtau)}\rvy\rangle\\ 
& \leq  \frac{(\qr^2+1)}{2}+\qr^2\sum_{x\in\gV_1\setminus\pruneV }\alpha_x \ervw_x^2+\sum_{x\in\gV_2\setminus\pruneV }\alpha_x \ervw_x^2 \\ 
&\,\quad +\sum_{x\in\pruneV }\alpha_x (\qr^2 \ervw_x^2\indi{x\in\gV_1} + \ervw_x^2\indi{x\in\gV_2})+ \frac{1}{8}\cdot\qr\error^{1/2}\\
&\leq \frac{\qr^2+1}{2}+\frac{\qr^2\tau_1+\tau_2^+}{2} + \qr^2(\max_{y\in[N]}\alpha_y)\cdot\sum_{x\in\gV^{(\vtau)} } \ervw_x^2+\frac{1}{8}\cdot\qr\error^{1/2}. \label{eqn:upperbpund_mu}
\end{align}
With Bennett's inequality in \Cref{lem:Bennett}, we know that for all $y\in[N]$ with very high probability, 
\[
\alpha_y\le c(1+\frac{\log N}{d})\le d.
\] 
Denote $\mu_0:=\qr+\qr^{-1}+\frac{1}{4}\error^{1/2}$. Recall the definition of $\error$ in \eqref{eqn:error_parameters} and radius $\radius$ in \eqref{eqn:radius}.  
Then, together with Proposition~\ref{prop:delocalization_vec},
since $\mu>\qr+\qr^{-1}+\frac{1}{2}\error^{1/2}> \mu_0$, we can obtain that
\begin{align}
 \qr(\max_{y\in[N]}\alpha_y)\cdot\sum_{x\in\gV^{(\vtau)} }w_x^2&\leq  \qr\cdot d \left(\frac{1}{1-\mu_0/\,u}\right)^{4}\left(1-\frac{\mu-\mu_0}{\mu}\right)^{2\radius}\\
&\leq \qr\cdot d \left(\frac{\qr+\qr^{-1}+\frac{1}{2}\error^{1/2}}{\frac{1}{4}\error^{1/2}}\right)^{4}\left(1-\frac{ \error^{1/2}}{12\qr}\right)^{2\radius}\\
&\leq\frac{ 256\qr(\qr+\qr^{-1}+1) d}{ \error^2}\exp\left(-\frac{ \radius}{6\qr}\error^{1/2}\right)\\
&\leq\frac{ 256\qr(\qr+\qr^{-1}+1)}{ \error^5}\cdot \exp\left(- \frac{1}{24\qr^2}\error^{-1/2}\right)\ll \frac{1}{8}\error^{1/2},
\end{align}
where the third step is followed by the definition of \eqref{eqn:error_parameters} and $\error^{-1/2} \gg \log^2(d)$.
Plugging this estimate into (\ref{eqn:upperbpund_mu}) and dividing both sides by $\qr$, we have
\begin{align}
    \mu \leq \frac{\qr+\qr^{-1}}{2}+\frac{\qr\tau_1+\qr^{-1}\tau_2^+}{2}+\frac{1}{4}\error^{1/4} = \qr  +\qr^{-1} +\frac{1}{2}\error^{1/4},
\end{align}
because of the choice of $\tau_1$ and $\tau_2^+$ in \eqref{eqn:tau_constrain}. However, this is a contradiction to the assumption \eqref{eq:lem612_contradiction}. Hence, we conclude the final result of this proposition.
\end{proof}

\begin{proof}[Proof of \Cref{thm:right_edge_behavior} (2)]
   Recall the definitions of $\sR_2$ and $\sR_1$ in \eqref{eqn:R_2} and \eqref{eqn:R_1}. Note that \[\frac{d}{dt}\Lambda_{\qr}(t)|_{t = \qr^{2}+1} = 0,  \quad \frac{d}{dt}\Lambda_{\qr^{-1}}(t)|_{t = \qr^{-2}+1} = 0,\] and $\frac{d^{2}}{dt^{2}}\Lambda_{\qr}(t)|_{t = \qr^{2}+1} \neq 0$, $\frac{d^{2}}{dt^{2}}\Lambda_{\qr^{-1}}(t)|_{t = \qr^{-2}+1} \neq 0$. Then by Taylor expansion, there exists a universal constant $\gC>0$ such that
   \begin{subequations}
    \begin{align}
            \sR_2 = &\, \{x\in \gV_2| \alpha_x  \geq 1+\qr^{2}+\gC\error^{1/8} \},\\
        \sR_1 = &\, \{x\in \gV_1| \alpha_x  \geq  1 + \qr^{-2} + \gC\error^{1/8} \}.
    \end{align}
\end{subequations}
Then, $\sR_2\cup \sR_1\subset\gV^{(\vtau^*)}$ from \eqref{def:W_set}. Hence, all the elements in $\Lambda^{\sR_1} \cup \Lambda^{\sR_2}$ in a non-increasing order,  
\begin{align}
    \Lambda_1 \geq \Lambda_2\geq \cdots \geq \Lambda_{|\sR_1| + |\sR_2|} \geq \qr + \qr^{-1} + \error^{1/4}\,,
\end{align}
are the positive eigenvalues of $\hatH$ defined in \eqref{eqn:hatH}.
As a convention, we denote $\sR_2=\emptyset$ if $\alpha_{x}<1+\qr^{2}+\gC\error^{1/4}$ for all $x\in\gV_2$; and $\sR_1=\emptyset$ if $\alpha_{x}<1+\qr^{-2}+\gC\error^{1/4}$ for all $x\in\gV_1$. 
Moreover, we know that except for positive eigenvalues $\Lambda_1\geq \Lambda_2\geq \cdots \geq \Lambda_{|\sR_1| + |\sR_2|}$, all other eigenvalues of $\hatH$ are smaller than $\qr^{-1}+\qr +\frac{1}{2}\error^{1/4}$ due to Lemma~\ref{lem:project_norm} with $\vtau$ defined in \eqref{eqn:tau_constrain}.
Then, combining \Cref{lem:hatHtau_norm} and Weyl's inequality \Cref{lem:weyl}, we can conclude that with high probability, $\lambda_{|\sR_1| + |\sR_2|+1}(\rmH) \leq \qr+\qr^{-1} + \error^{1/4}$.
\end{proof}

\subsection{Proof of \Cref{thm:left_edge_behavior} (2)} 
\begin{lemma}\label{lem:project_lower}
 Assume that $\vtau$ satisfies \eqref{eqn:tau_constrain}. Let $(\lambda,\rvw)$ be an eigenpair of $\hatH$ in \eqref{eqn:hatH}, such that $\lambda>0$ and $\rvw\perp \textnormal{span}\{\hatvtau_{\sigma}(x):x\in\gV^{(\vtau^{\star})}\}$, then with high probability
    \[
    \lambda\ge  \qr  -\qr^{-1}-\frac{1}{2}\error^{1/4}. 
    \]
\end{lemma}

\begin{proof}[Proof of \Cref{lem:project_lower}]
Assume that $0<\lambda <\qr  - \qr^{-1} - \frac{1}{2}\error^{1/4}$. We will derive a contradiction below using \Cref{prop:lower_bound_H} and \Cref{prop:delocalization_vec} (2).

Recall that $\rmH^{(\vtau)}$ in \eqref{eqn:def_H_tau} admits the block structure and $\rmX^{(\vtau)}$ denotes the $n\times m$ block. We first derive the following Loewner order \eqref{eqn:lowerboundX_1} in terms of $\rmX^{(\vtau)}$.
\begin{align}
    {\rmX}^{(\vtau)*}(\rmI - d^{-1}\rmD^{(1)})^{-1} {\rmX}^{(\vtau)}\succeq d^{-1}\rmD^{(2)}- \rmI - \rmE \label{eqn:lowerboundX_1}
\end{align}
Note that, with high probability, $\|\rmX^{(\vtau)}\|\lesssim 1$ due to \eqref{eqn:M1_inverse_bounds}, $\|(\rmI - d^{-1}\rmD^{(1)})^{-1}\|\le \epsilon^{-1}$ for some small $\epsilon>0$ due to \Cref{lem:invertiblility_D1}, and $\|\rmX-\rmX^{(\vtau)}\|\le \gC\error^{1/2}$ due to Lemma~\ref{lem:hatHtau_norm}. Then with high probability,
\begin{align}
   &\,\| {\rmX}^{(\vtau)*}(\rmI - d^{-1}\rmD^{(1)})^{-1} {\rmX}^{(\vtau)}-{\rmX}^*(\rmI - d^{-1}\rmD^{(1)})^{-1} {\rmX}\|\\
    \leq &\, \| (\rmI - d^{-1}\rmD^{(1)})^{-1} \|\cdot\|{\rmX}^{(\vtau)}-{\rmX}\|\cdot(\|{\rmX}^{(\vtau)}\|+\|{\rmX}\|) \lesssim \error^{1/2}.
\end{align}
We conclude \eqref{eqn:lowerboundX_1} from \eqref{eqn:lower_bound_H}, where the error matrix $\|\rmE\| \ll \error^{1/4}$ for sufficiently large $N$.

Consider the block structure of $\rvw$ in \eqref{eqn:w_block_structure}. We use $\rvw^{(2)}$ as a test vector on both sides of \eqref{eqn:lowerboundX_1}, then 
\begin{align}
    \<  {\rmX^{(\vtau)}}\rvw^{(2)}\, , \, (\rmI - d^{-1}\rmD^{(1)})^{-1} {\rmX^{(\vtau)}}\rvw^{(2)} \>  \geq \<\rvw^{(2)}\, , \, (d^{-1}\rmD^{(2)}- \rmI - \rmE)\rvw^{(2)}\> \label{eqn:test_w2_inequality}
\end{align}

For the left-hand side of \eqref{eqn:test_w2_inequality}, note that $(\lambda,\rvw)$ is an eigenpair of $\hatH$ with $\rvw\perp \textnormal{span}\{\hatvtau_{\sigma}(x):x\in\gV^{(\vtau^*)}\}$. By construction, $(\lambda,\rvw)$ is also an eigenpair of $\rmH^{(\vtau)}$ in \eqref{eqn:def_H_tau}, hence $\lambda \rvw^{(1)} = \rmX^{(\vtau)}\rvw^{(2)}$. Furthermore, \Cref{lem:invertiblility_D1} shows that with high probability, $0 \leq \alpha_x \leq \qr^{-2} - \epsilon$ for all $x\in\gV_1$. Then
\begin{align}
    &\<\rvw^{(1)}, (\rmI - d^{-1}\rmD^{(1)})^{-1} \rvw^{(1)}\> =   \sum_{x\in \gV_{1}} \ervw_{x}^{2} (1 - \alpha_{x})^{-1}\\
     = ~& \sum_{0 \leq \alpha_{x} < \tau_1} \ervw_{x}^{2} (1 -\alpha_{x})^{-1}  
    + \sum_{\tau_1\leq \alpha_{x} < 1 - \epsilon} \ervw_{x}^{2} (1 -\alpha_{x})^{-1}.\label{eqn:left_bound_1}
\end{align}
Using $\|\rvw^{(1)}\|^{2}_{2} = 1/2$ in \eqref{eqn:norm_half} and $\tau_{1} \geq \qr^{-2} + \frac{1}{4\qr}\error^{1/4}$ in \eqref{eqn:tau_constrain}, the first term is bounded by
\begin{align}
    \sum_{0\leq \alpha_{x} < \tau_1} \ervw_{x}^{2} (1 -\alpha_{x})^{-1} \leq \frac{1}{2} (1 - \qr^{-2} - \frac{1}{4\qr}\error^{1/4})^{-1}.\label{eqn:left_bound_2}
\end{align}
Note that \eqref{eqn:wx_ratio_left} holds for $x\in \pruneV \setminus \gV^{(\vtau^*)}$. By applying \eqref{eqn:left_delocalization}, for sufficiently large $d$, we have
\begin{align}
    &\, \sum_{\tau_1 \leq \alpha_{x} < 1 -  \epsilon} \ervw_{x}^{2} (1 -\alpha_{x})^{-1} \\
    \leq &\, \epsilon^{-1}    \sum_{x\in\pruneV} \ervw_{x}^{2} \leq  \epsilon^{-1} \Big(\frac{ \qr - \qr^{-1} - \frac{1}{2}\error^{1/2} }{ \qr - \qr^{-1} - \frac{1}{2}\error^{1/2} - \lambda} \Big)^{4} \Big( \frac{\lambda}{\qr - \qr^{-1} - \frac{1}{2}\error^{1/2}} \Big)^{ 2\radius}\\
    \leq &\,\epsilon^{-1} \Big(\frac{ \qr - \qr^{-1} - \frac{1}{2}\error^{1/2} }{\frac{1}{4}\error^{1/2} } \Big)^{4} \Big( 1 - \frac{\frac{1}{4}\error^{1/2}}{\qr - \qr^{-1} - \frac{1}{2}\error^{1/2}} \Big)^{ 2\radius}\\
    \leq &\,\epsilon^{-1} \exp \bigg( - \frac{\radius\error^{1/2}}{2\qr} -4\log\Big(\frac{\error^{1/2}}{4\qr} \Big) \bigg)\ll \frac{1}{4( \qr-\qr^{-1}) }\error^{1/4}.\label{eqn:left_bound_3}
\end{align} 
where the second line holds since $\qr - \qr^{-1} - \frac{1}{2}\error^{1/2} - \lambda > \frac{1}{2}\error^{1/4} - \frac{1}{2}\error^{1/2} > \frac{1}{4}\error^{1/2}$ due to the assumption, and in the last line, we use the facts  $1-x\leq e^{-x}$, $\epsilon\ge \error^{1/2}$,  $\error= \sqrt{\log(d)/d}$ and $\radius \asymp \sqrt{d/\log(d)}$ in \eqref{eqn:radius}.

For the right-hand side of \eqref{eqn:test_w2_inequality}, we consider a slightly different decomposition 
\begin{align}
    &\,\<\rvw^{(2)}, d^{-1}\rmD^{(2)}\rvw^{(2)}\> \\
    = &\,\sum_{x\in \gV_{2}} \ervw_{x}^{2} \alpha_{x} 
    = \sum_{x\in \gV_{2}} \ervw_{x}^{2} \alpha_{x} \cdot \Big( \indi{x \in \gV^{(\vtau^*)}} + \indi{x\in \pruneV \setminus \gV^{(\vtau^*)}} + \indi{x\in \gV_{2} \setminus \pruneV } \Big),\label{eqn:left_bound_4}
\end{align}
where $\vtau$ and $\vtau^{\star}$ satisfy \eqref{eqn:tau_constrain} and \eqref{def:W_set}, respectively. Firstly, eigenvectors $\{\hatvtau_{\sigma}(x)\}_{x\in\gV^{(\vtau^{\star})}}$ have pairwise disjoint supports due to \Cref{lem:existence_pruned_graph}, and $\rvw\perp \textnormal{span}\{\hatvtau_{\sigma}(x):x\in\gV^{(\vtau^{\star})}\}$, then $\ervw_{x} = 0$ for all $x\in \gV^{(\vtau^*)}$. Secondly, for $x\in \pruneV \setminus \gV^{(\vtau^{\star})}$, since $\lambda <\qr  -\qr^{-1}-\frac{1}{2}\error^{1/4}$, then by \eqref{eqn:wx_ratio_left} and calculations in \eqref{eqn:left_bound_3},
\begin{align}
    0\leq \sum_{x\in \gV_{2}} \ervw_{x}^{2} \alpha_{x} \indi{x\in \pruneV \setminus \gV^{(\vtau^*)}} \leq \tau_{2}^{-} \sum_{x\in \gV_{2}} \ervw_{x}^{2} \indi{x\in \pruneV \setminus \gV^{(\vtau^*)}} \ll \error^{1/2}. \label{eqn:left_bound_5}
\end{align}
Finally, vertices in $\gV_{2} \setminus \pruneV$ satisfy $\alpha_{x} \geq \qr^{2} - \frac{\qr}{4}\error^{1/4}$. Using calculations in \eqref{eqn:left_bound_5}, we then have
\begin{align}
    \sum_{x\in \gV_{2}} \ervw_{x}^{2} \alpha_{x} \indi{ x\in \gV_{2} \setminus \pruneV } \geq \Big(\qr^{2} - \frac{\qr}{4}\error^{1/4} \Big) \cdot \Big( \frac{1}{2} - \error^{1/2} \Big).\label{eqn:left_bound_6}
\end{align}

According to the upper bound of \eqref{eqn:left_bound_1} in \eqref{eqn:left_bound_2}, \eqref{eqn:left_bound_3}, and lower bound of \eqref{eqn:left_bound_4} in \eqref{eqn:left_bound_5} and \eqref{eqn:left_bound_6}, \eqref{eqn:test_w2_inequality} can be further expanded as
\begin{align}
&\, \frac{\lambda^{2}}{2} \cdot [ (1 - \qr^{-2} - \frac{1}{4\qr}\error^{1/4})^{-1}  +\frac{1}{4( \qr-\qr^{-1}) }\error^{1/4} ]\\
\geq &\, \lambda^{2}\<\rvw^{(1)}, (\rmI - d^{-1}\rmD^{(1)})^{-1} \rvw^{(1)}\> = \<  {\rmX^{(\vtau)}}\rvw^{(2)}\, , \, (\rmI - d^{-1}\rmD^{(1)})^{-1} {\rmX^{(\vtau)}}\rvw^{(2)} \>  \\
\geq &\, \<\rvw^{(2)},(d^{-1}\rmD^{(2)}- \rmI - \rmE)\rvw^{(2)}\> \geq \frac{1}{2}(\qr^{2} - 1 - \frac{\qr}{8}\error^{1/4}),
\end{align}
Consequently, the following holds
\begin{align}
     \lambda^{2} \geq \frac{(\qr - \qr^{-1})^2(1- \frac{\error^{1/4}}{8(\qr-\qr^{-1})})(1- \frac{\error^{1/4}}{4(\qr-\qr^{-1})})}{1+\frac{1}{4\qr}(1 - \frac{1}{4(\qr - \qr^{-1})}\error^{1/4}) \error^{1/4}} \geq (\qr - \qr^{-1})^2 \Big(1- \frac{\error^{1/4}}{4(\qr-\qr^{-1})}\Big)^{2}
\end{align}
which then implies $\lambda \geq \qr - \qr^{-1} - \frac{1}{2}\error^{1/4}$. This contradicts the initial assumption that $0< \lambda <\qr - \qr^{-1} - \frac{1}{2}\error^{1/4}$. Therefore, the proof is completed.
\end{proof}

\begin{proof}[Proof of \Cref{thm:left_edge_behavior} (2).]
 Recall the set $\sL_2$ in \eqref{eqn:L2_set} and the permutation $\pi$ in $[m]$ which arranges the normalized degrees of $\gV_2$ in a non-increasing order. Note that $\frac{d}{dt}\Lambda_{\qr}(t)|_{t = \qr^{2}-1} = 0$, and $\frac{d^{2}}{dt^{2}}\Lambda_{\qr}(t)|_{t = \qr^{2}-1} \neq 0$. Then equivalently, $\sL_2$ can be rewritten in terms of the normalized degree using Taylor expansion
    \begin{align}
     \sL_2  &= \{x\in\gV_2 | 0<\alpha_{\pi(l)}\leq \qr^{2} -1 - \gC \error^{1/8}\}.
    \end{align}
Note that $\sL_2\subseteq \gV^{(\vtau^*)}$ due to \eqref{def:W_set}. 
By the definition of $\sL_2$, we have
\begin{align}
    \Lambda_\qr(\alpha_{\pi(m)})\leq \Lambda_\qr(\alpha_{\pi(m-1)})\leq \cdots \leq \Lambda_\qr(\alpha_{\pi(m-|\sL_2|+1)})\leq \qr - \qr^{-1} - \error^{1/4}.
\end{align}
According to \Cref{lem:project_lower}, except for $|\sL_2|$ positive eigenvalues  above, all the other positive eigenvalues of $\hatH$ in \eqref{eqn:hatH} are larger than $\qr  -\qr^{-1}-\frac{1}{2}\error^{1/4}$. Then, combining \Cref{lem:hatHtau_norm} and Weyl's inequality, we can conclude that with high probability, $\lambda_{l}(\rmH) \geq \qr-\qr^{-1} -\error^{1/4}$  for any $1\le l\le m-|\sL_2|$.
\end{proof}
\section*{Acknowledgments}
H.X.W. would like to thank Hanbaek Lyu and Ruth Williams for their valuable feedback when he presented this paper at probability seminars at UW-Madison and UC San Diego, which lead to a better understanding of \Cref{thm:number_outliers}. Part of this manuscript was finished when H.X.W. was attending the \emph{Second Conference on Random Matrix Theory and Numerical Linear Algebra} at the University of Washington, Seattle in June 2025. H.X.W. thanks Tom Trogdon and Xiucai Ding for their generous hospitality, as well as comments from audience during the research talk.
I.D. and Y.Z. thank Antti Knowles for many helpful discussions.

I.D. was partially supported by NSF DMS-2154099. This material is based upon work supported by NSF DMS-1928930, while the I.D. and Y.Z. were in residence at the Simons Laufer Mathematical Sciences Institute in Berkeley, California, during the Universality and Integrability in Random Matrix Theory and Interacting Particle Systems Program (Fall 2021) and I.D. and Z.W. were in residence during the Probability and Statistics of Discrete Structures Program (Spring 2025). Y.Z. was partially supported by an AMS-Simons Travel Grant. This material is based upon work supported by the Swedish Research Council under grant no. 202106594 while Y.Z. was in residence at Institut Mittag-Leffler in Djursholm, Sweden during the Fall of 2024.

\printbibliography
\newpage
\appendix
\addcontentsline{toc}{section}{Appendices}
\section{Facts about bipartite Erd\H{o}s-R\'enyi graphs}
\begin{lemma}[Connectivity of random bipartite graphs]\label{lem:connectivity_bipartite}
    Let $\gG$ be a random bipartite graph, where each edge with one vertex in $\gV_1$ and the other in $\gV_2$ is sampled with probability $p$ with $|\gV_1| = n$ and $|\gV_2| = m$. Let $p \coloneqq (b\log(N) + c)/\sqrt{mn}$ where $N = n + m$ and $0\leq c\leq \log(N)$. Define the ratio $\ratio \coloneqq n/m \geq 1$ and $\qr \coloneqq \ratio^{1/4}$. If $b\geq \qr^2$, with probability at least $1 - 2e^{-\qr^{-2}c}- O(N^{-2b})$, there is only one connected component in $\gG$ as $n, m \to \infty$. Consequently, $\gG$ is connected with probability tending to $1$ as $c \to \infty$.
\end{lemma}
\begin{proof}[Proof of \Cref{lem:connectivity_bipartite}]
    Let $\rX_k \coloneqq \rX_{k, N}$ denote the number of connected components of $k$ vertices in $\gG$, where those components are pairwise disconnected, then
    \begin{align}
        \P(\gG\textnormal{ is not connected}) = \P\Big( \bigcup_{k=1}^{\lfloor N/2 \rfloor} \gG \textnormal{ has a component of } k \textnormal{ vertices} \Big) = \P\Big( \bigcup_{k=1}^{\lfloor N/2 \rfloor} \{\rX_k > 0\} \Big).
    \end{align}
    Note that $\rX_1$ counts the number of isolated vertices, then
    \begin{align}
        \P(\rX_1 > 0) \leq \P(\gG \textnormal{ is not connected}) \leq \P(\rX_1 > 0) + \sum_{k=2}^{\lfloor N/2 \rfloor}  \P(\rX_k > 0).
    \end{align}
    With the proof later deferred, we claim that the following holds when $b \geq \qr^{2}$
    \begin{align}
        \P(\gG \textnormal{ is connected}) \geq 1 - \P(\rX_1 \geq 1) - O(N^{-2b + o(1)}). \label{eqn:X1_connectivity}
    \end{align}
    Let $\Xi_{v}$ denote the event that vertex $v$ is isolated. Note that $p = (b\log(N) + c)/\sqrt{mn}$ and $\qr^{4} = \ratio \geq 1$, then
    \begin{align}
        \E \rX_1 =&\, \sum_{v\in \gV_1} \P(\Xi_{v}) + \sum_{v\in \gV_2}\P(\Xi_{v}) \\
        =&\, n(1-p)^{m} + m(1 - p)^{n} \\
        =&\, n \exp\big(m\log(1 - p)\big) + m\exp\big(n\log(1-p)\big)\\
        =&\, n \exp\big(-mp + O(mp^2)\big) + m\exp\big(-np + O(np^2)\big) \quad \Big(\log(1 + x) = x + O(x^2) \textnormal{ when } x = o(1) \Big)\\
        =&\, n \exp\Big(-b\qr^{-2}\log(N) - \qr^{-2}c + O((\log N)^2/n)\Big) + m\exp\Big(-b\qr^{2}\log(N) - \qr^{2}c + O((\log N)^2/m)\Big)\\
        =&\, \exp\Big( (1 - b\qr^{-2})\log(n) - \qr^{-2}c - b\qr^{-2}\log(1 + \qr^{-4}) + O((\log N)^2/n)\Big) \\
        &\, + \exp\Big( (1 - b\qr^{2})\log(m) - \qr^{2}c - b\qr^{2}\log(1 + \qr^{4}) + O((\log N)^2/m) \Big)\\
        \leq &\, 2\exp\big(-\qr^{-2}c + O((\log N)^2/n) \big)
    \end{align}
    where the last inequality holds since $1 - b\qr^{2} \leq 1 - b\qr^{-2} \leq 0$, $b\qr^{-2}\log(1 + \qr^{-4}) \geq 0$ and $b\qr^{2}\log(1 + \qr^{4})\geq 0$ when $b\geq \qr^{2}$. The desired result follows directly by Markov's inequality, since
    \begin{align}
        \P(\rX_1 \geq 1) \leq \E \rX_1 \leq 2\exp\big( - \qr^{-2}c + O((\log N)^2/n)\big).
    \end{align}
What remains, therefore, is the proof of~\eqref{eqn:X1_connectivity}. By Cayley's formula~\eqref{lem:Cayley_formula}, the number of trees on $k$ labeled vertices is $k^{k-2}$ for $k \geq 2$, then by Markov
\begin{align}
    \P(\rX_k > 0) \leq \E \rX_k \leq \binom{N}{k} k^{k-2} p^{k -1} (1 - p)^{k(N-k)},
\end{align}
where the factor $p^{k-1}$ ensures that the $k-1$ edges in the tree is connected, while the factor $(1 - p)^{k(N-k)}$ ensures that those $k$ vertices are disconnected from the rest $N-k$ nodes. We use the fact $1 - p = e^{-p}$ since $p = (b\log(N) + c)/\sqrt{mn} = o(1)$, and the Stirling's approximation \Cref{lem:stirling}, then
\begin{align}
    \log \P(\rX_k > 0) \leq &\, \frac{1}{2}\log \frac{N}{2\pi k (N-k)} + N\log(N) - k\log(k) - (N-k) \log(N-k) + O(N^{-1}) +  O(k^{-1}) \\
    &\, + (k-2)\log(k) + (k-1)\log (p) - k(N-k)p\\
    \leq &\, \Big(N + \frac{1}{2}\Big) \log \Big( 1 + \frac{k}{N-k} \Big) - \frac{5}{2}\log(k) - \frac{1}{2}\log(2\pi) + O(N^{-1}) +  O(k^{-1})\\
     &\, + 2k\log\Big((N-k)\frac{b\log(N) + c}{\sqrt{mn}} \Big) - k (N-k)\frac{b\log(N) + c}{\sqrt{mn}}\\
     \leq &\, -\frac{k(N-k)}{2N} (1 + \qr^{4})\qr^{-2}(b\log(N) + c)\cdot (1 + o(1)),
\end{align}
where in the second to last line, the last term dominates. For $k\geq 2$, we split the sum and bound the two terms separately, when $b \geq \qr^{2}$,
\begin{align}
    &\, \sum_{k=2}^{\lfloor N/2 \rfloor} \P( \rX_k > 0) =  \sum_{k=2}^{\lfloor 11 \rfloor} \P( \rX_k > 0) + \sum_{k=12}^{\lfloor N/2 \rfloor} \P( \rX_k > 0)\\
    \leq &\, 10 \cdot e^{-(b\log(N) + c)\cdot (\qr^{2} + \qr^{-2})\cdot (1 + o(1))} + \sum_{k=12}^{\lfloor N/2 \rfloor} e^{-(b\log(N) + c)\cdot k/2\cdot (\qr^{2} + \qr^{-2})\cdot (1 + o(1))}, \\
    \leq &\,10 \cdot e^{-(b\log(N) + c)\cdot (\qr^{2} + \qr^{-2})\cdot (1 + o(1))} + \frac{N}{2} \cdot e^{-6(b\log(N) + c)\cdot (\qr^{2} + \qr^{-2})\cdot (1 + o(1))} \\
    \leq &\, O(N^{-2b + o(1)}),
\end{align}
where the factor $2b$ on the exponent comes from the fact that $\qr^{2} + \qr^{-2} \geq 2$. Consequently,
\begin{align}
    &\, \P(\gG \textnormal{ is connected}) = 1 - \P(\gG \textnormal{ is not connected}) \leq 1 - \P(\rX_{1} \geq 1) = \P(\rX_{1} = 0),\\
    &\, \P(\gG \textnormal{ is connected}) \geq 1 - \P(\rX_1 > 0) - \sum_{k=2}^{\lfloor N/2 \rfloor}  \P(\rX_k > 0) = 1 - \P(\rX_1 \geq 1) - O(N^{-2b + o(1)}).
\end{align}
which completes the proof of \eqref{eqn:X1_connectivity}.
\end{proof}

\begin{lemma}[Bounds on degrees] \label{lem:deg_bound}
Recall $\Vonehigh, \Vtwohigh, \Vtwolow$ in \eqref{eqn:atypical_vertices}. Then for any $\nu >0$, there is a universal constant $\const >0$ such that
\begin{align}
    \P\bigg( x\in \Big( \gV_{1} \setminus \Vonehigh \Big) \cup \Big( \gV_{2} \setminus (\Vtwohigh \cup \Vtwolow ) \Big) \bigg) \geq 1 - \const N^{-\nu},
\end{align}
where $\tau_{1} = \qr^{-2} + 1$, $\tau^{+}_{2} = \qr^{2} + 1$ and $\tau^{-}_{2} = \qr^{2} - 1$. Furthermore, $1\leq j \leq r_{x}$ with $r_x$ in \eqref{eqn:defrx}, then the event
\begin{align}
     \{\rD_{x} \leq \sqrt{N}(2d)^{-\lfloor j/2 \rfloor}\}.
\end{align}
holds with probability at least $1 - \exp(-N^{1/8}\log(N))$.
\end{lemma}
\begin{proof}[Proof of \Cref{lem:deg_bound}]
Since $d_{2} = \qr^{2}d$ and $d_{1} = \qr^{-2}d$, by Bennett's inequality \Cref{lem:Bennett}, we have
    \begin{align}
        &\, \P(\rD_{x} \geq d_2 + u d) = \P(\alpha_{x} \geq \qr^{2} + u)  \leq \exp\big(-d \qr^{2} \benrate(u \qr^{-2})\big), \quad \text{for} \quad x\in \gV_2,\\
        &\, \P(\rD_{x} \leq d_2 - u d) = \P(\alpha_{x} \leq \qr^{2} - u) \leq \exp\big(-d \qr^{2} \benrate(u \qr^{-2})\big), \quad \text{for} \quad x\in \gV_2,\\
        &\, \P(\rD_{x} \geq d_1 + u d ) = \P(\alpha_{x} \geq \qr^{-2} + u) \leq \exp\big(-d \qr^{-2} \benrate(u \qr^{2})\big), \quad \text{for} \quad x\in \gV_1.
    \end{align}
    Note that $\benrate(u)= (1 + u)\log(1 + u) - u \geq 0$ for all $u\geq 0$, consequently,
    \begin{align}
        &\, \P(x\in \Vtwohigh) \leq \exp\Big(-d \qr^{2} \benrate\big( (\tau^{+}_{2} - \qr^{2})\qr^{-2} \big)\Big),\\
        &\, \P(x\in \Vtwolow) \leq \exp\Big(-d \qr^{2} \benrate\big( (\qr^{2} - \tau^{-}_{2})\qr^{-2} \big)\Big),\\
        &\, \P(x\in \Vonehigh) \leq \exp\Big(-d \qr^{-2} \benrate\big( (\tau_{1} - \qr^{-2})\qr^{2} \big)\Big).
    \end{align}
    The first argument follows since the three terms above are dominated by $N^{-\nu}$ for some constant $\nu >0$ since $\benrate(1) = O(1)$. 

    The second argument follows directly by taking $\tau_1 = \tau^{+}_{2} = \sqrt{N}(2d)^{-r_{x} - 1}$.
\end{proof}

\section{Deferred proofs in \Cref{sec:proof_emergence}}\label{sec:proof_emergence_lemma}
\begin{proof}[Proof of \Cref{lem:approxBinomProb}]
We first prove the case where $x \in \gV_2$. By definition,
\begin{align}
     \P(\rD_x^{(2)} \geq \alpha d) = \sum_{t = \lceil \alpha d \rceil}^{n} \P(\rD_x^{(2)} = t) = \sum_{t = \lceil \alpha d \rceil}^{n} \binom{n}{t} p^{t}(1 - p)^{n-t}.
\end{align}
We define the following ratio function
\begin{align}
    r_2(t) \coloneqq \frac{\P(\rD_x^{(2)} = t) }{\P(\rD_x^{(2)} = t - 1)} = \frac{(n-t+1)}{t} \cdot \frac{p}{1 - p}\,\,,
\end{align}
which decreases with respect to $t$. Note that the critical ratio $r_2(t^{\star}) = 1$ when $t^{\star} = p(n+1) = \qr^{2} d + o(1) = \qr^{2} d + o(1)$, thus $\P(\rD_x^{(2)} = t)$ achieves its maximum around $t^{\star} = \qr^{2} d$, since
\begin{align}
    &\, r_2(t) \geq r_2(\qr^{2} d) = \frac{n - \qr^{2} d + 1}{\qr^{2} d} \cdot \frac{d/\sqrt{mn}}{1 - d/{\sqrt{mn}}} = \frac{1 - d/\sqrt{mn} + 1/n  }{1 - d/\sqrt{mn}} > 1,\quad \forall t\leq \qr^{2} d,\\
    &\, r_2(t) \leq r_2(\qr^{2} d + 1) = \frac{n - \qr^{2} d}{\qr^{2} d + 1} \cdot \frac{d/\sqrt{mn}}{1 - d/{\sqrt{mn}}} = \frac{\qr^{2} d}{\qr^{2} d + 1} < 1,\quad \forall t\geq \qr^{2} d + 1.
\end{align}
Note that $\alpha \geq \qr^{2} + 1$, then for $t \geq \alpha d \geq \qr^{2} d + d > \qr^{d}$, we have
\begin{align}
    r_2(t) \leq r_2(\alpha d) = \frac{\qr^{2}}{\alpha} \cdot \frac{1 - \alpha d/n + 1/n}{1 - d/\sqrt{mn}}\leq \frac{\qr^{2} }{\alpha} < 1.
\end{align}
Consequently, using the fact $r_2(t) \leq \qr^{2}\alpha^{-1}$ for all $t \geq \alpha d$, we have
\begin{align}
     \P(\rD_x^{(2)} = \alpha d) \leq \P(\rD_x^{(2)} \geq \alpha d) \leq &\, \P(\rD_x^{(2)} = \alpha d)\cdot \sum_{l = 0}^{n - \lceil \alpha d \rceil} (\qr^{2} \alpha^{-1})^{l} \leq \frac{\alpha }{\alpha - \qr^{2} } \cdot \P(\rD_x^{(2)} = \alpha d),
\end{align}
which completes the proof of $(1)$. Similarly, when $\alpha \leq \qr^{2} - 1$, for $t < \alpha d$, we have
\begin{align}
    r_2(t) \geq r_2(\alpha d - 1) = \frac{\qr^{2}}{\alpha} \cdot \frac{1 - \alpha d/n + 2/n}{\big(1 - 1/(\alpha d) \big) \cdot (1 - d/\sqrt{mn})} \geq \frac{\qr^{2}}{\alpha} > 1,
\end{align}
Then by following the same technique, the proof of $(1)$ is finished by
\begin{align}
     \P(\rD_x^{(2)} = \alpha d) \leq \P(\rD_x^{(2)} \leq \alpha d) \leq &\, \P(\rD_x^{(2)} = \alpha d)\cdot \sum_{l = 0}^{\lceil \alpha d \rceil} (\qr^{-2} \alpha)^{l} \leq \frac{\qr^{2}}{\qr^{2} - \alpha} \cdot \P(\rD_x^{(2)} = \alpha d).
\end{align}
For the case $v \in \gV_1$, denote the ratio
\begin{align}
    r_1(t) = \frac{\P(\rD_x^{(1)} = t) }{\P(\rD_x^{(1)} = t - 1)} = \frac{(m-t+1)}{t} \cdot \frac{p}{1 - p},
\end{align}
which decreases with respect to $t$. Note that the critical ratio $r_1(t^{\star}) = 1$ for $t^{\star} = p(m+1) = d/\qr^{2} + o(1) = \qr^{-2} d + o(1)$. Under the condition that $\alpha \geq \qr^{-2} + 1$, for any $t \geq \lceil \alpha d\rceil$, the ratio
\begin{align}
    r_1(t) \leq r_1(\alpha d) = \frac{m-\alpha d + 1}{\alpha d} \cdot \frac{d/\sqrt{mn}}{1 - d/{\sqrt{mn}}} \leq \frac{1}{\alpha \qr^{2}} < 1,
\end{align}
leading to the fact $\P(\rD_x^{(1)} = t)$ decreases monotonically when $t\geq \lceil \alpha d\rceil$. Therefore,
\begin{align}
   \P(\rD_x^{(1)} = \alpha d) \leq \P(\rD_x^{(1)} \geq \alpha d) \leq  &\, \P(\rD_x^{(1)} = \alpha d) \sum_{l = 0}^{m - \lceil \alpha d\rceil } (\alpha \qr^{2})^{-l} \leq  \frac{\alpha \qr^{2}}{\alpha \qr^{2} - 1}  \P(\rD_x^{(1)} = \alpha d),
\end{align}
which completes the proof of $(3)$.
\end{proof}

\begin{proof}[Proof of \Cref{lem:degreeProbApprox}]
Using the facts $p = d/\sqrt{mn}$, $1 \ll \alpha d \asymp \log(n) \ll n$, $\log(1 + x )= x$ for $x = o(1)$, together with the Stirling's approximation \Cref{lem:stirling}, we have
\begin{align}
    &\, \log\P(\rD_x^{(1)} = \alpha d) = \log\binom{m}{\alpha d} + \alpha d\log p + (m - \alpha d) \log(1 - p)\\
    = &\, \frac{1}{2}\log \frac{m}{2\pi \alpha d(m-\alpha d)} + m\log(m) - \alpha d\log(\alpha d) - (m-\alpha d)\log(m-\alpha d) \\
    &\, + \alpha d \log \frac{d}{\sqrt{mn}} + (m - \alpha d) \log(1 - d/\sqrt{mn} ) + O\big( (\alpha d)^{-1}\big)\\
    = &\, -\frac{1}{2}\log(2\pi \alpha d) - \alpha d\log(\qr^{2} \alpha) -  m \log(1 -\alpha d/m) + m \log(1 - d/\sqrt{mn} ) + O\big( (\alpha d)^{-1}\big)\\
    = &\, -d\Big( \alpha \log(\qr^{2}\alpha) - \alpha + \ratio^{-1/2}\Big) - \frac{1}{2}\log(2\pi \alpha d) = -f_{\qr, d} (\alpha), 
\end{align}
where in the last line, we used the fact $\ratio = \qr^{4}$. Similarly,
\begin{align}
    &\,\log\P(\rD_x^{(2)} = \alpha d ) = \log\binom{n}{\alpha d} + \alpha d\log p + (n - \alpha d) \log(1 - p)\\
    = &\, -\frac{1}{2}\log(2\pi \alpha d) - \alpha d\log(\alpha/\qr^{2}) -  n \log(1 -\alpha d/n) + n \log(1 - d/\sqrt{mn} ) + O\big( (\alpha d)^{-1}\big)\\
    = &\, -d\Big( \alpha \log(\qr^{-2}\alpha) - \alpha + \qr^{2}\Big) - \frac{1}{2}\log(2\pi \alpha d) = - f_{\qr, d}(\alpha),
\end{align}
which completes the proof.
\end{proof}

\section{Spectral properties of biregular trees}\label{sec:tri_regular}
\begin{quote}
	\center{\emph{A computational trick can also be a theoretical trick.}}
		\begin{flushright}
			\small{-- Alan Edelman's 60th birthday\\ July 2023 at MIT}
		\end{flushright}
\end{quote}
Tridiagonalization is a powerful technique in numerical linear algebra that reduces a matrix to a simpler tridiagonal form while preserving its essential spectral properties. This method is particularly useful for analyzing the eigenvalues and eigenvectors of matrices, and plays a key role in random matrix theory, for example, in the construction of $\beta$-ensembles~\cite{dumitriua2002matrix}. In this section, we apply the tridiagonalization technique to study the spectral properties of the idealized tree $\tree$ in \Cref{def:biregular_tree}. For convenience, let $\rmA$ denote the adjacency matrix of $\tree$, and let $r \equiv r_x$ in \eqref{eqn:defrx}.

\subsection{Tridiagonal representation}
We study the spectral properties of $\tree$ by computing its tridiagonal representation, following the standard construction in \cite{trotter1984eigenvalue}.

\begin{lemma}[Tridiagonal representation]\label{lem:tridiagrep}
Let $\rvs_0, \rvs_1, \rvs_2, \ldots, \rvs_{r}$ denote the vectors obtained from $\ones_{x}$, $\rmA\ones_x$, $\rmA^2\ones_x$, $\ldots$, $\rmA^{r}\ones_x$ after Gram-Schmidt orthonormalization. Then the following holds.
\begin{enumerate}
    \item For all $j = 0, 1, \ldots, r$, we have
        \begin{align}
            \rvs_j = |\rvs_j(x)|^{-1/2} \ones_{\layer_j(x)}.
        \end{align}
    \item Let $\rvs_{r + 1}, \ldots, \rvs_{N-1}$ be a completion of $\rvs_0, \ldots, \rvs_{r}$ such that $\{\rvs_j\}_{j=0}^{N-1}$ form an orthonormal basis of $\R^{N}$. Define the basis $\rmS \coloneqq [\rvs_0, \rvs_1, \ldots, \rvs_{N-1}]$ and let
    \begin{align}
        \rmZ \coloneqq \frac{1}{\sqrt{d}} \rmS^{\sT} \rmA \rmS
    \end{align}
    denote the representation of $\rmA/\sqrt{d}$ in the basis $\rmS$. Then the upper-left $(r +1)\times (r +1)$ block $\rmZ_{\llbracket r \rrbracket}$ of $\rmZ$ has the following tridiagonal form
    \begin{align}
        \rmZ_{\llbracket r \rrbracket}(\alpha) = \indi{x\in \gV_1} \rmZ^{(1)}(\alpha_x ) + \indi{x\in \gV_2} \rmZ^{(2)}(\alpha_x),
    \end{align}
    where for $r$ being even, $\rmZ^{(1)}(\alpha)$ and $\rmZ^{(2)}(\alpha)$ are $(r+1)\times (r+1)$ matrices defined through
\begin{align}\label{eqn:Z1Z2}
    \rmZ^{(1)}(\alpha) \coloneqq \begin{bmatrix}
   0  & \sqrt{\alpha} & \\
   \sqrt{\alpha} & 0 & \qr\\
       &  \qr &  0 &  \qr^{-1}\\
        & &\qr^{-1} & 0 & \ddots\\
        & &  &  \ddots &  \ddots & \qr \\
        & &  & &   \qr & 0
    \end{bmatrix}, \quad
    \rmZ^{(2)}(\alpha)=\begin{bmatrix}
   0  & \sqrt{\alpha} & \\
   \sqrt{\alpha} & 0 & \qr^{-1}\\
       &  \qr^{-1} &  0 &  \qr\\
        & &\qr & 0 & \ddots\\
        & &  &  \ddots &  \ddots & \qr^{-1} \\
        & &  & &   \qr^{-1} & 0
    \end{bmatrix}.
\end{align}
The bottom right corner of $\rmZ^{(1)}(\alpha)$ (resp. $\rmZ^{(2)}(\alpha)$) changes to $\qr^{-1}$ (resp. $\qr$) when $r$ is odd. 
\end{enumerate}
\end{lemma}

\begin{proof}[Proof of \Cref{lem:tridiagrep}]
\emph{(1)}. Let $\rmQ_j$ denote the orthogonal projection onto the orthogonal complement of the space spanned by $\{\rmA^l \ones_x\}_{l=0}^{j-1}$ for each $1 \leq j \leq r$ and let $\rmQ_0 = \id_N$. The proof of \emph{(1)} is then established by inductively showing that
\begin{align}\label{eqn:inductionSi}
    \ones_{\layer_{j}(x)} = \rmQ_j \rmA^j \ones_x, \quad j=0, 1, \ldots, r. 
\end{align}
The base case is trivial since $\layer_0(x) = \{x\}$.
The following recurrence relation is established by considering the tree structure of ${\tree}$ in \Cref{def:biregular_tree}, 
\begin{align}\label{eqn:AoneSl}
    \rmA \ones_{\layer_l(x)} = \left\{
    \begin{aligned}
        &\ones_{\layer_1(x)}, &\, x\in \gV_{1} \cup \gV_{2}, \,\,l=0,\\
        &\ones_{\layer_2(x)} + \rD_x \ones_x, &\, x\in \gV_{1} \cup \gV_{2}, \,\, l=1,\\
        &\ones_{\layer_{l+1}(x)} + \indi{ l \textnormal{ even} } d_2\ones_{\layer_{l-1}(x)} + \indi{ l \textnormal{ odd} }d_1\ones_{\layer_{l-1}(x)}, &\,  x\in \gV_1, \,\, l\in [r -1]\setminus \{1\}\\
        &\ones_{\layer_{l+1}(x)} + \indi{ l \textnormal{ even} } d_1\ones_{\layer_{l-1}(x)} + \indi{ l \textnormal{ odd} }d_2\ones_{\layer_{l-1}(x)}, &\,  x\in \gV_2, \,\, l\in [r -1]\setminus \{1\}.
    \end{aligned}
    \right.,
\end{align}
where $d_1 \coloneqq \qr^{-2}d$ and $d_2 \coloneqq \qr^{2} d$ denote the degrees of vertices in $\gV_1$ and $\gV_2$ respectively.

 As for the induction step, suppose that \eqref{eqn:inductionSi} holds for the case $j-1$ when $j\geq 1$, then the orthogonality between $\rmQ_{j}$ and $\{\rmA^l \ones_x\}_{l=0}^{j-1}$ leads to
    \begin{align}
        \rmQ_{j} \rmA^{j} \ones_x = \rmQ_{j} \rmA \ones_{\layer_{j-1}(x)} = \rmQ_{j}\ones_{\layer_j(x)} = \ones_{\layer_j(x)},
    \end{align}
where $\ones_{\layer_j(x)}^{\sT}\ones_{\layer_l(x)} = 0$ for $l \neq j$ since $|\layer_j(x) \cap \layer_l(x)| = \emptyset$.

\noindent \emph{(2)} The proof of $\emph{(2)}$ is established by finishing the following computation
\begin{align}
    \<\rvs_i, \rvs_j\> = |\layer_i(x)|^{-1/2} |\layer_j(x)|^{-1/2} \<\ones_{\layer_i(x)}, \rmA \ones_{\layer_j(x)}\>.
\end{align}
According to \eqref{eqn:AoneSl}, $\<\rvs_i, \rvs_j\> = 0$ for $|i - j| \neq 1$. Moreover, $|\layer_0(x)| = 1$, and for $0\leq l \leq \lfloor r/2 \rfloor - 1$, $|\layer_{2l + 1}(x)| = \rD_x (d_1 d_2)^l$, while $|\layer_{2l + 2 }(x)| = \rD_x d_2(d_1 d_2)^l$ for $x\in \gV_1$ and $|\layer_{2l + 2 }(x)| = \rD_x d_1(d_1 d_2)^l$ for $x\in \gV_2$. We then have for $x \in \gV_1$
\begin{align}
    \<\rvs_i, \rmA \rvs_j\> = \left\{ 
    \begin{aligned}
        & \sqrt{\rD_x}, &\, |i - j| = 1 \textnormal{ and one of } i, j \textnormal{ is } 0,\\
        & \sqrt{d_1}, &\, i = 2l = j - 1, \textnormal{ or }  i = 2l + 1 = j + 1, \\
        & \sqrt{d_2}, &\, i = 2l = j + 1, \textnormal{ or }  i = 2l + 1 = j - 1,\\
        &0, &\, |i - j| \neq 1,
    \end{aligned}
    \right.
\end{align}
while for $x\in \gV_2$, we have
\begin{align}
    \<\rvs_i, \rmA \rvs_j\> = \left\{ 
    \begin{aligned}
        & \sqrt{\rD_x}, &\, |i - j| = 1 \textnormal{ and one of } i, j \textnormal{ is } 0,\\
        & \sqrt{d_1}, &\, i = 2l = j + 1, \textnormal{ or }  i = 2l + 1 = j - 1, \\
        & \sqrt{d_2}, &\, i = 2l = j - 1, \textnormal{ or }  i = 2l + 1 = j + 1,\\
        &0, &\, |i - j| \neq 1.
    \end{aligned}
    \right.
\end{align}
By factoring out $\sqrt{d}$, for $r$ even, the upper-left $(r +1)\times (r +1)$ block $\rmZ_{\llbracket r \rrbracket}$ of $\rmZ$ becomes
\begin{align}
    \indi{x\in \gV_1} \begin{bmatrix}
   0  & \sqrt{\alpha_x} & \\
   \sqrt{\alpha_x} & 0 & \qr\\
      & \qr &  0 &  \qr^{-1}\\
      &  &\qr^{-1} & 0 & \ddots\\
      &  & &   \ddots & \ddots & \qr\\
      &  & &    &  \qr & 0
    \end{bmatrix} + \indi{ x\in \gV_2} \begin{bmatrix}
   0  & \sqrt{\alpha_x} & \\
   \sqrt{\alpha_x} & 0 & \qr^{-1}\\
      & \qr^{-1} &  0 &  \qr\\
      &  &\qr & 0 & \ddots\\
      &  & &   \ddots & \ddots & \qr^{-1}\\
      &  & &    &  \qr^{-1} & 0
    \end{bmatrix},
\end{align}
where we use the facts $\alpha_x = \rD_x/d$, $d_1 = d/\qr^2$ and $d_2 = d\qr^2$. The bottom corner of the first (resp. second) matrix changes to $\qr^{-1}$ (resp. $\qr$) for the case $r$ being odd.
\end{proof}



\subsection{Spectral properties}
Note that $\rmZ$ and $\rmA /\sqrt{d}$ share the same spectrum. However, only $\rmZ_{\llbracket r \rrbracket}$ plays an important role in our spectral analysis.

\begin{lemma}[Approximate eigenpairs of $\rmZ^{(1)}(\alpha)$ and $ \rmZ^{(2)}(\alpha)$]\label{lem:eigenZ1Z2}
Let $\lambda^{(1)}_{1} \geq \cdots \geq \lambda^{(1)}_{r + 1}$ and $\lambda^{(2)}_{1} \geq \cdots \geq \lambda^{(2)}_{r + 1}$ denote the eigenvalues of $\rmZ^{(1)}(\alpha)$ and $\rmZ^{(2)}(\alpha)$, respectively.  The following holds.
\begin{enumerate}
    \item Bulk eigenvalues:
    \begin{align}
        \big\{ \{\lambda^{(1)}_{l}\}_{l=2}^{r} \cup \{\lambda^{(2)}_{l}\}_{l=2}^{r} \big\} \subset [-\qr -\qr^{-1},-\qr +\qr^{-1}]\cup [\qr-\qr^{-1}, \qr +\qr^{-1}]
    \end{align}
    \item Extreme eigenvalues:
    \begin{enumerate}
        \item If $\alpha>\qr^{-2}+1$, $\lambda^{(1)}_{1}$ and $\lambda^{(1)}_{r + 1}$ converge to $\Lambda_{\qr^{-1}}(\alpha)$ and $ -\Lambda_{\qr^{-1}}(\alpha)$, respectively, as $r\to\infty$.
        \item If $\alpha>\qr^2+1$, $\lambda^{(2)}_{1}$ and $\lambda^{(2)}_{r + 1}$ converge to $ \Lambda_{\qr}(\alpha)$ and $ -\Lambda_{\qr}(\alpha)$, respectively, as $r\to\infty$.
        \item If $0<\alpha<\qr^2-1$, the smallest positive eigenvalue and largest negative eigenvalue of $\rmZ^{(2)}(\alpha)$ converge to $\Lambda_{\qr}(\alpha)$ and $-\Lambda_{\qr}(\alpha)$, respectively, as $r\to\infty$.
    \end{enumerate}
        
    \item Approximate eigenvectors:
        \begin{enumerate}
            \item For $x\in \gV_{1}$ and $\alpha > \qr^{-2} + 1$, let $\rvu_{+} = \{\ervu_{j}\}_{j=0}^{r}$, $\rvu_{-} = \{(-1)^{j}\ervu_{j}\}_{j=0}^{r}$ have the components
            \begin{align}\label{eqn:Z1Eigenvec}
                     \ervu_0 \in \R \setminus \{0\}, \quad \ervu_1 &\,= \alpha^{-1/2} \Lambda_{\qr^{-1}} (\alpha)\ervu_0, \quad \ervu_2 = \alpha^{1/2}\qr(\alpha - \qr^{-2})^{-1} \ervu_0\,,\\
                      \ervu_{2j+1} &\, = (\alpha - \qr^{-2})^{-j} \ervu_1, \quad \ervu_{2j+2}= (\alpha -\qr^{-2})^{-j} \ervu_2, \quad 1 \leq j \leq \lfloor \frac{r-3}{2} \rfloor.
            \end{align}
            Here, $\ervu_0$ is chosen to ensure $\|\rvu_{+}\|_2 = 1$. Then as $r\to \infty$, $\rvu_{+}$ and $\rvu_{-}$ are the approximate eigenvectors of $\rmZ^{(1)}(\alpha)$ corresponding to $\Lambda_{\qr^{-1}}(\alpha)$ and $ -\Lambda_{\qr^{-1}}(\alpha)$, respectively. 
            \item For $x\in \gV_{2}$, and $\alpha > \qr^{2} + 1$ or $0<\alpha < \qr^{2} - 1$, define $\rvu_{+} = \{\ervu_{j}\}_{j=0}^{r}$ and $\rvu_{-} = \{(-1)^{j}\ervu_{j}\}_{j=0}^{r}$ by
                \begin{align}\label{eqn:Z2Eigenvec}
                    \ervu_0 \in \R \setminus \{0\}, \quad \ervu_1 &\,= \alpha^{-1/2} \Lambda_{\qr} (\alpha)\ervu_0, \quad \ervu_2 = \alpha^{1/2}\qr^{-1}(\alpha - \qr^2)^{-1} \ervu_0\,,\\
                     \ervu_{2j+1} &\,= (\alpha - \qr^2)^{-j} \ervu_1, \quad \ervu_{2j+2} = (\alpha -\qr^2)^{-j} \ervu_2, \quad 1 \leq j \leq \lfloor (r-3)/2 \rfloor.
                \end{align}
            Then as $r\to \infty$, $\rvu_{+}$ and $\rvu_{-}$ are the approximate eigenvectors of $\rmZ^{(2)}(\alpha)$ corresponding to $\Lambda_{\qr}(\alpha)$ and $ -\Lambda_{\qr}(\alpha)$, respectively.    
        \end{enumerate}
\end{enumerate}
\end{lemma}

\begin{proof} [Proof of \Cref{lem:eigenZ1Z2} (1)]
Let $\rmZ^{(1)}_{n}(\alpha)$ and $\rmZ^{(2)}_{n}(\alpha)$ denote the $n \times n$ matrices of the forms \eqref{eqn:Z1Z2}. Consider $\rmZ^{(1)}_{n}(\qr^{-2})$ and $\rmZ^{(2)}_{n}(\qr^{2})$. We claim that
\begin{claim}\label{clm:ZkappaSpectra}
    $\mathrm{Spec}\big(\rmZ^{(1)}_{n}(\qr^{-2})\big) \cup \mathrm{Spec}\big(\rmZ^{(2)}_{n}(\qr^{2})\big) \subset [-\qr -\qr^{-1},-\qr +\qr^{-1}]\cup [\qr-\qr^{-1}, \qr +\qr^{-1}]$.
\end{claim}
Since the matrix $\rmZ^{(1)}_{n}(\alpha) - \rmZ^{(1)}_{n}(\qr^{-2})$ (resp. $\rmZ^{(2)}_{n}(\alpha) - \rmZ^{(2)}_{n}(\qr^{2})$) has only two nontrivial eigenvalues, one positive and one negative, the proof is then completed by Weyl's interlacing inequality \Cref{lem:weyl}.
\end{proof}

\begin{proof} [Proof of the Claim \ref{clm:ZkappaSpectra}]
By calculating the determinants, we obtain the following recurrence relationship
\begin{align}
    \det(\rmZ^{(1)}_{n}(\qr^{-2}) - x\id) =&\, -x\det(\rmZ^{(2)}_{n-1}(\qr^{2}) - x \id) - \qr^{-2} \det(\rmZ^{(1)}_{n-2}(\qr^{-2}) - x\id),\\
    \det(\rmZ^{(2)}_{n}(\qr^{2}) - x\id) =&\, -x\det(\rmZ^{(1)}_{n-1}(\qr^{-2}) - x \id) - \qr^{2} \det(\rmZ^{(2)}_{n-2}(\qr^{2}) - x\id).
\end{align}
Denote $f_{n}(x) \coloneqq \det(\rmZ^{(1)}_{n}(\qr^{-2}) - x\id)$ and $g_{n}(x) \coloneqq \det(\rmZ^{(2)}_{n}(\qr^{2}) - x\id)$. We will analyze the roots of $f_n(x) = 0$ and $g_n(x) = 0$, which in turn allows us to understand the spectrum of $\rmZ^{(1)}_{n}(\qr^{-2})$ and $\rmZ^{(2)}_{n}(\qr^{2})$. By simply re-indexing, we have
\begin{subequations}
\begin{align}
    f_{2n}(x) =&\, -x g_{2n-1}(x) - \qr^{-2}f_{2n-2}(x),\label{eqn:f-rec}\\
    g_{2n-1}(x) =&\, -x f_{2n-2}(x) - \qr^{2}g_{2n-3}(x).\label{eqn:g-rec}
\end{align}
\end{subequations}
From \eqref{eqn:g-rec}, we then have
\begin{align}
    f_{2n-2}(x) = -\frac{1}{x}\big[ g_{2n-1}(x) + \qr^{2}g_{2n-3}(x) \big], \quad f_{2n}(x) = -\frac{1}{x}\big[ g_{2n+1}(x) + \qr^{2}g_{2n-1}(x) \big],
\end{align}
By substituting the results above into \eqref{eqn:f-rec}, it implies that
\begin{align}
    &\, g_{2n+1}(x) + (\qr^2 + \qr^{-2} - x^2) g_{2n-1}(x) + g_{2n-3}(x) = 0.
\end{align}
Similarly, from the recursive equations \eqref{eqn:f-rec} and \eqref{eqn:g-rec}, we can deduce that
\begin{align}
     &\, g_{2n+2}(x) + (\qr^2 + \qr^{-2} - x^2) g_{2n}(x) + g_{2n-2}(x) = 0.
\end{align}


Notice that $g_2(x)=x^2-\qr^2$, $g_1(x)=-x$ and $f_2(x)=x^2-\qr^{-2}$, $f_1(x)=-x$. Let $g_0(x) = f_0(x) = 1$. Define $\phi_n(y) \coloneqq g_{2n}(x)$ where we substitute $y \coloneqq (x^2 - \qr^2 - \qr^{-2})/2$, then
\begin{align}
    \phi_{n+1}(y) &\,= 2y\, \phi_{n}(y)-\phi_{n-1}(y), \quad n\geq 1\\
    \phi_1(y) &\,= 2y + \qr^{-2},\quad \phi_0(y) =1, 
\end{align}
 which is a generalization of the Chebyshev polynomials of the fourth kind, see \cite[Section 1.3.2]{mason2002chebyshev}. In fact, the fourth kind requires $\phi_1(y) = 2y + 1$, however $\qr^{-2} \leq 1$ in our scenario.
 
In the following, we will analyze the roots of polynomials $\phi_{n}(y) = 0$. Let $\Phi(y, t) = \sum_{n = 0}^{\infty} \phi_n(y) t^n$ denote the generating function. By plugging in the recurrence formula, one obtain
\begin{align}
    \Phi(y, t) = 1 + (2y + \qr^{-2})t + 2yt (\Phi(y, t) - 1) - t^{2} \Phi(y, t).
\end{align}
We substitute $y = \cos \theta$. By Euler's formula $e^{\ii \theta} = \cos \theta + \ii \sin \theta$ where $\ii = \sqrt{-1}$, we have $2y = e^{\ii \theta} + e^{-\ii \theta}$, then by geometric series expansion,
\begin{align}
     \Phi(y, t) =&\, \frac{k^{-2}t + 1}{t^2 - 2yt + 1} = \frac{1}{e^{-\ii \theta} - e^{\ii \theta}} \Big( \frac{1 + \qr^{-2}e^{-\ii \theta}}{t - e^{-\ii \theta}} - \frac{\qr^{-2}e^{\ii \theta} + 1}{t - e^{\ii \theta}} \Big)\\
     =&\,\sum_{n=0}^{\infty} \frac{1}{e^{\ii \theta} - e^{-\ii \theta}}\Big( \frac{1 + \qr^{-2}e^{-\ii \theta}}{e^{-\ii (n+1)\theta}} - \frac{\qr^{-2}e^{\ii \theta} + 1}{e^{\ii (n+1)\theta}} \Big) t^n\\
     =&\, \sum_{n=0}^{\infty} \Big( e^{-\ii n\theta}\sum_{l=0}^{n} e^{\ii 2l \theta}
     + \qr^{-2} e^{-\ii (n-1)\theta}\sum_{l=0}^{n-1} e^{\ii 2l \theta} \Big) t^n\\
     =&\, \sum_{n=0}^{\infty} \Big( e^{-\ii n\theta} (1 + \qr^{-2} e^{\ii \theta}) \cdot \big(\sum_{l=0}^{n-1} e^{\ii 2l \theta} \big) + e^{\ii n \theta} \Big) t^n,
 \end{align}
where in the second to last line we use the fact $a^n - b^n = (a - b)(a^{n-1} + a^{n-2}b + \ldots + b^{n-1})$ twice, which gives rise to the explicit formula
 \begin{align}
     \phi_n(y) = (1 + \qr^{-2}e^{\ii \theta})e^{-\ii n\theta}\frac{1 - e^{\ii 2n\theta}}{1 - e^{\ii 2\theta}} + e^{\ii n\theta}, \quad y = \cos \theta. 
 \end{align}
Applying Euler's formula again, the roots for $\phi_n(y) = 0$ can be obtained by solving
 \begin{align}
      1 + \qr^{-2} e^{\ii \theta} = e^{\ii (2n+2)\theta}(1 + \qr^{-2}e^{-\ii \theta}) \iff \sin[(n+1)\theta ] + \qr^{-2}\sin (n\theta) = 0.
 \end{align}
Furthermore, the function $h_1(\theta) = \sin[(n+1)\theta ]$ intersects with $h_2(\theta) = -\qr^{-2}\sin (n\theta)$ exactly once in $(l\pi/n, (l+1)\pi/n)$ for each $l = 0, \ldots, 2n-1$, while those two functions differ at the endpoints of each interval. Thus, the roots $\theta_0, \ldots, \theta_{2n-1}$ for $\phi_n(y) = 0$ are all real, indicating $|y_l| = |\cos \theta_l| \leq 1$ for each $0 \leq l \leq 2n-1$. Recall that $y \coloneqq (x^2 - \qr^2 - \qr^{-2})/2$, then each eigenvalue $x$ of $\rmZ^{(2)}_{2n}(\qr^2)$ satisfies $x \in [-\qr - \qr^{-1}, -\qr + \qr^{-1}] \cup  [\qr - \qr^{-1}, \qr + \qr^{-1}]$. The analysis for the spectrum of $\rmZ^{(1)}_{2n}(\qr^{-2})$ follows similarly.
\end{proof}

\begin{proof}[Proof of \Cref{lem:eigenZ1Z2} (2) and (3)]
The eigenpairs of $\rmZ^{(1)}(\alpha)$ and $\rmZ^{(2)}(\alpha)$ can be analyzed via the transfer matrix approach. 

    We focus on $\rmZ^{(2)}(\alpha)$ first. Let $(\eta, \rvu)$ be an eigenpair of $\rmZ^{(2)}(\alpha)\in \R^{(r+1)\times (r+1)}$ with $\rvu = [ \ervu_0, \cdots, \ervu_{r}]^{\sT} \in \R^{r + 1}$.  The recurrence relation between entries in $\rvu$ can be derived from $\rmZ^{(2)} \rvu = \eta \rvu$ and~\eqref{eqn:Z1Z2}, 
\begin{align}
    \sqrt{\alpha} \ervu_1  =&\, \eta \ervu_0, \quad  \sqrt{\alpha} \ervu_0 +\qr^{-1} \ervu_2 = \eta \ervu_1, \label{eqn:Z2recurrence}
    \\ 
    \qr^{-1} \ervu_{2j-1} + \qr \ervu_{2j+1} =&\, \eta \ervu_{2j},\quad \qr \ervu_{2j}+\qr^{-1} \ervu_{2j+2} = \eta \ervu_{2j+1}, \quad j = 1, \ldots, \lfloor (r-3)/2 \rfloor,
\end{align}
where $\qr^{-1}\ervu_{r-1} = \eta u_{r}$ for $r$ even, and $\qr \ervu_{r-1} = \eta u_{r}$ for $r$ odd. Then the recurrence can be written as
\begin{subequations}
    \begin{align}
  \begin{bmatrix}
      \ervu_{2j+1}      \\
      \ervu_{2j} 
   \end{bmatrix} = \,& \begin{bmatrix}
     \qr^{-1}\eta  &   -\qr^{-2} \\
     1  & 0
   \end{bmatrix}  \begin{bmatrix}
    \ervu_{2j}      \\
       \ervu_{2j-1} 
   \end{bmatrix}= \rmT^{(1)}(\eta) \begin{bmatrix}
    \ervu_{2j}      \\
       \ervu_{2j-1} 
   \end{bmatrix},\\
   \begin{bmatrix}
    \ervu_{2j+2}      \\
       \ervu_{2j+1} 
   \end{bmatrix} = \,& 
   \begin{bmatrix}
     \qr \eta  &    -\qr^{2} \\
     1  & 0
   \end{bmatrix}  
   \begin{bmatrix}
    \ervu_{2j+1}      \\
       \ervu_{2j} 
   \end{bmatrix}= \rmT^{(2)}(\eta) \begin{bmatrix}
    \ervu_{2j+1}      \\
       \ervu_{2j} 
   \end{bmatrix}.
\end{align}
\end{subequations}
Therefore, the transfer matrix $\rmT(\eta)$ is defined as
\begin{align}
   \rmT(\eta)\coloneqq \rmT^{(2)}(\eta) \rmT^{(1)}(\eta)=
   \begin{bmatrix}
     \eta^2-\qr^{2}  &   - \qr^{-1}\eta \\
     \qr^{-1}\eta  & -\qr^{-2}
   \end{bmatrix}. \label{eqn:T(eta)}
\end{align}
Then for $j = 1, \ldots, \lfloor (r-3)/2 \rfloor$,
\begin{align}
  \begin{bmatrix}
    \ervu_{2j+2}      \\
    \ervu_{2j+1} 
   \end{bmatrix}  = \rmT(\eta) \begin{bmatrix}
    \ervu_{2j}      \\
    \ervu_{2j-1} 
   \end{bmatrix}  =  [\rmT(\eta)]^{j} \begin{bmatrix}
        \ervu_{2}      \\
        \ervu_{1} 
   \end{bmatrix}. \label{eqn:T(eta)Recurrence}
\end{align}
The eigenvalues of $\rmT(\eta)$ are then given by
\begin{align}\label{eqn:lambda}
    \lambda_{\pm}(\eta) = &\, \frac{1}{2}\Big(  \eta^2 - (\qr^2 + \qr^{-2} ) \pm \sqrt{[\eta^{2} -(\qr^{-1} - \qr)^{2}][\eta^{2} - (\qr^{-1} + \qr)^{2}] }\Big). 
\end{align}
Note that $\eta \in \R$ since $\rmZ^{(2)}(\alpha)$ is Hermitian, then $\rvu \in \R^{r + 1}$. Thus, we also deduce that $\lambda_{\pm}(\eta) \in \R$ and the following condition has to be satisfied
\begin{align}
    |\eta|\geq \qr+\qr^{-1}\quad \textnormal{ or } \quad |\eta|\leq \qr-\qr^{-1}.
\end{align}
The eigenvectors corresponding to $\lambda_{+}(\eta)$ and $\lambda_{-}(\eta)$ are given by
\begin{align}
    \rvv_{+}(\eta) =\begin{bmatrix}
    \frac{\qr\lambda_{+}(\eta)}{\eta}+\frac{1}{\qr\eta}      \\
       1 
   \end{bmatrix} ,\quad  \rvv_{-}(\eta) =\begin{bmatrix}
    \frac{\qr\lambda_{-}(\eta)}{\eta}+\frac{1}{\qr\eta}      \\
       1 
   \end{bmatrix}.
\end{align}

We take $[\ervu_2, \ervu_1]^{\sT}$ which also satisfies \eqref{eqn:T(eta)Recurrence}. Then, by solving $[\ervu_2, \ervu_1]^{\sT}\in \mathrm{Ker}(\rmT(\eta) - \lambda_{-}(\eta)\id_2)$, we have
\begin{align}
    \frac{\ervu_2}{\ervu_1} = \frac{\qr\lambda_{-}}{\eta}+\frac{1}{\qr\eta} = \frac{\qr^2\lambda_{-} + 1}{\qr\eta} = \frac{\qr(\eta^2 - \alpha)}{\eta},
\end{align}
where the last equality comes from \eqref{eqn:Z2recurrence}, leading to the following explicit formula of $\eta$
\begin{align}
    \eta^2 = \frac{\alpha^2\qr^2 + \alpha(1 - \qr^4) }{\qr^2(\alpha - \qr^2)} = \alpha + \frac{1}{\qr^2} + \frac{1}{\alpha - \qr^2} = [\Lambda_{\qr}(\alpha)]^2,
\end{align}
and the conditions $|\eta| \geq \qr + \qr^{-1}$ and $|\eta| \leq \qr - \qr^{-1}$ are satisfied since $(\alpha - \qr^2) + (\alpha - \qr^2)^{-1} > 2$ (resp. $< -2$) when $\alpha > \qr^2 + 1$ (resp. $0 < \alpha < \qr^2 - 1$). 

We now construct the approximate eigenvector $\rvu$ corresponding to $\eta$. Note that $\lambda_{-}(\eta) = (\alpha - \qr^{2})^{-1}$ when $\eta = \Lambda_{\qr} (\alpha)$. Let $\ervu_0 \in \R \setminus \{0\}$, then \eqref{eqn:Z2Eigenvec} is obtained by applying \eqref{eqn:Z2recurrence}. Let $\{\rve_{l}\}_{l = 0}^{r}$ denote the set of standard basis vectors of $\R^{r + 1}$. According to \eqref{eqn:Z2recurrence} and \eqref{eqn:T(eta)Recurrence},
\begin{align}
    (\rmZ^{(2)}(\alpha) - \Lambda_{\qr} (\alpha) \rmI_{r + 1}) \rvu = \Big((\qr^{-1} \indi{r \textnormal{ even}} + \qr \indi{r \textnormal{ odd}}) \cdot \ervu_{r-1} - \Lambda_{\qr} (\alpha) \ervu_r \Big) \rve_{r}.
\end{align}
Hence $\|(\rmZ^{(2)}(\alpha) - \Lambda_{\qr} (\alpha) \rmI_{r + 1}) \rvu\| \to 0$ as $r \to \infty$ since $\alpha > \qr^{2} + 1$ or $\alpha < \qr^{2} - 1$ for $\rmZ^{(2)}(\alpha)$, as well as the exponential decay in \eqref{eqn:Z2Eigenvec}. Meanwhile, combining the conclusion in (1), as $r \to \infty$,
\begin{enumerate}
    \item $\Lambda_{\qr} (\alpha)$ is the limit of $\lambda^{(2)}_{1}$, since $\Lambda_{\qr} (\alpha) > \qr + \qr^{-1}$ when $\alpha > \qr^{2} + 1$.
    \item $\Lambda_{\qr} (\alpha)$ is the limit of the smallest positive eigenvalue, since $\Lambda_{\qr} (\alpha) < \qr - \qr^{-1}$ when $\alpha < \qr^{2} - 1$.
\end{enumerate}
This completes the proofs of (b) (c) in (2) and (b) in (3).

~\\
We now turn to $\rmZ^{(1)}(\alpha)$. Similarly, let $(\zeta,\rvu)$ be an eigenpair of $\rmZ^{(1)}(\alpha)$. The recurrence relation now is
\begin{align}
    \sqrt{\alpha} \ervu_1  =&\, \zeta \ervu_0, \quad  \sqrt{\alpha} \ervu_0 +\qr \ervu_2 = \zeta \ervu_1, \label{eqn:Z1recurrence}
    \\ 
    \qr \ervu_{2i-1} + \qr^{-1}\ervu_{2i+1} =&\, \zeta \ervu_{2i},\quad \qr^{-1}\ervu_{2i}+\qr \ervu_{2i+2} = \zeta \ervu_{2i+1}, \quad i = 1, \ldots, \lfloor (r-3)/2 \rfloor,
\end{align}
where $\qr \ervu_{r-1} = \eta u_{r}$ for $r$ even, and $\qr^{-1} \ervu_{r-1} = \eta u_{r}$ for $r$ odd. Then
\begin{align}
  \begin{bmatrix}
      \ervu_{2i+1}      \\
      \ervu_{2i} 
   \end{bmatrix} = \rmT^{(2)}(\zeta) \begin{bmatrix}
    \ervu_{2i}      \\
       \ervu_{2i-1} 
   \end{bmatrix},\quad 
   \begin{bmatrix}
    \ervu_{2i+2}      \\
       \ervu_{2i+1} 
   \end{bmatrix} = \rmT^{(1)}(\zeta) 
   \begin{bmatrix}
    \ervu_{2i+1}  \\
    \ervu_{2i} 
   \end{bmatrix}.
\end{align}
Consequently,
\begin{align}
  \begin{bmatrix}
    \ervu_{2i+2}      \\
    \ervu_{2i+1} 
   \end{bmatrix}  = \widetilde{\rmT}(\zeta) \begin{bmatrix}
    \ervu_{2i}      \\
    \ervu_{2i-1} 
   \end{bmatrix}  =  [\widetilde{\rmT}(\zeta)]^{i} \begin{bmatrix}
        u_{2}      \\
        u_{1} 
   \end{bmatrix},
\end{align}
where we define
\begin{align}
   \widetilde{\rmT}(\zeta)\coloneqq \rmT^{(1)}(\zeta) \rmT^{(2)}(\zeta)=
   \begin{bmatrix}
     \zeta^2-\qr^{-2}  &   - \qr\zeta \\
     \qr\zeta  & -\qr^{2}
   \end{bmatrix}, \label{eqn:Ttilde(eta)Recurrence}
\end{align}
with the eigenvalues given by $\lambda_{\pm}(\zeta)$ in \eqref{eqn:lambda}. Again, the condition $|\zeta|>\qr+\qr^{-1}$ or $|\zeta|<\qr-\qr^{-1}$ should be satisfied to ensure $\lambda_{\pm}(\zeta)$ being real. The corresponding eigenvectors are
\begin{align}
    \rvv_{+}(\zeta) =\begin{bmatrix}
    \frac{\lambda_{+}(\zeta)}{\qr\zeta}+\frac{\qr}{\zeta}      \\
       1 
   \end{bmatrix} ,\quad  \rvv_{-}(\zeta) =\begin{bmatrix}
    \frac{\lambda_{-}(\zeta)}{\qr\zeta}+\frac{\qr}{\zeta}      \\
       1 
   \end{bmatrix}
\end{align}
Let $[\ervu_2, \ervu_1]^{\sT} \in \mathrm{ker}(\widetilde{\rmT}(\zeta) - \lambda_{-}(\zeta) \id_2)$, then
\begin{align}
    \frac{\ervu_2}{\ervu_1} = \frac{\lambda_{-}(\zeta)}{\qr \zeta}+\frac{\qr}{\zeta} = \frac{\zeta^2 - \alpha}{\qr\zeta},
\end{align}
where the last equality is obtained from \eqref{eqn:Z1recurrence}, leading to the solution
\begin{align}
    \zeta^2 = \frac{\alpha^2\qr^{-2} + \alpha(1 - \qr^{-4}) }{\qr^{-2}(\alpha - \qr^{-2})} = \alpha + \qr^2 + \frac{1}{\alpha - \qr^{-2}} = [\Lambda_{\qr^{-1}}(\alpha)]^2.
\end{align}

We now construct the approximate eigenvector $\rvu$ corresponding to $\zeta$. Note that $\lambda_{-}(\zeta) = (\alpha - \qr^{-2})^{-1}$ when $\zeta = \Lambda_{\qr^{-1}} (\alpha)$. Let $\ervu_0 \in \R \setminus \{0\}$, then \eqref{eqn:Z1Eigenvec} is obtained by applying \eqref{eqn:Z1recurrence}. According to \eqref{eqn:Z1recurrence} and \eqref{eqn:Ttilde(eta)Recurrence},
\begin{align}
    (\rmZ^{(1)}(\alpha) - \Lambda_{\qr^{-1}}(\alpha)\rmI_{r + 1}) \rvu = \Big((\qr \indi{r \textnormal{ even}} + \qr^{-1} \indi{r \textnormal{ odd}}) \cdot \ervu_{r-1} - \Lambda_{\qr^{-1}}(\alpha)\ervu_r \Big) \rve_{r}.
\end{align}
Hence $\|(\rmZ^{(1)}(\alpha) - \Lambda_{\qr^{-1}}(\alpha) \rmI_{r + 1}) \rvu\| \to 0$ as $r \to \infty$ since $\alpha > \qr^{-2} + 1$ for $\rmZ^{(1)}(\alpha)$. Meanwhile, combining the conclusion in (1), as $r \to \infty$, $\Lambda_{\qr^{-1}}(\alpha)$ is the limit of $\lambda^{(1)}_{1}$ since $\Lambda_{\qr^{-1}}(\alpha) > \qr + \qr^{-1}$ when $\alpha > \qr^{-2} + 1$. This completes the proofs of (a) in (2) and (a) in (3).
\end{proof}

\subsection{Bounds on the adjacency matrices of biregular trees}
\begin{lemma}\label{lem:tree_norm_rough}
    Let $\mathbb{T}$ be a tree whose root has at most $q$ children, and all the other vertices in an even layer have at most
    $q$ children and all the other vertices in the odd layer have at most $p$ children. Then its adjacency matrix satisfies
    $$\|\rmA^{\mathbb{T}}\| \leq \sqrt{p}+\sqrt{q}.$$
\end{lemma} 
\begin{proof}
Denote the layers of $\mathbb{T}$ by $0,1,2,\dots$, with layer $0$ being the root of this tree.  
Since
\[
   \|\rmA^{\mathbb{T}}\|
   \;=\; 
   \sup_{\vw \neq 0}
   \frac{\bigl|\langle \vw,\,\rmA \vw\rangle\bigr|}{\|\vw\|_2^2},
\]
we focus on control the quadratic form $\langle \vw,\,\rmA \vw\rangle$, for any $\vw$.
Notice that
\[
  \langle \vw, \,\rmA \vw\rangle
  \;=\;
  \sum_{x,y} w_x \,\rmA_{xy}\, w_y
  \;=\;
  \sum_{\{x,y\} \in E(\mathbb{T})} 2\,w_x \,w_y.
\]
where each undirected edge $\{x,y\}$ contributes $2w_x\,w_y$.  
 
If $x$ is in an even layer and $y$ is in the odd layer directly below it, 
then $x$ has at most $q\ge 1$ children. For each such edge $(x,y)$, we use the elementary Young's inequality
\[
   2\,w_x\,w_y
   \;\le\;
   \frac{1}{\sqrt{q}}\,w_x^2 
   \;+\;
   \sqrt{q}\,w_y^2.
\]
Summing over all children $y$ of a fixed even-layer vertex $x$ gives
\[
   \sum_{y\in C_x} 2\,w_x\,w_y
   \;\le\;
   \sum_{y\in C_x}
      \Bigl(\tfrac{1}{\sqrt{q}}\,w_x^2 + \sqrt{q}\,w_y^2\Bigr)
   \;=\;
   \frac{\#C_x}{\sqrt{q}}\,w_x^2 + \sqrt{q}\,\sum_{y\in C_x} w_y^2.
\]
Since $x$ has at most $q$ children, $\#C_x\le q$, we have
\[
   \sum_{y\in C_x} 2\,w_x\,w_y
   \;\le\;
   \sqrt{q}\,w_x^2 \;+\;\sqrt{q}\,\sum_{y\in C_x} w_y^2.
\]
Summing this over all even-layer vertices $x$, we can get
\[
   \sum_{\substack{x\in\text{even}\\y\in C_x}}
   2\,w_x\,w_y
   \;\le\;
   \sqrt{q}\sum_{x\in\text{even}} w_x^2 
   \;+\;\sqrt{q}\sum_{y\in\text{odd}} w_y^2,
\]
because each odd-layer vertex appears \emph{exactly once} as a child of some even vertex 
in a tree.

A symmetric argument handles edges from an odd layer to the even layer below it:
if $x$ is odd and $y$ its child in the next even layer, then $x$ has at most $p\ge 1$ children, and
\[
2\,w_x\,w_y
\;\le\;
\frac{1}{\sqrt{p}}\,w_x^2 
\;+\;
\sqrt{p}\,w_y^2.
\]
Hence for each odd‐layer vertex $x$,
\[
\sum_{y\in C_x} 2\,w_x\,w_y
\;\le\;
\sqrt{p}\,w_x^2 
\;+\;\sqrt{p}\sum_{y\in C_x} w_y^2,
\]
and summing over odd‐layer vertices $x$ yields
\[
\sum_{\substack{x\in\text{odd}\\y\in C_x}}
2\,w_x\,w_y
\;\le\;
\sqrt{p}\sum_{x\in\text{odd}} w_x^2
\;+\;
\sqrt{p}\sum_{y\in\text{even}} w_y^2.
\]

All edges of $\mathbb{T}$ are accounted for in exactly one of the two sums above.  Therefore,
\begin{align}
\sum_{\{x,y\}\in E}
    2\,w_x\,w_y
\;\le\;&
\sqrt{q}\sum_{x\in\text{even}} w_x^2
    + \sqrt{q}\sum_{y\in\text{odd}} w_y^2
\;+\;
\sqrt{p}\sum_{x\in\text{odd}} w_x^2
    + \sqrt{p}\sum_{y\in\text{even}} w_y^2\\
=\;&
(\sqrt{p} + \sqrt{q})
\,\Bigl(\,\sum_{\text{even}}w_x^2 \;+\;\sum_{\text{odd}}w_x^2\Bigr)
\;=\;
(\sqrt{p}+\sqrt{q})\,
\sum_{x \in V(\mathbb{T})} w_x^2.
\end{align}
This completes the proof.
\end{proof}

\begin{lemma}\label{lem:tree_norm}
Let $s,p,q \in \mathbb{N}^*$. Let $\mathbb{T}$ be a tree whose root has $s$ children, all the other vertices in even layers have at most $q$ children, 
and all the other vertices in odd layers have at most $p$ children. Then the adjacency matrix $\rmA^{\mathbb{T}}$ of $\mathbb{T}$ satisfies
\[
\|\rmA^{\mathbb{T}}\| \le (qp)^{1/4}~\Lambda_{\big(\frac{p}{q}\big)^{1/4}}(s/\sqrt{pq}\vee (\sqrt{p/q}+1)),
\]where $\Lambda$ is the function defined by \eqref{eqn:Lambda_qr} and \eqref{eqn:Lambda_qr_inverse}.
\end{lemma}
\begin{proof} 
Let $r \in \mathbb{N}$ and denote by $\mathbb{T}_{s,p,q}(r)$ be the tree of depth $2r$ whose root $x$ has $s$ children, all vertices at distance $2i$ have at most $q$ children, 
all the other vertices at distance $2i+1$ have at most $p$ children, for $1\le i\le 2r$,
and all vertices at distance $2r+1$ from $x$ are leaves. For large enough $r$,
we can exhibit $\mathbb{T}$ as a subgraph of $\mathbb{T}_{s,p,q}(r)$. By the Perron–Frobenius theorem,
\begin{equation}\label{eqn:Perron–Frobenius}
\|\rmA^{\mathbb{T}}\| = \langle \rvw, \rmA^{\mathbb{T}} \rvw\rangle
\end{equation}
for some normalized eigenvector $\rvw$ whose entries are nonnegative. 
We can extend $\rvw$ to a vector indexed by the vertex set of $\mathbb{T}_{s,p,q}(r)$ by setting $\rvw_y = 0$ for $y$ not in the vertex set of $\mathbb{T}$. 
Then
\begin{equation}\label{eqn:quadratic_tree}
\langle \rvw, \rmA^{\mathbb{T}} \rvw\rangle \le \langle \rvw, \rmA^{\mathbb{T}_{s,p,q}(r)}\rvw\rangle.
\end{equation}
For simplicity, we denote by $\rmA \equiv \rmA^{\mathbb{T}_{s,p,q}(r)}$, it therefore remains to estimate the right-hand side of \eqref{eqn:quadratic_tree} for large enough $r$. 
By  Lemma~\ref{lem:tridiagrep} in Appendix~\ref{sec:tri_regular}, we can define $\rmZ$ as the tridiagonalization of $\rmA$ around the root up to radius $2r$. 
The associated orthonormal set $\rvs_0,\rvs_1,\dots,\rvs_r$ is given by $\rvs_i = \mathbf{1}_{S_i(x)}/\|\mathbf{1}_{S_i(x)}\|$, and 
$$\rmZ=(qp)^{1/4}~\rmZ_{\big(\frac{p}{q}\big)^{1/4}}(s/\sqrt{pq}),$$
where $\rmZ_{\big(\frac{p}{q}\big)^{1/4}}(\alpha)$ is the upper-left $(r+1)\times(r+1)$ block of \eqref{eqn:Z1Z2}. 
We introduce the orthogonal projections $P_0:=\rvs_0\rvs_0^*$ and $P:=\sum_{i=0}^{r}\rvs_i\rvs_i^*$. 
Notice that $P_0P=P_0$ and $(1-P)(1-P_0)=1-P$. 
For large enough $r$ the vectors $\rvs_r$ and $\rvw$ have disjoint support, and hence $(1-P)\rmA P \rvw=(1-P)A\sum_{i=0}^{r-1}\rvs_i\langle \rvs_i,\rvw\rangle=0$, 
since $\rmA\rvs_i \subset \text{Span}\{\rvs_{i-1},\rvs_{i+1}\}$ for $i<r$. Thus we have
\begin{align}
\langle \rvw, \rmA\rvw\rangle = ~&\langle \rvw, P\rmA P\rvw\rangle + \langle \rvw, (1-P)\rmA(1-P)\rvw\rangle\\
=~& \langle \rvw, P\rmA P\rvw\rangle + \langle \rvw, (1-P)(1-P_0)\rmA(1-P_0)(1-P)\rvw\rangle.
\end{align}
Then, Lemma~\ref{lem:eigenZ1Z2} shows that
\begin{equation}
\lim_{r \to \infty} \|P\rmA P\| = \lim_{r \to \infty} \|\rmZ\| = (qp)^{1/4}~\Lambda_{\big(\frac{p}{q}\big)^{1/4}}(s/\sqrt{pq}\vee (\sqrt{p/q}+1)).
\end{equation}
Next, notice that $(1-P_0)\rmA(1-P_0)$ is the adjacency matrix of a disjoint union of $s$ bipartite trees whose vertices in even layers have $p$ children and odd layers have $q$ children. Hence, by Lemma~\ref{lem:tree_norm_rough}, we therefore obtain $\|(1-P_0)\rmA(1-P_0)\| \le  \sqrt{q}+\sqrt{p}$.
Thus,
\begin{align}
\limsup_{r \to \infty} \langle \rvw, \rmA\rvw\rangle \le~& (qp)^{1/4}~\Lambda_{\big(\frac{p}{q}\big)^{1/4}}(s/\sqrt{pq}\vee (\sqrt{p/q}+1))\|P\rvw\|^2 + (\sqrt{p}+\sqrt{q})\|(1-P)\rvw\|^2\\
\le~& (qp)^{1/4}~\Lambda_{\big(\frac{p}{q}\big)^{1/4}}(s/\sqrt{pq}\vee (\sqrt{p/q}+1))\|\rvw\|^2.
\end{align}
We conclude this lemma by \eqref{eqn:Perron–Frobenius} and \eqref{eqn:quadratic_tree}.
\end{proof}

\section{Deferred proofs in \Cref{sec:approx_eigenpairs}}\label{sec:proofs_of_approximate_eigen_pairs}
\subsection{Proof of \Cref{lem:five_term_decomp}}\label{sec:five_term_decomp}
\begin{proof}[Proof of \Cref{lem:five_term_decomp}(i)] We focus on the proof of case (i). First, it follows directly that $\scA \rvv = \rvw_0 + d^{-1/2}\rmA \rvv$ since $\scA = (\rmA - \E \rmA)/\sqrt{d}$ in \eqref{eqn:scAdj}. Note that the following holds  
\begin{align}
    \rmA \ones_{\layer_{j}} = \indi{ j\geq 1 } \cdot \sum_{y\in \layer_{j-1}}|\NS{j}{y}| \ones_y + \sum_{y\in \layer_{j}}|\NS{j}{y}| \ones_y + \sum_{y\in \layer_{j+1}}|\NS{j}{y}| \ones_y,
\end{align}
Then by plugging $\rvv$ in \eqref{eqn:approx_eigenvector}, we have
\begin{align}
    \scA \rvv - \sum_{k=0}^{2}\rvw_k = \frac{1}{\sqrt{d}} \ervu_0 \ones_{\layer_1} + \frac{1}{\sqrt{d}} \sum_{j=1}^{r}\frac{\ervu_{j}}{\sqrt{|\layer_{j}|}}\Big[ \sum_{y\in \layer_{j-1}} \frac{|\layer_{j}|}{|\layer_{j-1}|} \ones_y + \ones_{\layer_{j+1}} \Big].
\end{align}
By rearranging the above terms with respect to $\rvs_j = |\layer_{j}|^{-1/2} \ones_{ \layer_{j} }$, it follows that
\begin{align}
    \scA \rvv - \sum_{k=0}^{2}\rvw_k =&\, \ervu_0 \sqrt{ \frac{|\layer_1|}{d}}\rvs_1 + \ervu_1 \sqrt{ \frac{|\layer_1|}{d}}\rvs_0 + \ervu_2 \frac{1}{\sqrt{d}} \cdot \frac{\sqrt{|\layer_2|}}{\sqrt{|\layer_1|}} \rvs_1 \\
    &\, + \frac{1}{\sqrt{d}} \sum_{j=2}^{r-1}\Big( \ervu_{j-1}\frac{\sqrt{|\layer_{j}|}}{\sqrt{|\layer_{j-1}|}} + \ervu_{j+1} \frac{\sqrt{|\layer_{j+1}|}}{\sqrt{|\layer_{j}|}} \Big)\rvs_j \notag\\
    &\, + \ervu_{r-1}\frac{1}{\sqrt{d}} \frac{\sqrt{|\layer_{r}|}}{\sqrt{|\layer_{r-1}|}} \rvs_{r}  + \ervu_{r} \frac{1}{\sqrt{d}} \frac{\sqrt{|\layer_{r+1}|}}{\sqrt{|\layer_r|}} \rvs_{r+1}.
\end{align}
We consider the case that $r$ is even. Recall $\Lambda_{\qr}(t)$ in \eqref{eqn:Lambda_qr} and the recurrence in \eqref{eqn:u012V2}, \eqref{eqn:ujsV2}, then
\begin{subequations}
\begin{align}
    \Lambda_{\qr}(\alpha_x) \ervu_0 = \sqrt{\alpha_x} \ervu_1, \quad &\, \Lambda_{\qr}(\alpha_x) \ervu_1 = \sqrt{\alpha_x}\ervu_0 + \qr^{-1} \ervu_2,\\
    \Lambda_{\qr}(\alpha_x) \ervu_{2j} = \qr^{-1} \ervu_{2j-1} + \qr \ervu_{2j+1}, \quad &\, \Lambda_{\qr}(\alpha_x) \ervu_{2j+1} = \qr \ervu_{2j} + \qr^{-1} \ervu_{2j+2},
\end{align}
\end{subequations}
for all $j = 2, \ldots, \lfloor r/2\rfloor$. Note that $\rD_x = |\layer_1(x)|$; it then implies
\begin{subequations}
\begin{align}
    \Lambda_{\qr}(\alpha_x) \rvv = &\, \ervu_1 \sqrt{|\layer_1|/d} \,\,\rvs_0 + (\ervu_0 \sqrt{|\layer_1|/d} + \ervu_2 \qr^{-1})\rvs_1\\
    &\, + \sum_{j=1}^{\lfloor r/2\rfloor} (\qr^{-1} \ervu_{2j-1} + \qr \ervu_{2j+1}) \rvs_{2j} + \sum_{j=1}^{\lfloor r/2\rfloor-1} (\qr \ervu_{2j} + \qr^{-1} \ervu_{2j+2}) \rvs_{2j+1}.
\end{align}
\end{subequations}
Together with previous results, we have
\begin{subequations}
\begin{align}
    &\,\big(\scA - \Lambda(\alpha_{x}) \big) \rvv - \sum_{k=0}^{2}\rvw_k\notag\\
    = &\, \ervu_2 \Big( \frac{\sqrt{|\layer_2|}}{\sqrt{d}\cdot \sqrt{|\layer_1|}} - \qr^{-1} \Big)\rvs_1 + \sum_{j=1}^{\lfloor r/2\rfloor - 1 } \Bigg( \ervu_{2j-1}\Big(\frac{\sqrt{|\layer_{2j}|}}{\sqrt{d} \cdot \sqrt{|\layer_{2j-1}|}} - \qr^{-1} \Big) + \ervu_{2j+1}\bigg(\frac{\sqrt{|\layer_{2j+1}|}}{\sqrt{d}\cdot \sqrt{|\layer_{2j}|}} - \qr \bigg) \Bigg) \rvs_{2j}\\
    &\, + \sum_{j=1}^{\lfloor r/2\rfloor - 1} \Bigg( \ervu_{2j}\Big(\frac{\sqrt{|\layer_{2j+1}|}}{\sqrt{d} \cdot \sqrt{|\layer_{2j}|}} - \qr \Big) + \ervu_{2j+2}\bigg(\frac{\sqrt{|\layer_{2j+2}|}}{\sqrt{d} \cdot \sqrt{|\layer_{2j+1}|}} - \qr^{-1} \bigg) \Bigg) \rvs_{2j+1}\notag\\
    &\, + \Bigg( \frac{\ervu_{r-1}}{\sqrt{d}} \frac{\sqrt{|\layer_{r}|}}{\sqrt{|\layer_{r-1}|}} - \ervu_{r-1}\qr^{-1} - u_{r+1}\qr \Bigg) \rvs_{r} + \frac{\ervu_{r}}{\sqrt{d}} \frac{\sqrt{|\layer_{r+1}|}}{ \sqrt{|\layer_r|}} \,\,\rvs_{r+1}\\
    = &\, \rvw_3 + \rvw_4,
\end{align}
\end{subequations}
which completes the proof. The case of odd $r$ can be verified similarly by checking the boundary conditions.
\end{proof}

\begin{proof}[Proof of \Cref{lem:five_term_decomp}(ii)] The proof follows similarly by using the recurrence relationship \eqref{eqn:u012V1}, \eqref{eqn:ujsV1}.
\end{proof}

\subsection{Proof of \Cref{lem:w0tow5}}\label{sec:w0tow5}
Analogous to {\cite[Lemma 5.4]{alt2021extremal}}, we establish the following concentration results on the neighborhood growth rates in $\gG$.
\begin{lemma}[Concentration of neighborhood growth rate]\label{lem:concentrationSi}
 Suppose $1\leq d\leq N^{1/4}$. Conditioned on the event \eqref{eqn:Dx_Event}, the following holds.
 \begin{enumerate}
     \item When $x\in \gV_2$, for $j \in \{ 1, \ldots, \lceil r/2 \rceil\}$, with very high probability, we have
     \begin{subequations}
        \begin{align}
       & \bigg|\frac{1}{d_1}\cdot\frac{|\layer_{2j}|}{|\layer_{2j - 1}|} -1 \bigg| \lesssim \bigg( \frac{\log(N)}{d_{1}|\layer_{2j-1}|}\bigg)^{1/2},    \label{eqn:S2j_ratio_bound}\\
           &\bigg| \frac{1}{d_2}\cdot \frac{|\layer_{2j+1}|}{|\layer_{2j}|} - 1 \bigg|  \lesssim \bigg( \frac{\log(N)}{d_{2}|\layer_{2j}|}\bigg)^{1/2}, \label{eqn:S2j+1_ratio_bound}
        \end{align}
    \end{subequations}
    where $d_{1} = \qr^{-2}d$, $d_{2} = \qr^{2}d$, and
    \begin{subequations}
    \begin{align}
        |\layer_{2j}|&=\rD_x d_1 d^{2j-2}\Big(1+ O\Big( \Big(\frac{\log(N)}{d \rD_{x}}\Big)^{1/2} \Big) \Big), \label{eqn:S2j}\\
        |\layer_{2j+1}|&=\rD_xd^{2j}\Big(1+ O\Big( \Big(\frac{\log(N)}{d \rD_{x}}\Big)^{1/2} \Big) \Big). \label{eqn:S2j+1}
    \end{align}
    \end{subequations}
    \item When $x\in \gV_1,$ the same statement holds by swapping $d_2$ with $d_1$ in (1).
\end{enumerate}
\end{lemma}

\begin{lemma}[Few cycles in small balls]\label{lem:few_even_cycle_prob}
    For $k,r\in \N$, we have
    \begin{align}
        \P \big( |\gE(\gG|_{\ball_r})|-|\ball_r|+1\geq k \mid \layer_1 \big) \lesssim \frac{d^{k}}{(2\sqrt{nm})^k} (d+\rD_{x})^{2kr}(2kr)^{2k},\label{eqn:few_even_cycle_prob}
    \end{align}
    where $\lesssim$ only hides some constants irrelevant to $m$ and $n$.

    Consequently, for $r\leq r_{x}$, conditioned on $\layer_{1}$, the following holds with very high probability
    \begin{align}
        \# \textnormal{ of cycles in } \gG|_{\ball_r} \leq C,\label{eqn:few_even_cycle}
    \end{align}
    where $\const >0$ is some absolute constant. 
\end{lemma}

\begin{corollary}\label{cor:few_cycle_counts}
For $x\in \gV$, conditioned on $\layer_{1}$, with very high probability, for all $1\leq j\leq r_x + 1$ 
\begin{subequations}
    \begin{align}
        &\, |\layer_{j}|=\sum_{y\in \layer_{j}} |\NS{j-1}{y}|+ \const_{1} = \sum_{y\in\layer_{j-1}} |\NS{j}{y}|+ \const_{1}\,, \label{eqn:few_cycle_different_layer}\\
        &\, \sum_{y\in \layer_{j}} |\NS{j}{y}| =\const_{2}\,,\label{eqn:few_cycle_same_layer}
    \end{align}
\end{subequations}
 where $\const_{1}, \const_{2}>0$ are some absolute constants.
\end{corollary}

\begin{proof}[Proof of \eqref{eqn:norm_w0}]
Let $\ones_{\gV_1}$ (resp. $\ones_{\gV_2}$) denote the vector where the entry is $1$ if $x \in \gV_1$ (resp. $x \in \gV_2$), otherwise $0$, then $\E \rmA = d/\sqrt{mn} \cdot ( \ones_{\gV_1}\ones_{\gV_2}^{\sT} + \ones_{\gV_2}\ones_{\gV_1}^{\sT})$. By plugging in $\rvv$ \eqref{eqn:approx_eigenvector}
    \begin{align}
        \|\rvw_0\|_2^2 =&\, \frac{1}{d} \rvv^{\sT}\E\rmA^{\sT} \E\rmA \rvv = \frac{d}{m} \rvv^{\sT} \ones_{\gV_2} \ones_{\gV_2}^{\sT} \rvv + \frac{d}{n}\rvv^{\sT}\ones_{\gV_1} \ones_{\gV_1}^{\sT}\rvv\\
        = &\, \frac{d}{m} \sum_{y\in \ball_{r} \cap \gV_2} |\<\ones_{y}, \rvv\>|^2 + \frac{d}{n} \sum_{y\in \ball_{r} \cap \gV_1} |\<\ones_{y}, \rvv\>|^2\\
        \leq &\, d\cdot (|\ball_{r} \cap \gV_2|/m + |\ball_{r} \cap \gV_1|/n) \leq d\cdot |\ball_{r}|/m
    \end{align}
 where the last line holds due to Cauchy–Schwarz inequality, and the facts $\|\rvv\|_2 = 1$ and $n \geq m$. 
 Consequently, \eqref{eqn:norm_w0} follows by \eqref{eqn:S2j} \eqref{eqn:S2j+1} and the event \eqref{eqn:Dx_Event}, since
 \begin{align}
     |\ball_r| \lesssim \sum_{j=0}^{r-1} \rD_{x} d^{j}\lesssim \rD_{x} d^{r-1} \leq \sqrt{N}.
 \end{align}
\end{proof}

\begin{proof}[Proof of \eqref{eqn:norm_w1}]
Since $\layer_0=\{x\}$, we have $|\NS{0}{y}| =1$ for all $y\in \layer_1$ and $|\NS{0}{x}| = 0$. According to the triangle inequality, \eqref{eqn:few_cycle_different_layer}, and \eqref{eqn:few_cycle_same_layer}, the following holds with very high probability
\begin{align}
    \|\rvw_1\|\lesssim  \sum_{j=1}^r \frac{1}{\sqrt{d|\layer_{j}|}}\lesssim (d\rD_{x})^{-1/2} \sum_{j=1}^{r} d^{-j+1} \lesssim (d\rD_{x})^{-1/2},
\end{align}
 where we applied \eqref{eqn:S2j} and \eqref{eqn:S2j+1} in the second inequality.
\end{proof}

\begin{proof}[Proof of \eqref{eqn:norm_w3}]
According to \eqref{eqn:S2j_ratio_bound}, \eqref{eqn:S2j+1_ratio_bound} and the facts $d_{1} = \qr^{-2}d$, $d_{2} = \qr^{2}d$, $\|\rvs_{j}\|_{2}^{2} = 1$, then by Cauchy–Schwarz inequality, the following holds with very high probability
 \begin{align}
     \|\rvw_3\|^2\lesssim \frac{\log N}{d\rD_x} \sum_{j=0}^{r} \ervu_{j}^2 \leq  \frac{\log N}{d\rD_x}.
 \end{align}
\end{proof}

\begin{proof}[Proof of \eqref{eqn:norm_w4}]
According to triangle inequality and the fact $\|\rvs_{j}\| = 1$, as well as \eqref{eqn:S2j_ratio_bound} \eqref{eqn:S2j+1_ratio_bound}, the following holds with very high probability when $r$ is odd
\begin{align}\label{eqn:trianglew4}
    \|\rvw_4\| \leq &\, |\ervu_{r-1}| \bigg(\sqrt{\frac{|\layer_r|}{d|\layer_{r-1}|}}-\qr \bigg)+\qr^{-1}|\ervu_{r+1}|+ \frac{|\ervu_r|}{\sqrt d}\sqrt{\frac{|\layer_{r+1}|}{|\layer_r|}}\\
     \lesssim &\, |\ervu_{r-1}| \bigg( \frac{\log(N)}{d_{2}|\layer_{r-1}|}\bigg)^{1/4} + \qr^{-1} |\ervu_{r+1}| + \qr^{-1}|\ervu_r| \bigg(1  + \bigg( \frac{\log(N)}{d_{1}|\layer_{r}|}\bigg)^{1/4} \bigg)\\
    \lesssim &\, |\ervu_{r-1}| \bigg( \frac{\log(N)}{\rD_{x} d^{(r-1)/2}}\bigg)^{1/4} + \qr^{-1} (|\ervu_{r+1}| + |\ervu_r|)\quad \textnormal{(first term is relatively negligible)}\\
    \leq &\, |\alpha_{x} - \qr^{2}|^{-\lceil(r-1)/2\rceil} \alpha_{x}^{-1/2} \qr^{-1} \bigg( \qr^{-1}\alpha_{x}|\alpha_{x} - \qr^{2}|^{-1} + \Lambda_{\qr}(\alpha_{x}) \bigg) |\ervu_{0}|\\
    \lesssim &\, |\alpha_{x} - \qr^{2}|^{-\lceil(r-1)/2\rceil},
\end{align}
where in the last two lines, we apply the definitions of $\ervu_{r-1}, \ervu_r$ in \eqref{eqn:ujsV2} and the fact $|\ervu_0| \leq 1$. The case of $r$ being even can be proved by a similar argument.
\end{proof}

 \begin{proof}[Proof of \eqref{eqn:norm_w2}]
We present the proof for the case $x\in \gV_{2}$. The proof for $x\in \gV_{1}$ follows in a similar way. Note that $\layer_{0} =\{x\}$ and $|\NS{1}{x}| = |\layer_{1}| =  |\layer_{1}| / |\layer_{0}|$ and $\sum_{j=2}^{r} \ervu_{j}^{2} \leq 1$, then
\begin{align}
    \|\rvw_{2}\|_{2}^{2} = &\, \frac{1}{d}\sum_{j=2}^{r} \frac{|\ervu_{j}|^{2}}{|\layer_{j}|} \sum_{y\in \layer_{j-1}}\bigg(|\NS{j}{y}|-\frac{|\layer_{j}|}{|\layer_{j-1}|} \bigg)^{2}\\
    \leq &\, \frac{3}{d}\sum_{j=2}^{r} \frac{|\ervu_{j}|^{2}}{|\layer_{j}|} \sum_{y\in \layer_{j-1}}\bigg[ \Big(|\NS{j}{y}|- \E\big[|\NS{j}{y}| \big| \ball_{j-1}\big] \Big)^{2} \\
    &\, \quad \quad \quad \quad \quad \quad + \Big(\E\big[|\NS{j}{y}| \big| \ball_{j-1}\big] - d_{1}\indi{j\textnormal{ even}} - d_{2}\indi{j\textnormal{ odd}}\Big)^2 \\
    &\, \quad \quad \quad \quad \quad \quad + \Big( d_{1}\indi{j\textnormal{ even}} + d_{2}\indi{j\textnormal{ odd}} - \frac{|\layer_{j}|}{|\layer_{j-1}|} \Big)^{2}\bigg]\\
    \leq &\, \frac{4}{d^2} \max_{2\leq j\leq r} \big(\qr^{2} \cdot \indi{j \textnormal{ odd}} + \qr^{-2} \cdot \indi{j \textnormal{ even}}\big) \bigg[ \rY_{j} + \const_{1}\log(N) + \const_{2} \frac{d\log(N)}{\rD_{x}} \bigg],
\end{align}
where in the second to last inequality, the last two terms are bounded with very high probability for some constants $\const_{1}, \const_{2} >0$ according to \Cref{lem:conditional_layer_concentration} and \Cref{lem:concentrationSi} respectively. Then \eqref{eqn:norm_w2} follows by the following claim
\begin{align}
    \rY_j = &\, \frac{1}{|\layer_{j-1}|}\sum_{y\in \layer_{j-1}} \Big( |\NS{j}{y}| - \E\big[|\NS{j}{y}| \mid \ball_{j-1}\big] \Big)^2 \leq \const d (1 + \log(N)/\rD_{x})\log(d), \quad 2\leq j \leq r.
\end{align}

We now present the proof of the claim above. For convenience, we abbreviate
\begin{align}
    \rE_{j}(y) \coloneqq |\NS{j}{y}| - \E\big[|\NS{j}{y}| \mid \ball_{j-1}\big],\quad y\in \layer_{j-1}.
\end{align}
 Equivalently, we obtain
\begin{align}
    \rY_{j} \cdot |\layer_{j-1}| = &\, \sum_{y \in \layer_{j-1}} (\indi{|\rE_{j}(y)|^{2} \leq d} + \indi{|\rE_{j}(y)|^{2} > d} )\cdot |\rE_{j}(y)|^{2} \notag \\
    \leq &\, d|\layer_{j-1}| + \sum_{y \in \layer_{j-1}} \indi{|\rE_{j}(y)|^{2} > d} |\rE_{j}(y)|^{2}. \label{eqn:Ey_decomposition}
\end{align}
The second term of \eqref{eqn:Ey_decomposition} adapts the following dyadic decomposition with very high probability
\begin{align}
    \sum_{y \in \layer_{j}} \indi{|\rE_{j}(y)|^{2} > d} |\rE_{j}(y)|^{2} \leq \sum_{k = k_{\min}}^{k_{\max}} \sum_{y\in \sN_{j}^{k}} |\rE_{j}(y)|^2 \leq d^{2} \sum_{k = k_{\min}}^{k_{\max}} e^{k+1} |\sN_{j}^{k}|, \label{eqn:Ey_dyadic_decomposition}
\end{align}
where we introduced
\begin{align}
    \sN_{j}^{k} \coloneqq &\, \{y\in \layer_{j-1}: d^{2} e^{k} \leq  |\rE_{j}(y)|^2 \leq d^{2} e^{k+1} \},\\
    k_{\min} \coloneqq &\, - \lceil \log(d)\rceil, \quad k_{\max} \coloneqq 0.
\end{align}
Together with \eqref{eqn:Ey_decomposition} and \eqref{eqn:Ey_dyadic_decomposition}, we conclude that
\begin{align}
     \rY_{j} \leq &\, d + \frac{d^{2}}{|\layer_{j-1}|}\sum_{k = k_{\min}}^{k_{\max}} e^{k+1} |\sN_{j}^{k}|\\
     \leq &\, d + \frac{d^{2}}{|\layer_{j-1}|} \cdot \frac{\const}{4d} (|\layer_{j-1}| +\log(N)) \cdot \Big(\sum_{k = k_{\min}}^{0} e^{k+1} e^{-k} \Big)\\
     \leq &\, d + \frac{\const d}{4} \Big(1 + \frac{\log(N)}{|\layer_{j-1}|} \Big) \cdot e|k_{\min}|\\
     \leq &\, \rY_{2} \leq \const d (1 + \log (N)/\rD_{x}) \log(d).
\end{align}

{\noindent \textit{Proof of \eqref{eqn:Ey_dyadic_decomposition}.}} We will show that $|\rE_{j}(y)|^{2} \leq d^{2}e^{k_{\max}+1}$ for all $y\in \layer_{j}$ with very high probability. To that end, we first introduce the following level set
\begin{align}
    \sL_{j}^{s} \coloneqq \{y\in \layer_{j}: |\rE_{j}(y)|^{2} \geq s^{2} d^{2}\},
\end{align}
for any $s>0$. The problem is then transferred to obtain a probabilistic tail bound on the size of $\sL_{j}^{s}$, since
\begin{align}
    \P(\exists y\in \layer_{j}, |\rE_{j}(y)|^{2} > s^{2}_{\max}d^{2} | \ball_{j-1}) = \P(|\sL_{j}^{s_{\max}}| \geq 1 | \ball_{j-1}), \quad s_{\max} = e^{(k_{\max}+1)/2}. 
\end{align}
For each $l \leq |\layer_{j}|$, we use a union bound to sum over all subsets $\setT \subset \layer_{j}$ with $|\setT| = l$, hence
\begin{align}
    &\, \P(|\sL_{j}^{s_{\max}}| \geq 1 | \ball_{j-1}) \label{eqn:level_set_size}\\
     \leq &\, \sum_{l = 1}^{|\layer_{j}|} \sum_{\setT \subset \layer_{j}, |\setT | = l} \P(|\rE_{j}(y)|^{2} > s_{\max}^{2}d^{2} \textnormal{ for all } y \in \setT |\ball_{j-1})\\
    \leq &\,\sum_{l = 1}^{|\layer_{j}|} \binom{|\layer_{j}|}{l} \max_{\setT \subset \layer_{j}, |\setT | = l} \P(|\rE_{j}(y)|^{2} > s_{\max}^{2}d^{2} \textnormal{ for all } y \in \setT |\ball_{j-1})\\
    \leq &\,2\sum_{l = 1}^{|\layer_{j}|} \binom{|\layer_{j}|}{l} \Big[ \P(\rE_{j}(y) > s_{\max} d )\Big]^{l} \quad (\textnormal{Bennett \Cref{lem:Bennett}})\\
    \leq &\,2 \sum_{l = 1}^{|\layer_{j}|} \exp\Big( l\log(|\layer_{j}|) -ld \cdot (\indi{j \textnormal{ even}}\qr^{2}\benrate(\qr^{-2}s_{\max}) + \indi{j \textnormal{ odd}}\qr^{-2}\benrate(\qr^{2}s_{\max}) \Big),\\
    \leq &\, 2 |\layer_{j}| \exp\Big( \log(|\layer_{j}|) -d \cdot (\qr^{2}\benrate(\qr^{-2} s_{\max}) \wedge \qr^{-2}\benrate(\qr^{2}s_{\max}) \Big) \leq N^{-\nu},
\end{align}
where the last inequality holds some constant $\nu >0$ since $|\layer_{j}| \leq N$ and $d = b\log(N)$, while $s_{\max}$ is bounded below by some sufficiently large constant.

We now estimate the upper bound of $|\sN_{j}^{k}|$. When conditioned on $\Xi_{x}(\epsilon)$ in \eqref{eqn:Dx_Event}, we have
\begin{align}
    \P\Big(|\sN_{j}^{k}| \geq \ell_{j}^{k} | \ball_{j-1} \Big) \leq \P\Big(|\sL_{j}^{\exp(k/2)}| \geq \ell_{j}^{k} | \ball_{j-1} \Big) \leq N^{-\nu},
\end{align}
where the last inequality follows by an argument similar to \eqref{eqn:level_set_size}, with 
\begin{align}
    \ell_{j}^{k} \coloneqq \frac{\const}{4d} (|\layer_{j}| +\log(N)) \cdot \Big(e^{-k/2}\indi{k > 0} + e^{-k}\indi{k \leq 0} \Big),
\end{align}
for some constant $\const >0$. Therefore, $|\sN_{j}^{k}| \leq \ell_{j}^{k}$ with very high probability.
\end{proof}

\subsection{Deferred proofs of Lemmas}
\begin{proof}[Proof of Lemma~\ref{lem:concentrationSi}]
We present the proof for the case $x\in \gV_2$. The proof for $x\in \gV_1$ follows similarly. 

First, the proofs of \eqref{eqn:S2j_ratio_bound} and \cref{eqn:S2j+1_ratio_bound} are established by inducting \Cref{lem:conditional_layer_concentration}. We conclude \cref{eqn:S2j_ratio_bound}, \cref{eqn:S2j+1_ratio_bound} and
\begin{subequations}
    \begin{align}
        \rD_x \left(d/2\right)^{2j-2}\leq |\layer_{2j-1}|\leq \rD_x (2d)^{2j-2}\label{eqn:ratio2j-1}\\
        \rD_x d_1\left(d/2\right)^{2j-2}\leq |\layer_{2j}|\leq \rD_x d_1 (2d)^{2j-2} \label{eqn:ratio2j}
    \end{align}
\end{subequations}
for $1\leq j\leq \lceil r/2 \rceil$ simultaneously.

(i) Base case $j=1$. First, $\eqref{eqn:ratio2j-1}_{j=1}$ holds trivially. We choose $\epsilon^2 = \gK \log(N)(9d\rD_x)^{-1}$. Conditioned on the event \eqref{eqn:Dx_Event}, $N^{-1/4} \lesssim \epsilon^{2} \leq 1/9$ and $\rE_{1} \asymp d|\rD_{x}|/N\lesssim N^{-1/4}$ in \eqref{eqn:error_layer_l}, then $\rE_{1} \leq \epsilon^{2} \leq \epsilon \leq 1/3$, thus $\eqref{eqn:S2j_ratio_bound}_{j=1}$ holds with very high probability by applying $\eqref{eqn:l_odd_and_x_V2}_{j=1}$, where $\epsilon$ dominates the error. Then $\eqref{eqn:ratio2j}_{j=1}$ is proved by $\eqref{eqn:S2j_ratio_bound}_{j=1}$, which finally leads to $\eqref{eqn:S2j+1_ratio_bound}_{j=1}$ by applying $\eqref{eqn:l_even_and_x_V2}_{j=1}$.

(ii) Induction step. Suppose that $\eqref{eqn:S2j_ratio_bound}_{j} \eqref{eqn:S2j+1_ratio_bound}_{j}$ and $\eqref{eqn:ratio2j-1}_{j} \eqref{eqn:ratio2j}_{j}$ hold up to $l \geq 2$ with very high probability. First, $\eqref{eqn:ratio2j-1}_{j=l+1}$ can be concluded from $\eqref{eqn:S2j+1_ratio_bound}_{j=l}$ and $\eqref{eqn:ratio2j}_{j=l}$. Then from $\eqref{eqn:ratio2j-1}_{j=l+1}$, we conclude that $|\ball_{2l+1}(x)| \leq \rD_{x}(2d)^{r}\leq N^{1/2}$ due to event \eqref{eqn:Dx_Event}. We choose \[\epsilon^2 = \gK \log(N)(9d|\layer_{2l+1}|)^{-1}.\] Similarly, $N^{-1/4} \lesssim \epsilon^{2} \leq 1/9$ and $\rE_{2l + 1} \lesssim N^{-1/4}$ in \eqref{eqn:error_layer_l}, then $\rE_{2l + 1} \leq \epsilon^{2} \leq \epsilon \leq 1/3$, thus $\eqref{eqn:S2j_ratio_bound}_{j=l+1}$ holds with very high probability by applying $\eqref{eqn:l_odd_and_x_V2}_{j=l+1}$, where $\epsilon$ dominates the error. 
Similarly, inequality $\eqref{eqn:ratio2j}_{j=l+1}$ follows from $\eqref{eqn:S2j_ratio_bound}_{j=l+1}$. Then $\eqref{eqn:S2j+1_ratio_bound}_{j=l+1}$ is concluded by applying $\eqref{eqn:l_even_and_x_V2}_{j=l+1}$. Note that the necessary union bounds are affordable for the induction argument to reach $j = \lceil r/2 \rceil$, since the right-hand sides of \eqref{eqn:S2j_ratio_bound} and \eqref{eqn:S2j+1_ratio_bound} are always at least $1 - 2N^{-\const}$ for some constant $\const >0$.

We now turn to \cref{eqn:S2j} and \cref{eqn:S2j+1}, which can be established by
    \begin{align}
        \rD_x d^{2j-2} (1 - \epsilon_{2j-1})\leq &\, |\layer_{2j-1}|\leq \rD_x d^{2j-2} (1 + \epsilon_{2j-1}), \\ 
         \epsilon_{2j-1} = &\, 2\const \Big(\frac{\log(N)}{d\rD_{x}}\Big)^{1/2}\sum_{l=0}^{j-1} d^{-(2l-1)/2},\label{eqn:S2j-1_error}
    \end{align}
    as well as
    \begin{align}
        \rD_x d_1 d^{2j-2} (1 - \epsilon_{2j}) \leq &\, |\layer_{2j}| \leq \rD_x d_1 d^{2j-2} (1 + \epsilon_{2j}), \\ 
        \epsilon_{2j} = &\, 2\const \Big(\frac{\log(N)}{d\rD_{x}}\Big)^{1/2}\sum_{l=0}^{j-1} d^{-l}, \label{eqn:S2j_error}
    \end{align}
for $j\in \{1,\ldots, \lceil r/2 \rceil\}$ and some constant $\const >0$. Here, \eqref{eqn:S2j_error} \eqref{eqn:S2j-1_error} are proved by inducting \eqref{eqn:S2j_ratio_bound} \eqref{eqn:S2j+1_ratio_bound} and \eqref{eqn:ratio2j-1} \eqref{eqn:ratio2j}.

(i) Base case $j=1$. First, $\eqref{eqn:S2j-1_error}_{j=1}$ has trivially the same validity as $\epsilon_{1} = 0$. $\eqref{eqn:S2j_error}_{j=1}$ is a direct consequence of $\eqref{eqn:S2j_ratio_bound}_{j=1}$ with $\epsilon_{2} = (\gK \log(N)/(9\qr^{2}d_{1}\rD_x) \big)^{-1/2}$.

(ii) Induction step. Suppose that $\eqref{eqn:S2j-1_error}_{j}$ and $\eqref{eqn:S2j_error}_{j}$ hold up to $l \geq 2$. We can conclude from $\eqref{eqn:S2j+1_ratio_bound}_{j=l}$ $\eqref{eqn:S2j_error}_{j=l}$ that for some $\const >0$, with very high probability
\begin{align}
    |\layer_{2l+1}| \geq d_2|\layer_{2l}| \bigg( 1 - \const \Big( \frac{\log(N)}{d_{2}|\layer_{2l}|}\Big)^{1/2} \bigg) \geq \rD_{x} d^{2l} \bigg( 1 - \epsilon_{2l} - \const \Big( \frac{\log(N)}{\rD_{x} d^{2l}(1 - \epsilon_{2l}) }\Big)^{1/2} \bigg),
\end{align}
where we use the facts $d_{2} = \qr^2 d$, $d_{1} = \qr^{-2}d$ and $\epsilon_{2l} \leq 3/4$ for sufficiently large $\gK$ in event \eqref{eqn:Dx_Event}, finishing the proof of lower bound. The proof for the upper bound can be finished analogously, thus ending the induction of $\eqref{eqn:S2j-1_error}_{j=l+1}$. The induction for $\eqref{eqn:S2j_error}_{j=l+1}$ can be finished through a similar argument. Note that the necessary union bounds are affordable for the induction argument to reach $j = \lceil r/2 \rceil$.
\end{proof}

\begin{lemma}\label{lem:conditional_layer_concentration}
Assume $1 \leq d\leq N^{1/4}$. For $x\in \gV$, assume $|\ball_l(x)| \leq N^{1/2}$ for some $1\leq l \leq r$. Define 
\begin{align}
    \rE_l \coloneqq  \frac{d|\layer_{l}|}{2\sqrt{mn}} + (mn)^{-1/4}. \label{eqn:error_layer_l}
\end{align}
Then for $x\in \gV_2$, there exists some constant $\const > 0$ such that the following holds for all $\epsilon \in [0, 1]$
\begin{subequations}
    \begin{align}
        &\, \P\big( \big| |\layer_{2j}|- d_{1}|\layer_{2j-1}| \big| \leq (\epsilon + \const \rE_{2j-1})d_{1}|\layer_{2j-1}| \big| \ball_{2j-1} \big) \geq 1 - 2 e^{-\frac{3}{8}\epsilon^2 d_1 |\layer_{2j-1}|},\label{eqn:l_odd_and_x_V2}\\
        &\, \P\big( \big| |\layer_{2j+1}|- d_{2}|\layer_{2j}| \big| \leq (\epsilon + \const \rE_{2j})d_{2}|\layer_{2j}| \big| \ball_{2j} \big) \geq 1 - 2 e^{-\frac{3}{8}\epsilon^2 d_2 |\layer_{2j}|}.\label{eqn:l_even_and_x_V2}
    \end{align}
\end{subequations}
On the other hand, for $x \in \gV_1$, \eqref{eqn:l_odd_and_x_V2} and \eqref{eqn:l_even_and_x_V2} hold by swapping $d_2$ and $d_1$.
\end{lemma}
\begin{proof}[Proof of \Cref{lem:conditional_layer_concentration}]
Consider the case $x\in \gV_2$ first. Recall $\layer_{0} = \{x\}$. Due to the constraint of the bipartite structure, we have $\layer_{2j} \subset \gV_2$ and $\layer_{2j + 1} \subset \gV_1$ for all $j\geq 0$. 
We now explore the law of $\layer_{l+1}$ conditioning on $\ball_{l}$, where the randomness of $\layer_{l+1}$ comes from the edges with at least one node in $\ball_l^{\complement} \coloneqq \gV \setminus \ball_l$. For each vertex $y \in \ball_l^{\complement}$, define the random variable $\rY(y) \coloneqq \indi{ \layer_1(y) \cap \layer_l \neq \emptyset }$, that is, $\rY(y)  = 1$ if $y$ is connected to $\layer_l$ on graph $\gG$, where $\P( \rY(y) = 1) = 1 - \P(\rY(y) = 0) = 1 - (1 - p)^{|\layer_{l}|}$. Conditioned on $\ball_{l}$, $\{\rY(y)\}_{y\in \ball_l^{\complement}}$ are i.i.d. random variables with its summation satisfying the law
    \begin{align}
        |\layer_{l+1}| \Big|_{\ball_{l}} =&\, \sum_{y\in \ball_l^{\complement}} \rY(y)  \sim \mathrm{Binom}(1 - (1 - p)^{|\layer_l|}, \,\, \rN_{l+1})
    \end{align}
where the capacity of such random variable is defined through
\begin{align}
    \rN_{l+1} \coloneqq \indi{ l \textnormal{ even} } \cdot \Big( n - \sum_{j = 0}^{\lfloor l/2 \rfloor - 1} |\layer_{2j + 1}|\Big) + \indi{ l \textnormal{ odd} } \cdot \Big( m - \sum_{j = 0}^{ \lfloor l/2 \rfloor } |\layer_{2j}| \Big).
\end{align}
Here, $|\layer_l|$ and $\rN_{l+1}$ are deterministic when conditioned on $\ball_{l}$. Then the conditional expectation of $|\layer_{l+1}|$ can be evaluated as
    \begin{align}
        \overline{\rS}_{l + 1} \coloneqq &\,\E[|\layer_{l+1}| \mid \ball_{l}] = \big( 1 - (1 - p)^{|\layer_l|} \big) \cdot \rN_{l+1} = (1 + o(1))\cdot\big( p |\layer_l| + p^2 |\layer_l|^2 /2 \big)\cdot\rN_{l+1}\\
        =&\, d|\layer_{l}| \cdot ( 1 + \rE_l )\cdot ( \indi{l \textnormal{ even} } \qr^{2} +  \indi{l \textnormal{ odd} }  \qr^{-2}) \label{eqn:ESl+1}
    \end{align}
 where in the first line we applied $1 - p = e^{-p + O(p^2)}$ and $1 - e^{-p} = p + p^2/2 + O(p^3)$ when $p = o(1)$, and in the second line we leverage the fact $p = d/\sqrt{mn}$ and the assumptions that $d \leq N^{1/4}$ and $|\ball_l | = |\cup_{j=0}^{l}\layer_{j}| \leq N^{1/2}$. Note that $\rE_l \lesssim N^{-1/4}$ in \eqref{eqn:error_layer_l}, then by Bernstein \Cref{lem:Bernstein}, for some $\epsilon \in [0, 1]$, the following holds
 \begin{align}
     \P( \big| |\layer_{l+1}| - \overline{\rS}_{l + 1} \big| \geq \epsilon \overline{\rS}_{l + 1} \mid \ball_{l}) \leq 2 \exp\Big( - \frac{\epsilon^2 \cdot \overline{\rS}_{l + 1}^{2}/2 }{ \overline{\rS}_{l + 1} + \epsilon \overline{\rS}_{l + 1}/3 } \Big) \leq 2 e^{-\frac{3}{8}\epsilon^2 \overline{\rS}_{l + 1}}.
 \end{align}
Consequently, \eqref{eqn:l_odd_and_x_V2} and \eqref{eqn:l_even_and_x_V2} can be proved by leveraging the definition of $\overline{\rS}_{l + 1}$ in \eqref{eqn:ESl+1} together with the parity of $l$ and the fact $d_1 = d\qr^{-2}$, $d_2 = d\qr^{2}$.

The proof for case $x\in \gV_1$ can be concluded similarly.
\end{proof}

\begin{proof}[Proof of \Cref{lem:few_even_cycle_prob}]
The proof idea follows in the same way as in \cite[Lemma 5.5]{alt2021extremal}. We include details for the sake of completeness. 

For $x\in \gV$, $r\in [N]$, $k\geq 1$, let $\setH_k$ be the collection of all connected graphs $\gH$ by
\begin{align}
    \setH_{k}\coloneqq \{\gH: x\in \gV(\gH), \layer_{1}^{\gH}\subset \layer_{1}, |\gE(\gH)|=|\gV(\gH)| - 1 + k, |\gV(\gH)| \leq 2kr + 1\}.
\end{align}
Note that the requirement of $\setH_{k}$ is deterministic when conditioned on $\layer_{1}\coloneqq \layer_{1}^{\gG}(x) = \rD_{x}$.

The right-hand side of \eqref{eqn:few_even_cycle_prob} is derived from the probability of a certain graph's existence through the following inclusion relationship,
\begin{align}
    \{|\gE(\gG|_{\ball_r})|-|\ball_r|+1\geq k \} \subseteq \{\exists \gH \in \setH_{k}, \gE(\gH) \subseteq \gE(\gG)\}. \label{eqn:inclusion_few_cycle}
\end{align}
Consequently, when conditioned on $\layer_{1}$, by a union bound,
\begin{align}
    \P \big( |\gE(\gG|_{\ball_r})|-|\ball_r| + 1 \geq k \mid \layer_{1} \big) \leq \P(\exists \gH \in \setH_{k}, \gE(\gH) \subseteq \gE(\gG) \mid \layer_1) \leq \sum_{\gH\in \setH_{k}} \P(\gE(\gH) \subseteq \gE(\gG) \mid \layer_1).
\end{align}
Note that for any $\gH\in \setH_k$, we have 
\begin{align}
    \P(\gE(\gH)\subset \gE(\gG)|\layer_1)=\bigg(\frac{d}{\sqrt{mn}} \bigg)^{\gE(\gH)-|\layer_{1}^{\gH}|}=\bigg(\frac{d}{\sqrt{mn}}\bigg)^{|\gV(\gH)|-1+k-|\layer_1^{\gH}|}.
\end{align}
We denote $q_{1} \coloneqq |\layer_{1}^{\gH}|$ and $q_{2} \coloneqq |\gV(\gH) \setminus (\layer_{1}^{\gH} \cup \{x\})|$, thus $|\gV(\gH)| = q_{1} + q_{2} + 1$. Let $C_{q, k}$ denote the number of connected graphs on $q$ vertices with exactly $q -1 + k$ edges. Note that such a graph can be written as a union of a tree on $q$ vertices and $k$ additional edges. Then by Cayley's formula \Cref{lem:Cayley_formula}, we have
\begin{align}
    C_{q, k} \leq q^{q-2} \Big[\binom{q}{2}\Big]^{k} \leq 2^{-k} q^{q + 2k-2}.
\end{align}

To sum over $\gH \in \setH_{k}$, we first sum over all such $\gH \in \setH_{k}$ with exactly $q_{1} + q_{2} + 1$ vertices, then over all possible $q_{1} + q_{2}$. Therefore,
\begin{align}
    &\, \P(\exists \gH \in \setH_{k}, \gE(\gH)\subset \gE(\gG)|\layer_1) \\
    \leq &\, \sum_{j=1}^{2kr} \sum_{q_{1} + q_{2} = j} \binom{|\layer_{1}|}{q_{1}} \binom{N - |\layer_{1}| - 1}{q_{2}} C_{j+1, k} \bigg(\frac{d}{\sqrt{mn}}\bigg)^{q_{2} + k}\\
    \leq &\, \sum_{j=1}^{2kr} \sum_{q_{1} + q_{2} = j} \frac{|\layer_{1}|^{q_{1}}}{q_{1}!} \frac{N^{q_{2}}}{q_{2} !} 2^{-k} (j + 1)^{j+2k - 1} \bigg(\frac{d}{\sqrt{mn}}\bigg)^{q_{2} + k}\\
    = &\, \frac{d^{k}}{(2\sqrt{mn})^{k}} \sum_{j=0}^{2kr}\frac{1}{j!}(|\layer_{1}| + d)^{j} (j + 1)^{j+2k - 1} \quad \textnormal{(Binomial theorem)}\\
    \lesssim &\, \frac{d^{k}}{(2\sqrt{nm})^k} (d + \rD_{x})^{2kr}(2kr)^{2k},
\end{align}
where in the last step, we only keep the last term in the summation since it maintains the largest order according to Stirling \Cref{lem:stirling}.

What remains is the proof of \eqref{eqn:inclusion_few_cycle}, which is established by explicitly building a graph $\gH \in \setH_{k}$ from the requirement $|\gE(\gG|_{\ball_r})|-|\ball_r|+1\geq k$. Let $\tree$ be a spanning tree of $\ball_{r}$ such that $d^{\tree}(x,y) = d^{\gG}(x, y)$ in \eqref{eqn:graph_distance} for all $x, y\in \ball_{r}$. From the left-hand side of \eqref{eqn:inclusion_few_cycle}, there exist $k$ edges in $\gG|_{\ball_{r}}$, denoted by $\gE_{1}$, such that $\gE_{1}\not\subset \gE(\tree)$. Let $\sU_{1}$ denote the set of vertices incident to the edges $\gE_{1}$. For each $y\in \sU_{1}$, according to the tree property, there is a unique path $\ell(x, y)$ connecting $x$ and $y$. Let $\sU_{2} \coloneqq \{z \in \gV(\cup_{y\in \sU_{1}}\ell(x, y)), z\notin \sU_{1} \}$ denote the vertices in the unique paths of $\tree$ connecting $x$ and $\sU_{1}$, and $\gE_{2}$ denote the edges of those paths. Define the graph $\gH$ by $\gV({\gH}) = \{x\} \cup \sU_{1} \cup \sU_{2}$, $\gE(\gH) = \gE_{1} \cup \gE_{2}$. The only thing left is to show that $\gH \in \setH_{k}$, where the non-trivial property to verify is $\gV({\gH}) \leq 2kr + 1$. Obviously, $|\sU_{1}|\leq 2k$ since $|\gE_{1}| =k$, and each path $\ell(x, y)$ for $y\in \sU_{1}$ has at most $r-1$ vertices in $\gV\setminus (\{x\} \cup \sU_{1})$. Consequently, $|\gV(\gH)| \leq 1 + 2k + 2k(r-1) = 2kr + 1$.

We now prove \eqref{eqn:few_even_cycle}. Note that according to \eqref{eqn:few_even_cycle_prob}, for some constant $k\geq 1$,
\begin{align}
    \P \big( |\gE(\gG|_{\ball_r})|-|\ball_r|+1 \leq k \mid \layer_1 \big) \geq 1 - \frac{d^{k}}{(2\sqrt{nm})^k} (d+\rD_{x})^{2kr}(2kr)^{2k}.
\end{align}
Then $\gG|_{\ball_r}$ can be written as a union of a tree $\tree$ and $k$ additional edges denoted by $\gE_{1}$, since $|\gV(\tree)|= |\gE(\tree)| + 1$ for every tree $\tree$. We add those $k$ edges sequentially to the tree $\tree$ to create cycles. In each step, the number of newly created cycles is finite, since the new cycle should either be completely new, or preserve some edges in existing cycles. Therefore, \eqref{eqn:few_even_cycle} follows since the power set of $\gE_{1}$ is finite.
\end{proof}

\begin{proof}[Proof of \Cref{cor:few_cycle_counts}]
    Note that in $\gG|_{\ball_{r}}$, we have
    \begin{align}
        \sum_{y\in \layer_{j}} |\NS{j-1}{y}| = \<\rmA \ones_{\layer_{j-1}}, \ones_{\layer_{j}}\> = \<\ones_{\layer_{j-1}}, \rmA \ones_{\layer_{j}}\> = \sum_{y\in \layer_{j-1}} |\NS{j}{y}|. 
    \end{align}
    According to \Cref{lem:few_even_cycle_prob}, conditioned on $\layer_{1}$, with very high probability, $\gG|_{\ball_{r}}$ can be written as a union of a tree $\tree$ and $k$ additional edges. Let $\rmA^{\tree}$ denote the adjacency matrix of $\tree$, then for each $1\leq j \leq r_x$
    \begin{align}
         \<\rmA^{\tree} \ones_{\layer_{j-1}^{\tree}}, \ones_{\layer_{j}^{\tree}}\> = \<\ones_{\layer_{j-2}^{\tree}}, \ones_{\layer_{j}^{\tree}}\> + \<\ones_{\layer_{j}^{\tree}}, \ones_{\layer_{j}^{\tree}}\> = |\layer_{j}^{\tree}|=|\layer_{j}|,
    \end{align}
    where the second to last equality holds since $\layer_{j-2} \cap \layer_{j} = \emptyset$, otherwise there are cycles in the tree. Thus \eqref{eqn:few_cycle_different_layer} follows directly since $k$ is finite. Then \eqref{eqn:few_cycle_same_layer} follows immediately since only finitely many edges are added to the tree $\tree$, and the number of edges with both ends in the same layer should also be finite.
\end{proof}

\section{Deferred proofs in \Cref{sec:outlier_locations}}
\subsection{Proof of \Cref{lem:existence_pruned_graph}}\label{sec:proof_existence_prunded_graph}
We first state the algorithms to construct $\gH^{(1)}$ and $\gH^{(2)}$, then the pruned graph is obtained by $\pruneG = \gG \setminus ( \gH^{(1)} \cup \gH^{(2)})$. After that, we will prove that $\pruneG$ satisfies all the properties above.

First, $\gH^{(1)}$ is the subgraph such that $\ball^{\gG \setminus \gH^{(1)}}_{2\radius}(x)$ is a tree for each $x\in \pruneV$. For each $y\in \layer_{1}^{\gG}(x)$, let
\begin{align}
   \sW_{x, \, r}^{(1)} (y)\coloneqq \{ z \in \gV \mid d^{\gG}(y, z) \leq r, \, e_{x, y}\notin \ell_{y, z},\,\, \ell_{y, z} \textnormal{ is the path connecting $y$ and $z$}\}
\end{align}
denote the set of vertices reachable via $y$ with distance at most $r$ on $\gG$, but not traversing the edge $e_{x, y}$ between $x$ and $y$. If the branch $\gG|_{\sW_{x, \, \radius }^{(1)}(y)}$ is not a tree, then the edge $e_{x, y}$ is added to $\gH^{(1)}$.

Second, $\gH^{(2)}$ is constructed such that the distance between any two vertices $x, y\in \pruneV$ is at least $2\radius + 1$. For each $x\in \pruneV$, let 
\begin{align}
    \sW^{(2)}_{\tau,\, r}(x) \coloneqq \Big( \pruneV \cap \ball^{\gG \setminus \gH^{(1)}}_{r}(x) \Big) \setminus \{x\}
\end{align}
denote the set of vertices in $\pruneV$ that are connected with $x$ by paths in $\gG \setminus \gH^{(1)}$ with distance at most $r$. According to the previous step, $\ball^{\gG \setminus \gH^{(1)}}_{2\radius}(x)$ is a tree for each $x\in \pruneV$. Then for each vertex $y \in \sW^{(2)}_{\tau,\, 2\radius}(x)$, there must exist a unique vertex $z\in \layer^{\gG \setminus \gH^{(1)}}_{1}(x)$ such that the path $\ell_{x, y}$, which is of length at most $2\radius$, has to cross the edge $e_{x, z}$. All such edges $e_{x, z}$ are added to $\gH^{(2)}$.

Let $\pruneG = \gG \setminus (\gH^{(1)} \cup \gH^{(2)})$. We now prove that $\pruneG$ satisfies all the properties. 

\begin{proof}[Proof of \eqref{item:one}, \eqref{item:two} and \eqref{item:three} in \Cref{lem:existence_pruned_graph}]  By the second step construction, each path $\ell_{x, y}$ in $\pruneG$ with $x, y \in \pruneV$ and $x \neq y$ has length at least $2\radius + 1$. Moreover, $\pruneG \subset \gG \setminus  \gH^{(1)}$, where $\gG \setminus  \gH^{(1)} |_{\ball^{\gG \setminus \gH^{(1)}}_{\radius}(x)}$ is a tree by first step, thus $\pruneG$ is a tree since only branches with cycles were removed. This completes the proof of \eqref{item:one} and \eqref{item:two}. Note that edges incident to $x\in \pruneV$ are added to $\gH^{(1)}$ and $\gH^{(2)}$ during the construction, namely
\begin{align}
   \gE(\gG\setminus \pruneG) = \gE (\gH^{(1)}\cup \gH^{(2)}) \subset \bigcup_{x\in \pruneV} \bigcup_{y \in \ball^{\gG}_{1}(x)} e_{x, y}.
\end{align}
Consequently, there exists at least one vertex in $\pruneV$ incident to $e_{x, y}$, which completes the proof of \eqref{item:three}.
\end{proof}

\begin{proof}[Proof of \eqref{item:four} in \Cref{lem:existence_pruned_graph}] It follows directly, since all branches in $\ball^{(\vtau)}_{\radius}(x)$ are unchanged for $x\in \pruneV$.
\end{proof}

To prove \eqref{item:five} in \Cref{lem:existence_pruned_graph}, we need the following Lemma. For any $x \in \pruneV$, with very high probability, it provides an upper bound on the number of vertices in $\pruneV$ other than $x$, whose distance from $x$ is sufficiently small. 
\begin{lemma}\label{lem:few_high_degree_in_neighbour}
Recall $\benrate_{\qr^{-1}}(\tau)$, $\benrate_{\qr}(\tau)$ and $\mathrm{r}_{\qr^{-1}}(\tau)$, $\mathrm{r}_{\qr}(\tau)$ defined in \eqref{eqn:htau} and \eqref{eqn:rtau} respectively. For $\tau_1 > \qr^{-2}$, $\tau^{+}_{2} > \qr^{2}$ and $0 < \tau^{-}_{2}< \qr^{2}$, the following holds with very high probability.
\begin{enumerate}
    \item When $x\in \Vonehigh$, for any $r\in \N$ with $r \leq \mathrm{r}_{\qr^{-1}}(\tau_1)$, we have 
    \begin{align}
        |\ball_{r}^{\gG}(x) \cap \Vonehigh | \lesssim &\, d^{-1} \log(N) / \benrate_{\qr^{-1}}(\tau_1).
    \end{align}
    
    \item When $x\in \Vtwohigh$, for any $r\in \N$ with $r \leq \mathrm{r}_{\qr}(\tau^{+}_{2})$, we have
    \begin{align}
        |\ball_{r}^{\gG}(x) \cap \Vtwohigh| \lesssim &\, d^{-1}\log(N) / \benrate_{\qr}(\tau^{+}_{2}).
    \end{align}
    For any $r\in \N$ with $r \leq \mathrm{r}_{\qr^{-1}}(\tau_1) \wedge \mathrm{r}_{\qr}(\tau^{+}_{2})$, we have
    \begin{align}
        |\ball_{r}^{\gG}(x) \cap \Vonehigh | \lesssim &\, d^{-1} \log(N)  \cdot [\benrate_{\qr^{-1}}(\tau_1) \wedge \benrate_{\qr}(\tau^{+}_{2})]^{-1}.
    \end{align}
    
    \item When $x\in \Vtwolow$, for any $r\in \N$ with $r \leq \mathrm{r}_{\qr}(\tau^{-}_{2})$, we have
    \begin{align}
        |\ball_{r}^{\gG}(x) \cap \gV_{1}^{(\leq \tau^{-}_{2})} | \lesssim &\, d^{-1} \log(N) / \benrate_{\qr}(\tau^{-}_{2}).
    \end{align}
    For any $r\in \N$ with $r \leq \mathrm{r}_{\qr}(\tau^{+}_{2}) \wedge \mathrm{r}_{\qr}(\tau^{-}_{2})$, we have
    \begin{align}
        |\ball_{r}^{\gG}(x) \cap \Vtwohigh | \lesssim &\, d^{-1} \log(N)  \cdot [\benrate_{\qr}(\tau^{+}_{2}) \wedge \benrate_{\qr}(\tau^{-}_{2})]^{-1}.
    \end{align}
    For any $r\in \N$ with $r \leq \mathrm{r}_{\qr^{-1}}(\tau_1) \wedge \mathrm{r}_{\qr}(\tau^{-}_{2})$, we have
    \begin{align}
        |\ball_{r}^{\gG}(x) \cap \Vonehigh | \lesssim &\, d^{-1} \log(N)  \cdot [\benrate_{\qr^{-1}}(\tau_1) \wedge \benrate_{\qr}(\tau^{-}_{2})]^{-1}.
    \end{align}
\end{enumerate}
\end{lemma}
The proof of \Cref{lem:few_high_degree_in_neighbour} is deferred to \Cref{sec:few_high_degree_in_neighbour}.

\begin{proof}[Proof of \eqref{item:five} in \Cref{lem:existence_pruned_graph}]
Note that the degree on $\gG \setminus \pruneG$ adapts the following decomposition
\begin{align}
    \rD_{x}^{\gG \setminus \pruneG} = \rD_{x}^{\gH^{(1)}} + \rD_{x}^{\gH^{(2)}}.
\end{align}
Let $q_{j}$ denote the maximal number of vertices in $\pruneV$, where those vertices are in the ball of radius $j$ around a vertex in $\pruneV$, formally
\begin{align}
    q_{j} \coloneqq \max_{x\in \pruneV} |\pruneV \cap \ball_{j}^{\gG}(x) \setminus \{x\}|.
\end{align}

For $x\in \pruneV$, let $y \in \layer_{1}^{\gG}$, where the branch $\gG|_{\sW^{(1)}_{x, \, \radius }(y)}$ is not a tree. Then $|\gE(\sW^{(1)}_{x, \, \radius }(y))| - |\sW^{(1)}_{x, \, \radius }(y)| + 1 >0$. According to \Cref{lem:few_even_cycle_prob} and \Cref{cor:few_cycle_counts}, with very high probability, the number of edges in $\ball^{\gG}_{2\radius}$ that prevent it from being a tree is bounded from above, then $\rD_{x}^{\gH^{(1)}} \leq \const + q_{1}$. Furthermore, $\rD_{x}^{\gH^{(2)}} \leq q_{2\radius}$, hence
\begin{align}
   \rD_{x}^{\gG \setminus \pruneG} &\, = \rD_{x}^{\gH^{(1)}} + \rD_{x}^{\gH^{(2)}} \leq \const + q_{1} + q_{2\radius} \\
   &\, \lesssim \const + \frac{\log(N)}{d} \cdot \Big( \min\Big\{\benrate_{\qr^{-1}}(\tau_{1}),\,\, \benrate_{\qr}(\tau^{-}_{2}),\,\,\benrate_{\qr}(\tau^{+}_{2}) \Big\}\Big)^{-1},
\end{align}
where the last inequality follows by considering all scenarios in \Cref{lem:few_high_degree_in_neighbour}.

For $y \not\in \pruneV$, the only meaningful case is $y\in \layer_{1}^{\gG}(x)$ for some $x\in \pruneV$, then by \Cref{lem:few_high_degree_in_neighbour},
\begin{align}
    \rD_{y}^{\gG \setminus \pruneG} \leq |\layer_{1}^{\gG}(x) \cap \pruneV| \leq |\ball_{2}^{\gG}(x) \cap \pruneV|.
\end{align}
which completes the proof of the desired result.
\end{proof}

\begin{proof}[Proof of \eqref{item:six} in \Cref{lem:existence_pruned_graph}] We present the proof for case $x\in \gV_{2}$. The proof follows similarly for $x\in \gV_{1}$. When constructing $\pruneG$, several branches of $x$ are removed from $\gG$, so we need to count the number of vertices in those branches. Note that the ball $\ball_{2\radius}^{\radius}$ present a bipartite bi-regular structure for $j\geq 2$, then by \Cref{lem:concentrationSi}
\begin{align}
    |\layer_{j}^{\gG}(x) \setminus \layer_{j}^{(\vtau)}(x)| \leq \sum_{\substack{z\in \layer_{1}^{\gG}(x)\\ \{x, z\} \in \gH^{(1)} \cup \gH^{(2)}} }  |\layer_{j-1}^{\gG}(z)| \leq \rD_{x}^{\gG \setminus \pruneG} \cdot ( d^{j-1} \indi{j \textnormal{ odd}} + d^{j-2} d_{1} \indi{j \textnormal{ even}} ).
\end{align}
The proof is then completed by applying the fact $d_{1} = \qr^{-2}d$.
\end{proof}

\subsection{Proof of \Cref{lem:few_high_degree_in_neighbour}}\label{sec:few_high_degree_in_neighbour}
\begin{proof}[Proof of \Cref{lem:few_high_degree_in_neighbour} (1) and (2)]
Throughout the proof, we denote $\tau_2$ instead of $\tau^{+}_{2}$ for convenience. For $x \in \gV_{2}^{(\geq \tau_2)}$, we first derive the upper bound on the number of vertices in $\Vonehigh  \cap \ball_{r}(x)$. For each $l\in \N$ and $r \leq \mathrm{r}_{\qr^{-1}}(\tau_1) \wedge \mathrm{r}_{\qr}(\tau_2)$ with $\tau_1 > \qr^{-2}$ and $\tau_2 > \qr^{2}$, define the event
    \begin{align}
        \Xi^{(l, \geq)}_{2,1} \coloneqq \{ \exists x\in [N]: x\in \gV_{2}^{(\geq \tau_2)}, |\ball_{r}(x) \cap \Vonehigh | \geq l\},
    \end{align}
then there exists vertices $y_1, \ldots, y_{l} \in \Vonehigh $ such that $x$ is connected with $y_j$ by path \[z^{(j)} = (x, z_{1}^{(j)}, \ldots, z_{r_{j}}^{(j)}, y_j)\] of length $r_j + 1$ for each $j\in [l]$. Our goal is to show that $\P(\Xi^{(l, \geq)}_{2,1})$ is vanishing for some $l \ll \log(N)$.

For convenience, we denote the $l$-tuples by $\rvy = (y_1, \ldots, y_{l})$ and $\rvz = (z^{(1)},\ldots ,z^{(l)})$. It suffices to consider the case where those $l$ paths are pairwise disjoint. Otherwise, cycles would be created by intersecting paths, which occurs with a probability upper bounded by $N^{-\nu}$ for some $\nu >0$ according to \Cref{lem:few_even_cycle_prob} and \Cref{cor:few_cycle_counts}. For each fixed pair of $x, \rvy, \rvz$, define the following event
   \begin{align}
       \Xi^{(l, \geq)}_{x, \rvy, \rvz} \coloneqq &\, \{x\in \gV_{2}^{(\geq \tau_2)},\,\, \rvy \subset \Vonehigh , \,\, \rvz \subset \gE(\gG)\}.
   \end{align}
Then for each $l\in \N$, the event $\Xi^{(l, \geq)}_{2,1}$ admits the decomposition
\begin{align}
    \Xi^{(l, \geq)}_{2,1} = \bigcup_{x, \rvy, \rvz} \Xi^{(l, \geq)}_{x, \rvy, \rvz},
\end{align}
where the union is taken over all $x\in \gV_{2}^{(\geq \tau_2)}$, $\rvy \subset \Vonehigh $, and $\rvz$. Now, the task is transferred to establish an upper bound on $\P(\Xi^{(l, \geq)}_{x, \rvy, \rvz})$.  Furthermore, $\gG\setminus \rvz$ denotes the bipartite graph on vertices $\gV(\gG) = \gV_1 \cup \gV_2$ excluding the edges in $\rvz$ and define the following event
   \begin{align}
       \Xi^{(l, \geq)}_{x, \rvy} \coloneqq &\, \{\rD_x \geq \tau_2 d - l, \rD_{y_1} \geq \tau_1 d - 1, \ldots, \rD_{y_{l}} \geq \tau_1 d - 1 \textnormal{ on } \gG\setminus \rvz\} \label{eqn:Xi_xy^l>=}.
   \end{align}
Note that $\Xi^{(l, \geq)}_{x, \rvy}$ is independent of $\rvz$, then 
    \begin{align}\label{eqn:Xi_xyz^l+_prob}
        \P(\Xi^{(l, \geq)}_{x, \rvy, \rvz}) = &\, \P\Big(\Xi^{(l, \geq)}_{x, \rvy} \mid  \rvz \subset \gE(\gG) \Big) \cdot \P(\rvz \subset \gE(\gG) ) \leq \P(\Xi^{(l, \geq)}_{x, \rvy}) \cdot \prod_{j=1}^{l}p^{r_j + 1}.
    \end{align}
With the proof deferred later, for each $x, \rvy$, the following bound holds
\begin{align}
    \P(\Xi^{(l, \geq)}_{x, \rvy}) \leq &\, \exp\big[ - d(l + 1)\cdot \benrate_{\qr^{-1}}(\tau_1) \wedge \benrate_{\qr}(\tau_2) \big] + (l+1) \cdot [g_l((\tau_1 - \qr^{-2})/2) \vee g_l((\tau_2 - \qr^{2})/2)], \label{eqn:Xi_xy^l>=_prob}
\end{align}
where $g_{l}(t)$ is defined via
\begin{align}
    g_{l}(t) \coloneqq \frac{1}{\sqrt{2\pi td}} \Big(\frac{e(l+1)}{t\sqrt{mn}} \Big)^{td}\label{eqn:glt}.
\end{align}

We now finish the proof of (1) by combining previous estimates. Recall $N = n + m$, $p = d/{\sqrt{mn}}$, then
\begin{align}
    &\, \P(\Xi^{(l, \geq)}_{2,1} ) \coloneqq \P(\cup_{x, \rvy, \rvz} \Xi^{(l, \geq)}_{x, \rvy, \rvz}) \leq \sum_{x, \rvy, \rvz} \P(\Xi^{(l, \geq)}_{x, \rvy, \rvz})\\
    \leq &\, \binom{m}{1} \cdot \binom{n}{l}\sum_{r_{1} = 0}^{r-1} \cdots \sum_{r_{l} = 0}^{r-1} \binom{N - l - 1}{r_1} \cdots \binom{N - l - 1 - \sum_{j = 1}^{l-1} r_j}{r_l} \cdot p^{l + \sum_{t=1}^{l}r_t} \cdot \max_{x, \rvy}  \P(\Xi^{(l, \geq)}_{x, \rvy})\\
    \leq &\, N^{l+1}\cdot \sum_{r_1, \ldots, r_l = 0}^{r-1}N^{\sum_{t=1}^{l} r_t} \cdot p^{l + \sum_{t=1}^{l}r_t} \cdot \max_{x, \rvy}  \P(\Xi^{(l, \geq)}_{x, \rvy}) \\
    \leq &\, N \cdot \Big( \sum_{t=0}^{r-1}(\qr^{2} + \qr^{-2})\cdot d^{t + 1}\Big)^{l} \cdot \max_{x, \rvy}  \P(\Xi^{(l, \geq)}_{x, \rvy}) \\
    \leq &\, \circled{1} + \circled{2} \lesssim N^{-\nu},
\end{align}
where the second line relies on \eqref{eqn:Xi_xyz^l+_prob} and all possible choices of $x \in \gV_{2}^{(\geq \tau_2)}, \rvy \subset \Vonehigh $ are went through, and the last line relies on \eqref{eqn:Xi_xy^l>=_prob} with \circled{1}, \circled{2} defined as 
\begin{align}
    \circled{1} =&\, N(\qr^{2} + \qr^{-2})^{l} \Big( \frac{d^{r+1} - 1}{d - 1}\Big)^{l} \cdot \exp\Big( - d(l + 1)\cdot [\benrate_{\qr^{-1}}(\tau_1) \wedge \benrate_{\qr}(\tau_2)]  \Big), \\
    \circled{2} =&\, N(\qr^{2} + \qr^{-2})^{l} \Big( \frac{d^{r+1} - 1}{d - 1}\Big)^{l} \cdot (l+1) \cdot [g_l((\tau_1 - \qr^{-2})/2) \vee g_l((\tau_2 - \qr^{2})/2)].
\end{align}
The remaining thing is to verify that \circled{1} and \circled{2} are both dominated by $N^{-\nu}$ for some positive constant $\nu$. Firstly, $\circled{1} \leq N^{-\nu}$ is true if the following holds
\begin{align}
    \log(N) + (l + 1) \cdot \Big( \log(\qr^{2} + \qr^{-2}) + (r + 1)\log d - d\cdot [\benrate_{\qr^{-1}}(\tau_1) \wedge \benrate_{\qr}(\tau_2)] \Big) \leq -\nu \log(N),
\end{align}
which follows easily since $\qr \asymp 1$, $(r + 1)\log d \leq d h_{\min}(\tau_1, \tau_2)/2$ for $r \leq r(\tau_1, \tau_2)$ and $l \ll \log(N)$. 

\noindent Secondly, $\circled{2} \leq N^{-\nu}$ is true if the following is satisfied
\begin{align}
    &\, \log(N) + l \cdot \Big( \log(\qr^{2} + \qr^{-2}) + (r + 1)\log d \Big)  \\
    &\, - \frac{1}{2}(\tau_2 - \qr^{2})d \cdot \Big( \log((\tau_2 - \qr^{2})/2) + \log(N) - \log(l + 1) \Big) - \frac{1}{2}\log[(\tau_2 - \qr^{2})d/2] \leq -\nu \log(N),
\end{align}
which holds true under the previous conditions.

~\\
For the second part of (2), we define the event for $r \leq \mathrm{r}_{\qr}(\tau_2)$
\begin{align}
    \Xi^{(l, \geq)}_{2,2} \coloneqq \{ \exists x\in [N]: x\in \gV_{2}^{(\geq \tau_2)}, |\ball_{r}(x) \cap \gV_{2}^{(\geq \tau_2)}| \geq l\}.
\end{align} 
Similarly, $\P(\Xi^{(l, \geq)}_{2,2} ) \lesssim N^{-\nu}$ can be established by following the analysis above and slightly modifying the proof of \eqref{eqn:Xi_xy^l>=_prob}, where the $l$-tuple $\rvy = (y_1, \ldots, y_l)$ is now a collection of distinct elements $y_j$ from $\gV_{2}^{(\geq \tau_2)}$, and $\benrate_{\qr^{-1}}(\tau_1) \wedge \benrate_{\qr}(\tau_2)$ is replaced by $\benrate_{\qr}(\tau_2)$ throughout the concentration analysis.

~\\
For (1), we define the event for $r \leq \mathrm{r}_{\qr^{-1}}(\tau_1)$
\begin{align}
    \Xi^{(l, \geq)}_{1,1} \coloneqq \{ \exists x\in [N]: x\in \Vonehigh , |\ball_{r}(x) \cap \Vonehigh | \geq l\}.
\end{align} 
Then $\P(\Xi^{(l, \geq)}_{1,1} ) \lesssim N^{-\nu}$ can be established similarly, and $\benrate_{\qr^{-1}}(\tau_1) \wedge \benrate_{\qr}(\tau_2)$ is replaced by $\benrate_{\qr^{-1}}(\tau_1)$.
\end{proof}

\begin{proof}[Proof of \eqref{eqn:Xi_xy^l>=_prob}]
Recall the definition of $\Xi^{(l, \geq)}_{x, \rvy}$ in \eqref{eqn:Xi_xy^l>=}. Based on the following observation
\begin{align}
       \Xi^{(l, \geq)}_{x, \rvy}
       \subseteq \{\rD_x \geq \tau_2 d - l, \rD_{y_1} \geq \tau_1 d - 1, \ldots, \rD_{y_{l}} \geq \tau_1 d - 1 \textnormal{ on } \gG\},
\end{align}
the remaining discussion will be focusing on $\gG$ instead of $\gG \setminus \rvz$, since only an upper bound of $\P(\Xi^{(l, \geq)}_{x, \rvy})$ is needed. We should point out that the involvement of the extra edges in $\rvz$ would not lead to a trivial upper bound, since only $l \ll \log N \ll N$ more edges are considered for each vertex, which will not make any non-trivial difference to the concentration arguments below.

We are only interested in the case $\tau_1 > \qr^{-2}$ and $\tau_2 > \qr^{2}$. First, $\rD_{y_0}, \ldots, \rD_{y_{l}}$ are conditionally independent. Denote the set of vertices by $\sY = \{y_0, \ldots, y_{l}\}$ where $y_0 \coloneqq x$, then $\rD_{y_j} \coloneqq \rD^{\gG}_{y_j}  = \rD^{\gG \setminus \sY}_{y_j} + \rD^{\sY}_{y_j}$ where
\begin{align}
    \rD^{\sY}_{y_j} \coloneqq \sum_{z\in \sY} \ermA_{z y_i}, \quad \rD^{\gG \setminus \sY}_{y_j} \coloneqq \sum_{z\in \gV \setminus \sY} \ermA_{z y_i}, \quad 0 \leq j \leq l.
\end{align}
Thus $\{\rD^{\gG \setminus \sY}_{y_j}\}_{j=0}^{l}$ are independent since $\rD^{\sY}_{y_j}$ is measurable when conditioned on $\rmA_{\sY}$, where $\rmA_{\sY}$ denotes the set of edges with both ends in $\sY$. Then by conditional independence and Tonelli's theorem, \eqref{eqn:Xi_xy^l>=_prob} admits the following decomposition
    \begin{align}
        \P(\Xi^{(l,\geq)}_{x, \rvy}) = &\, \E[\P(\rD_{y_0} \geq \tau_2 d - l, \ldots ,\rD_{y_{l}} \geq \tau_1 d - 1 \mid \rmA_{\sY})]\\
        = &\, \E \bigg[ \prod_{j=0}^{l} \P(\rD_{y_j} \geq a_j \mid \rmA_{\sY})\bigg],
    \end{align}
where $a_0 = \tau_2 d - l$, $a_1 = \ldots = a_l = \tau_1 d - 1$. Note that $\E\rD^{\gG}_{y_0} = \qr^{2}d$, $\E\rD^{\gG}_{y_j} = \qr^{-2}d$ for $j\in [l]$ since 
$y_0 \in \gV_2$ and $\sY \setminus \{y_0\} \subset \gV_1$, then $\rD_{y_j} \geq a_j$ if and only if the following holds
\begin{align}
    \rD^{\gG \setminus \sY}_{y_j} - \E\rD^{\gG \setminus \sY}_{y_j} \geq (a_j - \E\rD^{\gG}_{y_j}) - (\rD^{\sY}_{y_j} - \E\rD^{\sY}_{y_j}),
\end{align}
where $\E\rD^{\gG \setminus \sY}_{y_0} = (1-l/n)\qr^{2}d = (1 -o(1))\qr^{2}d$ and $\E\rD^{\gG \setminus \sY}_{y_j} = (1-1/m)\qr^{-2}d = (1 -o(1))\qr^{-2}d$ for $j\in [l]$ since $l \ll d \ll n$, $m \gg 1$. By Bennett inequality in \Cref{lem:Bennett},
\begin{align}
    \P(\rD_{y_j} \geq a_j \mid \rmA_{\sY}) \leq &\, \exp\bigg( -\Var(\rD^{\gG \setminus \sY}_{y_j}) \cdot \benrate\Big(\frac{(a_j - \E\rD^{\gG}_{y_j}) - (\rD^{\sY}_{y_j} - \E\rD^{\sY}_{y_j})}{\Var(\rD^{\gG \setminus \sY}_{y_j})} \Big) \bigg).
\end{align}
Note that $\Var(\rD^{\gG \setminus \sY}_{y_j}) = (1 -o(1)) \cdot \E\rD^{\gG \setminus \sY}_{y_j}$ since $p=o(1)$. Define $\overline{\rD}^{\sY}_{y_j} \coloneqq \rD^{\sY}_{y_j} - \E\rD^{\sY}_{y_j}$, then by conditioning on the events $\overline{\sD}_{0} \coloneqq \{\overline{\rD}^{\sY}_{y_0} \leq (\tau_2 - \qr^{2})d/2\}$ and $\overline{\sD}_{\max} \coloneqq \{\max_{j\in [l]}\overline{\rD}^{\sY}_{y_j} \leq (\tau_1 - \qr^{-2})d/2\}$, we have
\begin{align}
    \P(\rD_{y_0} \geq a_0 \mid \overline{\sD}_{0}) \leq &\, \exp\big( -\qr^{2}d \cdot h((\qr^{-2}\tau_2 -1)/2)\big),\\
    \P(\rD_{y_j} \geq a_j \mid \overline{\sD}_{\max}) \leq &\, \exp\big( -\qr^{-2}d \cdot h((\qr^{2}\tau_1 -1)/2)\big),
    \quad \forall j\in [l].
\end{align} 
Therefore, combining the arguments above, we have
\begin{align}
    &\, \P(\Xi^{(l,\geq)}_{x, \rvy}) \\
    = &\, \P(\Xi^{(l,\geq)}_{x, \rvy} \mid \overline{\sD}_{0} \cap \overline{\sD}_{\max}) \cdot \P(\overline{\sD}_{0} \cap \overline{\sD}_{\max}) + \P(\Xi^{(l,\geq)}_{x, \rvy} \mid (\overline{\sD}_{0} \cap \overline{\sD}_{\max})^{\complement}) \cdot  \P( (\overline{\sD}_{0} \cap \overline{\sD}_{\max})^{\complement})\\
    \leq &\, \P(\Xi^{(l,\geq)}_{x, \rvy} \mid \overline{\sD}_{0} \cap \overline{\sD}_{\max}) + \P( \overline{\sD}_{0}^{\complement} \cup \overline{\sD}_{\max}^{\complement})\\
    \leq &\, \exp\bigg( - d\Big[\qr^{2} \cdot h\Big( \frac{\qr^{-2}\tau_2 -1}{2} \Big) + l\qr^{-2} \cdot h \Big(\frac{\qr^{2}\tau_1 -1}{2}\Big) \Big] \bigg) + \P\big(\overline{\sD}_{0}^{\complement}\big) + \P\big(\overline{\sD}_{\max}^{\complement}\big)\\
    \leq &\, \exp\big( - d(l + 1) \cdot [\benrate_{\qr^{-1}}(\tau_1) \wedge \benrate_{\qr}(\tau_2)] \big) + (l+1) \cdot [g_l((\tau_1 - \qr^{-2})/2) \vee g_l((\tau_2 - \qr^{2})/2)],
\end{align}
where $\benrate_{\qr^{-1}}(\tau_1)$, $\benrate_{\qr}(\tau_2)$ and $g_{l}(t)$ in the last line are defined in  \eqref{eqn:htau} and \eqref{eqn:glt}, respectively. The upper bounds of $\P\big(\overline{\sD}_{0}^{\complement}\big)$ and $ \P\big(\overline{\sD}_{\max}^{\complement}\big)$ are derived as follows. First, recall $p = d/\sqrt{mn}$, then
\begin{align}
    &\, \P\big(\overline{\sD}_{0}^{\complement}\big) \coloneqq \P(\overline{\rD}^{\sY}_{y_0} > (\tau_2 - \qr^{2})d/2) \\
    \leq &\, \P(\rD^{\sY}_{y_0} > (\tau_2 - \qr^{2})d/2) \leq \binom{l+1}{(\tau_2 - \qr^{2})d/2} \Big(\frac{d}{\sqrt{mn}} \Big)^{(\tau_2 - \qr^{2})d/2}\\
    \leq &\, \frac{(l+1)^{(\tau_2 - \qr^{2})d/2}}{(\tau_2 d/2 - \qr^{2}d/2)!} \cdot \Big(\frac{d}{\sqrt{mn}} \Big)^{\frac{(\tau_2 - \qr^{2})d}{2}} = \frac{1}{\sqrt{\pi(\tau_2 - \qr^{2})d }}\Big(\frac{2e(l+1)d}{(\tau_2 - \qr^{2})d\sqrt{mn}} \Big)^{\frac{(\tau_2 - \qr^{2})d}{2}} = g_l((\tau_2 - \qr^{2})/2),
\end{align}
where the last line holds due to Stirling's approximation and binomial coefficients. At the same time,
\begin{align}
&\, \P\big(\overline{\sD}_{\max}^{\complement}\big) \coloneqq \P(\exists j \in [l] \textnormal{ such that } \overline{\rD}^{\sY}_{y_j} > (\tau_1 - \qr^{-2})d/2) \\
\leq &\, \sum_{j\in [l]} \P\big(\overline{\rD}^{\sY}_{y_j} > (\tau_1 - \qr^{-2})d/2 \big) \leq l \cdot \P(\rD^{\sY}_{y_j} > (\tau_1 - \qr^{-2})d/2 ) \\
    \leq &\, l \cdot \binom{l+1}{(\tau_1 - \qr^{-2})d/2} \Big(\frac{d}{\sqrt{mn}} \Big)^{(\tau_1 - \qr^{-2})d/2} \\
    \leq &\, \frac{l}{\sqrt{\pi(\tau_1 - \qr^{-2})d }}\Big(\frac{2e(l+1)d}{(\tau_1 - \qr^{-2})d\sqrt{mn}} \Big)^{(\tau_1 - \qr^{-2})d/2} = l\cdot g_l((\tau_1 - \qr^{-2})/2).
\end{align}
Therefore, it follows easily that $\P\big(\overline{\sD}_{0}^{\complement}\big) + \P\big(\overline{\sD}_{\max}^{\complement}\big) \leq (l+1) \cdot [g_l((\tau_1 - \qr^{-2})/2) \vee g_l((\tau_2 - \qr^{2})/2)]$.
\end{proof}

\begin{proof}[Proof of \Cref{lem:few_high_degree_in_neighbour} (3)]
For each $l \in \N$ and $r\leq \mathrm{r}_{\qr}(\tau^{-}_{2})$ with $\tau^{-}_{2} < \qr^{2}$, define the event
\begin{align}
    \Xi^{(l,\leq)}_{2,2} \coloneqq \{ \exists x\in [N]: x\in \Vtwolow, |\ball_{r}(x) \cap \Vtwolow | \geq l\}.
\end{align}
Similarly, for each $l$-tuples $\rvy = (y_1, \ldots, y_{l}) \subset \Vtwolow$ and $\rvz = (z^{(1)},\ldots ,z^{(l)})$, define
   \begin{align}
       \Xi^{(l,\leq)}_{x, \rvy, \rvz} \coloneqq &\, \{x\in \Vtwolow,\,\, \rvy \subset \Vtwolow , \,\, \rvz \subset \gE(\gG)\}.
   \end{align}
Then we can write $\Xi^{(l,\leq)}_{2,2} = \cup_{x, \rvy, \rvz} \Xi^{(l,\leq)}_{x, \rvy, \rvz}$. For each fixed $x, \rvy, \rvz$, define
   \begin{align}
       \Xi^{(l,\leq)}_{x, \rvy} \coloneqq &\, \{\rD_x \leq \tau^{-}_{2} d - l, \rD_{y_1} \leq \tau^{-}_{2} d - 1, \ldots, \rD_{y_{l}} \leq \tau^{-}_{2} d - 1 \textnormal{ on } \gG\setminus \rvz\} \label{eqn:Xi_xy^l<=}
   \end{align}
With the proof deferred later, for each $x, \rvy$, the following bound holds
\begin{align}
    \P(\Xi^{(l, \leq)}_{x, \rvy}) \leq &\, \exp\big[ - d(l + 1)\cdot \benrate_{\qr}(\tau^{-}_{2}) \big] + (l+1) \cdot g_l((\qr^{2} - \tau^{-}_{2})/2), \label{eqn:Xi_xy^l<=_prob}
\end{align}
where $\benrate_{\qr}(\tau)$ and $g_{l}(t)$ are defined in \eqref{eqn:htau} and \eqref{eqn:glt}, respectively. By taking the union bound and following the same analysis in the proof for (1) and (2), we obtain $\P(\Xi^{(l,\leq)}_{2,2}) \leq N^{-\nu}$ for some constant $\nu > 0$.

~\\
For the second part of (3), we define the event for $r \leq \mathrm{r}_{\qr}(\tau^{+}_{2}) \wedge \mathrm{r}_{\qr}(\tau^{-}_{2})$ and $l \in \N$
\begin{align}
    \Xi^{(l, \leq \geq)}_{2,2} \coloneqq \{ \exists x\in [N]: x\in \gV_{2}^{(\geq \tau^{-}_{2})}, |\ball_{r}(x) \cap \Vtwohigh| \geq l\}.
\end{align} 
Similarly, $\P(\Xi^{(l, \leq \geq)}_{2,2} ) \lesssim N^{-\nu}$ can be established by following the analysis above and combining the estimates in the proofs of \eqref{eqn:Xi_xy^l>=_prob} and \eqref{eqn:Xi_xy^l<=_prob}, where the $l$-tuple $\rvy = (y_1, \ldots, y_l)$ is now a collection of distinct elements $y_j$ from $\Vtwohigh$, and $\benrate_{\qr}(\tau^{+}_{2})$ is replaced by $\benrate_{\qr}(\tau^{+}_{2}) \wedge \benrate_{\qr}(\tau^{-}_{2})$ throughout the concentration analysis.

~\\
For the third part of (3), we define the event for $r \leq \mathrm{r}_{\qr^{-1}}(\tau_1) \wedge \mathrm{r}_{\qr}(\tau^{-}_{2})$ and $l \in \N$,
\begin{align}
    \Xi^{(l, \leq)}_{2,1} \coloneqq \{ \exists x\in [N]: x\in \gV_{2}^{(\geq \tau^{-}_{2})}, |\ball_{r}(x) \cap \Vonehigh | \geq l\}.
\end{align} 
Then $\P(\Xi^{(l, \leq)}_{2,1}) \lesssim N^{-\nu}$ can be established similarly, and $\benrate_{\qr}(\tau^{-}_{2})$ is replaced by $\benrate_{\qr^{-1}}(\tau_1) \wedge \benrate_{\qr}(\tau^{-}_{2})$.
\end{proof}

\begin{proof}[Proof of \eqref{eqn:Xi_xy^l<=_prob}]
    Recall the definition of $\Xi^{(l, \leq)}_{x, \rvy}$ in \eqref{eqn:Xi_xy^l<=} and the condition $\tau^{-}_{2} < \qr^{2}$. For convenience, we denote $\tau_2$ instead of $\tau^{-}_{2}$. Throughout the proof, we focus on $\gG$ instead of $\gG \setminus \rvz$ throughout the proof. We should point out that the involvement of the extra edges in $\rvz$ would not lead to a trivial upper bound, since at most $l \ll \log N$ more edges are considered for each vertex, which will not make any nontrivial difference to the concentration arguments below.

    Similar to the proof of \eqref{eqn:Xi_xy^l>=_prob}, denote $\sY = \{y_0, \ldots, y_{l}\}$ where $y_0 \coloneqq x$. When conditioned on $\rmA_{\sY}$, $\{\rD^{\gG \setminus \sY}_{y_j}\}_{j=0}^{l}$ are independent since each $\rD^{\sY}_{y_j}$ is measurable. Then \eqref{eqn:Xi_xy^l<=_prob} admits the decomposition
    \begin{align}
        \P(\Xi^{(l,\leq)}_{x, \rvy})  = &\, \E \Big[ \prod_{j=0}^{l} \P(\rD_{y_j} \leq a_j \mid \rmA_{\sY})\Big],
    \end{align}
where $a_0 = \tau_2 d - l$, $a_1 = \ldots = a_l = \tau_2 d - 1$. Since $\E\rD^{\gG}_{y_j} = \qr^{2}d$ for each $0 \leq j \leq l$, $\rD_{y_j} \leq a_j$ if and only if
\begin{align}
    \rD^{\gG \setminus \sY}_{y_j} - \E\rD^{\gG \setminus \sY}_{y_j} \leq (a_j - \E\rD^{\gG}_{y_j}) - (\rD^{\sY}_{y_j} - \E\rD^{\sY}_{y_j}),
\end{align}
where $\E\rD^{\gG \setminus \sY}_{y_j} = \Var(\rD^{\gG \setminus \sY}_{y_j}) = (1-l/n)\qr^{2}d = (1 -o(1))\qr^{2}d$ since $l \ll d \ll n$ and $p = o(1)$. Note that $a_j - \E\rD^{\gG}_{y_j} < 0$ since $\tau_2 < \qr^{2}$, then by second part of Bennett \Cref{lem:Bennett},
\begin{align}
    \P(\rD_{y_j} \leq a_j \mid \rmA_{\sY}) \leq &\, \exp\bigg( -\Var(\rD^{\gG \setminus \sY}_{y_j}) \cdot \benrate\Big(\frac{- (a_j - \E\rD^{\gG}_{y_j}) + (\rD^{\sY}_{y_j} - \E\rD^{\sY}_{y_j})}{\Var(\rD^{\gG \setminus \sY}_{y_j})} \Big) \bigg).
\end{align}
Define $\overline{\rD}^{\sY}_{y_j} \coloneqq \rD^{\sY}_{y_j} - \E\rD^{\sY}_{y_j}$, then by conditioning on the event $\overline{\sD}_{\min} \coloneqq \{\min_{0 \leq j \leq l}\,\, \overline{\rD}^{\sY}_{y_j} \geq (\tau_2 - \qr^{2})d/2\}$,
\begin{align}
    \P(\rD_{y_j} \leq a_j \mid \overline{\sD}_{\min}) \leq &\, \exp\big( -\qr^{2}d \cdot h((1 - \qr^{-2}\tau_2)/2)\big),
    \quad \forall \,\, 0 \leq j\leq l.
\end{align} 
Therefore, combining the arguments above, we have
\begin{align}
    &\, \P(\Xi^{(l,\leq)}_{x, \rvy}) = \P(\Xi^{(l,\leq)}_{x, \rvy} \mid \overline{\sD}_{\min}) \cdot \P(\overline{\sD}_{\min}) + \P(\Xi^{(l,\leq)}_{x, \rvy} \mid \overline{\sD}_{\min}^{\complement}) \cdot \P( \overline{\sD}_{\min}^{\complement})\\
    \leq &\, \P(\Xi^{(l,\leq)}_{x, \rvy} \mid \overline{\sD}_{\min}) + \P(\overline{\sD}_{\min}^{\complement})\\
    \leq &\, \exp\big( - d(l + 1) \cdot \benrate_{\qr}(\tau_2) \big) + (l+1) \cdot g_l((\qr^{2} - \tau_2)/2),
\end{align}
where $\benrate_{\qr}(\tau_2)$ and $g_{l}(t)$ in the last line are defined in  \eqref{eqn:htau} and \eqref{eqn:glt}, respectively. The remaining thing is to derive the upper bound of $ \P\big(\overline{\sD}_{\min}^{\complement}\big)$. Recall $p = d/\sqrt{mn}$, then
\begin{align}
&\, \P\big(\overline{\sD}_{\min}^{\complement}\big) \coloneqq \P( \,\, \exists 0 \leq j \leq l \textnormal{ such that }\overline{\rD}^{\sY}_{y_j} < ( \tau_2 - \qr^{2})d/2 ) \\
    \leq &\, \sum_{0 \leq j \leq l} \P\big(\overline{\rD}^{\sY}_{y_j} < ( \tau_2 - \qr^{2})d/2 \big) = \sum_{0 \leq j \leq l} \P\big(-\overline{\rD}^{\sY}_{y_j} < ( \tau_2 - \qr^{2})d/2 \big) \\
    \leq &\, (l + 1) \cdot \P(\rD^{\sY}_{y_j} > (\qr^{2} - \tau_2)d/2) \\
    \leq &\, (l + 1) \cdot \binom{l+1}{(\qr^{2} - \tau_2)d/2} \Big(\frac{d}{\sqrt{mn}} \Big)^{(\qr^{2} - \tau_2)d/2} \\
    \leq &\, \frac{l + 1}{\sqrt{(\qr^{2} - \tau_2)d }}\Big(\frac{2e(l+1)d}{(\qr^{2} - \tau_2)d\sqrt{mn}} \Big)^{(\qr^{2} - \tau_2)d/2} = (l + 1)\cdot g_l((\qr^{2} - \tau_2)/2).
\end{align}
where the second line holds since $\overline{\rD}^{\sY}_{y_j}$ is centered and symmetric, and the last line holds by Stirling's approximation and \Cref{lem:stirling}.
\end{proof}

\section{Deferred proofs in \Cref{sec:bulk_boundness}}\label{sec:proof_bulk_boundness}

\subsection{Proofs in \Cref{sec:Loewner-order}}\label{sec:Loewner-order-app}

\begin{proof}[Proof of \Cref{lem:invertiblility_D1}]
For each $x\in \gV_{1}$, $\E \rD_{x} = \qr^{-2} d \leq d$. Then by Bennett's inequality (\Cref{lem:Bennett}),
    \begin{align}
        \P(\rD_{x} - d \leq -\epsilon d) = \P(\rD_{x} - d_{1} \leq d - d_{1} - \epsilon d) \leq \exp\Big[-\qr^{-2} d \cdot \benrate \Big( \frac{(1-\qr^{-2})d - \epsilon d}{\qr^{-2} d} \Big) \Big].
    \end{align}
Since $d = b\log(N)$ for some constant $b>0$, the following holds through a union bound,
    \begin{align}
        \P(\exists x\in \gV_{1},\,\, \rD_{x} - d \leq -\epsilon d ) \leq n\cdot \exp\Big[- \log(N)\cdot b\cdot \qr^{-2} \cdot \benrate( \qr^{2} - 1 - \qr^{2}\epsilon ) \Big].
    \end{align}
Note that the function $f(\qr) \coloneqq \qr^{-2} \cdot \benrate( \qr^{2} - 1)$ decreases monotonically when $\qr \geq 1$.
    
Since \Cref{ass:Ihara_Bass} is satisfied, there exists a  constant $\epsilon$ such that 
    \begin{align}
        0< b \cdot \qr^{-2} \cdot \benrate(\qr^{2} - 1 - \qr^2 \epsilon) - 1.
    \end{align}
    Taking $\nu = \left (b \cdot \qr^{-2} \cdot \benrate(\qr^{2} - 1 - \qr^2 \epsilon) - 1 \right)/2$ and noting that $|\gV_{1}| = n \asymp N$, with probability at least $1 - N^{-\nu}$, the conclusion follows. 
\end{proof}

\begin{proof}[Proof of \Cref{prop:upper_bound_H}]
Following \Cref{lem:Ihara-Bass} and \Cref{lem:upper_bound_rhoB}, the $\rmH(\lambda)$ associated with $\scA$ in \eqref{eqn:scAdj} is
\begin{align}
   \rmH(\lambda) = \begin{bmatrix}
    \bzero & \rmX(\lambda) \\
    [\rmX(\lambda)]^{\sT} & \bzero
\end{bmatrix}, \quad \ermX_{jl}(\lambda) \coloneqq \frac{\lambda \ermX_{jl}}{\lambda^2 - \ermX_{jl}^2}, \quad j\in \gV_{1},\,\, l\in \gV_{2},
\end{align}
where $\gV_{1} = |n|$ and $|\gV_{2}| = m$. Similarly, the $\rmM(\lambda)$ associated with $\scA$ in \eqref{eqn:scAdj} is defined by
\begin{align}
   \rmM(\lambda) = \begin{bmatrix}
    \rmM^{(1)}(\lambda) & \bzero \\
    \bzero & \rmM^{(2)}(\lambda) 
\end{bmatrix}, \quad \ermM^{(1)}_{jj}(\lambda) \coloneqq 1 + \sum_{l\in \gV_{2}} \frac{\ermX_{jl}^2}{\lambda^2 - \ermX_{jl}^2}, \quad \ermM^{(2)}_{ll}(\lambda) \coloneqq 1 + \sum_{j\in \gV_{1}} \frac{\ermX_{jl}^2}{\lambda^2 - \ermX_{jl}^2}.
\end{align}
Note that both $\rmH(\lambda)$ and $\rmM(\lambda)$ are Hermitian for all $\lambda \in \R$, and $\rmM(\lambda) - \rmH(\lambda) \to \rmI_{N}$ when $\lambda\to \infty$. Hence, $\rmM(\lambda) - \rmH(\lambda)$ is strictly positive definite for sufficiently large $\lambda$. Define $\lambda_{\star}$ by
 \begin{align}
    \lambda_{\star} \coloneqq \inf\{t >0: \rmM(\lambda) - \rmH(\lambda) \succ \bzero \textnormal{ for all } \lambda > t\}.
\end{align}
By continuity, the smallest eigenvalue of $\rmM(\lambda_{\star}) - \rmH(\lambda_{\star})$ is zero. Let $\rmB$ denote the non-backtracking operator associated with $\scA$ in \eqref{eqn:scAdj}. Then by \Cref{lem:Ihara-Bass}, $\lambda_{\star} \in \mathrm{Spec}(\rmB)$ and $\lambda_{\star} \leq \rho(\rmB)$ since $\rho(\rmB)=\max_{\lambda \in \mathrm{Spec}(\rmB)}|\lambda|$. Thus, for any $\lambda \geq \rho(\rmB)$, it implies $\lambda \geq \lambda_{\star}$, then $\rmM(\lambda) \succeq \rmH(\lambda)$. By confirming that $\scA$ in \eqref{eqn:scAdj} satisfies the assumptions of \Cref{lem:upper_bound_rhoB} with $q = \sqrt{d}$, and choosing $\lambda = 1 + \epsilon$, one could conclude that
\begin{align}
    \P\big(\rmM(1 + \epsilon) - \scA(1 + \epsilon) \succeq \bzero \big) \geq \P(\rho(\rmB) \leq 1 + \epsilon ) \geq 1 - N^{3 - c\sqrt{d} \log(1 + \epsilon)}. \label{eqn:MH_lambda_PD}
\end{align}    
Then the following holds with probability at least $1 - N^{3 - c\sqrt{d} \log(1 + \epsilon)}$,
\begin{align}
    \scA \preceq &\, \lambda \id_N + d^{-1}\rmD + \lambda \big( \lambda^{-1} \scA - \scA(\lambda) + \rmM(\lambda) - \id_{N} - \lambda^{-1} d^{-1}\rmD \big)\\
        \preceq &\, \id_N + d^{-1} \rmD + (\lambda - 1)\id_N + \lambda \big( \| \lambda^{-1}\scA - \scA(\lambda)\| + \|\rmM(\lambda) - \id_{N} - \lambda^{-1}d^{-1}\rmD\| \big)\cdot \id_{N},
\end{align}
    where the second inequality follows from the triangle inequality. By choosing $\lambda = 1 + \epsilon$ with $\epsilon = C d^{-1/2}$ for some constant $C >0$, together with \eqref{eqn:MH_lambda_PD} and \eqref{eqn:upper_bound_HMlambda}, \Cref{prop:upper_bound_H} follows.
    
Deterministically, following similar steps as in \cite[Proposition 6.1]{alt2021extremal}, we will show that for $\lambda \geq 1$
\begin{align}
    \|\scA(\lambda) - \lambda^{-1} \scA\| \lesssim \frac{2\Delta}{(\lambda \sqrt{d})^{3}}, \quad \|\rmM(\lambda) - \id_N - \lambda^{-1}d^{-1}\rmD\| \lesssim \frac{\Delta + d}{\lambda^4 d^2}. \label{eqn:upper_bound_HMlambda}
\end{align}

For the first part, one can conclude from the Schur test \Cref{lem:schur_test} that
    \begin{align}
        \|\scA(\lambda) - \lambda^{-1}\scA\| 
        \leq &\, \bigg( \max_{j\in \gV_{1}} \sum_{l\in \gV_{2}}|\ermX_{jl}(\lambda) - \lambda^{-1}\ermX_{jl}| \bigg)^{1/2} \cdot \bigg( \max_{l\in \gV_{2}} \sum_{j\in \gV_{1}}|\ermX_{jl}(\lambda) - \lambda^{-1}\ermX_{jl}|  \bigg)^{1/2}.
    \end{align}
Recall that $\rmX = (\widetilde{\rmA} - \E \widetilde{\rmA})/\sqrt{d}$, then for $j\in \gV_{1}$ and $l \in \gV_{2}$, $|\ermX_{jl}|\leq 1/\sqrt{d}$, and
    \begin{align}
        \sum_{l \in \gV_{2}}|\ermX_{jl}|^2 \leq \frac{\rD_{j}}{d} + m\cdot \frac{d}{mn} \leq \frac{\Delta}{d} + \frac{d}{n}, \,\, \sum_{j\in \gV_{1}}|\ermX_{jl}|^2 \leq \frac{\rD_{l}}{d} + n\cdot \frac{d}{mn} \leq \frac{\Delta}{d} + \frac{d}{m}.
    \end{align}
 Since $4 \leq d \leq (mn)^{1/13}$, the first part of \eqref{eqn:upper_bound_HMlambda} follows due to
    \begin{align}
         \max_{j\in \gV_{1}} \sum_{l\in \gV_{2}}|\ermX_{jl}(\lambda) - \lambda^{-1}\ermX_{jl}| \leq \max_{j\in \gV_1} \sum_{l\in \gV_2} \frac{|\ermX_{jl}| \cdot |\ermX_{jl}|^{2}}{\lambda(\lambda^2 - \ermX_{jl}^2 )} \leq \frac{2}{\lambda^{3} \sqrt{d}} \max_{j\in \gV_1}  \sum_{l \in \gV_{2}}|\ermX_{jl}|^2 \leq \frac{2\Delta}{(\lambda \sqrt{d})^{3}},
    \end{align}
    where the second to last inequality holds since $|\ermX_{jl}|^{2} \leq \lambda^{2}$ for $\lambda \geq 1$.
    
    For the second part of \eqref{eqn:upper_bound_HMlambda}, note that $\rmM(\lambda)$ and $\rmD$ are diagonal matrices, then by definition
    \begin{align}
        &\, \|\rmM(\lambda) - \id_{N} - \frac{1}{\lambda d} \rmD\| = \max\bigg\{ \max_{j\in \gV_1}\bigg| \ermM^{(1)}_{jj}(\lambda) - 1 - \frac{1}{\lambda d}\sum_{l\in \gV_2} \ermA_{jl} \bigg|, \, \max_{l\in \gV_2}\bigg| \ermM^{(2)}_{ll}(\lambda) - 1 - \frac{1}{\lambda d}\sum_{j\in \gV_1}\ermA_{jl}  \bigg| \bigg\}.
    \end{align}
 For any $j\in \gV_1$, we apply the triangle inequality to the following identity,
\begin{align}
 \ermM^{(1)}_{jj}(\lambda) - 1 - \frac{1}{\lambda d}\sum_{l\in \gV_2} \ermA_{jl} \leq &\, \sum_{l\in \gV_2} \frac{\ermX^{2}_{jl}}{\lambda^{2} - \ermX^{2}_{jl}} - \frac{1}{\lambda^2 d}\sum_{l\in \gV_2}\ermA_{jl}^{2} \\
 =&\, \sum_{l\in \gV_2} \bigg( \frac{\ermX_{jl}^4}{\lambda^2 d(\lambda^2 - \ermX^{2}_{jl})} - \frac{2\cdot\ermX_{uv}}{\lambda^2 \sqrt{mn}} - \frac{d}{\lambda^2 mn} \bigg).
\end{align}
Note that $\sqrt{mn} \gtrsim d^{13/2} \asymp \lambda^{13}$, in the right-hand side of the above equation, the last two terms inside the summation are small compared with the first one. By applying $|\ermX_{jl}| \leq 1/\sqrt{d}$ and $|\ermX_{jl}|^2 \leq \lambda^2/2$, one can obtain the desired upper bound for the second inequality in \eqref{eqn:upper_bound_HMlambda}.
\end{proof}

\subsection{Proofs in \Cref{sec:delocalization}}\label{sec:delocalization_app}
\begin{proof}[Proof of \eqref{eqn:upperbound_Hhatx}]
For simplicity, we consider $x\in\gV_1$. The same argument can be applied for $x\in\gV_2$. Recall $\hatH$ in \eqref{eqn:hatH} and $\hatHx$ in \eqref{eqn:hatHx} for $x\in \gV^{(\vtau^{\star})}$. We separately consider two cases below.

\textbf{Case 1 for \eqref{eqn:upperbound_Hhatx}:} $x\in \gV^{(\vtau^{\star})}$ with $\rvw \perp \hatvtau_{\pm}(x)$.

Based on \Cref{lem:existence_pruned_graph} (1) and (2), since $x\in \gV^{(\vtau )}$, $\gG^{(\vtau)}$ restricted $\ball^{(\vtau)}_{2\radius}(x)$ is a tree and its root is  $x$ with $\alpha_x d$ children, and at most $\tau_2^+d$ children for odd layers and $\tau_1 d$ children for even layers. Then, \Cref{lem:tree_norm} and \Cref{lem:hatHtau_norm} show that 
\begin{equation}\label{eqn:upperbound_Hhatx0}
    \|\hatHx\|\le (\tau_1\tau_2^+)^{1/4}\cdot \Lambda_{(\tau_1/\tau_2^+)^{1/4}}\Big(\frac{\alpha_x}{\sqrt{\tau_1\tau_2^+}}\vee\big(1+\sqrt{\tau_1/\tau_2^+}\big)\Big) +\gC\xi.
\end{equation}
When $\alpha_x\le \tau_1+\sqrt{\tau_1\tau_2^+}$,  \eqref{eqn:upperbound_Hhatx0} implies $ \|\hatHx\|\le \sqrt{\tau_1}+\sqrt{\tau_2^+}+\gC\xi$. Then, \eqref{eqn:upperbound_Hhatx} is verified because 
\begin{equation}\label{eqn:tau_1_tau_2_bound}
    \sqrt{\tau_1}+\sqrt{\tau_2^+} \le \qr \tau_1+\qr^{-1}\tau_2^+. 
\end{equation}
When $\alpha_x> \tau_1+\sqrt{\tau_1\tau_2^+}$, we have
\[(\tau_1\tau_2^+)^{1/4}\Lambda_{(\tau_1/\tau_2^+)^{1/4}}\Big(\frac{\alpha_x}{\sqrt{\tau_1\tau_2^+}}\vee\big(1+\sqrt{\tau_1/\tau_2^+}\big)\Big) = \sqrt{\alpha_x+\tau_2^++\frac{1}{\alpha_x-\tau_1}}\ge \Lambda_{\qr^{-1}}(\alpha_x).\]
Hence, when $\Lambda_{\qr^{-1}}(\alpha_x)\le \sqrt{\tau_1}+\sqrt{\tau_2^+}+\gC\xi$, \eqref{eqn:upperbound_Hhatx} can be obtained by the above bound \eqref{eqn:tau_1_tau_2_bound}. The only case left is when $\alpha_x> \tau_1+\sqrt{\tau_1\tau_2^+}$ and $\Lambda_{\qr^{-1}}(\alpha_x)> \sqrt{\tau_1}+\sqrt{\tau_2^+}+\gC\xi$. In this case, by Lemma~\ref{lem:tree_norm_rough} and \eqref{eqn:tau_1_tau_2_bound}, we know that
\[
\|\rmH^{(\vtau)}|_{\ball^{(\vtau)}_{\radius }(x)\setminus \{x\}}\|\le \sqrt{\tau_1}+\sqrt{\tau_2^+} \leq \qr \tau_1+\qr^{-1}\tau_2^+.
\]
Notice that by Weyl's inequality, there exist at most one eigenvalue of $\rmH^{(\vtau)}|_{\ball^{(\vtau)}_{\radius }(x) }$ strictly larger than $\sqrt{\tau_1}+\sqrt{\tau_2^+}$ and at most one strictly smaller than $-\sqrt{\tau_1}-\sqrt{\tau_2^+}$. Then by Lemma~\ref{lem:hatHtau_norm}, there exist at most one eigenvalue of $\hatH|_{\ball^{(\vtau)}_{\radius }(x) }$ strictly larger than $\sqrt{\tau_1}+\sqrt{\tau_2^+}+\gC\xi$ and at most one strictly smaller than $-\sqrt{\tau_1}-\sqrt{\tau_2^+}-\gC\xi$. Meanwhile, $(\Lambda_{\qr^{-1}}(\alpha_x),\hatvtau_{+}(x))$ and $(-\Lambda_{\qr^{-1}}(\alpha_x),\hatvtau_{-}(x) )$ are the eigen-pairs for $\hatH|_{\ball^{(\vtau)}_{\radius }(x) }$. Recall \eqref{eqn:hatHx}, in $\hatHx$, we eliminate these two eigen-paris, hence all the other eigenvalues are between $\sqrt{\tau_1}+\sqrt{\tau_2^+}+\gC\xi$ and   $-\sqrt{\tau_1}-\sqrt{\tau_2^+}-\gC\xi$. Thus, combining \eqref{eqn:tau_1_tau_2_bound}, we can conclude \eqref{eqn:upperbound_Hhatx}.

\textbf{Case 2 for \eqref{eqn:upperbound_Hhatx}: $x\in\pruneV\setminus\gV^{(\vtau^*)}$.}

In this case, we denote $\hatHx:=\hatH |_{\ball^{(\vtau)}_{2\radius }(x) }$ and $\tau_1<\alpha_x\le 1+\qr^{-2}+\xi^{1/2}$. We can check that $\frac{\alpha_x}{\sqrt{\tau_1\tau_2^+}}\leq \big(1+\sqrt{\tau_1/\tau_2^+}\big)$. Then, with \eqref{eqn:upperbound_Hhatx0} and \eqref{eqn:tau_1_tau_2_bound}, we can conclude \eqref{eqn:upperbound_Hhatx} for this case.

The proof is then completed by applying \eqref{eqn:tau_1_tau_2_bound}, since $\qr \tau_1+\qr^{-1}\tau_2^+= \qr+\qr^{-1}+\frac{1}{2}\xi^{1/4}$ due to the values of $\vtau = (\tau_{1}, \tau_2^+, \tau_2^{-})$ in \eqref{eqn:tau_constrain}.
\end{proof}   

\begin{proof}[Proof of \eqref{eqn:supp_in_ball} ]
The proof is finished by induction. One can easily check the base case for $j=0$. Suppose that for $0\le j < r_*$, \eqref{eqn:supp_in_ball} holds, then we need to show
\begin{equation} 
    \mathrm{supp}\bigl(\hatHx\,\rvg_j \bigr) \subset \ball^{(\vtau)}_{\radius + j+1}(x).
\end{equation}
Recall that for $j+k\le 2\radius$ and for any $x\in \pruneV$ and vector $\rvv$, we have
\begin{equation}\label{eqn:containing_support}
\mathrm{supp} (\rvv) \subset \ball^{(\vtau)}_{j}(x)
\quad\Longrightarrow\quad
\mathrm{supp}\bigl((\rmH^{(\vtau)})^{k} \rvv\bigr) \subset \ball^{(\vtau)}_{j+k}(x) 
\end{equation}
and  $\mathrm{supp}(\hatvtau_{\sigma}(x)) \subset \ball^{(\vtau)}_{\radius}(x)$.
Then, by induction and \Cref{lem:existence_pruned_graph} (1), we have
    \[
        \hatHx \rvg_j = \Bigl(\id - \sum_{\sigma=\pm}\hatvtau_{\sigma}(x) [\hatvtau_{\sigma}(x)]^{*}\Bigr)
        \rmH^{(\vtau)}
    \Bigl(\id - \sum_{\sigma=\pm} \hatvtau_{\sigma}(x) [\hatvtau_{\sigma}(x)]^{*}\Bigr)
    \,\rvg_j,
    \]
    and we conclude \eqref{eqn:supp_in_ball}.
\end{proof}

\subsection{Proofs in \Cref{sec:proof_bulk_bound_maintext}}\label{sec:proof_bulk_bound_app}
\begin{proposition}\label{prop:approximate_eigenvalues_htau}
Suppose that $\vtau$ satisfies \eqref{eqn:tau_constrain} and let $\sigma = \pm$. Then, for any $x\in\gV^{(\vtau^*)}$ where $\gV^{(\vtau^*)}$ is defined in \eqref{def:W_set}, with very high probability, the following holds:
\begin{align}
    \| (\scA^{(\vtau)}-\sigma\cdot(\Lambda_\qr(\alpha_x)\indi{x\in\gV_2}+\Lambda_{\qr^{-1}}(\alpha_x)\indi{x\in\gV_1}))  \hatvtau_{\sigma}(x)\|\leq \const \error.
\end{align} 
\end{proposition}
\begin{proof}
    The proof is analogous to \cite[Proposition 3.9] {alt2021delocalization} by adjusting the proof of Proposition~\ref{prop:approximate_eigenvalues} to $\scA^{(\vtau)}$. 
\end{proof}

\begin{proof}[Proof of \Cref{lem:hatHtau_norm}]
The first inequality is proved by applying Lemma~\ref{lem:existence_pruned_graph}, where we use \eqref{eqn:approx_H_Htau} and the definition of $\error$ in \eqref{eqn:error_parameters}.

The second inequality follows the same idea as Lemma 3.11 of \cite{alt2021delocalization}. Notice that
\begin{align}
& \scA^{(\vtau)} - \hatH \label{eqn:decomp_hatHtau_norm}\\
 =~& \Pi^{(\vtau)} \scA^{(\vtau)} \Pi^{(\vtau)}
- \sum_{x\in \gV^{(\vtau^*)}}\sum_{\sigma = \pm} \sigma \cdot (\Lambda_\qr(\alpha_x)\indi{x\in\gV_2}+\Lambda_{\qr^{-1}}(\alpha_x)\indi{x\in\gV_1})  \hatvtau_{\sigma}(x) [\hatvtau_{\sigma}(x)]^{*}\\
& + \overline{\Pi^{(\vtau)}} \scA^{(\vtau)} \Pi^{(\vtau)} +(\overline{\Pi^{(\vtau)}} \scA^{(\vtau)} \Pi^{(\vtau)})^*.
\end{align}
Meanwhile, Proposition~\ref{prop:approximate_eigenvalues_htau} and \eqref{eqn:containing_support} indicate that  
for any $x \in \gV^{(\vtau^*)}$, these two cases imply that
\[
    \scA^{(\vtau)}\hatvtau_{\sigma}(x) =  \sigma\cdot(\Lambda_\qr(\alpha_x)\indi{x\in\gV_2}+\Lambda_{\qr^{-1}}(\alpha_x)\indi{x\in\gV_1})\hatvtau_{\sigma}(x) + \mathbf{e}_\sigma^{(\vtau)}(x), 
\]
where $\text{supp}\, \mathbf{e}_\sigma^{(\vtau)}(x) \subset \ball_{2\radius}^{(\vtau)}(x)$ and $ \|\mathbf{e}_\sigma^{(\vtau)}(x)\| \leq C\xi$ with very high probability. By Lemma~\ref{lem:existence_pruned_graph}, the balls $\ball_{2\radius}^{(\vtau)}(x)$ and $\ball_{2\radius}^{(\vtau)}(y)$
 are disjoint for $x,y \in \gV^{(\vtau)} $ with $x \neq y$. 
Therefore, $\hatvtau_{\sigma}(x), \mathbf{e}^{(\vtau)}_\sigma(x) \perp \hatvtau_{\sigma'}(y), \mathbf{e}^{(\vtau)}_{\sigma'}(y)$ for any $\sigma,\sigma'\in\{\pm \}$ and $x \neq y\in \gV^{(\vtau)} $. 
 For any vector $\mathbf{a} = \sum_{x \in \gV^{(\vtau^*)}} \sum_{\sigma=\pm} a_{x,\sigma}\hatvtau_{\sigma}(x)$, we obtain
\[
\overline{\Pi^{(\vtau)}} \scA^{(\vtau)} \Pi^{(\vtau)} \mathbf{a} = \sum_{x \in \gV^{(\vtau^*)}} \sum_{\sigma=\pm}   a_{x,\sigma}\overline{\Pi^{(\vtau)}} \scA^{(\vtau)}\hatvtau_{\sigma}(x)
= \overline{\Pi^{(\vtau)}} \sum_{x \in \gV^{(\vtau^*)}} \sum_{\sigma=\pm}  a_{x,\sigma}\mathbf{e}^{(\vtau)}_\sigma(x).
\]
By orthogonality, with very high probability,
we obtain $\|\overline{\Pi^{(\vtau)}} \scA^{(\vtau)} \Pi^{(\vtau)}\rva\|^2\leq \const^2\xi^2\|\rva\|^2 $ with very high probability. Hence, $\|\overline{\Pi^{(\vtau)}} \scA^{(\vtau)} \Pi^{(\vtau)}\|,\|(\overline{\Pi^{(\vtau)}} \scA^{(\vtau)} \Pi^{(\vtau)})^*\|\le C\error$ with very high probability. Similarly, in the rest term of \eqref{eqn:decomp_hatHtau_norm}, we know that
\begin{align}
&\left\|\left(\Pi^{(\vtau)} \scA^{(\vtau)} \Pi^{(\vtau)} - \sum_{x\in \gV^{(\vtau^*)}}\sum_{\sigma = \pm} \sigma \cdot (\Lambda_\qr(\alpha_x)\indi{x\in\gV_2}+\Lambda_{\qr^{-1}}(\alpha_x)\indi{x\in\gV_1})  \hatvtau_{\sigma}(x) [\hatvtau_{\sigma}(x)]^{*}\right)\mathbf{a} \right\|\\
=&\left\|\Pi^{(\vtau)}\sum_{x\in \gV^{(\vtau^*)}}\sum_{\sigma = \|\pm}a_{x, \sigma}\mathbf{e}^{(\vtau)}_\sigma(x)\right\|\leq C^2 \xi^2\|\rva\|^2.
\end{align}
Hence, we can conclude the proof of this lemma.  
\end{proof}

\section{Extension to general sparse random rectangular matrices.}\label{sec:sparse_rectangular}
This section is devoted to the proof of \Cref{thm:sparse_rectangular}. We only address necessary changes from the proof of \Cref{thm:right_edge_behavior} and \Cref{thm:left_edge_behavior}.

\subsection{Approximation error} 
Recall the definition of the sparse matrix $\rmM$ in \Cref{thm:sparse_rectangular}, which can be viewed as the adjacency matrix of an undirected weighted bipartite graph. We can define
\[\layer_j(x)=\{y\in[N]:\min\{t\ge 0: (\rmM^t)_{xy}\neq 0\}=j\}\]
and $\ball_{r}(x)=\cup_{j\in [r]}\layer_j(x)$. It turns out these sets are identical to the corresponding subsets defined before for the bipartite Erd\H{o}s-R\'{e}nyi graph $\gG$ sampled from $\gG(n, m,d/\sqrt{mn})$. 

Due to \Cref{ass:homogeneous_variance}, where entries of $\rmM$ maintain the same variance, we obtain the same tridiagonal matrix \eqref{eqn:Z1Z2} after Gram-Schmidt orthonormalization. The construction of approximate eigenvectors follows similarly by leveraging the new definition of the normalized degree $\alpha_{x}$ in \eqref{eqn:new_degree} and $\layer_j(x)$ above.

We now introduce the analogue of the approximate eigenvector $\rvv$ from \eqref{eqn:approx_eigenvector} for $x\in \gV_{2}$. The construction for $x\in \gV_{1}$ follows similarly. Define $\rvg_{0} \coloneqq \ones_{x}$. For $j\geq 1$, we define
\begin{align}\label{eqn:g_j_weighted}
    \rvg_{j} \coloneqq (\rmM \rvg_{j-1})|_{\layer_{j}(x)}.
\end{align}
With the choices of $\ervu_{j}$ in \eqref{eqn:u012V2} and \eqref{eqn:ujsV2}, we set
\begin{align}
    \rvv = \sum_{j=0}^{r}  \frac{\ervu_{j}}{\|\rvg_{j}\|} \rvg_{j}.
\end{align}where $\|\cdot\|$ is the $\ell_2$ norm.  
Below, we prove that  $\rvv$ is an approximate eigenvector for $\rmH$.  

Recall that $\rmH=\rmM/\sqrt{d}$. Similar to the proof of \Cref{lem:five_term_decomp}, we have
\begin{align}
    (\scA - \Lambda_{\qr}(\alpha_{x})) \rvv = \rvw_{1} + \rvw_{2} + \rvw_{3} + \rvw_{4},
\end{align}
where the error terms $\rvw_{1}, \rvw_{2}, \rvw_{3}, \rvw_{4}$ are defined through
\begin{subequations}\label{eqn:five_term_decomp_Wigner}
\begin{align}
    \rvw_1 &\, \coloneqq \frac{1}{\sqrt{d}}\sum_{j=0}^r \frac{\ervu_{j}}{\|\rvg_{j}\| }\Bigg( \sum_{y\in \layer_{j}} |\NS{j}{y}| \ones_y + \sum_{y\in \layer_{j+1}}(|\NS{j}{y}|- \<\ones_{y}, \rvg_{j+1}\> )\ones_y \Bigg)\,,\\
    \rvw_2 &\, \coloneqq \frac{1}{\sqrt{d}}\sum_{j=1}^r \frac{\ervu_{j}}{\|\rvg_{j}\|} \sum_{y\in \layer_{j-1}}\bigg(|\NS{j}{y}|\cdot \ones_y -\frac{\|\rvg_{j}\|^2}{\|\rvg_{j-1}\|^2} \cdot \rvg_{j-1} \bigg)\,,\\
    \rvw_3 &\, \coloneqq \frac{\ervu_2}{\|\rvg_1\|} \Big( \frac{\|\rvg_{2}\|}{\sqrt{d}\cdot \sqrt{\|\rvg_{1}\|}} - \qr^{-1} \Big)\rvg_1 \\
    &\, \quad \quad + \sum_{j=1}^{\lfloor r/2\rfloor - 1 } \frac{1}{\|\rvg_{2j}\|}\cdot \Bigg( \ervu_{2j-1}\Big(\frac{ \|\rvg_{2j}\|}{ \sqrt{d} \cdot \|\rvg_{2j-1}\|} - \qr^{-1} \Big) + \ervu_{2j+1}\bigg(\frac{\|\rvg_{2j+1}\|}{\sqrt{d}\cdot \|\rvg_{2j}\|} - \qr \bigg) \Bigg) \rvg_{2j}\\
    &\, \quad \quad + \sum_{j=1}^{\lfloor r/2\rfloor - 1} \frac{1}{\|\rvg_{2j+1}\|}\cdot \Bigg( \ervu_{2j}\Big(\frac{\|\rvg_{2j+1}\|}{\sqrt{d} \cdot \|\rvg_{2j}\|} - \qr \Big) + \ervu_{2j+2}\bigg(\frac{\|\rvg_{2j}\|}{\sqrt{d}\cdot \|\rvg_{2j}\|} - \qr^{-1} \bigg) \Bigg) \rvg_{2j+1}\,, \notag\\
    \rvw_4 &\, \coloneqq \frac{\ervu_{r-1}}{\sqrt{d}} \frac{1}{\|\rvg_{r-1}\|} \, \rvg_{r} + \frac{\ervu_{r}}{\sqrt{d}}\frac{1}{\|\rvg_{r}\|} \, \rvg_{r+1}  \\
    &\, \quad - \indi{ r\textnormal{ even}} \big( \qr^{-1} \ervu_{r-1} + \qr \ervu_{r+1} \big) \frac{1}{\|\rvg_{r}\|} \cdot \rvg_{r} - \indi{ r\textnormal{ odd}} \big( \qr \ervu_{r-1} + \qr^{-1} \ervu_{r+1}\big) \frac{1}{\|\rvg_{r}\|} \cdot \rvg_{r}.\notag
\end{align}
\end{subequations}
Here we ewdefine $|\NS{j}{y}|:=\langle \ones_y,\rmM\rvg_j\rangle$ for all $y\in [N]$.
Note that there is no $\rvw_{0}$ compared with \Cref{lem:five_term_decomp}, since the matrix $\rmM$ itself is already centered.

First, we prove that there is some constant $\gC>0$ such that
\begin{equation}\label{eqn:linfty_bound}
    \|\rvg_j\|_{\infty}\leq (\gC\gK)^j,
\end{equation}
where $\gK$ is the uniform bound for the entries of Wigner matrix $\widetilde \rmW$ in \Cref{thm:sparse_rectangular}. For $j=0$, we have $\|\rvg_0\|_{\infty} = 1$. Assume that \eqref{eqn:linfty_bound} holds for $j-1$. Then, by definition in \eqref{eqn:g_j_weighted}, we have for any $y\in \layer_{j}(x)$,
\begin{align}
 \left|\ones_y^\sT\rvg_j\right|= ~&\left|\ones_y^\sT\rmM \rvg_{j-1} \right|=  \left|\sum_{z\in \layer_{j-1}(x)} \ermM_{yz} \ones_z^\sT \rvg_{j-1} \right|  \\
 \leq  ~&(\gC\gK)^{j-1}\sum_{z\in \layer_{j-1}(x)} |\ermM_{yz}| = (\gC\gK)^{j-1}\gK\cdot\sum_{z\in \layer_{j-1}(x)}  \ermA_{yz}.
\end{align}
Hence, we return to the random bipartite graph $\gG$ and apply \eqref{eqn:few_cycle_different_layer} in \Cref{cor:few_cycle_counts} to obtain $\left|\ones_y^\sT\rvg_j\right|\le (\gC\gK)^j$ for all $y\in \layer_{j}(x)$. Thus, we have shown \eqref{eqn:linfty_bound} for all $j\ge 0$.

Next, we can mimic the proof of \Cref{lem:concentrationSi} to show the concentration on $\|\rvg_{j+1}\|/\|\rvg_{j}\|$. Without loss of generality, we consider $j$ is even, then $\layer_{j+1}(x)\subset \gV_1$. Let us define
\begin{equation}
    \rY_y \coloneq \frac{\rvg_{j+1}^\sT\ones_y}{\sqrt{d}\|\rvg_j\|} = \frac{1}{\sqrt{d}\|\rvg_j\|}\sum_{z\in \layer_j(x)} \ermM_{yz}\,\rvg_j^\sT\ones_z,
\end{equation}
where $\ermM_{yz}=\ermA_{yz}\ermW_{yz}$ with $\ermA_{yz}\sim \mathrm{Bernoulli}(p)$, and $\ermW_{yz}$ are independent entries with mean zero and bounded by $\gK$.  We aim to control
\begin{equation}
  \rS \coloneqq \sum_{y\in \gV_1} \bigl(\rY_{y}^{2} - \E[\rY_{y}^{2}\mid  \ball_{j}(x),\rvg_j]\bigr).
\end{equation}
Due to the distributions of $\rmA$ and $\rmW$, we have $\E[\rY_{y}^{2}\mid  \ball_{j}(x),\rvg_j]= p/d=1/\sqrt{nm}$. Hence, we can get
\[\rS = \frac{\|\rvg_{j+1}\|^2}{d\|\rvg_j\|^2}-\qr^2.\]
Notice that $|\rY_y|\le \frac{(\gC\gK)^{j+1}}{\sqrt{d}\|\rvg_j\|}$, and then $$\bigl|\rY_{y}^{2} - \E[\rY_{y}^{2}\mid  \ball_{j}(x),\rvg_j]\bigr|\le  \frac{2(\gC\gK)^{2j+2}}{d\|\rvg_j\|^2}$$
uniformly for all $y\in \gV_1$.
Meanwhile, we know that 
\begin{align}
    \sigma^2:=\sum_{y\in\gV_1}\Var (\rY_{y}^{2})\le \sum_{y\in\gV_1}\E[Y_y^4]\le \frac{\qr^2\|\rvg_j\|^4_4}{d\|\rvg_j\|^4}.
\end{align}
Then we can apply Bernstein inequality in Lemma~\ref{lem:Bernstein} for $S$ to obtain
\begin{align}
    \P\bigg(\bigg| \frac{\|\rvg_{j+1}\|^2}{d\|\rvg_j\|^2} - \qr^2\bigg| > \epsilon\Big| \rvg_j,\ball_j(x)\bigg) &\leq 2\exp\Big(-\frac{\epsilon^2d\|\rvg_j\|^2/2}{\frac{\qr^2\|\rvg_j\|^4_4}{\|\rvg_j\|^2} +  \epsilon\frac{2(\gC\gK)^{2j+2}}{3}}\Big)\\
    &\leq 2\exp\left(-\frac{c_1 d \|\rvg_j\|^2\epsilon^2}{(\gC\gK)^{2j}}\right),
\end{align}
for some constant $c_1>0$ and any sufficiently small $\epsilon>0$, where we use the fact that $\|\rvg_j\|_4^2\leq (\gC\gK)^j\|\rvg_j\|$. Hence we get a similar concentration inequality for $\|\rvg_{j+1}\|^2$ as \eqref{eqn:l_even_and_x_V2}. We can also apply the same argument for $j$ odd to get an analogous result as \eqref{eqn:l_odd_and_x_V2}. As a result, we have that 
for all small $\epsilon>0$,
\begin{subequations}
    \begin{align}
        &\, \P\Big( \Big| \,\|\rvg_{2j}\|^2- d_{1}\|\rvg_{2j-1}\|^2 \,\Big| \leq  \epsilon  d\|\rvg_{2j-1}\|^2 \Big| \ball_{2j-1} ,\rvg_{2j-1}\Big) \geq 1 - 2 e^{-c_1\epsilon^2 \frac{d \|\rvg_{2j}\|^2}{(\gC\gK)^{2(2j-1)}}},\label{eqn:g_odd_and_x_V2}\\
        &\, \P\Big( \Big| \,\|\rvg_{2j+1}\|^2- d_{2}\|\rvg_{2j}\|^2 \,\Big| \leq  \epsilon  d\| \rvg_{2j}\|^2 \Big| \ball_{2j},\rvg_{2j}\Big) \geq 1 - 2 e^{-c_1\epsilon^2 \frac{d \|\rvg_{2j}\|^2}{(\gC\gK)^{4j}}}.\label{eqn:g_even_and_x_V2}
    \end{align}
\end{subequations}
Notice that $\|\rvg_1\|^2=\rD_x := \sum_{y\in[m]}|\rmM_{xy}|^2$. Thus, we can also apply Bernstein inequality in Lemma~\ref{lem:Bernstein} to obtain $\rD_x\in(d_1/2,3d_1/2)$ with high probability. Therefore, we can apply the same induction argument as in the proof of Lemma~\ref{lem:concentrationSi} to simultaneously show that
\begin{subequations}
    \begin{align}
   & \bigg|\frac{1}{d_1}\cdot\frac{\|\rvg_{2j}\|^2}{\|\rvg_{2j-1}\|^2} -1 \bigg| \lesssim \bigg( \frac{\log(N)}{d_{1}\|\rvg_{2j-1}\|^2}\bigg)^{1/2},    \label{eqn:g2j_ratio_bound}\\
       &\bigg| \frac{1}{d_2}\cdot \frac{\|\rvg_{2j+1}\|^2}{\|\rvg_{2j}\|^2} - 1 \bigg|  \lesssim \bigg( \frac{\log(N)}{d_{2}\|\rvg_{2j}\|^2}\bigg)^{1/2}, \label{eqn:g2j+1_ratio_bound}
    \end{align}
\end{subequations}
and
\begin{subequations}
    \begin{align}
        \rD_x \left(d/2\right)^{2j-2}\leq \|\rvg_{2j-1}\|^2\leq \rD_x (2d)^{2j-2}\label{eqn:g_ratio2j-1}\\
        \frac{1}{2}\rD_x d_1\left(d/2\right)^{2j-2}\leq \|\rvg_{2j}\|^2\leq 2\rD_x d_1 (2d)^{2j-2} \label{eqn:g_ratio2j}
    \end{align}
\end{subequations}
for $1\leq j\leq \lceil r/2 \rceil$, with high probability. Thus, we can obtain a similar argument as in \Cref{lem:concentrationSi} for the concentration on $\|\rvg_{j+1}\|^2/\|\rvg_{j}\|^2$ and $\|\rvg_{j}\|^2$. We can have the control on $\|\rvw_3\|$ and $\|\rvw_4\|$ in \eqref{eqn:five_term_decomp_Wigner} analogously with the proof in \Cref{lem:w0tow5}.

The estimate of $\|\rvw_{1}\|$ is similar to the proof in Section 10 of \cite{alt2021extremal} by reproving \eqref{eqn:few_cycle_same_layer} in \Cref{cor:few_cycle_counts} as follows. For $j\ge 1$, we have
\begin{align}
\sum_{y \in \layer_j(x)} |\NS{j}{y}|^2 
= \sum_{y \in \layer_j(x)} \langle \ones_y, \rmM \rvg_j \rangle^2 
= \sum_{y \in \layer_j(x)} \left( \sum_{y_1 \in \layer_j(x)} \langle \rmM \mathbf{1}_y, \mathbf{1}_{y_1} \rangle \langle \mathbf{1}_{y_1}, \mathbf{g}_j \rangle \right)^2 
=0,
\end{align}
since $y_1 \in \layer_j(x)$ and $y \in \layer_j(x)$ imply that $y_1$ and $y$ are both in either $\gV_1$ or $\gV_2$, and $\langle \rmM \mathbf{1}_y, \mathbf{1}_{y_1} \rangle =0$ in this case. Moreover, $ |\NS{j}{y}| - \langle \mathbf{1}_y, \mathbf{g}_{j+1} \rangle = \langle \mathbf{1}_y, (\rmM\rvg_i)|_{[N] \setminus S_{j+1}(x)} \rangle = 0$, for any $ y \in S_{j+1} $. Hence,
\[
\|\rvw_1\|^2 
= \left\| \sum_{j=1}^r \frac{\ervu_{j}}{\|\rvg_i\|} \sum_{y \in S_i}|\NS{j}{y}| \mathbf{1}_y \right\|^2 
= \sum_{j=1}^r \frac{\ervu_{j}^2}{\|\rvg_j\|^2} \sum_{y \in S_j} |\NS{j}{y}| ^2=0.
\]

The estimate of $\|\rvw_{2}\|$ is similar to the proofs in Section 10 of \cite{alt2021extremal} combined with the proofs of \eqref{eqn:norm_w2} in \Cref{lem:w0tow5}. Hence, we ignore the details here.

\subsection{Pruned graph} In the pruned graph, each edge $\ermA_{xy}$ is assigned a weight $\ermW_{xy}$. The conclusions in \Cref{lem:existence_pruned_graph} follow directly for $\rmM$, since the tool (Bennett inequality in \Cref{lem:Bennett}) is still applicable for bounded random variables.

\subsection{Non-backtracking operators} Results for $\rmM$ analogous to \Cref{prop:upper_bound_H} and \Cref{prop:lower_bound_H} follow similarly since Lemmas \ref{lem:Ihara-Bass}, \ref{lem:block_Ihara} and \ref{lem:upper_bound_rhoB} hold for matrices with bounded entries.
This finishes the proof of \Cref{thm:sparse_rectangular}.

\section{Technical Lemmas}

\begin{lemma}[Weyl's inequality, \cite{Weyl1912DasAV}]\label{lem:weyl}
Let $\rmA, \rmE \in \R^{m \times n}$ be two real $m\times n$ matrices, then $|\sigma_i(\rmA + \rmE) - \sigma_i(\rmA)| \leq \|\rmE\|$ for every $1 \leq i \leq \min\{ m, n\}$. Furthermore, if $m = n$ and $\rmA, \rmE \in \R^{n \times n}$ are real symmetric, then $|\lambda_i(\rmA + \rmE) - \lambda_i(\rmA)| \leq \|\rmE\|$ for all $1 \leq i \leq n$.
\end{lemma}

\begin{lemma}[Bernstein's inequality, {\cite[Theorem $2.8.4$]{vershynin2018high}}]\label{lem:Bernstein}
    Let $\rX_1,\dots, \rX_n$ be independent mean zero random variables such that $|\rX_i|\leq K$ for all $i$. Let $\sigma^2 = \sum_{i=1}^{n}\E \rX_i^2$. Then for every $t \geq 0$,
    \begin{align}
        \P \Bigg( \Big|\sum_{i=1}^{n} \rX_i \Big| \geq t \Bigg) \leq 2 \exp \Bigg( - \frac{t^2/2}{\sigma^2 + Kt/3} \Bigg)\,.
    \end{align}
\end{lemma}

\begin{lemma}[Bennett's inequality, {\cite[Theorem $2.9.2$]{vershynin2018high} }]\label{lem:Bennett}
    Let $\rX_1,\dots, \rX_n$ be independent random variables. Assume that $|\rX_i - \E \rX_i| \leq K$ almost surely for every $i$. Then for any $t>0$, we have
    \begin{align}
        \P \Bigg( \sum_{i=1}^{n} (\rX_i - \E \rX_i) \geq t \Bigg) \leq \exp \Bigg( - \frac{\sigma^2}{K^2} \cdot \benrate \bigg( \frac{Kt}{\sigma^2} \bigg)\Bigg)\,, \notag 
    \end{align}
    where $\sigma^2 = \sum_{i=1}^{n}\Var(\rX_i)$ is the variance of the sum and $\benrate(u) \coloneqq (1 + u)\log(1 + u) - u$. Furthermore, define $\rY_i = - \rX_i$ and apply the inequality above, for any $t >0$, we then have
    \begin{align}
        \P \Bigg( \sum_{i=1}^{n} (\rX_i - \E \rX_i) \leq -t \Bigg) \leq \exp \Bigg( - \frac{\sigma^2}{K^2} \cdot \benrate \bigg( \frac{Kt}{\sigma^2} \bigg)\Bigg)\,. 
    \end{align}
    
\end{lemma}

\begin{lemma}[Cayley's formula]\label{lem:Cayley_formula}
For $n \geq 2$, the number of trees on $n$ labeled vertices is $n^{n-2}$.
\end{lemma}

\begin{lemma}\label{lem:stirling}
    For integers $n, k \geq 1$, we have
    \begin{align}
    \log(n!) = &\, n\log(n) - n + \frac{1}{2}\log(2\pi n) + O(n^{-1})\\
        \log \binom{n}{k} =&\, \frac{1}{2}\log \frac{n}{2\pi k (n-k)} + n\log(n) - k\log(k) - (n-k) \log(n-k) + O\Big( \frac{1}{n} + \frac{1}{k} + \frac{1}{n-k} \Big)\\
        \log \binom{n}{k} =&\, k\log\Big( \frac{n}{k} - 1 \Big) -\frac{1}{2}\log(2\pi k) + o(k^{-1})\,, \textnormal{ for } k = \omega(1) \textnormal{ and } \frac{k}{n} = o(1).
    \end{align}
    Moreover, for any $1\leq k \leq \sqrt{n}$, we have
    \begin{align}
        \frac{n^{k}}{4 \cdot k!} \leq \binom{n}{k} \leq \frac{n^{k}}{k!}\,,\quad \log \binom{n}{k} \geq k\log \Big(\frac{en}{k} \Big) - \frac{1}{2}\log(k)  - \frac{1}{12k} - \log(4\sqrt{2\pi})\,.
    \end{align}
\end{lemma}

\begin{lemma}[Schur test in matrix version, \cite{schur1911bemerkungen}]\label{lem:schur_test}
    Let $\{a_j\}_{j=1}^{n}$ and $\{b_l\}_{l=1}^{m}$ be two sequences of positive real numbers, and let $\lambda, \mu$ be positive real numbers. For matrix $\rmX \in \R^{n \times m}$, we have $\|\rmX\| \leq \sqrt{\lambda \mu}$ if
    \begin{align}
        \sum_{j=1}^{n} |\ermX_{jl}| \cdot a_{j} \leq \lambda \, b_{l}, \,\, \forall l\in [m], \quad \textnormal{ and }\quad \sum_{l=1}^{m} |\ermX_{jl}| \cdot b_{l} \leq \mu\, a_{j}, \,\, \forall j \in [n].
    \end{align}
    As a consequence, we take $a_{j} = b_{l} = 1$ and choose $\lambda$ (resp. $\mu$) as the maximum column (resp. row) sum, then
    \begin{align}
        \|\rmX\| \leq \bigg( \max_{j\in [n]} \sum_{l=1}^{m} |\ermX_{jl}| \bigg)^{1/2} \cdot \bigg( \max_{l\in [m]}\sum_{j=1}^{n} |\ermX_{jl}| \bigg)^{1/2}.
    \end{align}
\end{lemma}
\end{document}